\documentclass{article}
\usepackage{euscript,amsfonts,amssymb,amsmath,multicol,amsthm,
longtable,graphicx,epsfig,hhline,array,amscd,url,chngcntr,apptools,bm,mathabx}
\usepackage[all]{xy}
\usepackage[usenames, dvipsnames]{color}
\usepackage{tikz-cd}
\usepackage{hyperref}
\usepackage{dutchcal}
\tikzcdset{scale cd/.style={every label/.append style={scale=#1},
    cells={nodes={scale=#1}}}}

\author{A.~S.~Berdnikov, A.~G.~Gorinov\\ with an appendix by A.~G.~Gorinov and N.~S.~Konovalov}
\title{Conical resolutions and cohomology of the moduli spaces of nodal hypersurfaces}
\date{}
\newcolumntype{C}{>{$\displaystyle} c <{$}}
\makeatletter
\newcommand\level[1]{%
  \ifcase#1\relax\expandafter\chapter\or
    \expandafter\section\or
    \expandafter\subsection\or
    \expandafter\subsubsection\else
    \def\next{\@level{#1}}\expandafter\next
  \fi}
\newcommand{\@level}[1]{%
  \@startsection{level#1}
    {#1}
    {\z@}%
    {-3.25ex\@plus -1ex \@minus -.2ex}%
    {1.5ex \@plus .2ex}%
    {\normalfont\normalsize\bfseries}}

\newdimen\@leveldim
\newdimen\@dotsdim
{\normalfont\normalsize
 \sbox\z@{0}\global\@leveldim=\wd\z@
 \sbox\z@{.}\global\@dotsdim=\wd\z@
}

\newcounter{level4}[subsubsection]
\@namedef{thelevel4}{\thesubsubsection.\arabic{level4}}
\@namedef{level4mark}#1{}
\def\l@section{\@dottedtocline{1}{0pt}{\dimexpr\@leveldim*4+\@dotsdim*1+6pt\relax}}
\def\l@subsection{\@dottedtocline{2}{0pt}{\dimexpr\@leveldim*5+\@dotsdim*2+6pt\relax}}
\def\l@subsubsection{\@dottedtocline{3}{0pt}{\dimexpr\@leveldim*6+\@dotsdim*3+6pt\relax}}
\@namedef{l@level4}{\@dottedtocline{4}{0pt}{\dimexpr\@leveldim*7+\@dotsdim*4+6pt\relax}}

\count@=4
\def\@ncp#1{\number\numexpr\count@+#1\relax}
\loop\ifnum\count@<100
  \begingroup\edef\x{\endgroup
    \noexpand\newcounter{level\@ncp{1}}[level\number\count@]
    \noexpand\@namedef{thelevel\@ncp{1}}{%
      \noexpand\@nameuse{thelevel\@ncp{0}}.\noexpand\arabic{level\@ncp{1}}}
    \noexpand\@namedef{level\@ncp{1}mark}####1{}%
    \noexpand\@namedef{l@level\@ncp{1}}%
      {\noexpand\@dottedtocline{\@ncp{1}}{0pt}{\the\dimexpr\@leveldim*\@ncp{5}+\@dotsdim*\@ncp{0}\relax}}}%
  \x
  \advance\count@\@ne
\repeat
\makeatother
\newcommand{\Z}{{\mathbb Z}}
\newcommand{\p}{{\mathrm P}}
\newcommand{\C}{{\mathbb C}}
\newcommand{\R}{{\mathbb R}}
\newcommand{\Q}{{\mathbb Q}}
\newcommand{\OO}{{\mathcal{O}}}
\newcommand{\kk}{\mathbf{k}}
\newcommand{\eu}{\EuScript}
\newcommand{\N}{\eu N}
\newcommand{\X}{\mathcal{X}}

\newcommand{\z}{_{\setminus 0}}
\newcommand{\nod}{{\text{\upshape{Nod}}}}
\AtAppendix{\counterwithin{lemma}{section}}
\AtAppendix{\counterwithin{theorem}{section}}
\AtAppendix{\counterwithin{corollary}{section}}
\AtAppendix{\counterwithin{Prop}{section}}
\AtAppendix{\numberwithin{equation}{section}}
\newtheorem{theorem}{Theorem}
\newtheorem{cl}{Condition list}

\newtheorem{lemma}{Lemma}
\newtheorem{Prop}{Proposition}
\newtheorem{corollary}{Corollary}
\newtheorem*{assumptions}{Assumptions}
\newtheorem*{assumption}{Assumption}
\DeclareMathOperator\Gr{Gr}
\DeclareMathOperator\Grd{Grd}
\DeclareMathOperator\Geom{Geom}
\DeclareMathOperator\Sing{Sing}
\DeclareMathOperator\rk{rk}
\DeclareMathOperator\tot{tot}
\DeclareMathOperator\Hom{Hom}
\DeclareMathOperator\Hess{Hess}
\DeclareMathOperator\Lk{Lk}
\DeclareMathOperator\Gdm{Gdm}

\DeclareMathOperator\coker{coker}
\DeclareMathOperator\id{id}
\DeclareMathOperator\im{im}
\newcommand{\SL}{\mathrm{SL}}

\newcommand{\PGL}{\mathrm{PGL}}
\newcommand{\GL}{\mathrm{GL}}
\newcommand{\Hyp}{\mathrm {Hyp}}
\begin{document}
\maketitle
\begin{abstract}
A projective hypersurface is nodal if it does not have singularities worse than simple nodes. We calculate the rational cohomology of the spaces of equations of nodal cubic and quartic plane curves and also nodal cubic surfaces in the projective 3-space. We deduce the cohomology of the corresponding spaces of equations and moduli spaces. In the case of plane cubics we use the method of conical resolutions. We handle plane quartics by applying a more general version of this method, which we construct in this paper. For cubic surfaces we use an explicit description of the compact GIT quotient.
\end{abstract}

\addcontentsline{toc}{section}{Introduction}
\section*{Introduction}

A section of an algebraic vector bundle $E$ over a smooth complex algebraic variety $X$ is {\it smooth} if it is transversal to the zero section. Quite a lot of information is available about the topology of the space $U$ of smooth sections. Here are a few examples.

\begin{enumerate}
\item (Divisibility) If $E$ is equivariant with respect to an affine group $G$ acting on $X$, there is a divisibility phenomenon for rational cohomology. Namely, under some conditions on $E$ and $X$ we have $H^*(U,\Q)\cong H^*(U/G,\Q)\otimes H^*(G,\Q)$, see~\cite{ste_pet_division} for the case $X=\C\p^n,G=\GL_{n+1}(\C),E=\OO(d)$ and \cite{division} for the general case.
\item (Stability) If $X=\C\p^n,G=\GL_{n+1}(\C),E=\OO(d)$, then the rational cohomology of $U$ stabilises in each degree $i$ to $H^i(\GL_{n+1}(\C),\Q)$ as $d\to\infty$, see~\cite{tom_stab}. Moreover, it seems likely that there is an analogue of this for general bundles where stabilisation is taken with respect to tensoring with powers of a very ample line bundle.
\item In several cases the rational cohomology of $U$ is explicitly known, at least additively, see~\cite{vas3}, \cite{quintics} and \cite{tom}. In particular, it is known if $X=\C\p^n,G=\GL_{n+1}(\C),E=\OO(d)$, and $(d,n)=(3,2), (4,2), (5,2)$ and $(3,3)$.
\end{enumerate}

On the other hand, as far as we are aware, not as much is known about the topology of the spaces of sections with mild singularities. In this paper we present some experimental data in this direction and compare it with what happens in the smooth case. Namely, we calculate the rational cohomology rings of the spaces of cubic and quartic plane curves and cubic surfaces in $\C\p^3$ that are {\it nodal}, i.e.\ are only allowed to have Morse type singularities (see theorems~\ref{mhcubp2}, \ref{mhquartp2} and \ref{mhcubp3}). We also calculate the answer for the corresponding moduli and equation spaces, as well as the cohomology pullback maps induced by various natural maps of the varieties involved.

For curves we use the method of conical resolutions, see section~\ref{conres}. We describe the idea in the introduction to that section, but briefly, the main step consists in calculating the $E^1$ page and all differentials of a certain spectral sequence, see theorem~\ref{spseqcvector} for plane cubics and theorem~\ref{spseqqvector} for plane quartics, with the higher differentials being the trickiest part. The construction we give in this paper generalises the one described in~\cite{vas3},~\cite{quintics} and~\cite{tom}. 

For cubic curves the new construction gives the same result as the old one, but for quartic curves it allows one to reduce the number of columns in the spectral sequence from 22 to 11. The columns we are able to drop correspond to the configuration spaces given in table~\ref{geom_table_quartics}, middle column. Other applications of the construction of section~\ref{conres}, in particular to GIT-semi-stable and stable plane quartic curves, will be presented elsewhere, the present paper being long enough as it is.

For cubic surfaces in $\C\p^3$ we use a shortcut provided by the fact that in this case there is a description of the compact GIT quotient as a weighted projective space, see \cite[Chapter 10]{dolgachev} and references therein. So the arguments we use to prove theorem~\ref{mhcubp3} are largely standard, but we think it is useful to compare the result with theorems~\ref{mhcubp2} and~\ref{mhquartp2}. Let us also note that, to the best of our knowledge, apart from the case of cubic surfaces, an explicit description of the compact GIT quotient $\mathrm{Hyp}_{d,n}//\SL_{n+1}(\C),d>1$ is only available for binary forms of small degree and for cubic curves, see ibid. Here $\mathrm{Hyp}_{d,n}=\p\big(\Gamma(\C\p^n,\OO(d))\big)$ is the space of all degree $d$ hypersurfaces of $\C\p^n$; we  equip it with the line bundle $\OO(1)$ and the standard action of $\SL_{n+1}(\C)$ on it. For plane cubic curves the quotient is $\C\p^1$, but this does not help one much to calculate the cohomology of the space of curves or the space of equations because a non-smooth nodal cubic curve is only semi-stable and not stable.

\smallskip

{\bf Remark.} In~\cite{kolya} N.~Konovalov proves a divisibility theorem (i.e.\ an analogue of item 1 above) for the space of equations $\Pi_{d,n}\setminus\N_{d,n}$ of nodal hypersurfaces of $\C\p^n$ (see next paragraph). In particular, he shows that if $d>2,n>1$ and $(d,n)\neq (3,2)$, then any orbit map from $\SL_{n+1}(\C)\to \Pi_{d,n}\setminus\N_{d,n}$ induces a surjective map in rational cohomology.

\medskip

We will now state our main results. Let $\Pi_{d,n}=\Gamma(\C\p^n,\OO(d))$ be the space of homogeneous degree $d$ complex polynomials in $n+1$ variables, and set $\N_{d,n}$ to be the subvariety of $\Pi_{d,n}$ formed by all $f$ such that the kernel of the Hessian matrix of $f$ at a non-zero point $\mathrm{x}$ contains a 2-plane $L\ni\mathrm{x}$
(or, equivalently, the projective hypersurface defined by $f$ has a singularity worse than a simple node). We set $\nod_{d,n}$ to be the projectivisation $\p \Pi_{d,n}\setminus\p\N_{d,n}$ of $\Pi_{d,n}\setminus\N_{d,n}$, i.e.\ the space of all nodal hypersurfaces. We have $\nod_{d,n}=(\Pi_{d,n}\setminus\N_{d,n})/\C^*$. Elements of $\nod_{d,n}$ will be called {\it nodal} hypersurfaces. When the categorical quotient of $\nod_{d,n}$ by $\PGL_{n+1}(\C)$ (or equivalently, of $\Pi_{d,n}\setminus\N_{d,n}$ by $\GL_{n+1}(\C)$) exists, we will denote it by $\mathcal{M}^\nod_{d,n}$ and call it the {\it moduli space} of nodal degree $d$ hypersurfaces of $\C\p^n$.

If $X$ is a complex algebraic variety, we define the {\it mixed Hodge polynomial} of $X$ to be
$$P_{\mathrm{mH}}(X)=\sum_{n,p,q}t^nu^pv^q\dim_\C\Grd^p_F\Grd_{p+q}^W(H^n(X,\Q)\otimes\C).$$
By setting in this expression $v=u$, or $u=v=1$, we get the Poincar\'e-Serre polynomial, respectively the Poincar\'e polynomial of $X$. On the other hand, if we set $t=-1$, we get the Serre characteristic of $X$, also known as the Danilov-Khovansky characteristic. In the theorems below we set $w=u^{-1}v^{-1}$ for brevity.  

We will frequently consider {\it orbit maps} $g\mapsto g\cdot x$ where $g$ is an element of a group $G$ and $x$ is a (fixed) point of a variety $X$ on which $G$ acts. Whenever it is clear from the context which $G,X$ and $x$ are meant, we will denote the orbit map by $O$. Note that if $X$ is path connected, all orbit maps are homotopic. 

\medskip

\begin{theorem}\label{mhcubp2}
We have $$P_{\mathrm{mH}}(\Pi_{3,2}\setminus\N_{3,2})=1+t^3w^2+t^5w^3+t^8w^5+t^{10}w^6+t^{11}w^7\ \ \text{and}$$
$$P_{\mathrm{mH}}(\nod_{3,2})=1+t^2w^1+t^5w^3+t^7w^4+t^{10}w^6.$$

The ring structure is as follows. For the space $\Pi_{3,2}\setminus\N_{3,2}$, the cup product $H^3\otimes H^5\to H^8$ is non zero and the cup product $H^3\otimes H^8\to H^{11}$ is zero. For the space $\nod_{3,2}$, the cup product $H^2\otimes H^5\to H^7$ is non zero.

If we denote the projection $\Pi_{3,2}\setminus\N_{3,2}\to \nod_{3,2}$ by $p$, then in rational cohomology $p^*$ is an isomorphism in degrees $0, 5$ and $10$ and $0$ in all other degrees. The projection from $\nod_{3,2}$ to the categorical quotient $\mathcal{M}^\nod_{3,2}\cong\C\p^1$ induces an injective map of the rational cohomology groups.

The orbit maps for the action of $\GL_3(\C)$ on both $\Pi_{3,2}\setminus\N_{3,2}$ and $\nod_{3,2}$ induce zero map of rational cohomology in all positive degrees.
\end{theorem}

{\bf Remark.} Note that the cohomology of both $\Pi_{3,2}\setminus\N_{3,2}$ and $\nod_{3,2}$ contains classes that look like they might pull back to some of the standard generators of $H^*(\GL_3(\C),\Q)$ under an orbit map, but in fact they don't. This can be contrasted with what happens for smooth projective hypersurfaces (see~\cite{ste_pet_division},~\cite{division}), for nodal quartic curves and cubic surfaces (see below) and more generally, for nodal hypersurfaces of degree $d>2$ of $\C\p^n,n>1$ other than plane cubics (see~\cite{kolya}).

In the next two theorems we use $\Lambda$ to denote the exterior algebra, and we let $\xi_i, i=1,3,\ldots, 2m-1$ be the degree $i$ standard multiplicative generator of $H^*(\GL_m(\C),\Q)\cong\Lambda[\xi_1,\xi_3,\ldots,\xi_{2m-1}]$.

\begin{theorem}\label{mhquartp2}
There exists a geometric quotient $\mathcal{M}^\nod_{4,2}$. We have 
$$P_{\mathrm{mH}}(\Pi_{4,2}\setminus\N_{4,2})=(1+t^3w^2)(1+t^5w^3)(1+t^7w^4),$$
$$P_{\mathrm{mH}}(\nod_{4,2})=(1+t^2w^1)(1+t^5w^3)(1+t^7w^4),\ \ \text{and}$$
$$P_{\mathrm{mH}}(\mathcal{M}^\nod_{4,2})=1+t^2w^1+t^4w^2+t^6w^3.$$

The ring structure is as follows. There is an isomorphism
$$H^*(\mathcal{M}^\nod_{4,2},\Q)\cong\Q[a]/a^4,$$ 
where $a$ has degree 2 and is Tate of weight $2$.

Furthermore, we have $$H^*(\nod_{4,2},\Q)\cong\Lambda[\zeta_5,\zeta_7]\otimes\Q[q^*(a)]/(q^*(a))^2,$$ 
where $\zeta_i,i=5,7$ has degree $i$ and is Tate of weight $i+1$, and $q$ denotes the natural map $$q:\nod_{4,2}\rightarrow \mathcal{M}^\nod_{4,2}.$$

Under the orbit map $$O:\GL_3(\C)\rightarrow \nod_{4,2}$$ these classes are transformed as follows: $O^*(\zeta_5)=\xi_5$ and $O^*(\zeta_7)=0$. 

Finally, $$H^*(\Pi_{4,2}\setminus\N_{4,2},\Q)\cong\Lambda(\zeta_3,p^*(\zeta_5),p^*(\zeta_7))$$ where $p$ denotes the natural map $$p:\Pi_{4,2}\setminus\N_{4,2}\rightarrow \nod_{4,2}.$$ Furthermore, we have $O^*(\zeta_3)=\xi_3$.
 
\end{theorem}

\begin{theorem}\label{mhcubp3}
There exists a geometric quotient $\mathcal{M}^\nod_{3,3}$. We have 
$$P_{\mathrm{mH}}(\Pi_{3,3}\setminus\N_{3,3})=(1 + t^3 w^2) (1 + t^5 w^3) (1 + t^7 w^4)^2,$$
$$P_{\mathrm{mH}}(\nod_{3,3})	=
(1 + t^2 w^1) (1 + t^5 w^3) (1 + t^7 w^4)^2,\ \ \text{and}$$
$$P_{\mathrm{mH}}(\mathcal{M}^\nod_{3,3})=1+t^2w^1+t^4w^2+t^6w^3.$$

Moreover, $\mathcal{M}^\nod_{3,3}\cong \C\p(1,2,3,4,5)\setminus\{(0:0:0:0:1)\}$, so we have a ring isomorphism $$H^*(\mathcal{M}^\nod_{3,3},\Q)\cong\Q[a]/a^4,$$ 
where $a\in H^2(\mathcal{M}^\nod_{3,3},\Q)\cong H^2(\C\p(1,2,3,4),\Q)$ is Tate of weight $2$.

Furthermore, $$H^*(\nod_{3,3},\Q)\cong\Lambda[\zeta_5,\zeta_7,b]\otimes\Q[q^*(a)]/(q^*(a))^2,$$
where $\zeta_i,i=5,7$ lives in degree $i$ and is Tate of weight $i+1$, $b$ has degree 7 and is Tate of weight $8$, and $q$ denotes the natural map $$q:\nod_{3,3}\rightarrow \mathcal{M}^\nod_{3,3}.$$

Under the orbit map $$O:\GL_4(\C)\rightarrow \nod_{3,3}$$ these classes are transformed as follows: $O^*(\zeta_n)=\xi_n$ and $O^*(b)=0$. 

Finally, $$H^*(\Pi_{3,3}\setminus\N_{3,3},\Q)\cong\Lambda(\zeta_3,p^*(\zeta_5),p^*(\zeta_7),p^*(b)),$$
where $p$ denotes the natural map $$p:\Pi_{3,3}\setminus\N_{3,3}\rightarrow \nod_{3,3}.$$ Furthermore, we have $O^*(\zeta_3)=\xi_3$.
\end{theorem}

%
%

{\bf Remark.} The rational cohomology of all spaces involved in theorems~\ref{mhcubp2}-\ref{mhcubp3} is pure and Tate in each degree.

{\bf Remark.} The moduli space $\mathcal{M}^\nod_{3,3}=(\Pi_{3,3}\setminus\N_{3,3})/\GL_4(\C)$ of nodal cubic surfaces has an interpretation in terms of complex hyperbolic geometry, see~\cite{act}. More specifically, this space is the quotient of the complex hyperbolic 4-space by a certain discrete group $\Gamma$ acting with finite stabilisers. This means that $\mathcal{M}^\nod_{3,3}$ is a ``rational $K(\Gamma,1)$'', i.e. the rational cohomology of $\Gamma$ is isomorphic to the rational cohomology of $\mathcal{M}^\nod_{3,3}$.

\medskip

The paper is organised as follows: in section~\ref{conres} we construct conical resolutions of discriminants and prove their main properties (see theorems~\ref{thconical}, \ref{thconical_1}, \ref{thconical_2}, \ref{thhodge1}, \ref{thhodge2}, \ref{thhodge3}). We describe the motivation behind the construction and the contents of section~\ref{conres} in more detail in the introduction to this section. Section~\ref{nod} contains preliminary material on vector bundles and nodal projective hypersurfaces. In particular, we calculate the Borel-Moore homology map induced by the zero section embedding (proposition~\ref{thom_iso_zero_sec_prop}) and describe the transgression in the cohomology Leray spectral sequence of the spherisation of a vector bundle (proposition~\ref{euler-thom}). We will need these results later to calculate the differentials in the spectral sequences constructed using the natural filtrations on conical resolutions. In section~\ref{plane_cubics} we apply the method of section~\ref{conres} to the space of equations of plane nodal cubics and prove theorem~\ref{mhcubp2}. The proof of theorem~\ref{mhquartp2} is given in section~\ref{quartics}. Finally, in section~\ref{cubics} we prove theorem~\ref{mhcubp3}.

The paper has two appendices in which we summarise some background material and folklore. The main goal is to state the results we need in the main text, and to provide references when they are available and proofs otherwise. These results will hopefully also be useful in other examples similar to those that we consider in this paper. Specifically, in appendix~\ref{b} we consider group actions, quotients, slices and equivariant cohomology. For example, in theorem~\ref{groupactions} we give a sufficient condition for the direct images of a locally constant sheaf under the map into the quotient to be locally constant. We also equip equivariant cohomology with a mixed Hodge structure and show that the result is compatible with several natural maps (proposition~\ref{prop_equi}).

In appendix~\ref{app_mhs_spectral} we review Borel-Moore homology of semi-simplicial, in particular cubical varieties (section~\ref{borel_moore}), and introduce mixed Hodge structures on certain spectral sequences (propositions~\ref{mhs_sseq_filtred} and~\ref{mhs_sseq_leray}) that we use in the main text. We also show that for (semi-)simplicial complex algebraic varieties the K\"unneth decomposition and the relative cup product map are compatible with mixed Hodge structures (theorems~\ref{simplicial_kunneth_mhs} and~\ref{cup_product_mhs}), and prove the Thom isomorphism theorem for cubical affine bundles (theorem~\ref{thom_ccavb}). More details on the contents of appendices~\ref{b} and~\ref{app_mhs_spectral} can be found in the respective introductions.

This paper incorporates the 4th year undergraduate research project of the first named author written under supervision of the second named author at the Faculty of mathematics, HSE, Moscow, Russia.

\medskip

{\bf Acknowledgements:} We wish to thank Nikolay (Kolya) Konovalov and Alexander (Sasha) Kuznetsov for helpful conversations and, in the former case, for co-authoring appendix~\ref{app_mhs_spectral}.

\medskip

In the following table we summarise some of the notation used in the paper.
\pagebreak[3]

\begin{longtable}{|c|c|}
\hline
\rule{0pt}{11pt}$S_k$ & the symmetric group on $k$ letters\\
\hline
\rule{0pt}{11pt}$\Lambda[\mbox{odd variables}]$& exterior algebra generated by the variables\\
\hline
\rule{0pt}{11pt}$\p$& \begin{tabular}{c}projectivisation (of a vector bundle or\\ a $\C^*$-invariant subvariety of $\C^m$)\end{tabular}\\
\hline
\rule{0pt}{11pt}$\Q(n)$ & \begin{tabular}{c}Tate mixed Hodge structure on $\Q$
(the weight filtration on $\Q$ is\\ given by $W_{-2n-1}=0,W_{-2n}=\Q$, and the Hodge filtration on\\ $\C\cong\C\otimes_\Q\Q$ is given by $F^{-n}=\C,F^{-n+1}=0$)\end{tabular}\\
\hline
\rule{0pt}{11pt}$\bar H_*(X,\eu L)$&\begin{tabular}{c} Borel-Moore homology groups of a locally compact space $X$\\ with coefficients in a local system $\eu L$\end{tabular}\\
\hline
\rule{0pt}{11pt}$|K|$& the cardinality of a finite set $K$\\
\hline
\rule{0pt}{11pt}$\Pi_{d,n}$ & \begin{tabular}{c}$\Gamma(\C\p^n,\mathcal{O}(d))$, the space of complex homogeneous polynomials \\ of degree $d$ in $n+1$ variables\end{tabular}\\
\hline
\rule{0pt}{11pt}$\N_{d,n}$ & \begin{tabular}{c} the subvariety of $\Pi_{d,n}$ formed by all $f$ such that the kernel of the \\Hessian matrix of $f$ at a non-zero
point $\mathrm{x}$ contains a 2-plane $L\ni\mathrm{x}$\end{tabular}\\
\hline
\rule{0pt}{11pt}$\text{Nod}_{d,n}$ & \begin{tabular}{c} The space $\mathrm{P} \Pi_{d,n}\setminus \mathrm{P} \N_{d,n}$ of nodal hypersurfaces of degree $d$ in $\C \p^n$ \end{tabular}\\
\hline
\rule{0pt}{11pt}$\C \p^{n\vee}$ & the space of all lines in $\C \p^n$\\
\hline
\rule{0pt}{11pt}$\Gr(d,V)$ & the Grassmannian of $d$-planes in a vector space $V$\\
\hline
\rule{0pt}{11pt}$B(X,k)$ & the space of unordered configurations of $k$ (distinct) points in $X$\\
\hline
\rule{0pt}{11pt}$F(X,k)$ & the space of ordered configurations of $k$ (distinct) points in $X$\\
\hline
\rule{0pt}{11pt}$c(\xi)$ & the total Chern class of a vector bundle $\xi$\\
\hline
\rule{0pt}{11pt}$\xi_x$ & the fibre of a vector bundle $\xi$ over $x$\\
\hline
\rule{0pt}{11pt}$\tot(\xi)$ & the total space of a vector bundle $\xi$\\
\hline
\rule{0pt}{11pt}$Cone(X)$ & the cone over a topological space $X$\\
\hline
\rule{0pt}{11pt}$Cone^\circ(X)$ & the open cone over a topological space $X$, i.e.\ $Cone(X)\setminus X$\\
\hline
\rule{0pt}{11pt}$\Lambda(K)$& see p. \pageref{xlambda}\\
\hline
\rule{0pt}{11pt}$\tilde\Lambda(K)$& see p. \pageref{tlambda}\\
\hline
\rule{0pt}{11pt}$\partial\Lambda(K)$& see p. \pageref{plambda}\\
\hline
\rule{0pt}{11pt}$\tilde\partial\Lambda(K)$& see p. \pageref{tplambda}\\
\hline
\end{longtable}

By a {\it complex algebraic variety} we will mean a separated reduced (but not necessarily irreducible) scheme of finite type over $\C$.

\tableofcontents

\section{Conical resolutions}\label{conres}

In this section we describe a generalised version of the conical resolution from \cite{quintics,tom}. Before we do that we briefly explain what this construction does and why we generalise it. Suppose we have a finite-dimensional vector space $V$ over a field $\mathbf{k}$, and elements of $V$ parametrise geometric objects of some kind. These objects are often subdivided into ``generic'' and ``singular'' (what exactly we call ``singular'' depends on the problem we consider, see section~\ref{sing_loc}). The subset $\Sigma \subset V$ formed by all singular elements is called a (generalised) {\it discriminant}. To give an example, one can take $V$ to be the space of polynomials $\in\mathbf{k}[x]$ of some fixed degree $d-1$, and~$\Sigma$ equal the subset of $V$ formed by the polynomials $f$ such that $x^d+f$ has a multiple root in the algebraic closure of $\mathbf{k}$. Note that in this example $\Sigma$ is the zero locus of a single polynomial in the coefficients of $f\in V$; this polynomial is the classical discriminant, which explains the name.

Now suppose $\mathbf{k}=\R$ or $\C$ and $\Sigma\subset V$ is closed and sufficiently nice, e.g.\ a sub-polyhedron of some triangulation of $V$. By Alexander duality we have then an additive isomorphism
\begin{equation}\label{alexander}
\tilde H^*(V\setminus\Sigma)\cong\bar H_{\dim_\R (V)-*-1}(\Sigma),
\end{equation}
where on the right we have the Borel-Moore homology groups of $\Sigma$. (Recall that Borel-Moore homology groups $\bar H_*(X)$ of a locally compact and locally contractible topological space $X$ can be defined as the homology groups of the complex of locally finite singular chains; if $X=Y\setminus Z$ where $(Y,Z)$ is a finite CW-pair, then we have $\bar H_*(X)\cong H_*(Y,Z)$; for generalisations and more details see section~\ref{sec_bm_hom}.)

The space $\Sigma$ is more tractable than $V\setminus\Sigma$ but still quite complicated in all interesting examples. We construct below (see p.~\pageref{def_con_res}) another space $\sigma$, which we call a {\it resolution} of $\Sigma$, and a proper map $\sigma\to \Sigma$ that induces an isomorphism of the Borel-Moore homology groups. The resolution $\sigma$ comes equipped with a natural filtration
\begin{equation}\label{filtration}
\varnothing\subset\sigma_1\subset\sigma_2\subset\cdots\subset \sigma_N=\sigma
\end{equation}
such that the difference of any two consecutive terms is a fibre bundle over a configuration space of some kind. We use these fibre bundles to calculate all $\bar H_*(\sigma_i\setminus\sigma_{i-1})$, and then combine the results to get the page $E^1$ of the spectral sequence that corresponds to (\ref{filtration}).

Both $\sigma$ and the filtration on it were constructed in \cite{quintics,tom} starting from a family $X_1,\ldots X_N$ of configuration spaces that satisfy certain conditions, see list \ref{conditions}. One of these, namely condition \ref{five}$+$, requires that every space $\bigcup_{j\leq i} X_i$ must be closed with respect to the Hausdorff metric. As a result, the filtrations typically contain a lot of ``spurious'' terms that correspond to subvarieties which are limits of singular loci but are not singular loci themselves. This does not present a problem in~\cite{quintics,tom} as most of these terms do not contribute over $\Q$ to the spectral sequence.

Nevertheless, in other examples, and in particular for nodal quartic curves, the situation is different. While it is still true that most columns of the $E^\infty$ page of the spectral sequence are zero, many of the corresponding columns of $E^1$ are not, even rationally. To remedy this we replace condition \ref{five}$+$ with weaker conditions \ref{five} and \ref{five}$-$. As a result, the list of configuration spaces that need to be taken into account becomes shorter. For example, in the case of plane nodal quartics we would need to calculate the cohomology of all configuration spaces listed in the middle column of table~\ref{geom_table_quartics} if we applied the construction of~\cite{quintics, tom} directly.
%
Informally, we single out the configuration spaces that really matter and gather together spaces whose contributions to the spectral sequence are destined to kill each other in the end. 
%

We note that in the case of discrete singular loci resolutions of the kind we consider in this paper were first constructed by V.~Vassiliev in \cite{vas0}, and in fact, one can trace the idea all the way back to Euler's inclusion-exclusion formula.

\medskip

We will now describe the contents of section~\ref{conres}. In~\ref{sing_loc} we review the Hausdorff metric on the space $2^X$ of all compact non-empty subsets of a compact metric space $X$, and summarise its properties (proposition~\ref{conf_space}). Then, in section~\ref{sect_chow_hilb}, for $X$ complex projective we relate $2^X$ to its counterparts from algebraic geometry. In particular, we show that the natural cycle maps from the Chow variety and Hilbert scheme to $2^{\C\p^n}$ are continuous (propositions~\ref{cycle_maps} and~\ref{cycle_maps_hilbert}). We also prove a technical statement (corollary~\ref{Bkl}) about the closure in the Hausdorff metric of the image of the Chow variety under the cycle map.

In section~\ref{sect_fam_conf} we introduce conditions on a family of configuration spaces (condition list~\ref{conditions}) which we will use later to construct resolutions of discriminants. Conditions~\ref{first}-\ref{xxx} are the same as in~\cite{quintics} except that what was condition~\ref{five} there has now become~\ref{five}$+$; conditions~\ref{five} and~\ref{five}$-$ are weaker forms of condition~\ref{five}$+$. We prove a statement about subdividing the subset of $2^{\C\p^n}$ which consists of all subvarieties of not necessarily pure dimension $k$ with all components of degree $\leq l$ (proposition~\ref{chow}). This statement is used to check that a family of configuration spaces that satisfies~\ref{first}-\ref{five}$-$ can be completed to a family that satisfies~\ref{first}-\ref{five}$+$, at least in the projective case (lemma~\ref{auxconf}).

In section~\ref{conical_sect} we review joins and self-joins (also known as symmetric joins), construct the generalised conical resolution and prove its main properties (theorems~\ref{thconical}, \ref{thconical_1} and~\ref{thconical_2}).

If $\mathbf{k}=\C$ and $\Sigma$ is affine algebraic, then both sides of (\ref{alexander}) are mixed Hodge structures, which are preserved under the isomorphism. It is reasonable to expect that in this case the spectral sequence of~(\ref{filtration}) and all auxiliary spectral sequences also have natural mixed Hodge structures. In section~\ref{hodge_modge} we show that this is indeed the case under some assumptions which are given in condition list~\ref{condition_list_3}. More precisely, we show that if these assumptions hold, then the cubical spaces introduced by O.~Tommasi in~\cite{tom}, as well as the analogues we construct in this paper, can be blown up to cubical algebraic varieties. As a result, the spectral sequence of (\ref{filtration}) is isomorphic to the spectral sequence of a cubical algebraic variety filtered by cubical subvarieties (theorem~\ref{thhodge1}). We handle the auxiliary spectral sequences in theorems~\ref{thhodge2} and~\ref{thhodge3}. We also observe that the conical resolutions constructed in~\cite{quintics} and~\cite{tom} are in fact filtered homeomorphic (proposition~\ref{homeo_quint_tom}).

The main difficulty in section~\ref{hodge_modge} is that although the configuration spaces in this case have a natural complex algebraic structure, their closures in the space of all compact subsets sometimes don't; we give an example in a remark on p.~\pageref{nocomplexstr}. 
The main technical result is lemma~\ref{fib_contractible}.

\smallskip

{\bf Remark.} See~\cite{vas3} for a somewhat different approach to the conical resolution. Also, in~\cite{tom_stab} O.\ Tommasi constructs a ``partial'' conical resolution along the same lines. We will now briefly compare this approach with ours. In both~\cite{vas3} and~\cite{tom_stab} the starting point is a family of algebraic subvarieties of the Grassmannians $G_k(V)$. Each of these subvarieties contains an open part, which is typically isomorphic to some configuration space and is easy to understand, and the boundary, which is not. As a result, one obtains a mixed Hodge structure on the main spectral sequence for free, but the page $E^1$ is in general hard to calculate, and the higher differentials presumably even harder. So both constructions involve a second step, which simplifies the resolution by collapsing certain simplices.

In this paper, as in~\cite{quintics} and~\cite{tom}, we proceed in the opposite direction: We start with something which is already relatively simple, namely a family of configuration spaces, i.e.\ subspaces of $2^X$ where $X$ is compact complex algebraic (see section~\ref{sect_fam_conf}). All closures are taken in $2^X$, and not in a Grassmannian or its variants such as the Chow variety or the Hilbert scheme, which makes the boundary more tractable. Then we verify that the resulting spectral sequences have mixed Hogde structures using some axioms (see condition list~\ref{condition_list_3}), which we check in the projective case using Chow varieties (proposition~\ref{realistic}). It takes some work to prove that mixed Hodge structures exist and have the expected compatibility properties, but once we know this, calculating the mixed Hodge structures becomes straightforward.

It is possible that the construction given in~\cite{quintics} and~\cite{tom} and those described in~\cite{vas3} and~\cite{tom_stab} are in fact proper filtered homotopy equivalent under some assumptions; we do not know if this is the case.  
Either way, configuration spaces of algebraic origin seem to contain just the right amount of information for our purposes: enough to determine the rational cohomology of the discriminant complement, but not enough to make the task unfeasible.

{\bf Remark.} The idea that some contributions to the main spectral sequence can be eliminated in a systematic way was inspired by V.~Vassiliev's construction of reduced conical resolutions \cite[section 3.2]{vas3}.

\subsection{Configuration spaces}\label{conf_spaces_sect}
\subsubsection{Singular loci and the Hausdorff metric}\label{sing_loc}
To begin with, we need to specify what we mean by a singularity of a given type. Let $\kk$ and $V$ be as above, and let $\mathbf{M}$ be a finite CW-complex. Suppose we have a $\kk$-vector bundle $\mathcal{V}$ over $\mathbf{M}$ and constant rank map of vector bundles $ev:V\times\mathbf{M}\to\mathcal{V}$. For $K\subset\mathbf{M}$ we set $L(K)=\bigcap_{x\in K} \ker ev_x$ where $ev_x,x\in\mathbf{M}$ is the restriction of $ev$ to $V\times\{x\}$. If $f\in V$, we set $$\Sing(f)=\{x\in \mathbf{M}\mid f\in\ker ev_x\}.$$ Note that $K\subset\Sing f$ if and only if $f\in L(K)$. We define $$\Sigma=\{f\in V\mid\Sing(f)\neq\varnothing\}.$$

{\bf Example.} Suppose $\mathbf{M}$ is a smooth compact complex analytic manifold and $L$ is a complex analytic line bundle over $\mathbf{M}$. We take $\mathcal{V}$ to be the first jet bundle $J(L)$ of $L$ and let $V\subset\Gamma(X,L)$ be a vector subspace such that for all $x\in X$ the first jets of $f,f\in V$ at $x$ span the fibre $\mathcal{V}_x$. We set $ev:V\times\mathbf{M}\to\mathcal{V}=J(L)$ to be the jet evaluation map. Then $\Sing(f),f\in V$ is the usual singular locus of the zero locus $Z(f)$ of $f$, and $\Sigma$ is the set of all $f\in V$ such that the hypersurface $Z(f)\subset\mathbf{M}$ is singular.

In the particular case $\mathbf{M}=\C\p^n,L=\OO(d),d>0$, and $V=\Gamma(\mathbf{M},L)=\Pi_{d,n}$ the jet bundle $J(L)$ is (non-equivariantly) isomorphic to a direct sum of $n+1$ copies of $\OO(d-1)$, and the map $ev$ takes a couple $(f,x)\in\Pi_{d,n}\times\C\p^n$ to $$\left(\frac{\partial f}{\partial\mathrm{x}_0}(x),\frac{\partial f}{\partial\mathrm{x}_1}(x),\ldots, \frac{\partial f}{\partial\mathrm{x}_n}(x)\right).$$

{\bf Remark.} Most constructions below can be modified to work without assuming that $ev$ should have constant rank or that $\mathbf{M}$ should be finite. We postpone this until future work. 

The following lemma is straightforward.
%
\begin{lemma}
\begin{enumerate}\label{firstlist}
\item If $f_1,f_2\in\Sigma$, and $\Sing(f_1)\cap \Sing(f_2)\neq\varnothing$, then $f_1+f_2\in\Sigma$ and
$\Sing(f_1)\cap \Sing(f_2)\subset\Sing(f_1+f_2)$.

\item If $f\in\Sigma$, then for any $\lambda\neq 0$ we have $\Sing(\lambda f)=\Sing(f)$.

\item The zero element $0\in V$ belongs to $\Sigma$, and $\Sing(0)=\mathbf{M}$.
\item For every $K\subset \mathbf{M}$ we have $L(K)=\{f\in V\mid K\subset\Sing(f)\}$.
\item\label{cont_grassm} Let $\mathbf{d}=\dim_\kk\ker ev_x,x\in\mathbf{M}$. The map $g:\mathbf{M}\to \Gr(\mathbf{d},V)$ given by $g(x)=L(\{x\})$ is continuous, and the vector bundle $\ker ev$ is isomorphic to the pullback of the tautological bundle on the Grassmannian $\Gr(\mathbf{d},V)$ under $g$.
\item\label{Sigma_is_closed} The subspace $\Sigma\subset V$ is closed.
\end{enumerate}
\end{lemma}

$\clubsuit$
%

In all our applications $\mathbf{M}$ and all $\Sing(f)$ will be projective varieties. We will call $\Sing(f)$ the {\it (generalised) singular locus} of $f$.

\smallskip

Now we will review a few facts about configuration spaces. Let $X$ be a topological space. For a positive integer $k$ we define the {\it ordered configuration space} $F(X,k)$ to be
$$\{(x_1,\ldots, x_k)\in X^k \mid x_i\neq x_j\mbox{ {\normalfont for }} i\neq j\},$$ and we let the {\it unordered configuration space} $B(X,k)$ be $F(X,k)/S_k$. If $G$ is an Abelian group, the {\it alternating local system} $\pm G$ on $B(X,k)$ has fibre $G$ and corresponds to the alternating representation of $\pi_1(B(X,k))$: $$\pi_1(B(X,k))\to S_k\to\{\pm 1\}.$$ In most cases below $G$ will be $\Z$ or $\Q$.

Suppose now the topology on $X$ is induced by a metric $\rho$. We will call compact subsets of $X$ {\it configurations}, and we denote the set of all {\it non-empty} configurations by $2^X$. We equip $2^X$ with the {\it Hausdorff metric} $\rho_H$ given by
$$\rho_{H}(K_1,K_2)=\mathop{\mathrm{max}}(\mathop{\mathrm{max}}_{x\in K_1}\rho(x,K_2),\mathop{\mathrm{max}}_{x\in K_2}\rho(x,K_1)).$$ The basic properties of this metric can be summarised as follows:
\begin{Prop}\label{conf_space}
\begin{enumerate}
\item The metric space $(2^X,\rho_{H})$ is complete if $\rho$ is complete and compact if $X$ is compact.
\item Let $(K_i)$ be a convergent sequence in $(2^X,\rho_{H})$. Then $x\in\lim_{i\to\infty}K_i$, if and only if there exists a sequence $(x_i)$ that converges to $x$ and such that every $x_i\in K_i$.
\item Let $(K'_i)$ and $(K_i'')$ be two convergent sequences in $(2^X,\rho_{H})$, and let $K'$ and $K''$ be the limits. Then the following is true.
\begin{itemize}
\item If $K'_i\subset K''_i$ for every $i$, then $K'\subset K''$.
\item We have $\lim_{i\to\infty} (K'_i\cup K''_i)=K'\cup K''$.
\item If the limit $\lim_{i\to\infty} (K'_i\cap K''_i)$ exists, then it is $\subset K'\cap K''$.
\end{itemize}
\item Set $U=\{(x,K)\in X\times 2^X\mid x\in K\}$. Then the projection $p_2:U\to 2^X$ is open and proper.
\item A continuous map $(X,\rho)\to (X',\rho')$ of metric spaces induces a continuous map $2^X\to 2^{X'}$. So $(X,\rho)\mapsto (2^X,\rho_H)$ defines an endofunctor of the category of metric spaces and continuous maps.
\end{enumerate}
\end{Prop}

{\bf Remark.} Part 5 implies in particular that the topology on $2^X$ depends only on the topology on $X$, not on the choice of the metric. In fact, there is a topology on the set of all nonempty compact subspaces  of an arbitrary space $Y$, called the {\it Vietoris topology}, which depends only on the topology of $Y$ and coincides with the topology one gets from the Hausdorff metric in the metric case, see e.g.~\cite[Proposition 2.4.14]{srivastava}. We do not use this construction here because metric spaces suffice for all examples that we have in mind in this paper.

{\bf Proof.} Part 1 follows from e.g.~\cite[Propositions 2.4.15 and 2.4.17]{srivastava}. Parts 2 and 3 are straightforward. 
Part 4 is Theorem 1 in \cite{gillam_karan}. 
Let us prove part 5. Let $(K_i)$ be a convergent sequence in $2^X$, and let $K$ be the limit. We need to show that $\lim_{i\to\infty} f(K_i)=f(K)$. Take an $\varepsilon>0$. Let $V'\subset X'$ be the $\varepsilon$-neighbourhood of $f(K)$. The preimage $f^{-1}(V')$ contains an $\varepsilon'$-neighbourhood $V$ of $K$ for some $\varepsilon'>0$, and $V$ contains all $K_i$ for $i$ large enough. So for these values of $i$ we have $\max_{x'\in f(K_i)}\rho'(x',f(K))<\varepsilon$.

Now let $N\subset K$ be a finite subset such that $f(N)\subset f(K)$ is a finite $\frac{\varepsilon}{2}$-net. For every $x\in N$ choose a sequence $(y_i^x)$ such that $y^x_i\in K_i$ for all $i$ and $y_i^x\to x$ as $i\to\infty$. For $i$ large enough we will have $\rho'(f(y_i^x),f(x))<\frac{\varepsilon}{2}$ for all $x\in N$. This shows that, again for $i$ large enough, we have $\max_{x'\in f(K)}\rho'(x',f(K_i))<\varepsilon$. Combining this with the above we obtain part 5 of the proposition.$\clubsuit$

\begin{Prop}\label{bk}
The space $B(X,k)$ is homeomorphic to $\{K\in 2^X\mid |K|=k\}$. Moreover, we have $\overline{B(X,k)} = \bigcup_{j\le k}B(X,j)$.
\end{Prop}

$\clubsuit$

{\bf Remark.} Note however that $\bigcup_{k\geq 1} B(X,k)$ is dense in $2^X$.

\smallskip

Now we take $X$ to be the CW-complex $\mathbf{M}$ as above and let $\rho$ be any metric that induces the topology of $\mathbf{M}$. Using the part~\ref{cont_grassm} of lemma~\ref{firstlist} we can now deduce the following proposition.

\begin{Prop}\label{semicont}
\begin{enumerate}
\item Suppose $(K_i)$ is a convergent sequence of elements of $2^{\mathbf{M}}$ such that $\dim_{\mathbf{k}} L(K_i)$ is constant; denote this dimension $d$ and let $K=\lim_{i\to\infty} K_i$. If $\dim L(K)=d$, then $\lim_{i\to\infty} L(K_i)=L(K)$.
\item The map $K\mapsto \dim_\mathbf{k}L(K)$ from $2^{\mathbf{M}}$ to $\Z$ is upper semicontinuous, i.e., for all $l\in\Z$ the set of all $K\in 2^{\mathbf{M}}$ such that $\dim_\mathbf{k}L(K)\geq l$ is closed in $2^{\mathbf{M}}$.
\end{enumerate}
\end{Prop}
{\bf Proof.} Let us prove the first part of the proposition. Let $\{y_1,\ldots, y_m\}$ be a finite subset of $K$ such that $L(K)=\bigcap_{j=1}^m L(\{y_j\})$. By proposition \ref{conf_space} there are sequences $(y_{j,i})_{i\in\Z_{>0}},j=1,\ldots,m$ such that every $y_{j,i}\in K_i$ and $\lim_{i\to\infty} y_{j,i}=y_j$ for all $j=1,\ldots, m$. Then
$L(K_i)\subset \bigcap_{j=1}^m L(\{y_{j,i}\})$ for every $i$. 
Applying pat 3 of proposition~\ref{conf_space} and part~\ref{cont_grassm} of lemma~\ref{firstlist} we see that for every convergent subsequence of $(L(K_i))$ the limit will be a subspace of $\bigcap_{j=1}^m L(\{y_j\})=L(K)$. (In more detail, in order to apply proposition~\ref{conf_space} to vector subspaces of $V$ we replace them with their intersections with the unit sphere in some Euclidean metric on $V$.)
Since $\dim_\kk L(K_i)=\dim_\kk L(K)$ for all $i$, we conclude that
$\lim_{l\to\infty}L(K_i)=L(K)$.

We now prove the second part. Suppose $(K_i)$ is a sequence in $2^{\mathbf{M}}$ that converges to $K\in 2^{\mathbf{M}}$, and let $W_i\subset L(K_i)$ be a vector subspace of dimension $l$. By choosing a subsequence if necessary we may assume that $W_i\to W\in \Gr(l,V)$. Let us prove that $W\subset L(K)$. Take an $f\in W$ and $x\in K$. We will now show that $f\in L(\{x\})$. Choose a sequence $(f_i)$ of elements of $V$ such that $f_i\to f$ and every $f_i\in W_i$, and a sequence $(x_i)$ of elements of $\mathbf{M}$ such that $x_i\to x$ and every $x_i\in K_i$; the latter sequence exists by proposition~\ref{conf_space}. We have then $f_i\in W_i\subset L(K_i)\subset L(\{x_i\})$ for all $i$, which by part~\ref{cont_grassm} of lemma~\ref{firstlist} implies $f\in L(\{x\})$.$\clubsuit$

%
%

\subsubsection{Chow varieties and Hilbert schemes}\label{sect_chow_hilb}

Let us now see how configuration spaces are related to the Chow varieties and Hilbert schemes. 

Recall that the {\it (closed) Chow variety} $\bar {\eu C}_{k,l,n}$ is a complex projective variety such that its points naturally correspond to formal linear combinations $\sum_{i=1}^k l_i V_i$ where $V_1,\ldots, V_k$ are pairwise distinct irreducible subvarieties of $\C \p^n$ of dimension~$k$, all $l_i\in\Z_{> 0}$, and $\sum l_i\deg V_i=l$, see e.g.~\cite[Chapter 4]{gkz}, \cite[Chapter 1]{kollar} or \cite[Lecture 21]{harris}. There is an open dense subvariety $\eu{C}_{k,l,n}\subset\bar{\eu{C}}_{k,l,n}$, called the {\it (open) Chow variety}, whose closed points correspond to $\sum l_i V_i$ as above with all $l_i=1$, or in other words, to maybe reducible subvarieties of $\C \p^n$ of pure dimension $k$ and degree~$l$.

Let $X\subset\C\p^n$ be a closed subvariety. We then have closed subvarieties ${\eu C}_{k,l,n}^X\subset{\eu C}_{k,l,n}$ obtained by requiring that all $V_i$ in the definition above should be subvarieties of $X$. 
We set $\bar{\eu C}_{k,l,n}^X$ to be the closure of ${\eu C}_{k,l,n}^X$ in $\bar{\eu C}_{k,l,n}$. There exists a closed subvariety $\bar{\eu{U}}_{k,l,n}\subset\C\p^n\times\bar{\eu{C}}_{k,l,n}$, called the {\it universal family} over $\bar{\eu{C}}_{k,l,n}$, such that the fibre of $\bar{\eu{U}}_{k,l,n}$ over $x=\sum l_i V_i\in \bar{\eu{C}}_{k,l,n}$ is $\bigcup_i V_i$. We let $\eu{U}_{k,l,n},\eu{U}_{k,l,n}^X$, respectively $\bar{\eu{U}}_{k,l,n}^X$ be the pullback of $\bar{\eu{U}}_{k,l,n}$ to $\eu{C}_{k,l,n}, \eu{C}_{k,l,n}^X$ and $\bar{\eu{C}}_{k,l,n}^X$ respectively.

Let $P\in\Q[x]$ be a polynomial. 
There exists a complex projective scheme $\eu{H}_{P,n}$, called the {\it Hilbert scheme} of $\C\p^n$, the closed points of which correspond to closed subschemes of $\C \p^n$ with Hilbert polynomial $P$. Moreover, there exists a closed subscheme $\eu{U}_{P,n}\subset \C\p^n\times\eu{H}_{P,n}$, called the {\it universal family} over $\eu{H}_{P,n}$, flat over $\eu{H}_{P,n}$ and such that the fibre of $\eu{U}_{P,n}$ over a closed point $x\in \eu{H}_{P,n}$ is the subscheme of $\C\p^n$ that corresponds to $x$. For details see~\cite[expos\'e 221]{fga} and~\cite[Chapter 1]{kollar}. More generally, if $X\subset\C\p^n$ is a closed subscheme, there is a closed subscheme $\eu{H}_{P,n}^X\subset \eu{H}_{P,n}$ such that the closed points of $\eu{H}_{P,n}^X$ correspond to closed subschemes of $X$ with Hilbert polynomial $P$. 
We set $\eu{U}_{P,n}^X=\eu{H}_{P,n}^X\times_{\eu{H}_{P,n}}\eu{U}_{P,n}$.

The points of the Chow variety correspond to simpler geometric objects than those of the Hilbert scheme. On the other hand, the Hilbert scheme has the advantage that the universal family over it is flat, while the universal family over the Chow variety in general is not.

\smallskip

If $X\subset\C\p^n$ is a closed subvariety, there is a {\it cycle map} $Z^{C,X}_{k,l,n}:\bar{\eu C}_{k,l,n}^X\to 2^X$ which takes a point $\sum l_i V_i$ of $\bar{\eu C}_{k,l,n}^X$ to $\bigcup V_i$. Similarly, if $X\subset\C\p^n$ is a closed subscheme, there is a {\it cycle map} $Z^{H,X}_{P,n}:\eu{H}^X_{P,n}(\C)\to 2^X$ which takes a closed subscheme $Y\subset X$ to the set of the complex points of the underlying topological space. We will often write $Z^C$ and $Z^H$ instead of $Z^{C,X}_{k,l,n}$, respectively $Z^{H,X}_{P,n}$.

\begin{Prop}\label{cycle_maps}
For every closed subvariety $X\subset\C\p^n$ the cycle map $Z^C:\bar{\eu{C}}_{k,l,n}^X\to 2^X$ is continuous with respect to the complex analytic topology on the source and the topology on the target induced by the Hausdorff metric (obtained from any metric on $X$ which induces its complex analytic topology).
\end{Prop}

Before proving this proposition we state an analogue for Hilbert schemes.

\begin{Prop}\label{cycle_maps_hilbert}
For a closed subscheme $X\subset\C\p^n$, the cycle map $Z^H:\eu{H}^X_{P,n}(\C)\to 2^{X(\C)}$ is continuous with respect to the complex analytic topology on the source and the topology on the target induced by the Hausdorff metric obtained from any metric on $X(\C)$ which induces its complex analytic topology.
\end{Prop}

{\bf Proof.} See~\cite{gillam_karan}, in particular the map given by formula (1) there.$\clubsuit$

In order to prove proposition~\ref{cycle_maps} we need a cycle map
from the Hilbert scheme to the Chow variety. Using~\cite[Theorem I.6.3]{kollar} we get a map
\begin{equation}\label{cycle_map_hilb_chow}
Z^{C\gets H}:\bigsqcup_{P\in\Q[x]} \eu{H}^{sn}_{P,n}(\C)\to\bigsqcup_{k,l}\bar{\eu{C}}_{k,l,n},
\end{equation}
Here $sn$ stands for semi-normalisation, see e.g.~\cite[I.7.2]{kollar} or~\cite[tag 0EUK]{stacks}. So the Hilbert to Chow cycle map is defined on the semi-normalisation of the Hilbert scheme rather than on the Hilbert scheme itself, but this is not really a problem:

\begin{lemma}
For every scheme $Y$, the semi-normalisation $Y^{sn}$ is reduced and the natural morphism $sn:Y^{sn}\to Y$ is a universal homeomorphism. If $Y$ is separated of finite type over $\C$, then so is $Y^{sn}$, and $sn$ induces a homeomorphism $sn(\C):Y^{sn}(\C)\to Y(\C)$ of the spaces of complex points equipped with the complex analytic topology.
\end{lemma}

{\bf Proof.} All assertions except the last follow from e.g.~\cite[ibid.]{stacks}.
%
%
%
To prove the last assertion we use the fact that $sn$ is universally open, which implies that $sn(\C)$ is open, see~\cite{moret_bailly} and references therein.$\clubsuit$
%

The map $Z^{C\gets H}$ takes a closed subscheme of $\C\p^n$ to the effective cycle that is the union of the top dimensional components with multiplicities, see~\cite[I.3.1.3]{kollar}. In particular, all components of lower dimension are discarded, so we have $Z^H\circ sn=Z^C\circ Z^{C\gets H}$ on the Hilbert point of a closed subscheme $Y\subset\C\p^n$ if and only if all components of $Y$ have the same dimension. An irreducible component of ${\eu H}_{P,n}$ or of ${\eu H}_{P,n}^{sn}$ will be called {\it pure} if the subschemes of $\C\p^n$ that correspond to all closed points of this component are pure dimensional. The purpose of the next lemma is to show that the Hilbert scheme has sufficiently many pure components.

%

\begin{lemma}\label{loc_constant}
Let $X\subset\C\p^n$ be a closed subscheme. Suppose $Y\subset X$ is a reduced closed subscheme of pure dimension with Hilbert polynomial $P$, and let $\eu{H}\subset\eu{H}^X_{P,n}$ be an irreducible component which contains the point $y\in \eu{H}^X_{P,n}$ that corresponds to $Y$. Then the subscheme of $X$ that corresponds to any closed point of $\eu{H}$ is pure dimensional (although not necessarily reduced).
\end{lemma}

{\bf Proof.} The following argument is based on a suggestion by Will Sawin~\cite{mo_question_on_flat_deformations}. 
Let $\eu{U}=\eu{H}\times_{\eu{H}^X_{P,n}}\eu{U}^X_{P,n}$ be the (flat) universal family over $\eu{H}$. 
Let $y'\in \eu{H}$ be a closed point such that not all irreducible components of $\eu{U}_{y'}$ have the same dimension, take a reduced irreducible curve $C\subset \eu{H}$ though $y$ and $y'$, 
and let $\tilde C\to C$ be the normalisation of $C$. Let $\tilde y$ and $\tilde y'$ be points of $\tilde C$ over $y$ and $y'$ respectively, and set
$\eu{U}'\subset\eu{U}_{\tilde C}={\tilde C}\times_{\eu{H}}\eu{U}$ to be the union of all components of maximal dimension with the reduced scheme structure. 

The schemes $\eu{U}_{\tilde C}$ and $\eu{U'}$ are flat over $\tilde C$: for $\eu{U}_{\tilde C}$ this is clear, 
and for $\eu{U'}$ we use e.g.~\cite[Proposition 9.7 of chapter III]{hartshorne}.
The fibre $\eu{U}_{\tilde C,\tilde y}$ of $\eu{U}_{\tilde C}$ over $\tilde y$ is $\cong \eu{U}_{y}=Y$. The fibre $\eu{U}'_{\tilde y}$ of $\eu{U}'$ over $\tilde y$ is closed 
and has the same underlying topological space as $Y$: it follows from e.g.~\cite[Proposition 9.5 of chapter III]{hartshorne}
and the assumption that $Y$ has pure dimension that the local dimension of $\eu{U}_{\tilde C}$ at every point of the fibre $\eu{U}_{\tilde C,\tilde y}\cong Y$ is $1+\dim Y$, so every point of $\eu{U}_{\tilde C,\tilde y}$ belongs to an irreducible component of $\eu{U}_{\tilde C}$ of maximum dimension. 

Using the fact that $Y$ is reduced we conclude that $\eu{U}'_{\tilde y}=Y$, 
in particular, these subschemes of $\C\p^n$ have the same Hilbert polynomial. Meanwhile, the Hilbert polynomials of the fibres $\eu{U}_{\tilde C,\tilde y'}\supset \eu{U}'_{\tilde y'}$ are not equal because by our assumptions the former scheme contains some components the latter does not.$\clubsuit$

{\bf Proof of Proposition~\ref{cycle_maps}.} Let $\eu{H}\subset \bigsqcup_{P\in\Q[x]} \eu{H}^{sn}_{P,n}$ be the union of the pure components of $\eu{H}^{sn}_{P,n}$ for $P$ the Hilbert polynomial of an irreducible reduced subscheme of $\C\p^n$ of dimension $k$ and degree $l$. Note that the union is finite. 
Set ${\eu H}^{X}$ to be ${\eu H}$ intersected with the preimage of $\bar{\eu C}_{k,l,n}^X$ under $Z^{C\gets H}$. By lemma~\ref{loc_constant}, every element of ${\eu C}_{k,l,n}$ is in the image under $Z^{C\gets H}$ of a pure component $\subset\eu{H}$. Recall also that $\bar{\eu C}_{k,l,n}^X$ was defined as the closure of ${\eu C}_{k,l,n}^X$ in $\bar{\eu C}_{k,l,n}$. So $Z^{C\gets H}$ induces a surjective map $Z:{\eu H}^{X}(\C)\to \bar{\eu C}_{k,l,n}^X$.

Let us equip $\eu{H}^X(\C)$ and $\bar{\eu{C}}^X_{k,l,n}$ with the complex analytic topology. Since $\eu{H}^X(\C)$ is contained in a union of pure components, on $\eu{H}^X(\C)$ we have $Z^H\circ sn(\C)
=Z^C\circ Z^{C\gets H}=Z^C\circ Z$, where $$sn(\C):\bigsqcup_{P\in\Q[x]} \eu{H}^{sn}_{P,n}(\C)\to \bigsqcup_{P\in\Q[x]} \eu{H}_{P,n}(\C)$$ is induced by the semi-normalisation morphism. 
The map $Z^H$ is continuous by proposition~\ref{cycle_maps_hilbert}, and so are $sn(\C)$ and $Z$, and the latter is moreover surjective. Both $\eu{H}^X(\C)$ and $\bar{\eu{C}}^X_{k,l,n}$ are compact and Hausdorff, so $Z^C$ is continuous.$\clubsuit$

We state for the record the following lemma, which we will use several times in the sequel.

\begin{lemma}\label{gen_topo}
Let $f:A\to B$ be a continuous map of compact Hausdorff spaces. Let $A_1\subset A$ be a subspace on which $f$ is injective. Then $f$ induces a homeomorphism $A_1\to f(A_1)$ if and only if $f(A_1)\cap f(\bar A_1\setminus A_1)=\varnothing$.
\end{lemma}
%
%
$\clubsuit$

\begin{corollary}\label{Bkl}
We use the notation of proposition~\ref{cycle_maps}. Let $\eu B_{k,l}^X$ be $Z^C({\eu C}_{k,l,n}^X)$, i.e.\ the subspace of $2^{X}$ which consists of all (maybe reducible) subvarieties of $2^X$ of pure dimension $k$ and degree $l$. Then $Z^C|_{{\eu C}_{k,l,n}^X}:{\eu C}_{k,l,n}^X\to \eu B_{k,l}^X$ is a homeomorphism, and we have $$\bar{\eu B}_{k,l}^X\subset\bigcup_{j\leq l}\eu B_{k,j}^X.$$
\end{corollary}
{\bf Proof.} The space $\bar{\eu C}_{k,l,n}^X$ is compact in the complex analytic topology, and $2^X$ is Hausdorff. 
Moreover, by definition $Z^C$ takes ${\eu C}_{k,l,n}^X$ and $\bar{\eu C}_{k,l,n}^X\setminus{\eu C}_{k,l,n}^X$ to disjoint subsets of $2^X$. By lemma~\ref{gen_topo}, this implies the first assertion. The second one follows from the fact that $\bar{\eu B}_{k,l}^X\subset Z^C(\bar{\eu C}_{k,l,n}^X)$.$\clubsuit$

When no confusion is likely, we will write $\eu{B}_{k,l}$ instead of $\eu{B}_{k,l}^X$.


\smallskip

{\bf Remark.} The real algebraic analogue of corollary~\ref{Bkl} is not true: it may well happen that the limit of a sequence of $k$-dimensional real algebraic varieties $\subset\R \p^n$ has dimension $<k$.

{\bf Remark.}\label{nocomplexstr} If $\mathbf{M}$ is complex algebraic and a subspace $X\subset 2^\mathbf{M}$ has a natural complex algebraic structure, the closure $\bar X$ may not have one. Take for example $\mathbf{M}= \C \p^2$ and let $L\subset \C \p^2$ be a fixed line. Set $X$ to be the set of all unions $A\cup L$ such that $A\cap L=\varnothing$ and $|A|=1$. We have $X\cong\C^2$. The closure $\bar X$ is the set of all $A\cup L$ such that $A$ is an arbitrary one-element subset of $\C \p^2$. There is a natural surjective map $f:\mathbb{C} \p^2\to \bar X$ that collapses $L$ to a point, but there is no complex analytic structure on $\bar X$ which makes $f$ holomorphic: if $U$ is a neighbourhood of $L$ in $\mathbb{C} \p^2$, then every holomorphic mapping $f:U\to\C$ is constant. 
Note also that more generally, there are no complex analytic surfaces (possibly singular and possibly non-algebraic) that are homeomorphic to the sphere $S^4\cong \C \p^2/L\cong \bar X$. This follows from the claim below and the fact that no smooth manifold homeomorphic to $S^4$ has an almost complex structure.

{\bf Claim.} Suppose a topological 4-manifold $M$ is homeomorphic to a complex analytic surface $S$. Then the normalisation of $S$ is smooth and homeomorphic to $S$.

{\bf Proof of the claim.} Let $\nu:\tilde S\to S$ be the normalisation of $S$, and let $x\in S$ be a point. The germ $\mathcal{S}$ of $S$ at $x$ is irreducible: otherwise the local homology group $H_3(S,S\setminus\{x\},\Z/2)$ would have rank $>1$.
The link $L_{S,x}$ of $S$ at $x$ 
is simply-connected: $x$ has a neighbourhood homeomorphic to the open cone $Cone^\circ(L_{S,x})$, but it also has a neighbourhood homeomorphic to $Cone^\circ(S^3)$ (the homeomorphisms take $x$ to the vertices of the cones).
It now follows from~\cite{michel} that the normalisation $\tilde{\mathcal{S}}$ of $\mathcal{S}$ is smooth and the map $\tilde{\mathcal{S}}\to\mathcal{S}$ is a homeomorphism. This implies that $\tilde S$ is smooth and $\nu$ is open. Since $\nu$ is birational, it follows that it is also injective, so it is a global homeomorphism. 
%

\subsubsection{Families of configuration spaces}\label{sect_fam_conf}

A conical resolution of $\Sigma$ will be constructed below from an (ordered) family $X_1,\ldots,X_N$ of configuration spaces that satisfy the conditions which we now state.

\begin{cl}\label{conditions}
\newcounter{ccc}
\newcounter{ddd}
\begin{enumerate}
\item\label{first} For every $f\in\Sigma$ the singular locus $\Sing(f)$ belongs
to some $X_i$.

\item\label{sec} If $K\in X_i, L\in X_j, K\subsetneq L$, then
$i<j$.

\item\label{three} Recall that $L(K)$ was defined above to be the $\kk$-vector space formed by all
$f\in V$ such that $K\subset \Sing(f)$. We require that for every $i$ the $\kk$-dimension of
$L(K)$ is
the same for all $K\in X_i$; this dimension will be denoted by $d_i$.

\item\label{four} $X_i\cap X_j=\varnothing$ if $i\neq j$.

\item\label{five}
\addtocounter{ccc}{\value{enumi}}
Let $K$ be a configuration from $\bigcup \bar X_i$. We require that among all configurations $\in \bigcup X_i$ that contain $K$ there should be a smallest one, $K'$. We will call $K'$ the \emph{geometrisation} of $K$ and denote it $\Geom(K)$. Note that if $K\in \bigcup X_i$, then $\Geom(K)=K$. We require moreover that $L(K)=L(\Geom(K))$ and $\Geom (K)\in\bigcup_{j<i} X_j$ for all $K\in \bar X_i\setminus X_i$.
%

\item\label{last}
\addtocounter{ddd}{\value{enumi}}

We set
\begin{equation}\label{defai}
A_i=\left\{K\in\bigcup\bar X_j\middle\vert K\subset K'\mbox{ for some }K'\in\bigcup_{k\leq i} X_k	\right\},
\end{equation}
$$\eu T_i=\left\{(x,K)\in \mathbf{M}\times X_i\mid x\in K\right\},$$
$$\eu R_i=\left\{(K,K')\in A_{i}\times X_i\mid K\subset K'\right\},$$
$$\eu S_i=\left\{(K,K')\in A_{i-1}\times X_i\mid K\subset K'\right\}.$$
We require the spaces $\eu T_i, \eu R_i$ and $\eu S_i$ to be fibred over $X_i$ with the obvious projections. Moreover, a local trivialisation of $\eu T_i\to X_i$ must induce local trivialisations of $\eu R_i\to X_i$ and $\eu S_i\to X_i$.

\item\label{xxx}
If $K$ is a finite configuration from
$X_i$, then every nonempty $K'\subset K$ belongs to some $X_j$ with
$j<i$.
\setcounter{enumi}{-1}
\addtocounter{enumi}{\value{ccc}}
\renewcommand{\labelenumi}{\theenumi-.}
\item Same as condition \ref{five}, except that we require that for every $K\in\bar X_i\setminus X_i$ the geometrisation should belong to $\bigcup_{j\leq i}X_j$ (and not necessarily to $\bigcup_{j<i}X_j$);
\setcounter{enumi}{-1}
\addtocounter{enumi}{\value{ccc}}
\renewcommand{\labelenumi}{\theenumi+.}
\item For every $i$ we have $\bar X_i\subset\bigcup_{j\leq i}X_j$.
\end{enumerate}
\end{cl}
\renewcommand{\labelenumi}{\theenumi.}

%




{\bf Remark.} Condition \ref{five}$+$ implies condition \ref{five}, which in turn implies \ref{five}$-$. Note also that condition \ref{five}$+$ only involves configuration spaces, while conditions \ref{five} and \ref{five}$-$ also involve the vector subspaces of $V$ that correspond to configurations.

{\bf Remark.} Constructing a sequence of configuration spaces that satisfies the conditions of list \ref{conditions} usually takes some trial and error. The first step is to classify all possible singular loci of the elements of $\Sigma$. One then takes the union of the closures of the resulting configuration spaces (and maybe some additional configuration spaces, see the remarks after propositions~\ref{sing_plane_cub} and \ref{sing_plane_quart}) and stratifies it e.g.~as described in lemma~\ref{auxconf}. Usually there are several natural stratifications to choose from, and one has to find the one that gives a spectral sequence that converges the fastest.

Before going further we consider a simple example. 

{\bf Example.}\label{simpleex} Suppose $V=\Pi_{3,2}$, and $\Sigma$ is the set of all $f$ such that the zero locus of $f$ in $\C \p^2$ is singular. We set $\mathbf{M}=\C \p^2$ and define $\Sing(f),f\in\Sigma,f\neq 0$ to be the singular locus of the curve given by $f$, and we set $\Sing(0)=\C \p^2$. Let $X_1,\ldots,X_5$ be the following subspaces of $2^{\C \p^2}$.
\begin{itemize}
\item $X_1=B(\C \p^2,1)$;
\item $X_2=B(\C \p^2,2)$;
\item $X_3=\{l\mid\mbox{$l$ is a line in $\C \p^2$}\}$;
\item $X_4$ is the set of all $\{x_1,x_2,x_3\}\subset\C \p^2$ such that $x_1,x_2$ and $x_3$ are not on a line;
\item $X_5=\{\C \p^2\}$.
\end{itemize}

These configuration spaces satisfy conditions \ref{first}-\ref{xxx} from the above list, but not condition~\ref{five}$+$. For every $K\in\bar X_i\setminus X_i,i=1,\ldots,5$ we have $\Geom(K)=K$, except when $i=4$ and $K$ consists of three points on a line $l$, in which case $\Geom(K)$ is the line $l$. If we set $X'_1=X_1,X'_2=X_2,X'_3=\{K\mid K=\mbox{three points on a line}\}$ and $X'_i=X_{i-1}$ for $i=4,5,6$, we get a sequence of configuration spaces $X'_i$ that satisfy conditions \ref{first}-\ref{four} and \ref{five}$+$.

\medskip

Recall that for a closed subvariety $\mathbf{M}\subset\C\p^n$ we defined above the spaces $\eu B_{k,l}^\mathbf{M}\subset 2^\mathbf{M}$, that consist of all (maybe reducible) subvarieties of $\mathbf{M}$ of pure dimension $k$ and degree $l$. In the proof of the next proposition we take $\mathbf{M}=\C\p^n$ and write $\eu{B}_{k,l}$ for $\eu{B}_{k,l}^{\C\p^n}$.

\begin{Prop}\label{chow}
Let $n\geq k\geq 0$ and $l\geq 1$ be integers. There exists a sequence $Z_1,\ldots,Z_R$ of configuration spaces, $Z_i\subset 2^{\C \p^n}$, that satisfies conditions \ref{sec}, \ref{four} and \ref{five}$+$ from list \ref{conditions} and such that
\begin{equation}\label{chow_union}
\bigcup_{i\leq R}Z_i=\bigcup_{s=1}^{k+1}\bigcup_{0\leq i_1< \cdots <i_s\leq k}\left\{ K_1\cup\cdots\cup K_s\middle\vert K_r\in\bigcup_{j\leq l}{\eu B}_{i_r,j},r=1,\ldots,s\right\}.
\end{equation}
\end{Prop}

Let $Z$ be the right hand side of~(\ref{chow_union}). The elements of $Z$ are all possible subvarieties of $\C \p^n$ of (not necessarily pure) dimension~$\leq k$ with all components of degree $\leq l$.

{\bf Proof.} We will call the union of the components of dimension $k$ of a complex algebraic variety $K$ the {\it $k$-dimensional part} of $K$, and we denote this union $K(k)$. To each subvariety $K\subset\C\p^n$ we associate its {\it multidegree} $\mathop{\mathrm{mdeg}}(K)=(\deg K(n),\deg K(n-1),\ldots,\deg K(0))$, with the convention that if $K(i)=\varnothing$, then $\deg K(i)=0$. We introduce the {\it lexicographic order} on multidegrees: if $m_1,m_2$ are two integer sequences of length $n+1$, then we declare that $m_1\leq m_2$ if and only if either $m_1=m_2$, or the first non-zero element of $m_2-m_1$ is positive. Observe that
\begin{itemize}
\item if $K_1\subset K_2\subset\C\p^n$ are subvarieties, then $\mathop{\mathrm{mdeg}}(K_1)\leq \mathop{\mathrm{mdeg}}(K_2)$;
\item $Z\subset 2^{\C\p^n}$ is closed, and the function $K\mapsto \mathop{\mathrm{mdeg}}(K)$ is lower semicontinuous on $Z$, i.e.\ if $(K_i)$ is a sequence of elements of $Z$ that tends to $K\in 2^{\C\p^n}$ in the Hausdorff metric, then $K\in Z$ and $\mathop{\mathrm{mdeg}}(K)\leq\min\mathop{\mathrm{mdeg}}(K_i)$.
\end{itemize}

(The second statement follows from proposition~\ref{conf_space} and corollary~\ref{Bkl}.) We can now obtain $Z_1,\ldots,Z_R$ by subdividing $Z$ by multidegree.$\clubsuit$

The following lemma shows that in the complex projective case, once one has a sequence of configuration spaces $X_1,\ldots, X_N$ that satisfy conditions \ref{first}-\ref{four} and \ref{five}$-$ from list \ref{conditions}, one can construct auxiliary spaces $Y_1,\ldots, Y_\N$ that satisfy conditions \ref{first}-\ref{four} and \ref{five}$+$. These auxiliary configuration spaces may be more complicated than the ones we started with, but luckily, they will only play an auxiliary role.

\begin{lemma}\label{auxconf} Suppose $\mathbf{M}\subset \C \p^n$ is a subvariety.
Let $X_1,\ldots,X_N$ be a sequence of configuration spaces that satisfy conditions \ref{first}-\ref{four} and \ref{five}$-$ from list \ref{conditions} and such that $\bigcup_{j\leq N}X_j$ is contained in a finite union $\eu{B}$ of the spaces $\eu B_{k,l}$ defined in corollary~\ref{Bkl}. We also suppose that for every $i$ and every $k$, the total degree of the $k$-dimensional part $K(k)$ is the same for all $K\in X_i$.

Then there exist a positive integer $\N$, a strictly increasing function $\alpha:\{1,\ldots,N\}\to\{1,\ldots,\eu N\}$ and a sequence $Y_1,\ldots,Y_\N$ of configuration spaces such that
\begin{itemize}
\item For every $i=1,\ldots,N$ we have $Y_{\alpha(i)}=X_i$;
\item $\bigcup_{j\leq\N}Y_j=\bigcup_{j\leq N}\bar X_i$;
\item If $\alpha(i-1)<j\leq\alpha(i)$, then the geometrisation of every $K\in Y_j$ belongs to $X_i$;
\item The spaces $Y_1,\ldots,Y_\N$ satisfy conditions \ref{first}-\ref{four} and \ref{five}$+$ from list \ref{conditions}.
\end{itemize}
\end{lemma}

{\bf Proof.} Recall that above we defined the space $A_i$, see formula (\ref{defai}). It follows from condition \ref{five}$-$ that a configuration $K\in\bigcup_{j\leq N}\bar X_j$ is a subset of some $K'\in\bigcup_{j\leq i}\bar X_j$ if and only if it is a subset of  some $K''\in\bigcup_{j\leq i} X_j$. So every $A_i$ is closed in $2^{\mathbf{M}}$. Note also that by corollary~\ref{Bkl} and our assumptions on the spaces $X_i$, a configuration $K\in X_i$ cannot be the limit of a sequence $(L_j)$ such that every $L_j\in \eu{B}$ and is $\subsetneq$ some $K_j\in X_i$. Indeed, let $k$ be the maximum integer such that the $k$-dimensional part of an element of $X_i$ is non-empty, and let $l$ be the degree of this part. Then the $k$-dimensional part of every $L_j$ will be of degree $<l$, or will be empty, and by corollary~\ref{Bkl} the same will be true for $\lim_{j\to\infty} L_j$. We conclude that every $X_i$ is an open subset of $A_i\setminus A_{i-1}$.

We now subdivide every space $A_i\setminus (A_{i-1}\cup X_i)$ by intersecting it with the spaces $Z_1,\ldots,Z_R$ obtained by applying proposition \ref{chow} to sufficiently large $k$ and $l$. By inserting for every $i=1,\ldots,N$ the resulting subdivision of $A_i\setminus (A_{i-1}\cup X_i)$ between $X_i$ and $X_{i-1}$ and renaming we get a sequence of configuration spaces $Y_1,\ldots,Y_\N$ and a map $\alpha:\{1,\ldots, N\}\to\{1,\ldots\N\}$ such that $Y_{\alpha(i)}=X_i,i=1,\ldots, N$. The spaces $Y_1,\ldots,Y_\N$ automatically satisfy conditions~\ref{first} and~\ref{four}. Condition~\ref{sec} follows from the fact that an element of $A_i\setminus (A_{i-1}\cup X_i)$ is a subset of a unique configuration from $X_i$ and is not included in any configuration from $X_j,j<i$. Checking condition~\ref{five}$+$ is straightforward. 

We have $\bigcup_{j\leq\N}Y_j=\bigcup_{j\leq N}\bar X_i=A_N$ by construction. Finally, for $K\in A_i\setminus (A_{i-1}\cup X_i)=\bigcup_{\alpha(i-1)<j\leq\alpha(i)} Y_j$, the geometrisation is the unique element of $X_i$ that contains $K$. Since $L(\Geom(K))=L(K)$ by condition~\ref{five}$-$ for $X_1,\ldots, X_N$, we have $\dim_{\mathbf{k}} L(K)=d_i$ for all $K\in \bigcup_{\alpha(i-1)<j\leq\alpha(i)} Y_j$. So $Y_1,\ldots,Y_\N$ satisfy condition~\ref{three} as well.$\clubsuit$

\subsection{Generalised conical resolutions}\label{conical_sect}

\subsubsection{Joins and self-joins}\label{sect_joins}

We will now describe a version of the conical resolution that generalises the one given in~\cite{quintics}. Let us first recall the notion of the {\it $k$-th self-join} of a topological space (see \cite{vas3}), also known as the {\it $k$-th symmetric join}. For simplicity all spaces in~\ref{sect_joins} are assumed to be compactly generated and Hausdorff (e.g.\ CW-complexes or metric spaces).

Let $A$ be a topological space. Given a positive integer $k$, the {\it $k$-th self-join} $A^{*k}$ of $A$ is the space $(A^k\times \triangle^{k-1})/\sim$, where
$$\triangle^{k-1}=\left\{(t_1,\ldots,t_k)\in\R^k\middle\vert \sum t_i=1, \mbox{all } t_i\geq 0\right\}$$ is the standard $(k-1)$-simplex, the product $A^k\times\triangle^{k-1}$ is given the compactly generated topology (which coincides with the product topology e.g.\ if $A$ is locally compact), and $\sim$ is the equivalence relation obtained as follows. An injective map $\imath$ from $A$ to a real vector space $W$ is {\it $k$-generic} if every $(k-1)$-simplex with vertices in $\imath(A)$ is non-degenerate and the intersection of any two such simplices is their common face; in particular, the intersection is empty if the sets of the vertices are disjoint. (Note that here we do not require $W$ to have any topology.) We then define $\sim$ to be the equivalence relation induced by the map
\begin{equation}\label{two_models_self_join}
A^k\times\triangle^{k-1}\to W, ((a_1,\ldots, a_k),(t_1,\ldots, t_k))\mapsto \sum t_i \imath(a_i)
\end{equation}
where $\imath:A\to W$ is an arbitrary $k$-generic injective map. 

Explicitly, $\sim$ is generated by the diagonal action of the symmetric group $S_k$ on $A^k\times \triangle^{k-1}$, and by identifying any two $(a,t),(a',t')\in A^{k}\times\triangle^{k-1}$ such that (1) $t=t'$, and (2) if $i$ is such that the $i$-th coordinate of $t=t'$ is non zero, then the $i$-th components of $a$ and $a'$ coincide.

Suppose $A$ is compact, $W$ is a Hausdorff topological real vector space, and $\imath$ is a $k$-generic continuous embedding. Then (\ref{two_models_self_join}) is a quotient map, so in this case $A^{*k}$ is homeomorphic to the union of all $(k-1)$-simplices in $W$ with all vertices in $\imath(A)$. This is the model for $A^{*k}$ that we will use in the sequel. If $A$ is not necessarily compact, we can identify $A^{*k}$ with the direct limit of $B^{*k}$ for $B\subset A$ compact: denoting the set of all compact subspaces of $A$ by $\mathcal{K}(A)$, both 
$$A^{*k}=\left(\underset{B\in\mathcal{K}(A)}\varinjlim B\right)^{*k}\mbox{ and } \underset{B\in\mathcal{K}(A)}\varinjlim B^{*k}$$
are quotient spaces of $\bigsqcup_{B\in\mathcal{K}(A)} B^k\times\triangle^{k-1}$ with respect the equivalence relation induced by the natural map $\bigsqcup_{B\in\mathcal{K}(A)} B^k\times\triangle^{k-1}\to A^{*k}$.

We will sometimes denote the equivalence class of $((a_1,\ldots, a_k),(t_1,\ldots, t_k))$ by $\sum_{i=1}^k t_i a_i$. 
Also, if $I=\{i_1,\ldots, i_l\}\subset\{1,\ldots,k\}$ is a non-empty subset, the point $\frac{1}{l}\sum_{j=1}^l a_{i_j}$ will be called the {\it barycentre} of the simplex with vertices $a_{i_1},\ldots, a_{i_l}$.

\smallskip

The {\it join} $A_1 * \cdots * A_k$ of non-empty topological spaces $A_1,\ldots, A_k$ is the quotient $(A_1\times\cdots\times A_k\times\triangle^{k-1})/\approx$ where $A_1\times\cdots\times A_k\times\triangle^{k-1}$ is equipped with the compactly generated topology and $\approx$ is the equivalence relation induced by the map
$$(A_1\times\cdots\times A_k)\times\triangle^{k-1}\to W, ((a_1,\ldots, a_k),(t_1,\ldots, t_k))\mapsto \sum t_i \imath(a_i)$$
for some $k$-generic injective map $\imath$ from $A=A_1\sqcup\cdots\sqcup A_k$ to a real vector space $W$. Note that $A_1 * \cdots * A_k$ can be identified with a closed subspace of $A^{*k}$, namely the space of all $\sum t_i a_i$ such that all $t_i\geq 0,\sum t_i=1$ and every $a_i\in A_i$. 
This interpretation allows one to extend the definition to the case when some $A_i$ are empty, in which case we simply set $A_1 *\cdots * A_k=\Asterisk_{A_i\neq\varnothing} A_i$.

If all $A_i$ are compact, $W$ has a Hausdorff vector space topology, and $\imath$ is continuous, then $A_1*\cdots *A_k$ is homeomorphic to the union of all $(k-1)$-simplices $\triangle\subset W$ such that vertices of $\triangle$ belong to pairwise distinct spaces $\imath(A_1),\ldots, \imath(A_k)$. 
For general $A_1,\ldots,A_k$ we have $A_1*\cdots * A_k\cong\varinjlim B_1*\cdots * B_k$ for $(B_1,\ldots, B_k)\in \mathcal{K}(A_1)\times \cdots\times\mathcal{K}(A_k)$.

{\bf Remark.} The self-join $A^{*2}$ is not the same as the join $A*A$: for example, note that $A*A$ naturally contains two copies of $A$, while $A^{*2}$ contains only one.

\subsubsection{The conical resolution}\label{constr_con_res}

\begin{assumptions}\label{assumptions_conical}
In the rest of section~\ref{constr_con_res} $X_1,\ldots, X_N$ will be a sequence of configuration spaces satisfying conditions \ref{first}-\ref{four} and \ref{five}$-$ from list \ref{conditions}. Similarly, $\N,Y_1,\ldots,Y_\N$ and $\alpha$ will be respectively an integer, a sequence of configuration spaces and a map $\{1,\ldots,N\}\to\{1,\ldots,\N\}$ that satisfy the conclusions of lemma \ref{auxconf}. We set $\mathbf{X}=\bigcup_{i=1}^N X_i, \mathbf{Y}=\bigcup_{i=1}^\N\bar Y_i=\bigcup_{i=1}^\N Y_i$ and we will assume that the latter space is homeomorphic to a finite CW-complex.
%
\end{assumptions}


Let us choose a continuous $\N$-generic embedding $\imath:\mathbf{Y}\to W$ where $W=\R^\Omega,\Omega<\infty$ with the standard topology, and use it to construct a model for $\mathbf{Y}^{*\N}$, see section~\ref{sect_joins}. We will not distinguish between elements of $\mathbf{Y}$ (which are configurations, i.e. compact subsets of $\mathbf{M}$) and vertices of simplices of the self-join $\mathbf{Y}^{*\N}$ (which are the images of elements of $\mathbf{Y}$ under $\imath$).

We will call a simplex $\triangle\subset \mathbf{Y}^{*\N}$ {\it coherent} if its vertices form a chain, i.e.\ one of any two of them contains the other. Among the vertices of a coherent simplex $\triangle$ there is a configuration that contains all the other; such a configuration will be called the {\it principal vertex} of $\triangle$. Let $\tilde\Lambda$ be the union of all coherent simplices. Note that since $\Sing(0)=\mathbf{M}$, every coherent simplex is contained in a coherent simplex with principal vertex $\mathbf{M}$, so $\tilde\Lambda$ is contractible.

For every $K\in \mathbf{Y}$ we denote by $\tilde\Lambda(K)$\label{tlambda} the union of all coherent simplices with principal vertex $K$. The space $\tilde\Lambda(K)$ is contractible. Note that $\tilde\Lambda=\tilde\Lambda(\mathbf{M})$. For every $K\in \mathbf{X}$ let $\Lambda(K)$\label{xlambda} be the union of all coherent simplices $\triangle$ such that the principal vertex of $\triangle$ is $K$, and for every vertex $K'\neq K$ of $\triangle$ we have $\Geom (K')\neq K$; note that by condition~\ref{five}$-$ we then have $\Geom (K')\subsetneq K$, so $\Geom(K')\in\bigcup_{j<i}X_j$ by condition~\ref{sec}. In other words, $\Lambda(K)$ is the subspace of $\tilde\Lambda(K)$ that is the union of all coherent simplices $\triangle$ with principal vertex $K$ such that the largest vertex of $\triangle$ different from $K$ is a subset of an element of some $X_{<i}$. Same as $\tilde\Lambda(K)$, the space $\Lambda(K)$ is contractible.

For every $K\in X_i=Y_{\alpha(i)}$ we set
$$\partial\Lambda(K)=\bigcup_{j<i}\bigcup_{\genfrac{}{}{0pt}{}{K'\in X_j=Y_{\alpha(j)}}{K'\subset K}}\tilde\Lambda (K').$$\label{plambda}
Similarly, for $K\in Y_i$ set
$$\tilde\partial\Lambda(K)=\bigcup_{j<i}\bigcup_{\genfrac{}{}{0pt}{}{K'\in Y_j}{K'\subset K}}\tilde\Lambda (K').$$\label{tplambda}
In other words, $\partial\Lambda(K)$ is the subspace of $\tilde\Lambda(K)$ that is the union of all coherent simplices with principal vertex in some $X_{<i}$. Equivalently, $\partial\Lambda(K)$ is the union of the spaces $\tilde\Lambda (K')$ for maximal proper $K'\subsetneq K$ such that $K'\in\bigcup_{j<i} X_j$. 

The space $\tilde\Lambda(K), K\in Y_i$ is the union of all segments that join $K$ with elements of $\tilde\partial\Lambda (K)$. As we will shortly see (in lemma~\ref{confclosed}), the latter space is compact, so $\tilde\Lambda(K)$ is homeomorphic to the cone $Cone(\tilde\partial\Lambda(K))$. 

Similarly, for every $K\in X_i$ we have $\partial\Lambda(K)\subset\Lambda(K)$, and for every coherent simplex $\triangle\subset\Lambda(K)$ with principal vertex $K$ there is a possibly larger coherent simplex $\triangle'\subset\Lambda(K)$ such that the next largest vertex after $K$ (if any) belongs to $\bigcup_{j<i} X_j$, cf.\ above.
So lemma~\ref{confclosed} will allow us to identify $\Lambda(K)$ with $Cone(\partial\Lambda(K))$. 

{\bf Remark.} The notation we use here differs from the notation of~\cite{quintics}: what was denoted there by $\Lambda(K)$ and $\partial\Lambda(K)$ we now denote by $\tilde\Lambda(K)$ and $\tilde\partial\Lambda(K)$ respectively. The reason for this is that in our examples we will work more often with the spaces $\Lambda(K)$ and $\partial\Lambda(K)$ as defined above.

The following lemma is a straightforward consequence of the definitions.

\begin{lemma}\label{cup_cap}
Suppose $K_1,K_2,K_1\cap K_2\in \mathbf{Y}$. Then we have $\tilde\Lambda(K_1)\cap\tilde\Lambda (K_2)=\tilde\Lambda(K_1\cap K_2)$. If moreover $K_1,K_2,K_1\cap K_2\in \mathbf{X}$ and $K_1\neq K_1\cap K_2\neq K_2$, then $\Lambda(K_1)\cap\Lambda (K_2)=\tilde\Lambda(K_1\cap K_2)$.
\end{lemma}
$\clubsuit$

For every simplex $\triangle\subset \mathbf{Y}^{*\N}$ let $\stackrel{\circ}\triangle$ be the {\it interior} of $\triangle$, i.e. the set of all points of $\triangle$ that do not belong to any of the faces of lower dimension. For every $x\in \mathbf{Y}^{*\N}$ there exists a unique simplex $\triangle$ such that $x\in\stackrel{\circ}\triangle$.

\begin{lemma}\label{sequences_self_join}
Let $(x_l)$ be a convergent sequence in $\mathbf{Y}^{*\N}$, and let $\triangle$ be the simplex such that $\stackrel{\circ}\triangle\ni x=\lim_{l\to\infty} x_l$. Suppose $\triangle_l\subset \mathbf{Y}^{*\N}$ is a simplex that contains $x_l$. Then every vertex of $\triangle$ is the limit of a sequence of vertices of $\triangle_l$.
\end{lemma}
$\clubsuit$

We set $\Phi_0=\varnothing$ and $\Phi_i=\bigcup_{j\leq i}\bigcup_{K\in Y_j}\tilde\Lambda(K)$ for $i\in\{1,\ldots, \N\}$. By our assumptions on p.~\pageref{assumptions_conical} and condition~\ref{five}$-$ of list \ref{conditions}, if $K'\in Y_{j'}$ and $\alpha(i'-1)<j'\leq\alpha(i')$, then $K'\subset K$ for some $K\in Y_{\alpha(i')}$, so we have
$\Phi_{\alpha(i)}=\bigcup_{j\leq i}\bigcup_{K\in X_j=Y_{\alpha(j)}}\tilde\Lambda(K)$.
There is a filtration on
$\tilde \Lambda$: $\varnothing=\Phi_0\subset\Phi_1\subset\cdots\subset\Phi_\N=\tilde\Lambda$.

%
\begin{lemma}\label{confclosed}
The following spaces are closed in $\mathbf{Y}^{*\N}$:
\begin{itemize}
\item the spaces $\Lambda(K)$ and $\partial\Lambda(K)$ for every $K\in \mathbf{X}$;
\item the spaces $\tilde\Lambda(K)$ and $\tilde\partial\Lambda(K)$ for every $K\in \mathbf{Y}$;
\item the spaces $\Phi_i$ for every $i=1,\ldots, \N$.
\end{itemize}

Moreover, if $X_1,\ldots, X_N$ satisfy condition \ref{five} in addition to conditions \ref{first}-\ref{five}$-$, then for every $i=1,\ldots, N$ the space $\Phi_{\alpha(i-1)}\cup \bigcup_{K\in X_i}\Lambda(K)$ is also closed in $\mathbf{Y}^{*\N}$.
\end{lemma}

{\bf Proof.} The statement of the lemma for $\Phi_i$, $\tilde\Lambda(K)$ and $\tilde\partial\Lambda(K), K\in \mathbf{Y}$ follows from \cite[proposition 2.7]{quintics}.

%
Let us prove the ``moreover'' part of the lemma. Let $x\in \Phi_{\alpha(i)}$ be the limit of a sequence $(x_l)$ such that every $x_l\in\Lambda(K_l),K_l\in X_i$. We need to show that $x\in\Phi_{\alpha(i-1)}\cup \bigcup_{K\in X_i}\Lambda(K)$. In the following proof it will be convenient for us to allow empty configurations: if $\triangle$ is a coherent simplex with just one vertex, $K$, then the next largest vertex of $\triangle$ after $K$ will be the empty set, which we consider as an element of the space $2^{\mathbf{M}}\sqcup\{\varnothing\}$. 

For every $l$ there is a coherent simplex $\triangle_l\ni x_l$ with principal vertex $K_l$ and next largest vertex $K'_l\in\bigcup_{j<i} X_j\sqcup\{\varnothing\}$. By choosing subsequences if necessary and using condition \ref{five}$+$ satisfied by the spaces $Y_1,\ldots, Y_\N$, we may assume that $K_l\to K$ and $K_l'\to K'$ as $l\to\infty$, for some $K\in\bigcup_{j\leq \alpha(i)} Y_j$ and $K'\in\bigcup_{j\leq \alpha(i-1)}Y_j\sqcup\{\varnothing\}$. The point $x$ is in the interior of some coherent simplex $\triangle$ with principal vertex $\tilde K$ and next largest vertex $\tilde K'$. By lemma~\ref{sequences_self_join}, these are limits of sequences of vertices of $\triangle_l$, so we have $\tilde K\subset K$ and $\tilde K'\subset K'$. (Note that if $\tilde K'\neq\varnothing$, then none of $\triangle_l$ except maybe finitely many are points, so $K_l'\neq\varnothing$ for all sufficiently large $l$.)

If $K\in\bar X_i\setminus X_i$, then by condition~\ref{five}, $K\subset\Geom(K)\in\bigcup_{j<i} X_j$, which implies $x\in \Phi_{\alpha(i-1)}$. Suppose now $K\in X_i$ and $\tilde K\subsetneq K$. Then $\tilde K$ is the limit of a sequence of non-principal vertices of $\triangle_l$. So we have $\tilde K\subset K'$ and hence $K'\neq\varnothing$, which means that $K'\in\bigcup_{j\leq \alpha(i-1)}Y_j$. By the assumptions on p.~\pageref{assumptions_conical},  the geometrisation of $K'$ is in $\bigcup_{j<i} X_j$, so we get $\Geom(\tilde K)\subset\Geom(K')\in\bigcup_{j<i} X_j$, so $\Geom(\tilde K)\in\bigcup_{j<i} X_j$ too by condition~\ref{sec}. Also, $\Geom(\tilde K)\subset K$ by the definition of the geometrisation. The simplex $\triangle'$ spanned by $\triangle, K$ and $\Geom(\tilde K)$ contains $x$, is coherent, has principal vertex $K$ and the next largest vertex $\Geom(\tilde K)\in\bigcup_{j<i} X_j$, which implies $x\in\Lambda(K)$.

Suppose $K\in X_i$ and $\tilde K=K$. If $\tilde K'=\varnothing$, then $\triangle$ is a point and we have $x=\tilde K=K\subset\Lambda(K)$. If $\tilde K'\neq\varnothing$, then using the assumptions on p.~\pageref{assumptions_conical} again we have $\Geom(\tilde K')\subset\Geom(K')\in\bigcup_{j<i} X_j$, so $\Geom(\tilde K')\in\bigcup_{j<i} X_j$. We also have $\Geom(\tilde K')\subset K$. We set $\triangle'$ to be the simplex spanned by $\Geom(\tilde K')$ and $\triangle$. The principal vertex of $\triangle'$ is $K=\tilde K$ and the next largest one is $\Geom(\tilde K')$. We conclude again that $x\in\Lambda(K)$.

A similar argument shows that for every $K\in X_i$ the space $\Lambda(K)$ is closed. (In this case we take $K_l=K$ for all $l$.) Finally, if $K\in X_i$ then $\partial\Lambda(K)=\Lambda (K)\cap \Phi_{\alpha(i-1)}$, which implies that $\partial\Lambda(K)$ is also closed.$\clubsuit$

\bigskip

The {\it conical resolution} $\sigma$\label{def_con_res} of $\Sigma$ is the subspace of
$\Sigma\times\tilde\Lambda$ that is formed by all pairs $(f,x)$ such that
$f\in\Sigma, x\in\tilde\Lambda(\Sing(f))$. There exist natural projections $\pi:\sigma\to\Sigma$ and $p:\sigma\to\Lambda$. We introduce a filtration on
$\sigma$ by setting $ F_i=p^{-1}(\Phi_i)$.

\begin{lemma}\label{sigma_is_closed}
\begin{enumerate} 
\item Let $(K_i)$, respectively $(f_i)$ be a sequence of elements of $2^{\mathbf{M}}$, respectively $V$ that converges to $K$, respectively $f$. If all $f_i\in L(K_i)$, then $f\in L(K)$.
\item The space $\sigma\subset V\times\tilde\Lambda$ is closed and the map $\pi:\sigma\to\Sigma$ is proper.
\end{enumerate}
\end{lemma}

{\bf Proof.} To prove the first part recall that $f_i\in L(K_i)$ if and only if $K_i\subset\Sing f_i$. For every $x\in K$ we choose a sequence $(x_i)$ in $\mathbf{M}$ that converges to $x$ and such that each $x_i\in K_i$ (proposition~\ref{conf_space}). We have $ev(f_i,x_i)=0$ for all $i$, so $ev(f,x)=0$ by the continiuty of $ev$. So $x\in\Sing f$. We see that $K\subset \Sing f$, i.e.\ $f\in L(K)$. To prove the second part of the lemma we use lemma~\ref{sequences_self_join} and the first part.$\clubsuit$

Each point $x\in\Phi_{\alpha(i)}\setminus\Phi_{\alpha(i-1)},i=1,\ldots,N$ belongs to the
interior of some coherent simplex $\triangle$; let $K$ be the principal vertex of $\triangle$. The geometrisation of $K$ belongs to
$Y_{\alpha(i)}=X_i$; we denote the map
$\Phi_{\alpha(i)}\setminus\Phi_{\alpha(i-1)}\ni x\mapsto \Geom(K)\in Y_{\alpha(i)}=X_i$ by $e_i$, and let $g_i(x)$ be the $\kk$-vector space $L(e_i(x))$.

Similarly, for $i=1,\ldots, \N$ any $x\in\Phi_i\setminus\Phi_{i-1}$ belongs to the interior of a coherent simplex with principal vertex in $Y_i$; we denote this vertex $f_i(x)$, and we set $\tilde g_i(x)=L(f_i(x))$.

Set
$$\Psi_i=(\Phi_{\alpha(i)}\setminus\Phi_{\alpha(i-1)})\cap\left(\bigcup_{K\in X_i}\Lambda(K)\right).$$

\begin{Prop}\label{lambda}
Suppose the spaces $X_1,\ldots,X_N$ satisfy condition \ref{five} from list \ref{conditions} in addition to \ref{first}-\ref{five}$-$. Then
\begin{enumerate}
\item the maps $e_i, i=1,\ldots,N$ and $f_i,i=1,\ldots, \N$ are continuous;
\item $\Psi_i$ is a closed subspace of $\Phi_{\alpha(i)}\setminus\Phi_{\alpha(i-1)}$, and for every $K\in X_i$ the inclusion $e_i^{-1}(K)\cap\Psi_i\subset e_i^{-1}(K)$ induces an isomorphism of the groups $H^*_c(-,M)$ for an arbitrary abelian group $M$.
\end{enumerate}
\end{Prop}

{\bf Proof.} Let us prove the continuity of $e_i,i=1,\ldots, N$. Let $(x_l)$ be a convergent sequence in $\Phi_{\alpha(i)}\setminus\Phi_{\alpha(i-1)}$, and set $x=\lim_{l\to\infty}x_l$. Let $\triangle_l$, respectively $\triangle$ be the coherent simplex such that $x_l\in\stackrel{\circ}\triangle_l$, respectively such that $x\in\stackrel{\circ}\triangle$. Let $K$ be the principal vertex of $\triangle$. By lemma~\ref{sequences_self_join}, there exists a sequence $(K_l)$ of configurations such that every $K_l$ is a vertex of $\triangle_l$, and $\lim_{l\to\infty} K_l=K$. Since $x, x_l\in\Phi_{\alpha(i)}\setminus\Phi_{\alpha(i-1)}$, we have $K, K_l\in \bigcup_{\alpha(i-1)<j\leq\alpha(i)} Y_j$.
It follows our assumptions on p.~\pageref{assumptions_conical}, in particular condition~\ref{sec} of list~\ref{conditions}, that $e_i(x_l)$, which is the geometrisation of the principal vertex of $\triangle_l$, coincides with $\Geom(K_l)$. Similarly, $e_i(x)=\Geom(K)$.

Let $K'$ be the limit in $\Phi_{\alpha(i)}$ of a convergent subsequence of $(\Geom(K_l))$. The configuration $K'\in X_i$: otherwise by condition \ref{five} it would be a subset of some $K''\in\bigcup_{j<i}X_j$, and so $x$ would belong to $\Phi_{\alpha(i-1)}$. Moreover, $K_l\subset\Geom(K_l)$ and $K_l\to K$ as $l\to\infty$, so $K\subset K'$, which implies $\Geom(K)=K'$ (using condition \ref{sec} and the fact that $\Geom(K)\in X_i$). We have shown that every subsequence of $(\Geom(K_l))$ that converges in $\Phi_{\alpha(i)}$ converges in fact to $\Geom(K)$, so $$\lim_{l\to\infty}e_i(x_l)=\lim_{l\to\infty}\Geom(K_l)=\Geom(K)=e_i(x).$$ This proves that $e_i$ is continuous. For $f_i,i=1,\ldots, \N$ the argument is similar.

It follows from the ``moreover'' part of lemma \ref{confclosed} that $\Psi_i$ is closed in $\Phi_{\alpha(i)}\setminus\Phi_{\alpha(i-1)}$. Finally, for every $K\in X_i$, we have $e_i^{-1}(K)=\tilde\Lambda(K)\setminus\Phi_{\alpha(i-1)}$ and $e_i^{-1}(K)\cap\Psi_i=\Lambda(K)\setminus\Phi_{\alpha(i-1)}$. (The second equality follows from the fact that set-theoretically, we have $\Psi_i=\bigsqcup_{K\in X_i}(\Lambda(K)\setminus\Phi_{\alpha(i-1)})$.) 
%
The difference $e_i^{-1}(K)\setminus (e_i^{-1}(K)\cap\Psi_i)$ is the same as $\tilde\Lambda(K)\setminus\Lambda(K)$ because $\tilde\Lambda(K)\cap \Phi_{\alpha(i-1)}=\Lambda(K)\cap \Phi_{\alpha(i-1)}=\partial\Lambda(K)$. Using e.g.~\cite[II.10.3]{bredon} we get 
$H^*_c(\tilde\Lambda(K)\setminus\Lambda(K),M)=0$
as both $\tilde\Lambda(K)$ and $\Lambda(K)$ are compact and contractible, 
so 
$H^*_c(e_i^{-1}(K)\setminus (e_i^{-1}(K)\cap\Psi_i),M)=0$, 
which implies the last statement of the proposition.$\clubsuit$

\begin{lemma}\label{pullback_taut}
Suppose $V'$ is a finite-dimensional $\mathbf{k}$-vector space, $g:A\to \Gr(d,V')$ is a continuous map, and $\gamma$ is the tautological $\mathbf{k}$-vector bundle over $\Gr(d,V')$. Set $B=\{(v,a)\in V'\times A\mid v\in f(a)\}$ and let $p:B\to A$ be the projection. Then $(B,A,p)$ is a $\mathbf{k}$-vector bundle over $A$ which is isomorphic to $g^*(\gamma)$.
\end{lemma}
$\clubsuit$

Let $M$ be a module over a principal ideal domain $R$.

\begin{lemma}\label{proper_contr_fib}
Let $f:A\to B$ be a proper continuous map where $A,B$ are locally compact Hausdorff spaces. Assume all fibres of $f$ are contractible. Then $f$ induces an isomorphism $\bar H_*(A,M)\to\bar H_*(B,M)$. 
\end{lemma}

{\bf Proof.} The lemma follows from~\cite[Theorem V.6.1]{bredon}.$\clubsuit$

\begin{theorem}
\label{thconical}
%
\begin{enumerate}
Suppose the assumptions on p.~\pageref{assumptions_conical} are satisfied. Then the following holds.

\item The map $\pi$ is proper and induces an isomorphism $\bar H_*(\sigma,M)\to\bar H_*(\Sigma,M)$ for arbitrary~$M$. 
\item For every $i=1,\ldots, N$ the map $g_i$ is continuous, and the triple $$(F_{\alpha(i)}\setminus F_{\alpha(i-1)},\Phi_{\alpha(i)}\setminus\Phi_{\alpha(i-1)},p|_{F_{\alpha(i)}\setminus F_{\alpha(i-1)}})$$ is a $\kk$-vector bundle of rank $d_i$ isomorphic to the pullback of the tautological vector bundle over the Grassmannian $\Gr(d_i,V)$ under $g_i$.
\item For every $i=1,\ldots, \N$ the map $\tilde g_i$ is continuous, and the triple $$(F_i\setminus F_{i-1},\Phi_i\setminus\Phi_{i-1},p|_{F_{i}\setminus F_{i-1}})$$ is a $\kk$-vector bundle of rank $d_i$ isomorphic to the pullback of the tautological vector bundle over the Grassmannian $\Gr(d_i,V)$ under $\tilde g_i$.
\end{enumerate}
\end{theorem}

{\bf Proof.}
The first part follows from lemmas~\ref{sigma_is_closed} and \ref{proper_contr_fib}. The continuity of the maps $g_i,i=1,\ldots, N$ and $\tilde g_i,i=1,\ldots, \N$ follows from part~1 of proposition~\ref{semicont} and part~1 of proposition~\ref{lambda}. 
%
%
To prove the second part of the theorem we use lemma~\ref{pullback_taut} and the fact that for every $x\in\Phi_{\alpha(i)}\setminus\Phi_{\alpha(i-1)}$ we have $p^{-1}(x)=g_i(x)$, which follows immediately from the definitions, cf.~\cite{quintics}, proof of theorem 2.8. The proof of the third part is similar.$\clubsuit$
%

\begin{theorem}\label{thconical_1}
Suppose the assumptions on p.~\pageref{assumptions_conical} are satisfied, and suppose that $X_1,\ldots,X_N$ satisfy conditions \ref{five} and \ref{last} from list \ref{conditions} (in addition to conditions \ref{first}-\ref{four} and \ref{five}$-$). Then the following holds.
\begin{enumerate}
\item For every $i$ the spaces $\Phi_{\alpha(i)}\setminus\Phi_{\alpha(i-1)}$ and $\Psi_i$ are fibred over
$X_i$ with projection $e_i$ and generic fibres
$\tilde\Lambda(K)\setminus\tilde\partial\Lambda(K),K\in X_i$ and $\Lambda(K)\setminus\partial\Lambda(K),K\in X_i$ respectively.
\item If $M$ is finitely generated, then for all $i=1,\ldots, N$ the inclusion $\Psi_i\subset (\Phi_{\alpha(i)}\setminus\Phi_{\alpha(i-1)})$ induces an isomorphism of the groups $\bar H_*(-,M)$.
\end{enumerate}
\end{theorem}

{\bf Proof.} The proof of the first assertion repeats, with minor modifications, the proof of part 3 of \cite[Theorem 2.8]{quintics}. To prove the second assertion we first apply proposition~\ref{lambda} and base change (see e.g.~\cite[Proposition IV.4.2]{bredon}) to show that $\Psi_i\subset (\Phi_{\alpha(i)}\setminus\Phi_{\alpha(i-1)})$ induces an isomorphism of compactly supported cohomology with arbitrary constant coefficients. (Note that a $\Psi$-closed map in ibid.\ is a map that takes every element of $\Psi$ to a closed set, see p.~74, and that we need condition \ref{five}, rather than \ref{five}$-$, to be able to apply proposition \ref{lambda}.) We now use exact sequences (9) and (8) from~\cite[V.3]{bredon} to pass from compactly supported cohomology to Borel-Moore homology.$\clubsuit$

The next theorem shows that if $X_i$ consists of finite configurations and condition \ref{xxx} of list \ref{conditions} is satisfied, then the spaces $\Phi_{\alpha(i)}\setminus\Phi_{\alpha(i-1)}$ and $\Lambda(K),\tilde\Lambda(K),\partial\Lambda(K),\tilde\partial\Lambda(K), K\in X_i$ are easy to describe.

\begin{theorem}\label{thconical_2}
Suppose the assumptions on p.~\pageref{assumptions_conical} are satisfied, and suppose that $X_1,\ldots,X_N$ satisfy conditions \ref{five}, \ref{last} and \ref{xxx} from list \ref{conditions} (in addition to conditions \ref{first}-\ref{four} and \ref{five}$-$). Let $i$ be an index such that $X_i$ consists of finite configurations, and set $k=|K|, K\in X_i$.

Then for all $K\in X_i=Y_{\alpha(i)}$ we have $\tilde\Lambda(K)=\Lambda(K), \tilde\partial\Lambda(K)=\partial\Lambda(K)$, and there exists a continuous map
$(\Phi_{\alpha(i)}\setminus\Phi_{\alpha(i-1)})\to \mathbf{M}^{*k}$ that takes $K\in\Lambda(K)\setminus\partial\Lambda(K)\subset\Phi_{\alpha(i)}\setminus\Phi_{\alpha(i-1)}$
to the barycentre of the simplex $\triangle (K)$
spanned by the points of $K$. This map is
a homeomorphism onto its image, and it takes $\Lambda(K)$
(respectively $\partial\Lambda(K)$) homeomorphically to
$\triangle(K)$ (respectively $\partial\triangle (K)$).
\end{theorem}

\begin{figure}[h]
\centering
\includegraphics[scale=.3]{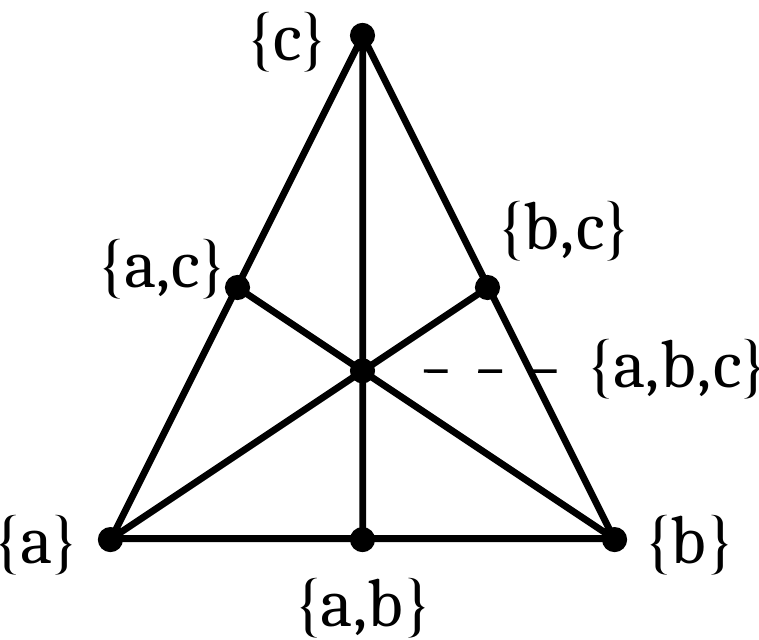}
\caption{$\Lambda(K)$ for $K=\{a,b,c\}$}
\label{fig:triangle_1}
\end{figure}

{\bf Proof.} The theorem follows from \cite[Lemmas 2.10 and 2.11]{quintics}; see figure~\ref{fig:triangle_1} for an illustration.$\clubsuit$

\subsection{Mixed Hodge structures}\label{hodge_modge}
From now on we take $\kk=\C$. So $V$ is a complex vector space, and we assume that the space $\mathbf{M}$ is a compact complex algebraic variety, the vector bundle $\mathcal{V}$ is algebraic and that the bundle map $ev$ (see the beginning of section~\ref{sing_loc}) is regular. This implies that $\Sigma\subset V$ is affine complex algebraic. We will now explain how the construction described in section~\ref{conical_sect} can be used to determine the mixed Hodge structure on the groups $\bar H_*(\Sigma,\Q)$. The conical resolution $\sigma$ as constructed above is only a topological space. Here we construct a cubical complex algebraic variety filtered by closed cubical subvarieties that (after passing to the geometric realisations) is filtered proper homotopy equivalent to $\sigma$ with its natural filtration.

\subsubsection{Cubical objects}

Let us recall basic facts and terminology related to cubical objects; for more details see e.g.~\cite[\S 5.1]{stepet} and references therein. Let ${\mathcal{I}}$ be a finite set. We set $\Box_{\mathcal{I}}$ to be the category of subsets of ${\mathcal{I}}$ with inclusions as morphisms, i.e.\ if $I,J\subset {\mathcal{I}}$ are subsets, then $\Box_{\mathcal{I}}(I,J)$ consists of one element if $I\subset J$ and is empty otherwise. 
Recall that an {\it ${\mathcal{I}}$-cubical object} of a category $\mathcal{C}$ is a contravariant functor $\Box_{\mathcal{I}}\to\mathcal{C}$, or explicitly, a collection $(C_I),I\subset {\mathcal{I}}$ of objects of $\mathcal{C}$, and a collection $\varphi_{IJ}:C_J\to C_I, I\subset J\subset {\mathcal{I}}$ of morphisms of $\mathcal{C}$, which we will call the {\it structure morphisms}, such that $\varphi_{II}=\id_{C_I}$ for all $I\subset {\mathcal{I}}$ and $\varphi_{IJ}\circ\varphi_{JK}=\varphi_{IK}$ for all $I\subset J\subset K\subset {\mathcal{I}}$. 
We say that an $\mathcal{I}$-cubical object is {\it constant} if it takes every object of $\Box_\mathcal{I}$ to some fixed $C\in\mathcal{C}$ and every morphism to $\id_C$.
If $\mathcal{I}=\{1,\ldots,k\}$, then $\mathcal{I}$-cubical objects will be referred to as {$k$-\it cubical objects}.

A {\it morphism} of $\mathcal{I}$-cubical objects is a natural transformation of the corresponding functors. Explicitly, if $C_\bullet=(C_I,\varphi_{IJ})$ and $C'_\bullet=(C'_I,\varphi'_{IJ})$ are $\mathcal{I}$-cubical objects, then a morphism $C_\bullet\to C'_\bullet$ is given by a collection $f_\bullet=(f_I:C_I\to C_I'),I\subset\mathcal{I}$ of morphisms of $\mathcal{C}$ that commute with the structure morphisms, i.e. we have $f_I\circ \varphi_{IJ}=\varphi'_{IJ}\circ f_J$ for all $I\subset J\subset\mathcal{I}$. If $\mathcal{C}$ is {\it concrete}, i.e.\ if it has a faithful functor into sets, then after choosing such a functor 
we define {\it $\mathcal{I}$-cubical subobjects} as subfunctors; also, in this case we will sometimes say ``maps of $\mathcal{I}$-cubical objects'' instead of ``morphisms''. If $\mathcal{C}$ has limits or colimits of a given shape, so does the category of $\mathcal{I}$-cubical objects of $\mathcal{C}$, and they are computed object-wise.

In the sequel we will mainly be interested in the cases when $\mathcal{C}$ is the category of complex algebraic varieties or the category of topological spaces. Cubical subspaces and subvarieties, their intersections and preimages are defined in a straightforward way.

Let $Y_\bullet=(Y_I,\varphi_{IJ})$ be a cubical topological space, and let $Y^*$ be the join of all $Y_I$ such that $I\subset {\mathcal{I}}$ is {\it non-empty}. A simplex $\triangle\subset Y^*$ with vertices $y_1\in Y_{I_1},\ldots, y_l\in Y_{I_l}$ is {\it coherent} with {\it principal vertex} $y_l$ if $I_1\subsetneq\cdots\subsetneq I_l$ and $y_j=\varphi_{I_j I_l}(y_l)$ for all $j\in\{1,\ldots, l\}$. (We will shortly see how this terminology is related to the one in section~\ref{conical_sect}.) The {\it geometric realisation} of $Y_\bullet$, denoted $|Y_\bullet|$, is the union of all coherent simplices. The space $Y_\varnothing$ does not participate in this construction, but we have a natural {\it augmentation map} $\varepsilon:|Y_\bullet|\to Y_{\varnothing}$ induced by the maps $\varphi_{\varnothing I}$ for $\varnothing\subsetneq I\subset {\mathcal{I}}$. Note that geometric realisation commutes with taking intersections and preimages.

If $Z\subset Y_\varnothing$ is a subspace, the preimages $\varphi_{\varnothing I}^{-1}(Z),I\subset\mathcal{I}$ form a cubical subspace of $Y_\bullet$, which we will denote $\varepsilon^{-1}(Z)_\bullet$. If $Y_\bullet$ is a cubical algebraic variety and $Z\subset Y_\varnothing$ is a subvariety, then $\varepsilon^{-1}(Z)_\bullet$ is a cubical subvariety of~$Y_\bullet$.

A continuous map $f_\bullet:Y_\bullet\to Y'_\bullet$ of $\mathcal{I}$-cubical topological spaces induces a continuous map $|f_\bullet|:|Y_\bullet|\to |Y'_\bullet|$ of the geometric realisations. We note for the sequel that if all components $f_I$ of $f_\bullet$ are proper, so is $|f_\bullet|$.

\subsubsection{O.\ Tommasi's construction}\label{tom_construction}

We will now review O.~Tommasi's construction from \cite{tom}.

\begin{assumption}
In the rest of section~\ref{hodge_modge} $Y_1,\ldots, Y_\N$ will be a sequence of spaces that satisfies conditions \ref{first}-\ref{five}$+$ of list~\ref{conditions}.
\end{assumption}

O.~Tommasi defines in~\cite{tom} $\eu N$-cubical topological spaces $\tilde\Lambda_\bullet$ and $\tilde{\mathcal{Y}}_\bullet$ by setting
\begin{align}\label{cubical_1ver}
\begin{split}
\tilde\Lambda_I=&\{(K_{i_1},\ldots,K_{i_l})\in \bar Y_{i_1}\times\cdots\times \bar Y_{i_l}\mid K_{i_1}\subset\cdots\subset K_{i_l}\},\\
\tilde{\mathcal{Y}}_I=&\{(f,(K_{i_1},\ldots,K_{i_l}))\in V\times\tilde\Lambda_I\mid K_{i_l}\subset \Sing f\}
\end{split}
\end{align}
for $\varnothing\subsetneq I=\{i_1<\ldots <i_l\}\subset\{1,\ldots,\N\}$. The structure maps $\tilde\Lambda_J\to \tilde\Lambda_I$ and $\tilde{\mathcal{Y}}_J\to\tilde{\mathcal{Y}}_I$ for $\varnothing\subsetneq I\subset J\subset\{1,\ldots,\N\}$ are the natural projections. Finally, set $\tilde\Lambda_{\varnothing}$ to be a point and $\tilde{\mathcal{Y}}_\varnothing=\Sigma$; the map $\tilde{\mathcal{Y}}_I\to \tilde{\mathcal Y}_\varnothing, I\subset\{1,\ldots,\N\}$ is induced by the projection to $V$.

\begin{figure}[h]
\centering
\includegraphics[scale=.3]{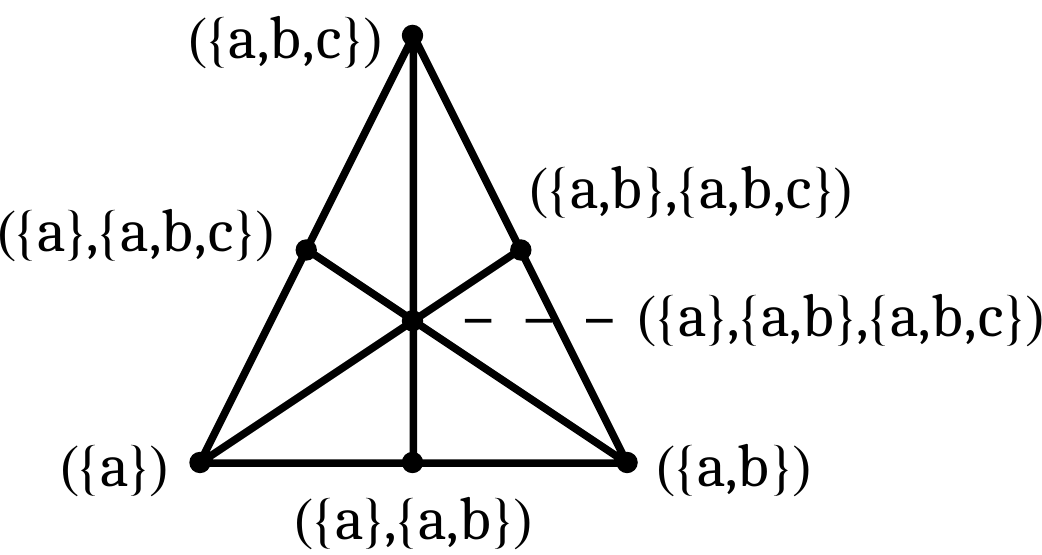}
\caption{The map $\beta$}
\label{fig:triangle_2}
\end{figure}

Every space $\tilde\Lambda_I$ is compact, and every $\tilde{\mathcal Y}_I$ is closed in $V\times \tilde\Lambda_I$, so all maps $\tilde{\mathcal{Y}}_I\to \tilde{\mathcal{Y}}_J,J\subset I$ are proper. In particular, every $\tilde{\mathcal{Y}}_I\to \tilde{\mathcal{Y}}_\varnothing$ is proper, which implies that so is the augmentation map $\varepsilon:|\tilde{\mathcal{Y}}_\bullet|\to\Sigma$. Moreover, for every $f\in\Sigma$ there is an $i\in\{1,\ldots,\N\}$ such that $\Sing f\in Y_i$, and the preimage $\varepsilon^{-1}(f)$ is deformation retracts onto $\{\Sing f\}\subset Y_i\subset\tilde\Lambda_{\{i\}}$, cf.\ figure~\ref{fig:triangle_2}. So $\varepsilon$ induces an isomorphism of Borel-Moore homology groups by lemma~\ref{proper_contr_fib}.

We now compare this construction with the conical resolution $\sigma$ from section~\ref{conical_sect}. 
Let $I=\{i_1<\ldots<i_l\}\subset\{1,\ldots,\N\}$ be a non-empty subset. Recall that above we set $\mathbf{Y}=\bigcup_{i=1}^\N Y_i$. The maps $(K_{i_1},\ldots, K_{i_l})\mapsto \frac{1}{l}\sum K_{i_j}$, the barycentre of the coherent simplex with vertices $K_{i_1},\ldots, K_{i_l}$, induce a continuous map $\Lambda_I\to\mathbf{Y}$. Extending by linearity, we get a continuous surjective map $\beta:|\tilde\Lambda_\bullet|\to\tilde\Lambda\subset\mathbf{Y}^{*\N}$, which takes the union of all coherent simplices in $|\tilde\Lambda_\bullet|$ with vertices equal to subsequences of $(K_{i_1},\ldots, K_{i_l})$ to the (barycentric subdivision of the) coherent simplex in $\tilde\Lambda$ with vertices $K_{i_1},\ldots, K_{i_l}$. See figure~\ref{fig:triangle_2} for an illustration for $K_{i_1}=\{a\}, K_{i_2}=\{a,b\}, K_{i_3}=\{a,b,c\}$; note that the part shown in this picture corresponds to the triangle with vertices $\{a\},\{a,b\}$ and $\{a,b,c\}$ in figure~\ref{fig:triangle_1}.

Note that 
we have a commutative diagram
\begin{equation}\label{cart_square_1}
\begin{tikzcd}
{{|}\tilde{\mathcal Y}}_{\bullet}{|} \ar[r,"\gamma"]\ar[d] &\sigma\ar[d,"p"]\\
{{|}\tilde{\Lambda}}_{\bullet}{|} \ar[r,"\beta"] & \tilde\Lambda
\end{tikzcd}
\end{equation}
where $\gamma$ is induced by
${{|}\tilde{\mathcal Y}}_{\bullet}{|}\to V\times{{|}\tilde{\Lambda}}_{\bullet}{|} \xrightarrow{\id_V\times \beta} V\times\tilde\Lambda$.

\begin{lemma}\label{lemma_cart}
Let
$$
\begin{tikzcd}
A\arrow[r]\arrow[d] & B\arrow[d]\\
C\arrow[r] & D
\end{tikzcd}
$$
be a commutative diagram in the topological category. Suppose the diagram is set-theoretically Cartesian (i.e.\ the map $A\to B\times_D C$ is bijective), one of the arrows that start at $A$ is proper and the spaces $A,B,C$ are locally compact and Hausdorff. Then the diagram is topologically Cartesian, i.e.\ $A\to B\times_D C$ is a homeomorphism.
\end{lemma}

$\clubsuit$

By inspection, (\ref{cart_square_1}) is set-theoretically Cartesian. 
Moreover, the following diagram commutes:
\begin{equation}\label{cub_spaces_diag1}
\begin{tikzcd}
{{|}\tilde{\mathcal Y}}_{\bullet}{|} \ar[r,"\gamma"]\ar[dr,"\varepsilon"] & \sigma\ar[d,"\pi"]\\
& \Sigma,
\end{tikzcd}
\end{equation}
so $\gamma$ is proper, as $\varepsilon$ is. We conclude using lemma~\ref{lemma_cart} that diagram~(\ref{cart_square_1}) is in fact topologically Cartesian.

The map $\beta$ is proper, its source and target being compact, so the topologies on $\sigma$ and $\tilde\Lambda$ coincide with those induced by the maps $\gamma$ and $\beta$ respectively.

O.~Tommasi also introduces cubical spaces $\Lambda_\bullet$ and ${\mathcal Y}_\bullet$, which are defined in the same way as $\tilde\Lambda_\bullet$ and $\tilde{\mathcal Y}_\bullet$ above, except that each closure $\bar Y_i$ is replaced by $Y_i$. The maps $\gamma$ and $\beta$ induce bijections (which are in general not homeomorphisms) $|\mathcal{Y}_\bullet|\to\sigma$ and $|\Lambda_\bullet|\to\tilde\Lambda$, and O.~Tommasi equips $|\mathcal{Y}_\bullet|$ and $|\Lambda_\bullet|$ with the topologies induced by these bijections.


Moreover, the spaces $|\mathcal{Y}_\bullet|$ and $|\Lambda_\bullet|$ have natural increasing filtrations $(\mathrm{Fil}_i)$ defined as follows: we set $\mathrm{Fil}_i |\Lambda_\bullet|$ to be the image of $|\Lambda^{\leq i}_\bullet|$ where $\Lambda_\bullet^{\leq i}$ denotes the restriction of the functor $\Lambda_\bullet$ from $\Box_{\{1,\ldots,\N\}}$ to $\Box_{\{1,\ldots,i\}}$, and similarly for $|\mathcal{Y}_\bullet|$. By inspection, the map $\beta$, respectively $\gamma$ identifies the resulting filtration with the filtration $(\Phi_i)$, respectively $(F_i)$ constructed in section~\ref{conical_sect}. Let us summarise these results.

\begin{Prop}\label{homeo_quint_tom}
Suppose $Y_1,\ldots, Y_\N$ is a sequence of spaces that satisfies conditions \ref{first}-\ref{five}$+$ of list~\ref{conditions}. Then the spaces $|\Lambda_\bullet|$ and $|\mathcal{Y}_\bullet|$ are filtered homeomorphic to $\tilde\Lambda$ and $\sigma$ respectively.
\end{Prop}

$\clubsuit$

\subsubsection{Algebraic configuration spaces}

The spaces $|\tilde{\mathcal{Y}}_\bullet|$ and $|\mathcal{Y}_\bullet|$ are still not quite what we are after since none of them is the geometric realisation of a cubical complex algebraic variety even if we assume that all $Y_i$ are complex algebraic: the spaces $\bar Y_i$ in general have no complex algebraic structure, see p.~\pageref{nocomplexstr}, and in the construction of $|\mathcal{Y}_\bullet|$ we replaced the natural topology on it by another one. One way to handle this problem is to ``blow up'' all $\bar Y_i$ to complex algebraic varieties in a compatible way. More precisely, in the rest of section~\ref{hodge_modge} we will assume the following.

\begin{cl}\label{condition_list_3}
There exist compact complex algebraic varieties $\bar{\pmb{Y}}_1,\ldots,\bar{\pmb{Y}}_\N$ and open dense subvarieties $\pmb{Y}_i\subset\bar{\pmb{Y}}_i,i=1,\ldots,\N$ such that
\begin{enumerate}
\item\label{cond_1_list_2} For every $i\in\{1,\ldots,\N\}$ there is a continuous map $h_i:\bar{\pmb{Y}}_i\to 2^{\mathbf{M}}$ that takes $\pmb{Y}_i$ homeomorphically to $Y_i$ and $\bar{\pmb{Y}}_i\setminus\pmb{Y}_i$ to $\bar Y_i\setminus Y_i$.
\item\label{cond_2_list_2} For every $i,j\in\{1,\ldots,\N\}$, the subset $\{(y_i,y_j)\in \bar{\pmb{Y}}_i\times \bar{\pmb{Y}}_j\mid h_i(y_i)\subset h_j(y_j)\}\subset \bar{\pmb{Y}}_i\times \bar{\pmb{Y}}_j$ is Zariski closed, and so is $\{(x,y)\in \mathbf{M}\times \bar{\pmb{Y}}_i\mid x\in h_i(y)\}\subset \mathbf{M}\times \bar{\pmb{Y}}_i$.
\item\label{cond_3_list_2} Recall that the vector bundle $(F_i\setminus F_{i-1},\Phi_i\setminus\Phi_{i-1},p|_{F_i\setminus F_{i-1}})$ from part 3 of Theorem~\ref{thconical} is the pullback of the tautological bundle over the Grassmannian $\Gr(d_i,V)$ along the map $\tilde g_i$. Recall also that $\tilde g_i$ restricted to $Y_i\subset (\Phi_i\setminus\Phi_{i-1})$ is given by $K\mapsto L(K)$.

We require that for every $i\in\{1,\ldots,\N\}$ the map $\pmb{Y}_i\to \Gr(d_i,V)$ given by $y\mapsto L(h_i(y))$ should be the restriction of a regular map $\bar{\pmb{Y}}_i\to\Gr(d_i,V)$.

Note that this implies that the pullback of $(F_i\setminus F_{i-1},\Phi_i\setminus\Phi_{i-1},p|_{F_i\setminus F_{i-1}})$ along $h_i$ extends to an algebraic vector bundle over $\bar{\pmb{Y}}_i$, which we will denote $\mathcal{E}^i$\label{bundle_ei} and which is moreover a subbundle of the trivial bundle $V\times \bar{\pmb{Y}}_i$.
\item\label{cond_4_list_2} For every $i,j\in\{1,\ldots,\N\}$ the preimage 
$h_i^{-1}(\bar Y_j)$ is Zariski closed in $\bar{\pmb{Y}}_i$.
\end{enumerate}
\end{cl}

{\bf Remark.} In most examples we have $\bar{\pmb{Y}}_1=\bar Y_1=Y_1=\mathbf{M}$, in which case the second part of condition~\ref{cond_2_list_2} is automatic.

{\bf Remark.} If $i<j$, then it may happen that there are $y_i\in\bar{\pmb{Y}}_i, y_j\in\bar{\pmb{Y}}_j$ such that $h_i(y_i)\supsetneq h_j(y_j)$. For example, this will be the case when $Y_i=B(\mathbf{M},2), Y_j=B(\mathbf{M},3)$ and $\dim\mathbf{M}>0$, as $\bar Y_j$ will then contain one-point configurations.


Let us check that these assumptions (i.e.\ the conditions in list \ref{condition_list_3}) are realistic. Suppose that $\mathbf{M}$ is projective and fix an embedding $\mathbf{M}\subset\C\p^n$. Recall that in section~\ref{conf_spaces_sect} we defined the cycle map $Z^C$ from the Chow variety $\bar{\eu{C}}_{k,l,n}^{\mathbf{M}}$ to $2^{\mathbf{M}}$. Recall that $Z^C$ is continuous (proposition~\ref{cycle_maps}) and takes the open Chow variety ${\eu{C}}_{k,l,n}^{\mathbf{M}}$ homeomorphically to a subspace of $2^\mathbf{M}$ which we denoted $\eu{B}_{k,l}^{\mathbf{M}}$ (corollary~\ref{Bkl}).

\begin{Prop}\label{realistic}
Suppose that for every $i\in\{1,\ldots,\N\}$ there exist a product $\bar{\eu{P}}$ of the Chow varieties $\bar{\eu{C}}_{k_j,l_j,n}^\mathbf{M}, j=1,\ldots, m$ such that $k_1<\cdots<k_m$ (i.e.\ the sequence of the dimensions of the corresponding subvarieties of $\mathbf{M}$ is strictly increasing), and a locally closed subvariety $\pmb{Z}_i$ of
$$\eu{U}=\{(x_1,\ldots,x_m)\in\bar{\eu{P}}\mid \mbox{{\normalfont each $x_j\in \eu{C}_{k_j,l_j,n}^{\mathbf{M}}$; no component of $x_{j_1}$ is contained in $x_{j_2}$ for $j_1<j_2$}}\}$$
such that $Y_i=u_i(\pmb{Z}_i)$ where $u_i:\bar{\eu{P}}\to 2^{\mathbf{M}}$ is given by $(x_1,\ldots, x_m)\mapsto \bigcup_{j=1}^m Z^C(x_j)$.

%
%
Then there exist varieties $\pmb{Y}_i, \bar{\pmb{Y}}_i$ and continuous maps $h_i:\bar{\pmb{Y}}_i\to 2^{\mathbf{M}},i=1,\ldots, N$ that have the properties given in condition list~\ref{condition_list_3}.
\end{Prop}

{\bf Proof.} Let us fix an $i\in\{1,\ldots,\N\}$. We first construct $\pmb{Y}_i, \bar{\pmb{Y}}_i$ and $h_i:\bar{\pmb{Y}}_i\to 2^{\mathbf{M}},i=1,\ldots, N$ such that condition~\ref{cond_3_list_2} is satisfied. We set
$$\pmb{Y}_i=\{(W,x)\in\Gr(d_i,V)\times\pmb{Z}_i\mid \mbox{every $f\in W$ is $0$ on $u_i(x)$}\}.$$
This is a Zariski closed subset of $\Gr(d_i,V)\times\pmb{Z}_i$, and we take $\bar{\pmb{Y}}_i$ to be the closure of $\pmb{Y}_i$ in $\Gr(d_i,V)\times\bar{\pmb{Z}}_i$ where $\bar{\pmb{Z}}_i$ is the closure of $\pmb{Z}_i$ in $\bar{\eu{P}}$. Set $h_i$ to be the restriction of $u_i\circ pr_{\bar{\pmb{Z}}_i}$ to $\bar{\pmb{Y}}_i$. The map $u_i$ is continuous by the continuity of $Z^C$ and part 3 of proposition~\ref{conf_space}, so $h_i$ is continuous.

Note that $\pmb{Y}_i$ is the graph of the map $\pmb{Z}_i\to \Gr(d_i,V)$ given by $x\mapsto L(u_i(x))$; note that this map is continuous by the continuity of $u_i$ and proposition~\ref{semicont}. So for every $y=(W,x)\in\pmb{Y}_i$ we have $L(h_i(y))=L(u_i(x))=W$, which implies that on $\pmb{Y}_i$ the map $y\mapsto L(h_i(y))$ coincides with $pr_{\Gr(d_i,V)}$, which is clearly regular and extends to a regular map $\bar{\pmb{Y}}_i\to\Gr(d_i,V)$. We have checked condition~\ref{cond_3_list_2} for $\pmb{Y}_i,\bar{\pmb{Y}}_i$ and $h_i:\bar{\pmb{Y}}_i\to 2^\mathbf{M}$.

Let us now do the same for condition~\ref{cond_1_list_2}. Since $\pmb{Y}_i$ is closed in $\Gr(d_i,V)\times \pmb{Z}_i$, we have $\bar{\pmb{Y}}_i\cap (\Gr(d_i,V)\times \pmb{Z}_i)=\pmb{Y}_i$, so the projection $pr_{\bar{\pmb{Z}}_i}:\Gr(d_i,V)\times\bar{\pmb{Z}}_i\to\bar{\pmb{Z}}_i$ takes $\bar{\pmb{Y}}_i\setminus\pmb{Y}_i$ to $\bar{\pmb{Z}}_i\setminus\pmb{Z}_i$. Note that $pr_{\bar{\pmb{Z}}_i}$ also induces a regular map of complex algebraic varieties $\pmb{Y}_i\to\pmb{Z}_i$ that is a homeomorphism in the complex analytic topology, the inverse (not necessarily regular) homeomorphism being given by $x\mapsto (L(u_i(x)),x)$.

Let us show that the map $u_i$ takes $\pmb{Z}_i$ and $\bar{\pmb{Z}}_i\setminus\pmb{Z}_i$ to disjoint subsets of $2^{\mathbf{M}}$. Recall that the cycle map $Z^C$ induces a bijection between the open Chow variety ${\eu{C}}_{k,l,n}^\mathbf{M}$ and $\eu{B}_{k,l}^{\mathbf{M}}$, and takes the boundary $\bar{\eu{C}}_{k,l,n}^\mathbf{M}\setminus {\eu{C}}_{k,l,n}^\mathbf{M}$ to the boundary $\bar{\eu{B}}_{k,l}^{\mathbf{M}}\setminus\eu{B}_{k,l}^{\mathbf{M}}$. Similarly, by the definition of $\eu{U}$, the map $u_i$ a.\ is injective on $\eu{U}$, and b.\ takes $\eu{U}$ and $\bar{\eu{U}}\setminus\eu{U}$ to disjoint subsets of $2^{\mathbf{M}}$. Let now $x\in \pmb{Z}_i,x'\in\bar{\pmb{Z}}_i\setminus\pmb{Z}_i$ be elements. If $x'\in\eu{U}$, we have $h_i(x)\neq h_i(x')$ using assertion a., and if $x'\in\bar{\eu{U}}\setminus\eu{U}$, we get the same result using assertion b. So we get $u_i(\pmb{Z}_i)\cap u_i(\bar{\pmb{Z}}_i\setminus\pmb{Z}_i)=\varnothing$. Note also that $u_i$ induces a bijection $\pmb{Z}_i\to Y_i$: surjectivity follows from our assumptions on $\pmb{Z}_i$, and injectivity from $\pmb{Z}_i\subset\eu{U}$.

Combining the observations from the last two paragraphs we see that $h_i:\bar{\pmb{Y}}_i\to 2^{\mathbf{M}}$ takes $\pmb{Y}_i$ bijectively to $Y_i$, and $\bar{\pmb{Y}}_i\setminus\pmb{Y}_i$ to $\bar Y_i\setminus Y_i$. 
Using lemma~\ref{gen_topo} we see that $h_i$ induces a homeomorphism $\pmb{Y}_i\to Y_i$. This proves that $\pmb{Y}_i,\bar{\pmb{Y}}_i$ and $h_i:\bar{\pmb{Y}}_i\to 2^\mathbf{M}$ satisfy condition~\ref{cond_1_list_2}.

Checking conditions 2 and 4 is straightforward.$\clubsuit$

{\bf Remark.} The map $\pmb{Y}_i\to \pmb{Z}_i$ from the proof of proposition~\ref{realistic} is a homeomorphism with respect to the complex analytic topology, so the corresponding morphism of schemes over $\C$ is a universal homeomorphism, 
which implies that there is a morphism from the absolute weak normalisation $\pmb{Z}_i^{awn}$ (\cite[tag 0EUK]{stacks}) of $\pmb{Z}_i$ to $\pmb{Y}_i$ that makes the following diagram commutative:
$$
\begin{tikzcd}
\pmb{Z}_i^{awn}\ar[dr]\ar[d] &\\
\pmb{Y}_i\ar[r]&\pmb{Z}_i
\end{tikzcd}
$$
where the diagonal arrow is the canonical morphism. Note also that for schemes over $\C$ absolute weak normalisations coincide with semi-normalisations (see ibid.)
\smallskip



%
%

\subsubsection{Cubical varieties constructed from algebraic configuration spaces}

As above (see (\ref{cubical_1ver})) we now define two $\N$-cubical complex algebraic varieties $\tilde{\pmb{\Lambda}}_\bullet$ and $\tilde{\pmb{\mathcal{Y}}}_\bullet$ by setting for $\varnothing\subsetneq I=\{i_1<\cdots<i_l\}\subset\{1,\ldots,\N\}$
\begin{align}\label{cubical_2ver}
\begin{split}
\tilde{\pmb{\Lambda}}_I=&\overline{\{(y_{i_1},\ldots,y_{i_l})\in \pmb{Y}_{i_1}\times\cdots\times \pmb{Y}_{i_l}\mid h_{i_1}(y_{i_1})\subset\cdots\subset h_{i_l}(y_{i_l})\}},\\
\tilde{\pmb{\mathcal{Y}}}_I=&\{(f,(y_{i_1},\ldots,y_{i_l}))\in V\times\tilde{\pmb{\Lambda}}_I\mid h_{i_l}(y_{i_l})\subset \Sing f\}
\end{split}
\end{align}
where the closure in the definition of $\tilde{\pmb{\Lambda}}_\bullet$ is taken in $\bar{\pmb{Y}}_{i_1}\times\cdots\times \bar{\pmb{Y}}_{i_l}$. As in the case of $\tilde\Lambda_\bullet$ and $\tilde{\mathcal{Y}}_\bullet$ above, we set $\tilde{\pmb{\Lambda}}_\varnothing=$ a point, $\tilde{\pmb{\mathcal{Y}}}_\varnothing=\Sigma$, and the structure maps of both $\tilde{\pmb{\Lambda}}_\bullet$ and $\tilde{\pmb{\mathcal{Y}}}_\bullet$ are defined to be the natural projections. All structure maps of $\tilde{\pmb{\Lambda}}_\bullet$ are proper.

\begin{lemma}\label{lem_prop}
Every $\tilde{\pmb{\mathcal{Y}}}_I,I=\{i_1<\cdots<i_l\}\subset\{1,\ldots,\N\}$ is Zariski closed in $V\times \tilde{\pmb{\Lambda}}_I$. All structure maps of $\tilde{\pmb{\mathcal{Y}}}_\bullet$ are proper, and so is the augmentation map $\pmb{\varepsilon}:|\tilde{\pmb{\mathcal{Y}}}_\bullet|\to\Sigma$.
\end{lemma}

{\bf Proof.} We will prove the first assertion of the lemma, and the rest will then follow. The subset
$$\{(f,x,y)\in V\times\mathbf{M}\times\bar{\pmb{Y}}_{i_l}\vert x\in\Sing f; x\in h_{i_l}(y)\}\subset V\times\mathbf{M}\times\bar{\pmb{Y}}_{i_l}$$
is Zariski closed since so are $\{(x,y)\in \mathbf{M}\times \bar{\pmb{Y}}_{i_l}\vert x\in h_{i_l}(y)\}$ in $\mathbf{M}\times \bar{\pmb{Y}}_{i_l}$ by condition~\ref{cond_2_list_2} of list~\ref{condition_list_3} and $\tot(\mathcal{E})=\{(f,x)\in V\times\mathbf{M}\vert x\in\Sing f\}$ in $V\times\mathbf{M}$. So $\{(f,y)\in V\times \bar{\pmb{Y}}_{i_l}\vert h_{i_l}(y)\subset\Sing f\}$ is Zariski closed in $V\times \bar{\pmb{Y}}_{i_l}$, which implies that $\tilde{\pmb{\mathcal{Y}}}_I$ is Zariski closed in $V\times \tilde{\pmb{\Lambda}}_I$.$\clubsuit$

{\bf Remark.} It follows from condition~\ref{cond_2_list_2} of list~\ref{condition_list_3} that if $(y_{i_1},\ldots,y_{i_l})\in\tilde{\pmb{\Lambda}}_I$, then $h_{i_1}(y_{i_1})\subset\cdots\subset h_{i_l}(y_{i_l})$.

Using the maps $(y_{i_1},\ldots,y_{i_l})\mapsto\frac{1}{l}\sum h_{i_j}(y_{i_j})$ we get a continuous surjective map $\pmb{\beta}:|\tilde{\pmb{\Lambda}}_\bullet|\to\tilde\Lambda$, which takes the union of all coherent simplices in $|\tilde{\pmb{\Lambda}}_\bullet|$ with vertices equal to subsequences of $(y_{i_1},\ldots,y_{i_l})$ to the barycentric subdivision of the coherent simplex in $\tilde\Lambda$ with vertices $(h_{i_1}(y_{i_1}),\ldots,h_{i_l}(y_{i_l}))$ (cf.\ figure~\ref{fig:triangle_2}). There is a commutative diagram (cf.\ (\ref{cart_square_1}))
\begin{equation}\label{cart_square_2}
\begin{tikzcd}
{{|}\tilde{\pmb{\mathcal Y}}}_{\bullet}{|} \ar[r,"\pmb{\gamma}"]\ar[d] &\sigma\ar[d,"p"]\\
{{|}\tilde{\pmb{\Lambda}}}_{\bullet}{|} \ar[r,"\pmb{\beta}"] & \tilde\Lambda
\end{tikzcd}
\end{equation}
where $\pmb{\gamma}$ is induced by ${{|}\tilde{\pmb{\mathcal Y}}}_{\bullet}{|}\to V\times{{|}\tilde{\pmb{\Lambda}}}_{\bullet}{|} \xrightarrow{\id_V\times \pmb{\beta}} V\times\tilde\Lambda$. 

We can extend diagram~(\ref{cub_spaces_diag1}) as follows:
\begin{equation}\label{cub_spaces_diag2}
\begin{tikzcd}
{{|}\tilde{\pmb{\mathcal Y}}}_{\bullet}{|} \ar[r]\ar[drr,"\pmb{\varepsilon}"]\ar[rr,bend left=30,"\pmb{\gamma}"] & {{|}\tilde{\mathcal Y}}_{\bullet}{|} \ar[r,"\gamma"]\ar[dr,"\varepsilon"] & \sigma\ar[d,"\pi"]\\
& & \Sigma.
\end{tikzcd}
\end{equation}
Here the map ${{|}\tilde{\pmb{\mathcal Y}}}_{\bullet}{|}\to {{|}\tilde{\mathcal Y}}_{\bullet}{|}$ is induced by the map ${\tilde{\pmb{\mathcal Y}}}_{\bullet}\to {\tilde{\mathcal Y}}_{\bullet}$ of cubical spaces constructed using $h_i:\bar{\pmb{Y}}_i\to \bar Y_i$. The map $\pmb{\gamma}$ is proper and~(\ref{cart_square_2}) is topologically Cartesian, cf.\ the proof of these properties for the map $\gamma$ and diagram~(\ref{cart_square_1}) in section~\ref{tom_construction} above.

\subsubsection{The map \texorpdfstring{$\pmb{\gamma}$}{} has contractible fibres}

\begin{lemma}\label{fib_contractible}
All fibres of $\pmb{\gamma}$ are contractible.
\end{lemma}

{\bf Proof.} Since~(\ref{cart_square_2}) is Cartesian, it suffices to prove that the fibres of $\pmb{\beta}$ are contractible. In order to do this, it will be convenient to use an alternative description of $|\tilde{\pmb{\Lambda}}_\bullet|$. Namely, we will now realise this space as a subspace of the join $\pmb{Y}^*=\Asterisk_{i=1}^\N \bar{\pmb{Y}}_i$ rather than $\Asterisk_{I\neq\varnothing}\tilde{\pmb{\Lambda}}_I$, cf.\ the proof of proposition~\ref{homeo_quint_tom} above, in particular figure~ \ref{fig:triangle_2}.

We call a simplex $\triangle\subset\pmb{Y}^*$ with vertices $y_{i_1}\in\bar{\pmb{Y}}_{i_1},\ldots, y_{i_l}\in\bar{\pmb{Y}}_{i_l}, 1\leq i_1<\cdots <i_l\leq\N$ {\it coherent} with {\it principal vertex} $y_{i_l}$ if $h_{i_1}(y_{i_1})\subset\cdots \subset h_{i_l}(y_{i_l})$, cf.\ similar definitions for the spaces $\tilde\Lambda, |\tilde\Lambda|_\bullet$ and $|\tilde{\pmb{\Lambda}}|_\bullet$. 
We define a map from $|\tilde{\pmb{\Lambda}}_\bullet|$ to the union of all coherent simplices in $\pmb{Y}^*$ by taking a sequence $(y_{i_1},\ldots, y_{i_l})\in\tilde{\pmb{\Lambda}}_{\{i_1,\ldots,i_l\}}$ to the barycentre of the coherent simplex in $\pmb{Y}^*$ with vertices $y_{i_1},\ldots, y_{i_l}$, and after that extending by linearity. The resulting map is injective, cf.\ figure~\ref{fig:triangle_2}, and we denote the image by $\pmb{L}$. We can identify $\pmb{\beta}:|\tilde{\pmb{\Lambda}}_\bullet|\to\tilde\Lambda$ with a map $\pmb{B}:\pmb{L}\to \tilde\Lambda$. Observe that unlike $\pmb{\beta}$, the map $\pmb{B}$ takes principal vertices of coherent simplices to principal vertices.

For an $x\in\tilde\Lambda$ the preimage $\pmb{B}^{-1}(x)$ can be described as follows. Let $K_1\subsetneq\cdots\subsetneq K_l$ be the vertices of the coherent simplex $\triangle\subset\tilde\Lambda$ such that $x\in\stackrel{\circ}{\triangle}$. Let $t_1,\ldots, t_l$ be the barycentric coordinates of $x$ in $\triangle$. For each $j\in\{1,\ldots, l\}$ there is a unique $i(j)\in\{1,\ldots,\N\}$ such that $K_j\in Y_{i(j)}$, but $K_j$ may also belong to one or several spaces $\bar Y_{i'}\setminus Y_{i'},i'>i(j)$. Let $\triangle'\subset \pmb{L}$ be a coherent simplex with vertices $y_{i_1}\in\bar{\pmb{Y}}_{i_1},\ldots, y_{i_m}\in\bar{\pmb{Y}}_{i_m}$. We say a vertex $y_{i_k}$ of $\triangle'$ is {\it inner} if $y_{i_k}\in\pmb{Y}_{i_k}$ and {\it boundary} if $y_{i_k}\in\bar{\pmb{Y}}_{i_k}\setminus\pmb{Y}_{i_k}$. Moreover, $\triangle'$ is {\it $x$-relevant} if
$$\{h_{i_1}(y_{i_1}),\ldots, h_{i_m}(y_{i_m})\}=\{K_1,\ldots, K_l\}.$$

The preimage $\pmb{B}^{-1}(x)$ is then
$$\bigcup_{\begin{array}{c}\triangle'\mbox{ }x\mbox{-relevant}\\ \mbox{with vertices }\\ y_{i_1},\ldots,y_{i_m}\end{array}}\left\{\sum_{k=1}^m s_k y_{i_k}\middle\vert \mbox{ all } s_k\geq 0; \sum_{k=1}^m s_k=1;\sum_{h_{i_k}(y_{i_k})= K_j} s_k=t_j\mbox{ for all }j=1,\ldots, l\right\}.$$

There does not seem to be a straightforward way to contract $\pmb{B}^{-1}(x)$ in one go: since the maps $h_i:\bar{\pmb{Y}}_i\to\bar Y_i$ need not be injective, it may happen that there are two $x$-relevant coherent simplices $\triangle'$ and $\triangle''$ and vertices $y'\in \triangle'$ and $y''\in \triangle''$ such that $y',y''\in\bar{\pmb{Y}}_i, y'\neq y''$ and $h_i(y')=h_i(y'')$. Then there are no coherent simplices that contain both $y'$ and $y''$.

So we proceed in several steps. Set $B_0=\pmb{B}^{-1}(x)$ and let $B_k,l\geq k\geq 1$ be the intersection of $\pmb{B}^{-1}(x)$ with the union of all $x$-relevant coherent simplices $\triangle'\subset\pmb{L}$ such that all vertices of $\triangle'$ that are mapped to $K_j,j\leq k$ are inner; note that if $\triangle'$ has this property, so does every face $\triangle''\subset\triangle'$. By condition~\ref{four} of list~\ref{conditions} and the fact that every map $h_k$ is assumed injective on $\pmb{Y}_k$ (condition~1 of list~\ref{condition_list_3}), the space $B_l$ is a point. So in order to prove lemma~\ref{fib_contractible} it would suffice to prove the following result.

\begin{lemma}\label{def_retr_bk}
Suppose $l>k\geq 0$. Then the following holds.
\begin{enumerate}
\item Each of the spaces $\bar{\pmb{Y}}_{i(j)},j=1,\ldots, k+1$ contains precisely one point that can be a vertex of a coherent simplex $\subset B_k$, namely the point $y_{i(j)}\in\pmb{Y}_{i(j)}$ that the map $h_{i(j)}$ takes to $K_j$.
\item The space $B_{k}$ deformation retracts onto $B_{k+1}$.
\end{enumerate}
\end{lemma}

{\bf Proof of lemma~\ref{def_retr_bk}.} Let us prove the first part of the lemma. Suppose $j\in\{1,\ldots,k+1\}$ and $y'_{i(j)}\in \bar{\pmb{Y}}_{i(j)},y'_{i(j)}\neq y_{i(j)}$ is a vertex of a coherent simplex $\triangle'\subset B_k$. Then $h_{i(j)}(y'_{i(j)})$ is a configuration $K_{j'}\neq K_j$, by our assumptions on $h_{i(j)}$, see condition~\ref{cond_1_list_2} of list~\ref{condition_list_3}. In particular, we have either $K_{j'}\subsetneq K_j$ or $K_j\subsetneq K_{j'}$. If $y'_{i(j)}\in {\pmb{Y}}_{i(j)}$, then $K_{j'}=h_{i(j)}(y'_{i(j)})\in Y_{i(j)}$, which is impossible as none of any two elements of $Y_{i(j)}$ can strictly contain the other by condition~\ref{sec} of list~\ref{conditions}.

Suppose now $y'_{i(j)}\in \bar{\pmb{Y}}_{i(j)}\setminus\pmb{Y}_{i(j)}$. Then $K_{j'}\in Y_{i'}$ for some $i'<i(j)$ by condition~\ref{five}$+$ from list~\ref{conditions} and condition~\ref{cond_1_list_2} of list~\ref{condition_list_3}. We have $K_{j'}\subsetneq K_j$, again using condition~\ref{sec} of list~\ref{conditions}. If $j=1$, we arrive at a contradiction, and the first part of the lemma follows. Otherwise, $K_{j'}$ is one of the configurations $K_1,\ldots, K_{j-1}$, i.e.\ $j'\leq j-1\leq k$. We conclude that the boundary vertex $y'_{i(j)}$ of the simplex $\triangle'$ maps to $K_{j'}$ with $j'\leq k$, which contradicts the definition of $B_k$. This completes the proof of the first part of lemma~\ref{def_retr_bk}.

We will now prove the second part. Suppose first that $1\leq k<l-1$. Consider the following three subspaces of $\pmb{Y}^*=\Asterisk_{i=1}^\N \bar{\pmb{Y}}_i$:
\begin{enumerate}
\item the simplex $\mathcal{A}_1$ with vertices $y_{i(j)},j=1,\ldots, k$;
\item the union $\mathcal{A}_2$ of coherent simplices $\triangle'$ such that all vertices of $\triangle'$ are mapped to $K_{k+1}$, one of the vertices is $y_{i(k+1)}\in \pmb{Y}_{i(k+1)}$ and all other vertices belong to $\bar{\pmb{Y}}_{i'},i'\not\in\{i(1),\ldots,i(k+1)\}$;
\item the union $\mathcal{A}_3$ of coherent simplices $\triangle''$ such that all vertices of $\triangle''$ are mapped to $K_{j},j>k+1$ and belong to $\bar{\pmb{Y}}_{i'},i'\not\in\{i(1),\ldots,i(k+1)\}$.
\end{enumerate}

By part 1 of the lemma, the vertices of each $x$-relevant coherent simplex $\triangle'\subset B_k$ are $y_{i(1)},\ldots, y_{i(k+1)}$, plus possibly some number of points $y_{i'}\in \bar{\pmb{Y}}_{i'},i'\not\in\{ i(1),\ldots, i(k+1)\}$ that are mapped to $K_{k+1}$, plus some points $y_{i'}\in \bar{\pmb{Y}}_{i'},i'\not\in\{i(1),\ldots, i(k+1)\}$ that are mapped to $K_j$ with $j>k+1$. So we can identify $B_k$ with a subspace of the join $\mathcal{A}_1 * \mathcal{A}_2 * \mathcal{A}_3$.
The space $\mathcal{A}_2$ deformation retracts onto $\{y_{i(k+1)}\}$. Let $H$ be the linear homotopy between the identity of $\mathcal{A}_2$ and the map $\mathcal{A}_2\to \{y_{i(k+1)}\}$. The join of $H$ and the identities of $\mathcal{A}_1$ and $\mathcal{A}_3$ induces a homotopy $B_k\times[0,1]\to B_k$ between the identity and a retraction $B_k\to B_{k+1}$. (During the homotopy, the vertices of each $x$-relevant coherent simplex $\triangle'\subset B_k$ that go to $K_{k+1}$ move linearly towards $y_{i(k+1)}$, and all other vertices stand still.) This proves the second part of lemma~\ref{def_retr_bk} in the case $1\leq k<l-1$.

The cases $k=0$ and $k=l-1$ are similar. The only difference is that if $k=0$, then the space $\mathcal{A}_1$ is absent from the join, and if $k=l-1$, then $\mathcal{A}_3$ is.$\clubsuit$

Lemma~\ref{fib_contractible} is now proved.$\clubsuit$

\subsubsection{The mixed Hodge structure on the main spectral sequence}\label{sec_main_spec_sec}

By construction we have $\pmb{\beta}^{-1}(\Phi_i)=\left|\tilde{\pmb{\Lambda}}_\bullet^{(i)}\right|$ and $\pmb{\gamma}^{-1}(F_i)=\left|\tilde{\pmb{\mathcal{Y}}}_\bullet^{(i)}\right|$ where $\tilde{\pmb{\Lambda}}_\bullet^{(i)}\subset \tilde{\pmb{\Lambda}}_\bullet$ and $\tilde{\pmb{\mathcal{Y}}}_\bullet^{(i)}\subset \tilde{\pmb{\mathcal{Y}}}_\bullet$ are cubical complex algebraic subvarieties defined by setting $\tilde{\pmb{\Lambda}}_\varnothing^{(i)}=pt, \tilde{\pmb{\mathcal{Y}}}^{(i)}_\varnothing=\Sigma$ and
\begin{align*}
\tilde{\pmb{\Lambda}}^{(i)}_I=&\left\{(y_{i_1},\ldots, y_{i_l})\in\tilde{\pmb{\Lambda}}_I\middle\vert h_{i_l}(y_{i_l})\in\bigcup_{j\leq i} Y_j\right\},\\
\tilde{\pmb{\mathcal{Y}}}^{(i)}_I=&\left\{(f,(y_{i_1},\ldots,y_{i_l}))\in V\times\tilde{\pmb{\Lambda}}^{(i)}_I\middle\vert h_{i_l}(y_{i_l})\subset\Sing(f)\right\}
\end{align*}
for $\varnothing\neq I=\{i_1<\cdots<i_l\}\subset\{1,\ldots,\N\}$. 
Similarly, the preimage $\pmb{\beta}^{-1}(\tilde\Lambda(K)), K\in Y_i$ is the geometric realisation of the cubical subvariety $\tilde{\pmb{\Lambda}}(K)_\bullet\subset \tilde{\pmb{\Lambda}}_\bullet$ with $\tilde{\pmb{\Lambda}}(K)_\varnothing=pt$ and
$$\tilde{\pmb{\Lambda}}(K)_I=\left\{(y_{i_1},\ldots, y_{i_l})\in\tilde{\pmb{\Lambda}}_I\middle\vert
h_{i_l}(y_{i_l})\subset K\right\}$$ for $\varnothing\neq I=\{i_1<\cdots<i_l\}\subset\{1,\ldots,\N\}$.

Let us generalise the last remark. Suppose $\pmb{Y}\subset\bar{\pmb{Y}}_i$ is a closed subvariety and set $Y=h_i(\pmb{Y})$. Then $$\pmb{\beta}^{-1}\left(\bigcup_{K'\in Y} \tilde\Lambda(K')\right)=|\tilde\Lambda(\pmb{Y})_\bullet|$$ where $\tilde{\pmb{\Lambda}}(\pmb{Y})_\bullet\subset \tilde{\pmb{\Lambda}}_\bullet$ is the cubical subvariety given by $\tilde{\pmb{\Lambda}}(\pmb{Y})_\varnothing=pt$ and $$\tilde{\pmb{\Lambda}}(\pmb{Y})_I=\left\{(y_{i_1},\ldots, y_{i_l})\in\tilde{\pmb{\Lambda}}_I\middle\vert
h_{i_l}(y_{i_l})\subset K\mbox{ \normalfont{for some }}K\in Y\right\}$$ for $\varnothing\neq I=\{i_1<\cdots<i_l\}\subset\{1,\ldots,\N\}$. (By condition~\ref{cond_2_list_2} of condition list~\ref{condition_list_3}, $\tilde{\pmb{\Lambda}}(\pmb{Y})_I\subset \tilde{\pmb{\Lambda}}_I$ is a closed subvariety.)

Finally, suppose we are given a sequence $X_1,\ldots, X_N$ of configuration spaces and a map $\alpha:\{1,\ldots, N\}\to \{1,\ldots,\N\}$ such that the assumptions on page~\pageref{assumptions_conical} are satisfied for $X_1,\ldots, X_N, Y_1, \ldots, Y_\N$ and $\alpha$. If $K\in X_i=Y_{\alpha(i)}$, then the spaces $\pmb{\beta}^{-1}(\Lambda(K))$ and $\pmb{\beta}^{-1}(\partial\Lambda(K))$ are the geometric realisations of the cubical subvarieties $\pmb{\Lambda}(K)_\bullet\subset\tilde{\pmb{\Lambda}}_\bullet$ and $\partial\pmb{\Lambda}(K)_\bullet\subset\tilde{\pmb{\Lambda}}_\bullet$ such that 
\begin{align}\label{lambdak}
\begin{split}
\pmb{\Lambda}(K)_I=&\left\{(y_{i_1},\ldots, y_{i_l})\in\tilde{\pmb{\Lambda}}_I\middle\vert h_{i_l}(y_{i_l})\subset K; 
\begin{array}{c}\mbox{{\normalfont the largest (i.e.\ $h_{i_l}(y_{i_l})$) or second largest configuration}}\\
\mbox{{\normalfont among }}h_{i_1}(y_{i_1}),\ldots, h_{i_l}(y_{i_l})\mbox{ {\normalfont belongs to }}\bigcup_{j\leq \alpha(i-1)}Y_j
\end{array}
\right\},\\
\partial\pmb{\Lambda}(K)_I=&\left\{(y_{i_1},\ldots, y_{i_l})\in\tilde{\pmb{\Lambda}}_I\middle\vert h_{i_l}(y_{i_l})\subset K; h_{i_{l}}(y_{i_{l}})\in\bigcup_{j\leq \alpha(i-1)}Y_j\right\}\\
\end{split}
\end{align}
for $\varnothing\neq I=\{i_1<\cdots<i_l\}\subset\{1,\ldots,\N\}$, and $\pmb{\Lambda}(K)_{\varnothing}=\partial\pmb{\Lambda}(K)_{\varnothing}=pt$. (To prove that ${\pmb{\Lambda}}(K)_I$ is a closed subvariety of $\tilde{\pmb{\Lambda}}_I$ one could e.g.\ observe that ${\pmb{\Lambda}}(K)_I$ is constructible and also closed in the complex analytic topology.) Note that we have $\partial\pmb{\Lambda}(K)_\bullet=\pmb{\Lambda}(K)_\bullet\cap\tilde{\pmb{\Lambda}}^{({\alpha(i-1)})}_\bullet$.

It follows from lemma~\ref{proper_contr_fib} that a proper continuous map between locally finite CW-complexes with fibres homeomorphic to contractible CW-complexes is a proper homotopy equivalence. So combining the above observations with lemma~\ref{fib_contractible} we get the following result.
\begin{theorem}\label{thhodge1}
Suppose $Y_1,\ldots, Y_\N$ is a sequence of configuration spaces that satisfies conditions~\ref{first}-\ref{five}$+$ of list~\ref{conditions} and also the conditions of list~\ref{condition_list_3}. Then the following holds.

\begin{enumerate}

\item The map $\pmb{\gamma}$ is a proper filtered homotopy equivalence between the geometric realisation of the cubical complex algebraic variety $\tilde{\pmb{\mathcal{Y}}}_\bullet$ filtered by the cubical subvarieties $\tilde{\pmb{\mathcal{Y}}}^{(i)}_\bullet, i=1,\ldots, \N$, and the conical resolution $\sigma$ with its natural filtration $\varnothing=F_0\subset F_1\subset\cdots\subset F_\N=\sigma$.

\item The map $\pmb{\beta}$ is a proper filtered homotopy equivalence between the geometric realisation of the cubical complex algebraic variety $\tilde{\pmb{\Lambda}}_\bullet$ filtered by the cubical subvarieties $\tilde{\pmb{\Lambda}}^{(i)}_\bullet, i=1,\ldots, \N$, and the space $\tilde\Lambda$ with its natural filtration $\varnothing=\Phi_0\subset\Phi_1\subset\cdots\subset\Phi_\N=\tilde\Lambda$.

\item Suppose $i\in\{1,\ldots,\N\}, K\in Y_i$, and $\pmb{Y}\subset\bar{\pmb{Y}}_i$ is a closed subvariety. Then $\pmb{\beta}$ induces proper filtered homotopy equivalences
\begin{align*}
|\tilde{\pmb{\Lambda}}(K)_\bullet|&=\pmb{\beta}^{-1}(\tilde\Lambda(K))\to \tilde\Lambda(K),\\
|\tilde{\pmb{\Lambda}}(K)_\bullet\cap\tilde{\pmb{\Lambda}}^{(i-1)}_\bullet|&=\pmb{\beta}^{-1}(\tilde\partial\Lambda(K))\to \tilde\partial\Lambda(K),\\
|\tilde{\pmb{\Lambda}}(\pmb{Y})_\bullet|&=\pmb{\beta}^{-1}\left(\bigcup_{K'\in Y} \tilde\Lambda(K')\right)\to \bigcup_{K'\in Y} \tilde\Lambda(K'),\\
\end{align*}
where $\tilde{\pmb{\Lambda}}(K)_\bullet$ and $\tilde{\pmb{\Lambda}}(\pmb{Y})_\bullet$ are cubical subvarieties of $\tilde{\pmb{\Lambda}}_\bullet$.

\item If $X_1,\ldots,X_N$ is a sequence of configuration spaces and $\alpha:\{1,\ldots,N\}\to\{1,\ldots,\N\}$ is a map such that the assumptions on page~\pageref{assumptions_conical} hold for $X_1,\ldots,X_n, Y_1, \ldots, Y_\N$ and $\alpha$, then for every $i\in\{1,\ldots, N\}$ and $K\in X_i=Y_{\alpha(i)}$ the map $\pmb{\beta}$ induces proper filtered homotopy equivalences $|\pmb{\Lambda}(K)_\bullet|=\pmb{\beta}^{-1}(\Lambda(K))\to\Lambda(K)$ and $|\partial\pmb{\Lambda}(K)_\bullet|=\pmb{\beta}^{-1}(\partial\Lambda(K))\to\partial\Lambda(K)$, where $\pmb{\Lambda}(K)_\bullet, \partial\pmb{\Lambda}(K)_\bullet\subset\tilde{\pmb{\Lambda}}_\bullet$ are cubical subvarieties.
\end{enumerate}
The filtrations in parts 3 and 4 are induced by the filtrations from part 2.

\end{theorem}

$\clubsuit$

Suppose $i,j\in\{1,\ldots,\N\},K\in Y_i$. We apply parts 1-3 of theorem~\ref{thhodge1} and proposition~\ref{mhs_compact_support_bm} to equip the groups $\bar H_*(F_i\setminus F_j,\Q)$, $\bar H_*(\Phi_i\setminus \Phi_j,\Q)$ and $\bar H_*(\tilde\Lambda(K)\setminus\tilde\partial\Lambda(K),\Q)$ with mixed Hodge structures using the isomorphisms
\begin{equation}\label{mhs_res}
H_*(F_i\setminus F_{j},\Q)\cong \bar H_*(|\tilde{\pmb{\mathcal{Y}}}^{(i)}_\bullet|\setminus|\tilde{\pmb{\mathcal{Y}}}^{(j)}_\bullet|, \Q)\cong
\bar H_*(|\tilde{\pmb{\mathcal{Y}}}^{(i)}_\bullet|,|\tilde{\pmb{\mathcal{Y}}}^{(j)}_\bullet|, \Q),
\end{equation}
\begin{equation}\label{mhs_order_compl}
\bar H_*(\Phi_i\setminus \Phi_{j},\Q)\cong \bar H_*(|\tilde{\pmb{\Lambda}}^{(i)}_\bullet|\setminus|\tilde{\pmb{\Lambda}}^{(j)}_\bullet|, \Q)\cong
H_*(|\tilde{\pmb{\Lambda}}^{(i)}_\bullet|,|\tilde{\pmb{\Lambda}}^{(j)}_\bullet|, \Q).
\end{equation}
\begin{equation}\label{mhs_tilde_lambdak}
\bar H_*(\tilde\Lambda(K)\setminus\tilde\partial\Lambda(K),\Q)\cong H_*(\tilde\Lambda(K),\tilde\partial\Lambda(K),\Q)\cong H_*(|\tilde{\pmb{\Lambda}}(K)_\bullet|, |\tilde{\pmb{\Lambda}}(K)_\bullet\cap\tilde{\pmb{\Lambda}}^{(i-1)}_\bullet|,\Q).
\end{equation}
Note that by diagram (\ref{cub_spaces_diag2}), the isomorphism $$\pmb{\gamma}_*:\bar H_*(|\tilde{\pmb{\mathcal{Y}}}_\bullet|,\Q)=\bar H_*(|\tilde{\pmb{\mathcal{Y}}}^{(\N)}_\bullet|,\Q)\to \bar H_*(F_\N,\Q)=\bar H_*(\sigma,\Q)$$ composed with the map $\pi_*:\bar H_*(\sigma,\Q)\to\bar H_*(\Sigma,\Q)$ from theorem~\ref{thconical} is the map induced by the augmentation $\pmb{\varepsilon}:|\tilde{\pmb{\mathcal{Y}}}_\bullet|\to \tilde{\pmb{\mathcal{Y}}}_\varnothing=\Sigma$, which is a map of mixed Hodge structures.

Suppose the assumptions of part 4 of theorem~\ref{thhodge1} hold, $i\in\{1,\ldots, N\}$ and $K\in X_i=Y_{\alpha(i)}$. We then introduce a mixed Hodge structure on $\bar H_*(\Lambda(K)\setminus\partial\Lambda(K),\Q)$ using proposition~\ref{mhs_compact_support_bm} again and the isomorphism
\begin{equation}\label{mhs_lambdak}
\bar H_*(\Lambda(K)\setminus\partial\Lambda(K),\Q)\cong H_*(\Lambda(K),\partial\Lambda(K),\Q)\cong H_*(|{\pmb{\Lambda}}(K)_\bullet|, |\partial{\pmb{\Lambda}}(K)_\bullet|,\Q)\cong H_*(|{\pmb{\Lambda}}(K)_\bullet|, |\pmb{\Lambda}(K)_\bullet\cap\tilde{\pmb{\Lambda}}^{({\alpha(i-1)})}_\bullet|,\Q).
\end{equation}

If in addition $X_1,\ldots, X_N$ satisfy condition~\ref{five} from list~\ref{conditions}, we get a mixed Hodge structure on  $\bar H_*(\Psi_i,\Q),i\in\{1,\ldots,\N\}$ via the isomorphism
\begin{equation}\label{mhs_psi}
\bar H_*(\Psi_i,\Q)\cong\bar H_*(\Phi_{\alpha(i)}\setminus \Phi_{{\alpha(i-1)}},\Q)\cong \bar H_*(|\tilde{\pmb{\Lambda}}^{({\alpha(i)})}_\bullet|\setminus|\tilde{\pmb{\Lambda}}^{({\alpha(i-1)})}_\bullet|, \Q)\cong
H_*(|\tilde{\pmb{\Lambda}}^{({\alpha(i)})}_\bullet|,|\tilde{\pmb{\Lambda}}^{({\alpha(i-1)})}_\bullet|, \Q).
\end{equation}
Here the first isomorphism follows from theorem~\ref{thconical_1}, and the second from theorem~\ref{thhodge1}. 

Recall that all structure maps of $\tilde{\pmb{\mathcal{Y}}}_\bullet$ are proper (lemma~\ref{lem_prop}). Applying proposition~\ref{mhs_sseq_filtred} we get the following result.

\begin{corollary}\label{cor_mhs_main_sseq}
The isomorphism $\bar H_*(\sigma,\Q)=\bar H_*(F_\N,\Q)\to\bar H_*(\Sigma,\Q)$ from theorem~\ref{thconical} is a map of mixed Hodge structures. Moreover, the Borel-Moore homology spectral sequence that corresponds to the filtration $\varnothing=F_0\subset F_1\subset\cdots\subset F_\N=\sigma$ has a natural mixed Hodge structure starting from $E^1$ (see section~\ref{sec_mhs_sseq}).
\end{corollary}
$\clubsuit$

\subsubsection{The mixed Hodge structure on auxiliary spectral sequences}

In this section all homology groups will be with rational coefficients, so we will not write the coefficients for homology except in theorems.

\level{4}{The spectral sequences of \texorpdfstring{$F_i\setminus F_{i-1}\to\Phi_i\setminus\Phi_{i-1}$}{} and \texorpdfstring{$\Phi_i\setminus\Phi_{i-1}\to Y_i$}{}}\label{aux_spec_seq_1}

In order to calculate the term $E^1$ of the Borel-Moore homology spectral sequence for $\sigma$ we will use theorem~\ref{thconical} (parts 2 and 3) and theorem~\ref{thconical_1}. We will now show that the corresponding spectral sequences also have natural mixed Hodge structures. The main obstacle is that the natural maps $\Phi_i\setminus\Phi_{i-1}\to Y_i$ and $|\tilde{\pmb{\Lambda}}^{(i)}_\bullet|\setminus |\tilde{\pmb{\Lambda}}^{(i-1)}_\bullet|\to Y_i$ need not extend to $\Phi_i$, respectively $|\tilde{\pmb{\Lambda}}^{(i)}_\bullet|$. In order to remedy this we fix an $i\in\{1,\ldots, \N\}$ and introduce auxiliary cubical complex algebraic varieties $\tilde M_\bullet\supset \tilde\partial M_\bullet$ and $\tilde Z_\bullet\supset \tilde\partial Z_\bullet$ as follows. For a non-empty $I=\{i_1<\cdots<i_l\}\subset\{1,\ldots,\N\}$ we set
\begin{align}\label{aux_var}
\begin{split}
\tilde{M}_I=& \left\{((y_{i_1},\ldots, y_{i_l}),y)\in \tilde{\pmb{\Lambda}}_I\times\bar{\pmb{Y}}_i\middle\vert h_{i_l}(y_{i_l})\subset h_i(y)\right\},\\
\tilde\partial{M}_I=& \left\{((y_{i_1},\ldots, y_{i_l}),y)\in \tilde{\pmb{\Lambda}}_I\times\bar{\pmb{Y}}_i\middle\vert h_{i_l}(y_{i_l})\subset h_i(y); h_{i_l}(y_{i_l})\in\bigcup_{j\leq i-1} Y_j\right\},\\
\tilde Z_I=&\left\{(f,(y_{i_1},\ldots, y_{i_l}),y)\in V\times\tilde M_I\middle\vert f\in\mathcal{E}^i_{h_i(y)}\right\},\\
\tilde\partial Z_I=&\left\{(f,(y_{i_1},\ldots, y_{i_l}),y)\in V\times\tilde \partial M_I\middle\vert f\in\mathcal{E}^i_{h_i(y)}\right\}
\end{split}
\end{align}
where $\mathcal{E}^i_{h_i(y)}$ is the fibre of the vector bundle $\mathcal{E}^i$ (see page~\pageref{bundle_ei}) over $h_i(y)\in\bar{\pmb{Y}}_i$. We also set $\tilde M_\varnothing=\tilde\partial M_\varnothing=\bar{\pmb{Y}}_i, \tilde Z_\varnothing=\tilde\partial Z_\varnothing=\tot(\mathcal{E}^i)\subset V\times\bar{\pmb{Y}}_i$, and we take all structure maps to be the natural projections. 
Using condition list~\ref{condition_list_3} we deduce that $\tilde M_\bullet,\tilde\partial M_\bullet,\tilde Z_\bullet$ and $\tilde\partial Z_\bullet$ are indeed cubical complex algebraic varieties, and moreover that all components of $\tilde M_\bullet,\tilde\partial M_\bullet$ are compact and $\tilde\partial Z_\bullet$ is closed in~$\tilde Z_\bullet$.

We have natural maps $\tilde Q_\bullet:\tilde Z_\bullet\to \tilde{\pmb{\mathcal{Y}}}_\bullet^{(i)}$ and $\tilde q_\bullet:\tilde M_\bullet\to \tilde{\pmb{\Lambda}}_\bullet^{(i)}$ of cubical complex algebraic varieties. All components of $\tilde q_\bullet$ are clearly proper. Let us show that so are the components of $\tilde Q_\bullet$. Suppose first $I\neq \varnothing$. The space $\tilde Z_I$ is closed in $V\times\tilde M_I$ as it is the total space of a vector subbundle of the trivial bundle $V\times\tilde M_I$. The map $\id_V\times \tilde q_I:V\times\tilde M_I\to\ V\times\tilde{\pmb{\Lambda}}_I^{(i)}$ is proper, hence so is $\tilde Q_I$. The map $\tilde Q_\varnothing:\tot(\mathcal{E}^i)\to \Sigma$ is the restriction of the proper projection map $V\times\bar{\pmb{Y}}_i\to V$ to the closed subset $\tot(\mathcal{E}^i)$, so $\tilde Q_\varnothing$ is also proper.

\smallskip

\begin{lemma}\label{excision}
\begin{enumerate}
\item Suppose $(A,B)$ is a locally finite CW-pair. If $U\subset A$ is an open subset such that $A\setminus B\subset U$, then we have $\bar H_*(A,B)\cong \bar H_*(U,U\cap B)$.
\item If $f:(A,B)\to (A',B')$ is a map of finite CW-pairs such that $f:A\to A'$ is surjective, $f^{-1}(B')=B$, and $f$ is injective on $A\setminus B$, then $f$ induces a homeomorphism $A\setminus B\to A'\setminus B'$, and $f_*:H_*(A,B)\to H_*(A',B')$ is an isomorphism.
\end{enumerate}
\end{lemma}
$\clubsuit$

Suppose $I=\{i_1<\cdots<i_l\}\subset\{1,\ldots,\N\}, (y_{i_1},\ldots,y_{i_l})\in\tilde{\pmb{\Lambda}}_I$ and $h_{i_l}(y_{i_l})\in Y_i$. Then there is a unique $y\in \pmb{Y}_i$ such that $h_i(y)=h_{i_l}(y_{i_l})$ (condition~\ref{cond_1_list_2} of list~\ref{condition_list_3}). This implies that $|\tilde q_\bullet|:|\tilde M_\bullet|\to|\tilde{\pmb{\Lambda}}_\bullet^{(i)}|$ is injective on $|\tilde M_\bullet|\setminus |\tilde \partial M_\bullet|$. Moreover, it is also surjective, and the preimage of $|\tilde{\pmb{\Lambda}}^{(i-1)}|$ under $|\tilde q_\bullet|$ is $|\tilde\partial M_\bullet|$.

Using lemma~\ref{excision} we see that $\tilde q_\bullet$ induces a homeomorphism
$$|\tilde M_\bullet|\setminus|\tilde\partial M_\bullet|\to |\tilde{\pmb{\Lambda}}^{(i)}_\bullet|\setminus|\tilde{\pmb{\Lambda}}^{(i-1)}_\bullet|$$ and an isomorphism
\begin{equation}\label{iso1}
H_*(|\tilde M_\bullet|,|\tilde\partial M_\bullet|)\cong\bar H_*(|\tilde M_\bullet|\setminus|\tilde\partial M_\bullet|)\to \bar H_*(|\tilde{\pmb{\Lambda}}^{(i)}_\bullet|\setminus|\tilde{\pmb{\Lambda}}^{(i-1)}_\bullet|)\cong H_*(|\tilde{\pmb{\Lambda}}^{(i)}_\bullet|,|\tilde{\pmb{\Lambda}}^{(i-1)}_\bullet|).
\end{equation}

\smallskip

Let $\varepsilon:|\tilde M_\bullet|\to \bar{\pmb{Y}}_i$ be the augmentation map. Recall that for a subvariety $\pmb{A}\subset \bar{\pmb{Y}}_i$ we defined in the beginning of section~\ref{hodge_modge} a cubical subvariety $\varepsilon^{-1}(\pmb{A})_\bullet$. 
If $I=\{i_i<\cdots<i_l\}\subset\{1,\ldots,N\}, (y_{i_1},\ldots, y_{i_l})\in\tilde{\pmb{\Lambda}}_I, y\in\bar{\pmb{Y}}_i\setminus\pmb{Y}_i$ and $h_{i_l}(y_{i_l})\subset h_i(y)$, then $h_{i_l}(y_{i_l})\in\bigcup_{j\leq i-1} Y_j$, so $\varepsilon^{-1}(\bar{\pmb{Y}}_i\setminus \pmb{Y}_i)_\bullet\subset\tilde\partial M_\bullet$, which implies $|\varepsilon^{-1}(\pmb{Y}_i)_\bullet|\supset |\tilde M_\bullet|\setminus |\tilde\partial M_\bullet|$. 

So applying of lemma~\ref{excision} we get an isomorphism
\begin{equation}\label{iso_aux_1}
H_*(|\tilde M_\bullet|,|\tilde\partial M_\bullet|)\cong  \bar H_*(|\varepsilon^{-1}(\pmb{Y}_i)_\bullet|,|\varepsilon^{-1}(\pmb{Y}_i)_\bullet\cap\tilde\partial M_\bullet|).
\end{equation}

\smallskip

Let $y\in\pmb{Y}_i$ be a point, and let $K=h_i(y)$. By construction, $\tilde q_\bullet$ takes $\varepsilon^{-1}(y)_\bullet$ surjectively to $\tilde{\pmb{\Lambda}}(K)_\bullet$. Since $\tilde q_\bullet^{-1}\left(\tilde{\pmb{\Lambda}}^{(i-1)}_\bullet\right)=\tilde\partial M_\bullet$, and $|\tilde q_\bullet|$ is injective on $|\tilde M_\bullet|\setminus |\tilde \partial M_\bullet|$, we obtain using part 2 of lemma~\ref{excision} that $\tilde q_\bullet$ induces an isomorphism
\begin{equation}\label{isofibre}
H_*(|\varepsilon^{-1}(y)_\bullet|,|\varepsilon^{-1}(y)_\bullet\cap\tilde \partial M_\bullet|)\to H_*(|\tilde{\pmb{\Lambda}}(K)_\bullet|, |\tilde{\pmb{\Lambda}}(K)_\bullet\cap \tilde{\pmb{\Lambda}}^{(i-1)}_\bullet|).
\end{equation}

\smallskip

Note that for $y\in\pmb{Y}_i$, we have $f\in\mathcal{E}^i_{h_i(y)}$ if and only if $h_i(y)\subset\Sing(f)$. This implies that $\tilde Q_\bullet$ induces a proper continuous bijection, hence a homeomorphism
$$|\tilde Z_\bullet|\setminus|\tilde\partial Z_\bullet|\to |\tilde{\pmb{\mathcal{Y}}}^{(i)}_\bullet|\setminus|\tilde{\pmb{\mathcal{Y}}}^{(i-1)}_\bullet|$$ and an isomorphism
\begin{equation}\label{iso2}
\bar H_*(|\tilde Z_\bullet|,|\tilde\partial Z_\bullet|)\cong\bar H_*(|\tilde Z_\bullet|\setminus|\tilde\partial Z_\bullet|)\to \bar H_*(|\tilde{\pmb{\mathcal{Y}}}^{(i)}_\bullet|\setminus|\tilde{\pmb{\mathcal{Y}}}^{(i-1)}_\bullet|)\cong \bar H_*(|\tilde{\pmb{\mathcal{Y}}}^{(i)}_\bullet|,|\tilde{\pmb{\mathcal{Y}}}^{(i-1)}_\bullet|).
\end{equation}

\smallskip

There is a map $\tilde p_\bullet:\tilde Z_\bullet\to\tilde M_\bullet$ of cubical complex algebraic varieties. For a subvariety $\pmb{A}\subset\bar{\pmb{Y}}_i$ set $\epsilon^{-1}(\pmb{A})_\bullet=\tilde p_\bullet^{-1}(\varepsilon^{-1}(\pmb{A})_\bullet)$. The preimage $|\tilde p_\bullet|^{-1}(|\varepsilon^{-1}(\pmb{Y}_i)_\bullet|)$ is by construction $|\epsilon^{-1}(\pmb{Y}_i)_\bullet|$, and $|\tilde p_\bullet|^{-1}(|\tilde\partial M_\bullet|)=|\tilde\partial Z_\bullet|$. The open set $|\epsilon^{-1}(\pmb{Y}_i)_\bullet|$ of $|\tilde Z_\bullet|$ contains $|\tilde Z_\bullet|\setminus |\tilde\partial Z_\bullet|$ because $|\varepsilon^{-1}(\pmb{Y}_i)_\bullet|\supset|\tilde M_\bullet|\setminus |\tilde\partial M_\bullet|$. So by lemma~\ref{excision}, part 1 we have an isomorphism 
\begin{equation}\label{iso_aux_2}
\bar H_*(|\tilde Z_\bullet|,|\tilde\partial Z_\bullet|)\cong \bar H_*(|\epsilon^{-1}(\pmb{Y}_i)_\bullet|,|\epsilon^{-1}(\pmb{Y}_i)_\bullet\cap\tilde\partial Z_\bullet|).
\end{equation}

We note that $\tilde Z_\bullet$ can be viewed as a {\it cubical complex algebraic vector bundle} (or CCAVB for brevity\label{ccavb}) over $\tilde M_\bullet$, meaning that each component $\tilde Z_I,I\subset\{1,\ldots,\N\}$ is an algebraic vector bundle over $\tilde M_I$ and the structure maps are fibrewise isomorphisms of algebraic vector bundles. Similarly, $\tilde\partial Z_\bullet$ is a CCAVB over $\tilde \partial M_\bullet$, and $\epsilon^{-1}(\pmb{A})_\bullet$ is a CCAVB over $\varepsilon^{-1}(\pmb{A})_\bullet$ for every subvariety $\pmb{A}\subset\bar{\pmb{Y}}_i$. Let us also note that all these bundles are obtained by pulling back the bundle $\mathcal{E}^i$ over $\bar{\pmb{Y}}_i$. Recall that the rank of $\mathcal{E}^i$ is denoted $d_i$, see condition list~\ref{conditions}.

The Thom isomorphism theorem extends to CCAVBs over cubical algebraic varieties (see theorem~\ref{thom_ccavb}), and we get
\begin{equation}\label{isothom}
\bar H_*(|\tilde Z_\bullet|,|\tilde\partial Z_\bullet|) \cong
H_*(|\tilde M_\bullet|,|\tilde\partial M_\bullet|)\otimes\Q(d_i)[-2d_i]
\end{equation}
where $\Q(d_i)[-2d_i]$ is the Tate Hodge structure $\Q(d_i)$ placed in degree $2d_i$.

Let us collect some of the maps we have constructed so far in a single commutative diagram:
\begin{equation}\label{big_comm_diag}
\begin{tikzcd}[scale cd=.9]
{|}{\tilde{Z}}_{\bullet}{|}\arrow{d}{{|}\tilde p_\bullet{|}} & {|}\epsilon^{-1}(\pmb{Y}_i)_\bullet{|}\arrow[swap]{l}{o}\arrow{d}{{|}\tilde p_\bullet{|}}& {|}\epsilon^{-1}(\pmb{Y}_i)_\bullet{|}\setminus {|}\epsilon^{-1}(\pmb{Y}_i)_\bullet\cap \tilde\partial Z_\bullet{|}=|\tilde Z_\bullet|\setminus|\tilde\partial Z_\bullet|\arrow[swap]{l}{o}\arrow{d}{{|}\tilde p_\bullet{|}}\arrow{r}{\cong}[swap]{|\tilde Q_\bullet|} & {|}\tilde{\pmb{\mathcal{Y}}}^{(i)}_\bullet|\setminus |\tilde{\pmb{\mathcal{Y}}}^{(i-1)}_\bullet{|}\arrow[d]\arrow{r}{\sim}[swap]{\pmb{\gamma}}& F_i\setminus F_{i-1}\arrow{d}{p}\\
{|}{\tilde{M}}_{\bullet}{|}\arrow{dr}{\varepsilon} & {|}\varepsilon^{-1}(\pmb{Y}_i)_\bullet{|}\arrow[swap]{l}{o}\arrow{dr}{\varepsilon}& {|}\varepsilon^{-1}(\pmb{Y}_i)_\bullet{|}\setminus {|}\varepsilon^{-1}(\pmb{Y}_i)_\bullet\cap \tilde\partial M_\bullet{|}=|\tilde M_\bullet|\setminus|\tilde\partial M_\bullet|\arrow[swap]{l}{o}\arrow{d}{\varepsilon}\arrow{r}{\cong}[swap]{|\tilde q_\bullet|} & {|}\tilde{\pmb{\Lambda}}^{(i)}_\bullet|\setminus |\tilde{\pmb{\Lambda}}^{(i-1)}_\bullet{|}\arrow[dl]\arrow{r}{\sim}[swap]{\pmb{\beta}}& \Phi_i\setminus\Phi_{i-1}\arrow{d}{f_i}\\
&\bar{\pmb{Y}}_i & \pmb{Y}_i\arrow[swap]{l}{o}\arrow{rr}{\cong}[swap]{h_i} & & Y_i.
\end{tikzcd}
\end{equation}
Here the arrows marked $\cong$ are homeomorphisms, those marked $o$ are open embeddings, and those marked $\sim$ are proper homotopy equivalences. 

\begin{theorem}\label{thhodge2}
We keep the assumptions of theorem~\ref{thhodge1}. For every $i\in\{1,\ldots,\N\}$ the following holds.
\begin{enumerate}
\item Recall that above we equipped $\bar H_*(F_i\setminus F_{i-1},\Q)$ and $\bar H_*(\Phi_i\setminus \Phi_{i-1},\Q)$ with mixed Hodge structures using isomorphisms~(\ref{mhs_res}) and~(\ref{mhs_order_compl}). We then have
$$\bar H_*(F_i\setminus F_{i-1},\Q)\cong\bar H_*(\Phi_i\setminus\Phi_{i-1},\Q)\otimes\Q(d_i)[-2d_i]$$ as mixed Hodge structures.

\item Both $\varepsilon^{-1}(\pmb{Y}_i)_\bullet$ and $\varepsilon^{-1}(\pmb{Y}_i)_\bullet\cap\tilde\partial M_\bullet$ are cubical complex algebraic varieties, and the latter is a closed subvariety of the former. The Borel-Moore homology Leray spectral sequence of $\Phi_{i}\setminus\Phi_{i-1}\to Y_i$ is isomorphic to the Leray spectral sequence for $\bar H_*(|\varepsilon^{-1}(\pmb{Y}_i)_\bullet|,|\varepsilon^{-1}(\pmb{Y}_i)_\bullet\cap\tilde\partial M_\bullet|,\Q)$ obtained from the map ${|}\varepsilon^{-1}(\pmb{Y}_i)_\bullet{|}\to\pmb{Y}_i$, which is the augmentation map for both $\varepsilon^{-1}(\pmb{Y}_i)_\bullet$ and $\varepsilon^{-1}(\pmb{Y}_i)_\bullet\cap\tilde\partial M_\bullet$. If $\pmb{Y}_i$ is quasi-projective, then this Leray spectral sequence has a natural mixed Hodge structure starting from $E_2$ (see section~\ref{sec_mhs_sseq}).
\item Recall that for $K\in Y_i$ the mixed Hodge structure on $\bar H_*(\tilde\Lambda(K)\setminus\tilde\partial\Lambda(K),\Q)$ was introduced using isomorphism~(\ref{mhs_tilde_lambdak}). This mixed Hodge structure agrees with the one on the fibre of the map ${|}\varepsilon^{-1}(\pmb{Y}_i)_\bullet{|}\to\pmb{Y}_i$ from part 2, i.e.\ if $y\in\pmb{Y}_i$ is such that $h_i(y)=K$, then we have
$$\bar H_*(\tilde\Lambda(K)\setminus\tilde\partial\Lambda(K),\Q)\cong H_*(|\varepsilon^{-1}(y)_\bullet|,|\varepsilon^{-1}(y)_\bullet\cap\tilde \partial M_\bullet|,\Q)$$ as mixed Hodge structures.
\end{enumerate}
\end{theorem}

{\bf Proof.} We conclude from formulas~(\ref{iso2}), (\ref{isothom}) and~(\ref{iso1}), that
\begin{multline*}
\bar H_*(|\tilde{\pmb{\mathcal{Y}}}^{(i)}_\bullet|\setminus|\tilde{\pmb{\mathcal{Y}}}^{(i-1)}_\bullet|)\cong
\bar H_*(|\tilde Z_\bullet|,|\tilde\partial Z_\bullet|) \cong
\bar H_*(|\tilde M_\bullet|,|\tilde\partial M_\bullet|)\otimes\Q(d_i)[-2d_i]
\cong
\bar H_*(|\tilde{\pmb{\Lambda}}^{(i)}_\bullet|,|\tilde{\pmb{\Lambda}}^{(i-1)}_\bullet|)
\otimes\Q(d_i)[-2d_i].
\end{multline*}
Moreover, all these isomorphisms are morphisms of mixed Hodge structures, which proves the first part of the theorem. Since all structure maps of $\tilde M_\bullet$ are proper, so are the structure maps of $\varepsilon^{-1}(\pmb{Y}_i)_\bullet$. The second part now follows from diagram (\ref{big_comm_diag}) and proposition~\ref{mhs_sseq_leray}, and the third part from formula (\ref{isofibre}).$\clubsuit$

{\bf Remark.} We expect that the Leray spectral sequence for $\bar H_*(|\varepsilon^{-1}(\pmb{Y}_i)_\bullet|,|\varepsilon^{-1}(\pmb{Y}_i)_\bullet\cap\tilde\partial M_\bullet|,\Q)$ from part 2 of the theorem has a mixed Hodge structure regardless whether $\pmb{Y}_i$ quasi-projective, see the remark after proposition~\ref{mhs_sseq_leray}.

\level{4}{The spectral sequences of \texorpdfstring{$F_{\alpha(i)}\setminus F_{\alpha(i-1)}\to\Phi_{\alpha(i)}\setminus\Phi_{\alpha(i-1)}$}{} and \texorpdfstring{$\Psi_i\to X_i$}{}}

Suppose now $X_1,\ldots,X_N$ is a sequence of configuration spaces and $\alpha:\{1,\ldots, N\}\to \{1,\ldots,\N\}$ is a map such that the assumptions on page~\pageref{assumptions_conical} hold for $X_1,\ldots,X_N, Y_1, \ldots, Y_\N$ and $\alpha$. Suppose moreover that $X_1,\ldots,X_N$ satisfy condition \ref{five} of list~\ref{conditions} (in addition to~\ref{five}$-$).

Let us now show that then for $i=1,\ldots, N$ the Borel-Moore homology Leray spectral sequences for the maps $F_{\alpha(i)}\setminus F_{\alpha(i-1)}\to \Phi_{\alpha(i)}\setminus \Phi_{\alpha(i-1)}$ (part 2 of theorem~\ref{thconical}), $\Phi_{\alpha(i)}\setminus \Phi_{\alpha(i-1)}\to X_i$ and $\Psi_i\to X_i$ (part 1 of theorem~\ref{thconical_1}) also have natural mixed Hodge structures. The argument is similar to the one we have just seen in section~\ref{aux_spec_seq_1}, but we sketch it for completeness with a focus on what is different. We fix an $i\in\{1,\ldots, N\}$ and define cubical complex algebraic varieties $M_\bullet \supset\partial M_\bullet$ and $Z_\bullet\supset\partial Z_\bullet$. Namely, we take $M_\bullet$ and $Z_\bullet$ to be $\tilde M_\bullet$, respectively $\tilde Z_\bullet$ constructed in section~\ref{aux_spec_seq_1} for $\alpha(i)$. For a non-empty $I=\{i_1<\cdots<i_l\}\subset\{1,\ldots,\N\}$ we set
\begin{align*}
\partial{M}_I=& \left\{((y_{i_1},\ldots, y_{i_l}),y)\in \tilde{\pmb{\Lambda}}_I\times\bar{\pmb{Y}}_{\alpha(i)}\middle\vert h_{i_l}(y_{i_l})\subset h_{\alpha(i)}(y); h_{i_l}(y_{i_l})\in\bigcup_{j\leq \alpha(i-1)}Y_j\right\},\\
\partial Z_I=&\left\{(f,(y_{i_1},\ldots, y_{i_l}),y)\in V\times\partial M_I\middle\vert f\in\mathcal{E}^{\alpha(i)}_{h_{\alpha(i)}(y)}\right\}.
\end{align*}
We define $\partial M_\varnothing=\bar{\pmb{Y}}_{\alpha(i)}, \partial Z_\varnothing=\tot(\mathcal{E}^{\alpha(i)})$ and take all structure maps of $\partial M_\bullet$ and $\partial Z_\bullet$ to be the natural projections. All components of $M_\bullet$ and $\partial M_\bullet$ are compact, and $\partial Z_\bullet$ is closed in $Z_{\bullet}$. Moreover, in~\ref{aux_spec_seq_1} we defined maps $\tilde q_\bullet:M_\bullet\to\tilde{\pmb{\Lambda}}_\bullet^{\alpha(i)}$ and $\tilde Q:Z_\bullet\to\tilde{\pmb{\mathcal{Y}}}_{\bullet}^{\alpha(i)}$ of cubical varieties and proved that all components of these maps are proper. Let us now denote these maps $q_\bullet$ and $Q_\bullet$ respectively.

Suppose now $I=\{i_1<\cdots<i_l\}\subset\{1,\ldots,\N\}, (y_{i_1},\ldots,y_{i_l})\in\tilde{\pmb{\Lambda}}_I$ and $h_{i_l}(y_{i_l})\in \bigcup_{\alpha(i-1)<j\leq \alpha(i)}Y_j$. Then by condition~\ref{five} of list~\ref{conditions} and condition~\ref{cond_1_list_2} of list~\ref{condition_list_3} there is a unique $y\in\pmb{Y}_{\alpha(i)}$ such that $h_{i_l}(y_{i_l})\subset h_{\alpha(i)}(y)$, namely $y$ is the unique point of $\pmb{Y}_{\alpha(i)}$ such that $h_{\alpha(i)}(y)=\Geom(h_{i_l}(y_{i_l}))$. So $|q_\bullet|:|M_\bullet|\to|\tilde{\pmb{\Lambda}}_\bullet^{\alpha(i)}|$ is injective on $|M_\bullet|\setminus |\partial M_\bullet|$. It is also surjective, and we have $|q_\bullet|^{-1}(|\tilde{\pmb{\Lambda}}^{(\alpha(i-1))}|)=|\partial M_\bullet|$. Using lemma~\ref{excision} we conclude that $|q_\bullet|$ induces a homeomorphism $$| M_\bullet|\setminus|\partial M_\bullet|\to |\tilde{\pmb{\Lambda}}^{(\alpha(i))}_\bullet|\setminus|\tilde{\pmb{\Lambda}}^{(\alpha(i-1))}_\bullet|,$$ and an isomorphism (cf.\ (\ref{iso1}))
$$H_*(|M_\bullet|,|\partial M_\bullet|)\cong\bar H_*(|M_\bullet|\setminus|\partial M_\bullet|)\to \bar H_*(|\tilde{\pmb{\Lambda}}^{(\alpha(i))}_\bullet|\setminus|\tilde{\pmb{\Lambda}}^{(\alpha(i-1))}_\bullet|)\cong H_*(|\tilde{\pmb{\Lambda}}^{(\alpha(i))}_\bullet|,|\tilde{\pmb{\Lambda}}^{(\alpha(i-1))}_\bullet|).$$ We also have a homeomorphism
$$|Z_\bullet|\setminus|\partial Z_\bullet|\to |\tilde{\pmb{\mathcal{Y}}}^{(\alpha(i))}_\bullet|\setminus|\tilde{\pmb{\mathcal{Y}}}^{(\alpha(i-1))}_\bullet|$$ induced by $|Q_\bullet|$, and an isomorphism (cf.\ (\ref{iso2}))
$$\bar H_*(|Z_\bullet|,|\partial Z_\bullet|)\cong\bar H_*(|Z_\bullet|\setminus|\partial Z_\bullet|)\to \bar H_*(|\tilde{\pmb{\mathcal{Y}}}^{(\alpha(i))}_\bullet|\setminus|\tilde{\pmb{\mathcal{Y}}}^{(\alpha(i-1))}_\bullet|)\cong \bar H_*(|\tilde{\pmb{\mathcal{Y}}}^{(\alpha(i))}_\bullet|,|\tilde{\pmb{\mathcal{Y}}}^{(\alpha(i-1))}_\bullet|).$$

Let $\varepsilon:|M_\bullet|\to \bar{\pmb{Y}}_{\alpha(i)}$ and $p_\bullet:Z_\bullet\to M_\bullet$ be the augmentation map and the natural projection respectively. 
As in section~\ref{aux_spec_seq_1}, starting from a subvariety $\pmb{A}\subset\bar{\pmb{Y}}_{\alpha(i)}$ we construct a cubical complex algebraic subvariety $\varepsilon^{-1}(\pmb{A})_\bullet\subset M_\bullet$ and set $\epsilon^{-1}(\pmb{A})_\bullet=p^{-1}_\bullet(\varepsilon^{-1}(A)_\bullet)\subset Z_\bullet$. Let us show that $\varepsilon^{-1}(\bar{\pmb{Y}}_{\alpha(i)}\setminus {\pmb{Y}}_{\alpha(i)})_\bullet\subset \partial M_\bullet$. Suppose $I=\{i_i<\cdots<i_l\}\subset\{1,\ldots,N\}, (y_{i_1},\ldots, y_{i_l})\in\tilde{\pmb{\Lambda}}_I, y\in\bar{\pmb{Y}}_{\alpha(i)}\setminus\pmb{Y}_{\alpha(i)}$ and $h_{i_l}(y_{i_l})\subset h_{\alpha(i)}(y)$. Since we assume that $X_1,\ldots, X_N$ satisfy condition~\ref{five} of list~\ref{conditions}, the geometrisation of $h_{\alpha(i)}(y)$ belongs to some $X_j=Y_{\alpha(j)},j<i$. So $h_{i_l}(y_{i_l})$, being a subset of $\Geom(h_{\alpha(i)}(y))$, belongs to $\bigcup_{j\leq \alpha(i-1)}Y_j$, i.e.\ $((y_{i_1},\ldots, y_{i_l}),y)\in \partial M_I$.

As in~\ref{aux_spec_seq_1}, we obtain isomorphisms
\begin{gather*}
H_*(|M_\bullet|,|\partial M_\bullet|)\cong  \bar H_*(|\varepsilon^{-1}(\pmb{Y}_{\alpha(i)})_\bullet|,|\varepsilon^{-1}(\pmb{Y}_{\alpha(i)})_\bullet\cap\partial M_\bullet|),\\
\bar H_*(|Z_\bullet|,|\partial Z_\bullet|)\cong \bar H_*(|\epsilon^{-1}(\pmb{Y}_{\alpha(i)})_\bullet|,|\epsilon^{-1}(\pmb{Y}_{\alpha(i)})_\bullet\cap\partial Z_\bullet|),\\
H_*(|\varepsilon^{-1}(y)_\bullet|,|\varepsilon^{-1}(y)_\bullet\cap \partial M_\bullet|)\cong H_*(|\tilde{\pmb{\Lambda}}(K)_\bullet|, |\tilde{\pmb{\Lambda}}(K)_\bullet\cap \tilde{\pmb{\Lambda}}^{(\alpha(i-1))}_\bullet|)
\end{gather*}
where $y\in\pmb{Y}_{\alpha(i)}$ and $K=h_{\alpha(i)}(y)$, cf.\ (\ref{iso_aux_1}), (\ref{iso_aux_2}) and (\ref{isofibre}). Let us now compare $H_*(|\tilde{\pmb{\Lambda}}(K)_\bullet|, |\tilde{\pmb{\Lambda}}(K)_\bullet\cap {\tilde{\pmb{\Lambda}}_\bullet}^{(\alpha(i-1))}|)$ with $H_*(|{\pmb{\Lambda}}(K)_\bullet|, |{\pmb{\Lambda}}(K)_\bullet\cap \tilde{\pmb{\Lambda}}^{(\alpha(i-1))}_\bullet|)$. As we observed in section~\ref{sec_main_spec_sec}, we have $$\pmb{\beta}^{-1}(\tilde\Lambda(K))=|\tilde{\pmb{\Lambda}}(K)_\bullet|, \pmb{\beta}^{-1}(\Lambda(K))=|{\pmb{\Lambda}}(K)_\bullet|\mbox{ and }\pmb{\beta}^{-1}(\Phi_{\alpha(i)})=\left|\tilde{\pmb{\Lambda}}_\bullet^{(\alpha(i))}\right|.$$
So $\pmb{\beta}$ induces isomorphisms
\begin{align*}
H_*(|\tilde{\pmb{\Lambda}}(K)_\bullet|, |\tilde{\pmb{\Lambda}}(K)_\bullet\cap {\tilde{\pmb{\Lambda}}_\bullet}^{(\alpha(i-1))}|)&\to H_*(\tilde\Lambda(K),\tilde\Lambda(K)\cap\Phi_{\alpha(i-1)}),\\
H_*(|{\pmb{\Lambda}}(K)_\bullet|, |{\pmb{\Lambda}}(K)_\bullet\cap {\tilde{\pmb{\Lambda}}_\bullet}^{(\alpha(i-1))}|)&\to H_*(\Lambda(K),\Lambda(K)\cap\Phi_{\alpha(i-1)}).
\end{align*}
By proposition~\ref{lambda} (see also theorem~\ref{thconical_1}), we have $H_*(\tilde\Lambda(K),\tilde\Lambda(K)\cap\Phi_{\alpha(i-1)})\cong H_*(\Lambda(K),\Lambda(K)\cap\Phi_{\alpha(i-1)})$, which implies that
$$H_*(|\varepsilon^{-1}(y)_\bullet|,|\varepsilon^{-1}(y)_\bullet\cap \partial M_\bullet|)\cong H_*(|\tilde{\pmb{\Lambda}}(K)_\bullet|, |\tilde{\pmb{\Lambda}}(K)_\bullet\cap \tilde{\pmb{\Lambda}}^{(\alpha(i-1))}_\bullet|)\cong H_*(|{\pmb{\Lambda}}(K)_\bullet|, |{\pmb{\Lambda}}(K)_\bullet\cap {\tilde{\pmb{\Lambda}}_\bullet}^{(\alpha(i-1))}|).$$
\smallskip

Moreover, we again have
$$\bar H_*(|Z_\bullet|,|\partial Z_\bullet|) \cong
H_*(|M_\bullet|,|\partial M_\bullet|)\otimes\Q(d_{\alpha(i)})[-2d_{\alpha(i)}]
$$
using the Thom isomorphism, see theorem~\ref{thom_ccavb}. The analogue of diagram (\ref{big_comm_diag}) reads
\begin{equation*}
\begin{tikzcd}[scale cd=.71]
{|}{{Z}}_{\bullet}{|}\arrow{d}{{|} p_\bullet{|}} & {|}\epsilon^{-1}(\pmb{Y}_{\alpha(i)})_\bullet{|}\arrow[swap]{l}{o}\arrow{d}{{|} p_\bullet{|}}& {|}\epsilon^{-1}(\pmb{Y}_{\alpha(i)})_\bullet{|}\setminus {|}\epsilon^{-1}(\pmb{Y}_{\alpha(i)})_\bullet\cap \partial Z_\bullet{|}={|}Z_\bullet{|}\setminus{|}\partial Z_\bullet{|}\arrow[swap]{l}{o}\arrow{d}{{|} p_\bullet{|}}\arrow{r}{\cong}[swap]{|Q_\bullet|} & {|}\tilde{\pmb{\mathcal{Y}}}^{({\alpha(i)})}_\bullet|\setminus |\tilde{\pmb{\mathcal{Y}}}^{({\alpha(i-1)})}_\bullet{|}\arrow[d]\arrow{r}{\sim}[swap]{\pmb{\gamma}}& F_{\alpha(i)}\setminus F_{{\alpha(i-1)}}\arrow{d}{p} & p^{-1}(\Psi_i)\arrow[swap]{l}{\approx}\arrow{d}{p}\\
{|}{{M}}_{\bullet}{|}\arrow{dr}{\varepsilon} & {|}\varepsilon^{-1}(\pmb{Y}_{\alpha(i)})_\bullet{|}\arrow[swap]{l}{o}\arrow{dr}{\varepsilon}& {|}\varepsilon^{-1}(\pmb{Y}_{\alpha(i)})_\bullet{|}\setminus {|}\varepsilon^{-1}(\pmb{Y}_{\alpha(i)})_\bullet\cap \partial M_\bullet{|}={|}M_\bullet{|}\setminus{|}\partial M_\bullet{|}\arrow[swap]{l}{o}\arrow{d}{\varepsilon}\arrow{r}{\cong}[swap]{|q_\bullet|} & {|}\tilde{\pmb{\Lambda}}^{({\alpha(i)})}_\bullet|\setminus |\tilde{\pmb{\Lambda}}^{({\alpha(i-1)})}_\bullet{|}\arrow[dl]\arrow{r}{\sim}[swap]{\pmb{\beta}}& \Phi_{\alpha(i)}\setminus\Phi_{{\alpha(i-1)}}\arrow{d}{e_i}& \Psi_i\arrow[swap]{l}{\approx}\arrow{dl}\\
&\bar{\pmb{Y}}_{\alpha(i)} & \pmb{Y}_{\alpha(i)}\arrow[swap]{l}{o}\arrow{rr}{\cong}[swap]{h_{\alpha(i)}} & & Y_{\alpha(i)}=X_i,&
\end{tikzcd}
\end{equation*}
where $\sim,\cong$ and $o$ have the same meaning as in (\ref{big_comm_diag}), and the arrows marked $\approx$ induce an isomorphism of Borel-Moore homology groups by theorem~\ref{thconical_1}. 

\smallskip

Let us summarise these results.

\begin{theorem}\label{thhodge3}
Suppose $X_1,\ldots,X_N$ are configuration spaces and $\alpha:\{1,\ldots, N\}\to \{1,\ldots,\N\}$ is a map such that the assumptions on page~\pageref{assumptions_conical} hold for $X_1,\ldots,X_N, Y_1, \ldots, Y_\N$ and $\alpha$. Suppose also that $X_1,\ldots,X_N$ satisfy condition \ref{five} of list~\ref{conditions} (in addition to~\ref{five}$-$). Then for every $i\in\{1,\ldots, N\}$ the following holds.

\begin{enumerate}
\item Recall that above we equipped $\bar H_*(F_{\alpha(i)}\setminus F_{{\alpha(i-1)}},\Q)$ and $\bar H_*(\Psi_i,\Q)$ with mixed Hodge structures using isomorphisms~(\ref{mhs_res}) and~(\ref{mhs_psi}). We then have
$$\bar H_*(F_{\alpha(i)}\setminus F_{{\alpha(i-1)}},\Q)\cong\bar H_*(\Phi_{\alpha(i)}\setminus \Phi_{{\alpha(i-1)}},\Q)[-2d_{\alpha(i)}]\cong \bar H_*(\Psi_i,\Q)\otimes\Q(d_{\alpha(i)})[-2d_{\alpha(i)}]$$ as mixed Hodge structures.

\item Both $\varepsilon^{-1}(\pmb{Y}_{\alpha(i)})_\bullet$ and $\varepsilon^{-1}(\pmb{Y}_{\alpha(i)})_\bullet\cap\partial M_\bullet$ are cubical complex algebraic varieties, and the latter is a closed subvariety of the former. The Borel-Moore homology Leray spectral sequence of $\Psi_i\to Y_{\alpha(i)}$ is isomorphic to the Leray spectral sequence for $\bar H_*(|\varepsilon^{-1}(\pmb{Y}_{\alpha(i)})_\bullet|,|\varepsilon^{-1}(\pmb{Y}_{\alpha(i)})_\bullet\cap\partial M_\bullet|,\Q)$ obtained from the map ${|}\varepsilon^{-1}(\pmb{Y}_{\alpha(i)})_\bullet{|}\to\pmb{Y}_{\alpha(i)}$, which is the augmentation map for both $\varepsilon^{-1}(\pmb{Y}_{\alpha(i)})_\bullet$ and $\varepsilon^{-1}(\pmb{Y}_{\alpha(i)})_\bullet\cap\partial M_\bullet$. If $\pmb{Y}_{\alpha(i)}$ is quasi-projective, then this Leray spectral sequence has a natural mixed Hodge structure starting from $E_2$ (see section~\ref{sec_mhs_sseq}).
\item Recall that for $K\in X_i=Y_{\alpha(i)}$ the mixed Hodge structure on $\bar H_*(\Lambda(K)\setminus\partial\Lambda(K),\Q)$ was introduced using isomorphism~(\ref{mhs_lambdak}). This mixed Hodge structure agrees with the one on the fibre of the map ${|}\varepsilon^{-1}(\pmb{Y}_{\alpha(i)})_\bullet{|}\to\pmb{Y}_{\alpha(i)}$ from part 2, i.e.\ if $y\in\pmb{Y}_{\alpha(i)}$ is such that $h_{\alpha(i)}(y)=K$, then we have
$$\bar H_*(\Lambda(K)\setminus\partial\Lambda(K),\Q)\cong H_*(|\varepsilon^{-1}(y)_\bullet|,|\varepsilon^{-1}(y)_\bullet\cap \partial M_\bullet|,\Q)$$ as mixed Hodge structures.
\end{enumerate}
\end{theorem}
$\clubsuit$

{\bf Remark.} As with theorem~\ref{thhodge2}, we expect that the assumption that $\pmb{Y}_{\alpha(i)}$ should be quasi-projective is not necessary in order for the last conclusion of part 2 to hold.

\section{Preliminaries on nodal hypersurfaces and vector bundles}\label{nod}

\subsection{Nodal hypersurfaces and equivariant vector bundles}

Recall that in the Introduction we have set $\Pi_{d,n}=\Gamma(\C\p^n,\mathcal{O}(d))$, or equivalently, the space of all complex homogeneous polynomials of degree $d$ in $n+1$ variables, which we will denote $\mathrm{x}_0,\ldots,\mathrm{x}_n$. For $f\in\Pi_{d,n}$ we set $\partial_i f=\frac{\partial f}{\partial\mathrm{x}_i}$ and we set $df$ to be the (column) $n+1$-tuple $(\partial_0 f,\ldots,\partial_n f)^T$. In the rest of the paper elements of $\C^{n+1}$ are denoted by Roman letters and are considered as column vectors unless stated otherwise. Elements of $\C \p^n$ will be denoted using italics.

Recall (see Introduction) that we set $\N_{d,n}$ to be the subvariety of $\Pi_{d,n}$ formed by all $f$ such that the kernel of the Hessian of $f$ at some $\mathrm{x}\in \C^{n+1}\setminus\{0\}$ contains $\mathrm{x}$ and has dimension $>1$. 
Recall also that we denote the space $\p\Pi_{d,n}\setminus \p\N_{d,n}$ of nodal degree $d$ hypersurfaces of $\C \p^n$ by $\nod_{d,n}$.

If $V=\Pi_{d,n}$ and $\Sigma$ is the set of all $f\in V$ that define singular hypersurfaces (cf.~\cite{vas3,quintics}), then it is natural to take $\mathbf{M}=\C \p^n$, and to define $\Sing(f)$ for $f\in \Sigma$ as the subvariety
$$\{x\in\C\p^n\mid f(x)=0, df|_x=0\}\subset\C\p^n.$$

In the examples we consider in the rest of the paper $V$ will again be $\Pi_{d,n}$, but the discriminant $\Sigma$ will be $\N_{d,n}$. 
In this situation we set $\mathbf{M}$ to be the flag variety $F(n+1;1,2)$ formed by all couples (a point $x\in\C \p^n$, a line through $x$), and define $\Sing(f), f\in \Sigma$ as the subvariety of $\mathbf{M}$ formed by all $(x,l)\in\mathbf{M}$ such that
the 2-plane in $\C^{n+1}$ that corresponds to $l$ is in the kernel of the Hessian matrix $\Hess_{x} f$ of $f$ at $x$. Recall that if $g$ is a homogeneous polynomial of degree $d$ and $\mathrm{x}=(\mathrm{x}_0,\ldots,\mathrm{x}_n)^T\in \C^{n+1}$, then
\begin{equation}\label{euler}
d\cdot g(\mathrm{x})=\sum_{i=0}^n\mathrm{x}_i \frac{\partial g}{\partial \mathrm{x}_i}(\mathrm{x}).
\end{equation}
It follows from this formula that if $(x,l)\in F(n+1; 1,2)$ belongs to $\Sing (f)$, then $df|_x=0$ and $f(x)=0$. It will sometimes be convenient to identify $F(n+1;1,2)$ with $\p (T \C \p^n)$, the projectivised tangent bundle of $\C \p^n$. 


Let us review the cohomology of $F(n+1;1,2)$ for $n=2$. We have $F(3;1,2)\cong \SL_3(\C)/B$ where $B$ is the Borel subgroup of all upper triangular matrices in $\SL_3(\C)$. Let $T\subset B$ be the maximal torus that consists of all diagonal matrices, and let $\mathfrak{t}$ be the tangent Lie algebra of $T$. We let $\epsilon_i:\mathfrak{t}\to\C,i=0,1,2$ be the function that takes a matrix to the $i$-th diagonal entry. These span the dual space $\mathfrak{t}^*=\mathrm{\mathop{Hom}}(\mathfrak{t},\C)$ and satisfy the relation $\sum\epsilon_i=0$. Let $\mathfrak{t}^*(\Z)\subset\mathfrak{t}^*$ be the Abelian subgroup generated by $\epsilon_0,\epsilon_1,\epsilon_2$ and set $\mathfrak{t}^*(\Q)=\mathfrak{t}^*(\Z)\otimes\Q$. 
The Weyl group $W=S_3$ of $\SL_3(\C)$ acts on $\mathfrak{t}^*$ by permuting the $\epsilon_i$'s.

Suppose $H$ is a Lie subgroup of a Lie group $G$. If $R:H\to \GL_m(\C)$ is a representation, then we set $G\times_R \C^m=(G\times \C^m)/H$, where the action of $H$ on $G\times \C^m$ is given by $h\cdot (g,v)=(gh^{-1},R(h)v), h\in H,g\in G,v\in V$. The map $G\times_R \C^m\to G/H$ is the projection of a rank $m$ equivariant complex vector bundle over $G/H$, and all such vector bundles are obtained in this way. Note that the representation of $H$ on the fibre of $G\times_R\C^m$ over $\{H\}\in G/H$ is isomorphic to $R$.

\begin{lemma}\label{alt_flag_variety}
\begin{enumerate}
\item We have $H^*(\SL_3(\C)/B,\Q)\cong \Q[\epsilon_0,\epsilon_1,\epsilon_2]/I$ where $I$ is the ideal generated by the elementary symmetric polynomials, and moreover, the isomorphism takes the first Chern class of an equivariant line bundle $\SL_3(\C)\times_\chi\C$ where $\chi:B\to\C^*$ is a character to $d\chi|_\mathfrak{t}\in\mathfrak{t}^*(\Q)$ written in terms of $\epsilon_i,i=0,1,2$.
\item Suppose $R:B\to \GL_m(\C)$ is a representation and $\vartheta=\SL_3(\C)\times_R\C^m$ is the corresponding homogeneous vector bundle. The restriction $R|_T$ decomposes as a direct sum of one-dimensional representations; let $\chi_1,\ldots , \chi_m:T\to\C^*$ be the corresponding characters of $T$. The total Chern class $c(\vartheta)$ is then
$$c(\vartheta)=\prod_{i=1}^m(1+d\chi_i).$$
\end{enumerate}
\end{lemma}

{\bf Proof.} The first part of the lemma is a particular case of the description of the rational cohomology of complete flag varieties, see~\cite[Chapter VI]{borel}; cf.\ also~\cite[Proposition 1.3]{bgg}. The second part follows from the first.$\clubsuit$

{\bf Remark.} In fact, we have an isomorphism $H^*(\SL_3(\C)/B,\Z)\cong \Z[\epsilon_0,\epsilon_1,\epsilon_2]/I_\Z$ where $I_\Z$ is again the ideal generated by the elementary symmetric polynomials. A similar result holds for $\SL_n(\C)$; for general complex semi-simple Lie groups the integral cohomology ring of the complete flag variety is more complicated, see~\cite{bgg}.
\subsection{Zero section embedding, the Euler class, and transgression}

We will need the following proposition, which is likely folklore material. We were not able to find a proof in the literature, so we give it here for completeness. All homology and cohomology groups in this section are with coefficients in a commutative ring $R$, and the orientations of all vector bundles are taken with respect to $H^*(-,R)$.

\begin{Prop}\label{thom_iso_zero_sec_prop}
Let $E$ be a real oriented vector bundle of rank $r$ over a locally finite and finite-dimensional CW-complex $X$, and let $i:X\to \tot(E)$ be the zero section embedding.  
There is an isomorphism $Th:\bar H_*(\tot(E))\to \bar H_{*-r}(X)$, called the \emph{Thom isomorphism}, and we have
\begin{equation}\label{thom_class_zero_sec}
Th(i_*(a))=a\frown e(E)
\end{equation}
for all $a\in\bar H_*(X)$, where $e(E)\in H^r(X)$ is the Euler class of $E$.
\end{Prop}

{\bf Proof.} We set $\tot^{\neq 0}(E)$ to be $\tot (E)$ minus the zero section. Let 
$p:\tot(E)\to X$ be the projection. We will use $\bar C_*$ to denote complexes of locally finite singular chains and $\#$'s to denote (co)chain maps induced by continuous maps of spaces. We introduce a metric on $E$ and set $\tot^{\leq 1}(E)\subset\tot(E)$ to be the space of all elements of length~$\leq 1$.

The Thom isomorphism $\bar H_*(\tot(E))\to \bar H_{*-r}(X)$ is constructed as follows. Let $\bar C_*^{\mathcal{U}}(\tot(E))\subset \bar C_*(\tot(E))$ be the subcomplex of locally finite chains subordinate to the cover $\mathcal{U}=\{\tot^{\neq 0}(E),\tot^{\leq 1}(E)\}$ of $\tot(E)$. The inclusion $\bar C_*^{\mathcal{U}}(\tot(E))\to \bar C_*(\tot(E))$ is a chain homotopy equivalence: the proof given e.g.\ in~\cite[2.1]{hatcher} for finite chains works verbatim. Let $u\in C^r(\tot(E),\tot^{\neq 0}(E))$ represent the Thom class of $E$. The Thom isomorphism is then induced by the chain map $\bar C_*^{\mathcal{U}}(\tot(E))\to \bar C_{*-r}(X)$ given by
\begin{equation}\label{thom_iso_chain}
s\mapsto p_\#(s\frown u).
\end{equation}

Note that the right hand side of this formula is a well-defined locally finite chain: for an $x\in X$ let $K\subset X$ be a compact neighbourhood (recall that we assume $X$ locally compact); for every $s\in \bar C_*^{\mathcal{U}}(\tot(E))$ the singular simplices that occur in $s\frown u$ have support in $\tot^{\leq 1}(E)$, so there are finitely many such simplices with support that intersects $p^{-1}(K)\cap \tot^{\leq 1}(E)$. In order to see that (\ref{thom_iso_chain}) is a quasi-isomorphism, one could e.g.\ consider the skeletal filtration on $X$ and its preimage under $p$, and observe that (\ref{thom_iso_chain}) preserves the resulting filtrations on locally finite singular chains and induces an isomorphism of the terms $E^1$ of the spectral sequences. Note that the skeletal filtration is finite as we assume $X$ finite-dimensional.

Using the projection formula we have
$$p_{\#}(i_{\#}(t)\frown u)=p_{\#}\big(i_{\#}(t\frown i^{\#}(u))\big)=t\frown i^{\#}(u)$$ for every $t\in\bar C_*(X)$. By definition, $i^{\#}(u)$ represents the Euler class $e(E)$. So if $t$ is closed and represents $a\in \bar H_*(X)$, then passing to homology we get (\ref{thom_class_zero_sec}).$\clubsuit$

\begin{corollary}\label{eulerclass}
Suppose $E$ is a real vector bundle of rank $r$ over a smooth manifold $X$, and let $E'\subset E$ be a vector subbundle of rank $r'$. If the quotient bundle $E/E'$ and the total space of $E'$ are oriented, then the following diagram commutes
$$
\begin{tikzcd}
\bar H_*(\tot(E'))\ar[d]\ar[r,"="] & \bar H_{*}(\tot (E'))\ar[r,"\cong"]\ar[d,"-\frown e(E/E')"] & H^{\dim X+r'-*}(X)\ar[d,"-\smile e(E/E')"]\\
\bar H_*(\tot(E))\ar[r,"Th"] & \bar H_{*-r+r'}(\tot (E'))\ar[r,"\cong"] & H^{\dim X+r-*}(X)
\end{tikzcd}
$$
where the left vertical arrow is induced by the inclusion $\tot (E')\subset \tot (E)$ and the right horizontal arrows are the Poincar\'e duality isomorphisms.
\end{corollary}

{\bf Proof.} If we identify $\bar H_*(\tot(E))$ with $\bar H_{*-r+r'}(\tot(E'))$ using the Thom isomorphism for $E/E'$, the map $\bar H_*(\tot(E'))\to \bar H_*(\tot(E))$ induced by the inclusion becomes $a\mapsto a\frown e(E/E')$ by proposition~\ref{thom_iso_zero_sec_prop}. Let us also identify $\bar H_*(\tot(E'))$ with $H^*(X)$ using the Poincar\'e isomorphism $H^{\dim X+r'-*}(X)\to \bar H_*(\tot(E'))$ given by capping with the fundamental class $[E']$. We have $[E']\frown (b\smile e(E/E'))=([E']\frown b)\frown e(E/E')$ for all $b\in H^*(X)$, which implies the corollary.$\clubsuit$

\smallskip

We will use the next lemma to calculate some of the differentials in the spectral sequences of the next two sections.

\begin{lemma}\label{aux_lemma_diff}
Suppose $(X,Y)$ is a CW-pair and $E_1$ and $E_2$ are real oriented vector bundles over $X$ and $Y$ respectively. Set $r_1=\rk E_1,r_2=\rk E_2$. Let $Z$ be the result of glueing $E_1|_Y$ and $E_2$ along an orientation preserving isomorphism $E_1|_Y\cong E_2'$ where $E_2'\subset E_2$ is an oriented vector subbundle. Then the following diagram commutes
$$
\begin{tikzcd}
\bar H_*(Z\setminus \tot(E_2))\arrow{r}{\cong} & \bar H_*(\tot(E_1)\setminus\tot(E_1|_Y))\arrow{d}{\cong}[swap]{Th}\arrow{r}{\partial} & \bar H_{*-1}(\tot(E_1|_Y))\arrow{d}{\cong}[swap]{Th}\arrow{r} &[7em] \bar H_{*-1}(\tot(E_2))\arrow{d}{\cong}[swap]{Th}\\
&\bar H_{*-r_1}(X\setminus Y)\arrow{r}{(-1)^{r_1}\partial}& \bar H_{*-r_1-1}(Y)\arrow{r}{-\frown e(E_2/E_2')\cdot (-1)^{r_1(r_2-r_1)}} & \bar H_{*-r_2-1}(Y).
\end{tikzcd}
$$
In this diagram, $\partial$ denotes the connecting homomorphisms, the arrows marked $Th$ are the Thom isomorphisms, and the upper right horizontal arrow is induced by the composite $E_1|_Y\cong E_2'\subset E_1$. Moreover, the composition of the top horizontal arrows is the connecting homomorphism of the pair $(Z,\tot(E_2))$.
\end{lemma}

{\bf Proof.} The left square commutes because capping chains with the Thom class of a rank $r$ vector bundle $(-1)^r$-commutes with the differential. To see that the right square commutes we first apply proposition~\ref{thom_iso_zero_sec_prop} to $E_2$ viewed as a vector bundle over $\tot(E_2')$ and then use the Thom isomorphism for $E_2'$, see formula~(\ref{thom_iso_chain}). The last statement follows from the fact that the proper map of pairs $(\tot(E_1),\tot(E_1|_Y))\to (Z,\tot (E_2))$ induces a map of the long exact sequences of Borel-Moore homology groups.$\clubsuit$

\medskip

In the following proposition we review the transgression in the Leray spectral sequence of the spherisation of a vector bundle.
\begin{Prop}
\label{euler-thom}
Let $E$ be a real rank $r\geq 2$ oriented vector bundle over a CW-complex $X$. As in the proof of proposition~\ref{thom_iso_zero_sec_prop}, we set $\tot^{\neq 0}(E)$ to be the total space of $E$ minus the zero section. For $x\in X$ we let $E_x$, respectively $E_{0,x}$ be the fibre of $p:\tot(E)\to X$, respectively $p|_{\tot^{\neq 0}(E)}:\tot^{\neq 0}(E)\to X$ over $x$.

Let $(E^{p,q}_r,d_r)$ be the Leray spectral sequence of $p|_{\tot^{\neq 0}(E)}$ for cohomology with coefficients in $R$. The group $E_2^{0,r-1}$ contains a unique element $\alpha$ that goes, for every $x\in X$, to the restriction of the Thom class under the map
$$E_2^{0,r-1}\to H^{r-1}(E_{0,x})\cong H^r(E_x,E_{0,x}).$$
This element transgresses to $d_r(\alpha)=-e(E)\in E^{r,0}_r\cong H^r(X)$.
\end{Prop}

{\bf Proof.} 
Let $s_0:X\to E$ be the zero section embedding and let $u\in C^r(\tot(E),\tot^{\neq 0}(E))$ be a singular cochain representing the Thom class of $E$. We may and will assume that for every $x\in X$ the pullback of $u$ under $\{x\}\subset X\xrightarrow{s_0} \tot(E)$ is zero. The Euler class $e(E)$ of $E$ is represented by the cochain $\varepsilon=s_0^*(u)$, so we have $s_0^*(u)-\varepsilon=0$. Since $p:\tot(E)\rightarrow X$ is a homotopy inverse of $s_0$, the cochain $u-p^*(\varepsilon)$ must be zero in cohomology, so $u-p^*(\varepsilon)=\delta\omega$ for some $\omega\in C^{r-1}(\tot(E))$. Note that for every $x\in X$, we have $\delta\omega|_{E_x}=u|_{E_x}$, which is zero when restricted to $E_{0,x}$ and represents the orientation class of $E_x$. So $\omega|_{E_{0,x}}$ is closed and represents the preferred generator $\alpha_x\in H^{r-1}(E_{0,x})$, namely the image of $\alpha$. And given that 
$$
\delta\omega|_{\tot^{\neq 0}(E)}=u|_{\tot^{\neq 0}(E)}-p^*(\varepsilon)|_{\tot^{\neq 0}(E)}=-p^*(\varepsilon)|_{\tot^{\neq 0}(E)},
$$
we conclude that the transgression $d_r(\alpha)$ is $-e(E)$.$\clubsuit$


\section{Nodal cubic curves in \texorpdfstring{$\C \p^2$}{}}\label{plane_cubics}

We now apply the method of section \ref{conres} to the simplest non-trivial example, namely the space of equations of cubic curves in $\C \p^2$.

\subsection{The space of equations}\label{eq32}

We set $V=\Pi_{3,2}$ and $\Sigma=\N_{3,2}$, and we take the flag variety $F(3;2,1)\cong \p T(\C \p^2)$ as $\mathbf{M}$. Let us define $\Sing(f),f\in\Sigma$ as described in the beginning of section~\ref{nod}: for each singular point $x\in \C \p^2$ of the curve given by $f=0$ we include in $\Sing(f)$ all lines $l\in P (T_x\C \p^2)$ such that $l\subset$ the kernel of $\Hess_x(f)$. 
Below we schematically represent subsets of $\p T(\C \p^2)$ as follows: the locus of points $x$ appearing in $\p (T_x \C \p^2)$ is shown in black, and the set of tangent directions taken at those points is shown in grey.  

\begin{Prop}\label{sing_plane_cub}
The spaces $X_1,\ldots,X_6$ whose generic elements are described below satisfy conditions \ref{first}-\ref{four}, \ref{five}$+$,\ref{last} and \ref{xxx} from list \ref{conditions} with the numbers in square brackets being the $d_i$'s from condition~\ref{three}.
\begin{enumerate}
\item A one-element subset of $\mathbf{M}$ [5]. 
\newline
\includegraphics[scale=0.3]{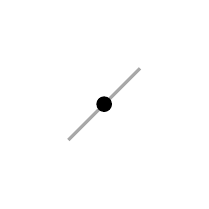}
\item A subset of $\mathbf{M}$ formed by all $(x_0,l)\in\mathbf{M}$ where $x_0$ is a fixed point [4].
\newline
\includegraphics[scale=0.3]{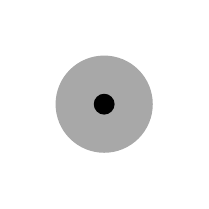}
\item A subset of $\mathbf{M}$ formed by all $(x,l_0)\in\mathbf{M}$ where $l_0$ is a fixed line [3].
\newline
\includegraphics[scale=0.3]{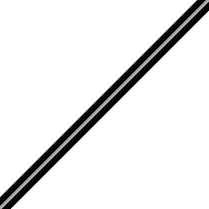}
\item A subset of $\mathbf{M}$ formed by all $(x,l)$ such that $x=x_0$ or $l=l_0$ where $x_0$ and $l_0\ni x_0$ are a fixed point and a fixed line respectively [2].
\newline
\includegraphics[scale=0.3]{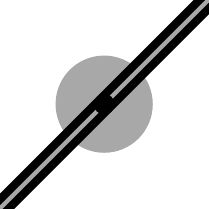}
\item A subset of $\mathbf{M}$ formed by all $(x,l)$ such that $x\in l_0$ where $l_0$ is a fixed line [1].
\newline
\includegraphics[scale=0.3]{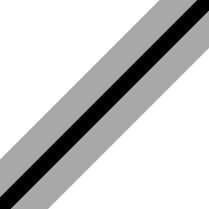}
\item The whole $\mathbf{M}$ [0].
\newline
\includegraphics[scale=0.3]{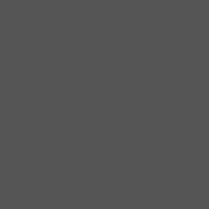}
\end{enumerate}
\end{Prop}

{\bf Remark.} Note that $X_3$ does not intersect the closure of $\{\Sing(f)\mid f\in\N_{3,2}\}$ in $2^{\mathbf{M}}$. Nevertheless, we include $X_3$ in our list, because otherwise the spectral sequence of theorem \ref{spseqcvector} would be more complicated.

{\bf Proof.} A plane cubic curve has singularities worse than simple nodes if and only if it is of one of the following types, illustrated below: an irreducible cuspidal cubic; the union of a conic and a tangent line; the union of three distinct lines through a point; the union of a double line $l'$ and a line $l''\neq l'$; a triple line. 

\includegraphics[scale=0.67]{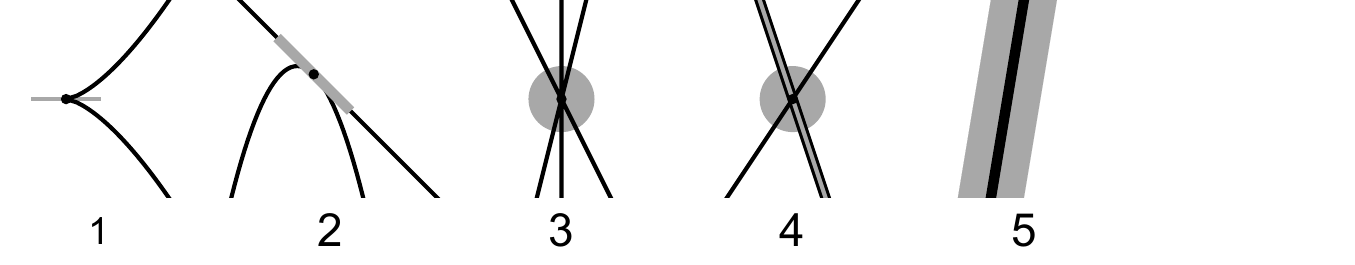}

The corresponding singular locus belongs to $X_1,X_1, X_2, X_4$ and $X_5$ respectively, which proves condition \ref{first}. The remaining conditions are easy to check directly.$\clubsuit$


Note that since the spaces $X_1,\ldots, X_6$ satisfy condition \ref{five}$+$, we have $\Geom (K)=K$ for every $K\in\bigcup X_i$. Set $Y_i=X_i,i=1,\ldots, 6.$ Using the results of section \ref{conres} we construct a conical resolution $\sigma$ of $\Sigma$ and a filtration on it. The first assertion of Theorem~\ref{mhcubp2} is a consequence of the following theorem.

\begin{theorem}\label{spseqcvector}
Let $(E^r_{p,q}, d^r_{p,q})$ be the resulting spectral sequence that converges to $\bar H_*(\sigma,\Q)\cong\bar H_*(\Sigma,\Q)$. The groups $(E^1_{p,q})$ are given by table \ref{spseqcub}.
\begin{table}
\caption{Spectral sequence for $\sigma$ in the case of plane nodal cubics}
\label{spseqcub}
$$\begin{array}{c|c|c|c|c|c|c|}
&1&2&3&4&5&6\\
\hline
15&\Q(8)&&&&&\\
\hline
14&&&&&&\\
\hline
13&\Q(7)^2&\Q(7)&&&&\\
\hline
12&&&&&&\\
\hline
11&\Q(6)^2&\Q(6)&&&&\\
\hline
10&&&\Q(6)&&&\\
\hline
9&\Q(5)&\Q(5)&&&&\\
\hline
8&&&\Q(5)&&&\\
\hline
7&&&&&&\\
\hline
6&&&\Q(4)&&&\\
\hline
5&&&&&&\\
\hline
4&&&&&&\\
\hline
3&&&&&&\\
\hline
2&&&&&&\Q(3)\\
\hline
\end{array}$$
\end{table}	
The differentials $d^1_{2,9},d^1_{2,11},d^1_{2,13}$ and $d^2_{3,10}$ are non-zero.
\end{theorem}

Note that for every $K\in\bigcup X_i$ the spaces $\Lambda(K)$ and $\tilde\Lambda (K)$ coincide, as do $\partial\Lambda(K)$ and $\tilde\partial\Lambda(K)$.

{\bf Proof.}  We first introduce some notation we will need in the proof. Let $X$ be a topological space. By writing $X\cong Y\ltimes Z$ we will mean that $X$ is homeomorphic to a locally trivial bundle with base $Y$ and fibre $Z$. Recall that we denote the cone over a space $X$ by $Cone(X)$. Note that if $X$ is a compact CW-complex, then we have
$$\bar H_*(Cone(X)\setminus X)\cong H_*(Cone(X),X)\cong \tilde H_{*-1}(X).$$
If the top non-zero group $\bar H_i(X,\Q)$ is $\Q$ has a preferred generator, then we will denote it by $[X]$. If $A_*$ is a graded Abelian group and $k\in\Z$, then we define another graded Abelian group $A_*[k]$ by $A_i[k]=A_{k+i}$.

We will now calculate in turn all columns and differentials of the spectral sequence.

{\bf Column 1.} We have $\Phi_1=X_1=\mathbf{M}\cong \p T \C \p^2$ and $F_1$ is a complex vector bundle over $\Phi_1$ of rank 5. Hence the groups $E^1_{1,*} = \bar{H}_{1+*}(F_1)$ are as shown in table~\ref{spseqcub}.

\medskip

{\bf Column 2.}\label{column2planecubbegin} The space $\Phi_2\setminus \Phi_1$ is a fibre bundle over $X_2\cong\C\p^2$ with the fibre over $x\in\C\p^2$ being an open cone $Cone^\circ(\p T_x\C\p^2\cong\C\p^1)$. Moreover, $F_2\setminus F_1$ is a complex vector bundle over $\Phi_2 \setminus \Phi_1$ of rank 4, so the groups $E^1_{2,*}$ are as shown in table~\ref{spseqcub}.

\begin{lemma}\label{diff_2nd_column_plane_cubics}
None of the differentials that originate in the second column is zero.
\end{lemma}

{\bf Proof.} Let $p:\p T\C \p^2\to\C \p^2$ be the bundle projection. 
Let $\xi_{d,n}$ and $\vartheta_{d,n}$ be the vector bundles over over $F(n+1;1,2)$ with total spaces $$\tot\xi_{d,n}=\{(f,(y,l))\in\Pi_{d,n}\times F(n+1;2,1)\mid\Hess_y f=0\}\mbox{ and }$$
$$\tot(\vartheta_{d,n})=\left\{ (f,(x,l))\in\Pi_{d,n}\times F(n+1;1,2)\mid l\subset \ker(\Hess_x(f))\right\}.$$\label{thetabundle}
respectively. Note that the fibre of $\xi_{3,2}$ over $(y,l)\in F(3; 1,2)$ is the vector subspace $L(K)\subset\Pi_{3,2}$ where $K\in X_2$ is the configuration that corresponds to $y$.  Note also that $\xi_{d,n}$ is a subbundle of $\vartheta_{d,n}$, and if $n=2$, the corank is~1. 
%
The differentials $d^1_{2,13},d^1_{2,11}$ and $d^1_{2,9}$ of the spectral sequence are the connecting homomorphisms in the long exact sequence of the pair $(F_2,F_1)$. 
Applying lemma~\ref{aux_lemma_diff} we see that to calculate these it would suffice to determine the map 
\begin{equation}\label{cubics_aux_comp}
H_*(\C\p^2,\Q)\xrightarrow{\cong} \bar H_{*+3}(\Phi_2\setminus\Phi_1)\rightarrow H_{*+2}(\Phi_1)\xrightarrow{\cong} H_{*+2}(\p T\C \p^2,\Q) \rightarrow H_*(\p T\C \p^2,\Q)
\end{equation}
where the last arrow is capping with $e(\vartheta_{3,2}/\xi_{3,2})$. A straightforward check shows that the composite of the first three arrows is $p^!$, the extraordinary pullback induced by $p$. 
After passing to cohomology and applying corollary~\ref{eulerclass} map (\ref{cubics_aux_comp}) becomes
\begin{equation}\label{cubics_aux_comp_1}
H^*(\C \p^2,\Q)\to H^*(\p T\C \p^2,\Q)\to H^{*+2}(\p T\C \p^2,\Q)
\end{equation}
where the first map is $p^*$ and the second is multiplication by $e(\vartheta_{3,2}/\xi_{3,2})=c_1(\vartheta_{3,2}/\xi_{3,2})$.

So we need to calculate $c_1(\vartheta_{3,2}/\xi_{3,2})$. We use lemma~\ref{alt_flag_variety} and the notation of this lemma. The Borel subgroup~$B$ is the stabiliser of the couple $(x_0,l_0)$ where $x_0=(1:0:0)$ and $l_0$ is the projectivisation of the plane in $\C^3$ spanned by $(1,0,0)$ and $(0,1,0)$. For $s,t\in \C^*$ set $T(s,t)=\mathbf{diag}((st)^{-1},s,t)$. The fibre of $\vartheta_{3,2}/\xi_{3,2}$ is spanned by $\mathrm{x}_0\mathrm{x}_2^2$, and we have $T(s,t)\cdot \mathrm{x}_0\mathrm{x}_2^2=(st\mathrm{x}_0)(t^{-1}\mathrm{x}_2)^2=st^{-1}\mathrm{x}_0\mathrm{x}_2^2$. The differential of the resulting character of $T$ is $\epsilon_1-\epsilon_2$.

Recall that $\C\p^2\cong\SL_3(\C)/P$ where $P\subset\SL_3(\C)$ is the stabiliser of the line spanned by $(1,0,0)$ in the standard representation. The tautological bundle $\OO(-1)$ on $\C\p^2$ is obtained from the character of $P$ that takes a matrix to its top left entry. If we restrict this character to $T$ and take the differential, we get $\epsilon_0$, which corresponds to $c_1\big(p^*(\OO(-1))\big)$ under the isomorphism of lemma~\ref{alt_flag_variety}.

So under the identification $H^*(\p T\C \p^2,\Q)\cong H^*(\SL_3(\C)/B,\Q)\cong\Q[\epsilon_0,\epsilon_1,\epsilon_2]/I$, see lemma~\ref{alt_flag_variety}, the image of the first arrow of (\ref{cubics_aux_comp}) is spanned by $1,\epsilon_0,\epsilon_0^2$, and the second arrow is multiplication by $\epsilon_1-\epsilon_2$. An easy check using e.g.\ Singular~\cite{singular} shows that none of the elements
\begin{equation}\label{image_diff_2nd_col}
\epsilon_1-\epsilon_2,\epsilon_0(\epsilon_1-\epsilon_2),\epsilon_0^2(\epsilon_1-\epsilon_2)\in\Q[\epsilon_0,\epsilon_1,\epsilon_2]/I
\end{equation}
is zero, which implies lemma~\ref{diff_2nd_column_plane_cubics}.$\clubsuit$
\medskip

{\bf Column 3.} The third column of the spectral sequence can be calculated in a similar way to the second one, the only significant difference being that $F_3\setminus F_1$ is a complex vector bundle over $\Phi_3\setminus \Phi_1$ of rank 3 rather than 4.

\begin{lemma}\label{diff_3d_column_plane_cubics}
The differential $d_{3,10}$ is non zero.
\end{lemma}

{\bf Proof.} Let $p^\vee:\p T\C\p^2\cong F(3;1,2)\to\C\p^{2\vee}$ be projection given by $(x,l)\mapsto l$. To calculate the differential $d_{3,10}$ we note that $F=F_1\sqcup (F_3\setminus F_2)$ is closed in $\sigma$. The connecting homomorphisms in the long exact sequence of the pair $(F,F_1)$ can be calculated in the same way as for the pair $(F_2,F_1)$ in the proof of lemma~\ref{diff_2nd_column_plane_cubics}, with the following differences. Instead of the vector bundle $\xi_{3,2}$ we use another bundle, $\xi'$, with total space equal the set of all couples  $(f,(y,l))\in\Pi_{3,2}\times F(3;1,2)$ such that $f$ is divisible by the square of an equation of $l$. Same as $\xi_{3,2}$, this is a subbundle of~$\vartheta_{3,2}$. 

We also replace (\ref{cubics_aux_comp}) by
$$H_*(\C \p^{2\vee},\Q)\to H_{*+2}(\p T\C \p^2,\Q)\to H_{*-2}(\p T\C \p^2,\Q),$$
where the first arrow is $(p^\vee)^!$ and the second is capping with $e(\vartheta_{3,2}/\xi')$. So (\ref{cubics_aux_comp_1}) gets replaced by
\begin{equation}\label{aux_cubics_col3}
H^*(\C \p^{2\vee},\Q)\to H^*(\p T\C \p^2,\Q)\to H^{*+4}(\p T\C \p^2,\Q)
\end{equation}
where the first map is $(p^\vee)^*$ and the second is multiplication by $e(\vartheta_{3,2}/\xi')=c_2(\vartheta_{3,2}/\xi')$.

Let $(x_0,l_0)$ be as in the proof of lemma~\ref{diff_2nd_column_plane_cubics}. The fibre of $\vartheta_{3,2}/\xi'$ over $(x_0,l_0)$ is spanned by $\mathrm{x}_0\mathrm{x_2}^2$ and $\mathrm{x}_0\mathrm{x_2}^2$. The lines spanned by these monomials are invariant under $T$, and the differentials of the corresponding characters of $T$ are $\epsilon_1-\epsilon_2$ and $-\epsilon_1-2\epsilon_2$ respectively, where we use the notation of lemma~\ref{alt_flag_variety}. Using the isomorphism of part 1 of this lemma we can write the second arrow of (\ref{aux_cubics_col3}) as multiplication by $(\epsilon_1-\epsilon_2)(-\epsilon_1-2\epsilon_2)$. Moreover, the image of the first arrow is spanned by $1,c,c^2$ where $c=c_1\big((p^{\vee})^*(\OO_{\C\p^2}(-1))\big)$. To prove lemma~\ref{diff_3d_column_plane_cubics} it remains to check that the elements
$$(\epsilon_1-\epsilon_2)(-\epsilon_1-2\epsilon_2)\mbox{ and } \epsilon_0(\epsilon_1-\epsilon_2)\in\Q[\epsilon_0,\epsilon_1,\epsilon_2]/I$$
(cf.\ (\ref{image_diff_2nd_col})) are linearly independent, which is again easy to do using Singular~\cite{singular}.$\clubsuit$

\medskip

{\bf Column 4.} A configuration $K\in X_4$ is a union $K'\cup K''$ where $K'\in X_2$ and $K''\in X_3$. The intersection $\Lambda(K')\cap\Lambda(K'')=\Lambda(K'\cap K'')$ is a point, so the space $\partial\Lambda(K)=\Lambda(K')\cup\Lambda (K'')$ is contractible, and we have $\tilde H_*(\partial\Lambda(K),\Q)=0$, which implies that the fourth column of the spectral sequence is zero.

\medskip

{\bf Column 5.} For a configuration $K\in X_5$ let $K_1$ be the unique configuration in $X_3$ such that $K_1\subset K$. A configuration $K'\in\bigcup X_i$ that is strictly between $K$ and $K_1$ must belong to $X_4$. Moreover, for any two configurations $K',K''\in X_4$ such that $K_1\subsetneq K'\subsetneq K, K_1\subsetneq K''\subsetneq K$ and $K'\neq K''$ we have $\Lambda(K')\cap \Lambda(K'')=\Lambda(K_1)$. So
$$\partial\Lambda(K)\setminus\Lambda(K_1)=\bigsqcup_{\mbox{\tiny $\begin{array}{c}K_1\subsetneq K'\subsetneq K\\K'\in X_4\end{array}$}}\Lambda(K')\setminus\Lambda(K_1).$$ Every space in the disjoint union in the right hand side of this equality has zero Borel-Moore homology, which implies that the fifth column of the spectral sequence is zero.

\medskip

{\bf Column 6.} The space $X_6$ contains just one element, $\mathbf{M}$. The space $\partial\Lambda(\mathbf{M})$ is the union of $\Lambda(K)$ for all $K\in X'=\bigcup_{1\leq i\leq 5} X_i$, or, in other words, the union of all coherent simplices with all vertices in $X'$. Let $\triangle$ be such a simplex. There is a smallest coherent simplex $\triangle'$ such that $\triangle'\supset\triangle$, all vertices of $\triangle'$ are in $X'$ and at least one vertex is in $X_2\cup X_5$. This implies that $\partial\Lambda(\mathbf{M})$ is homotopy equivalent to the {\it coherent join} of $X_2\cong\C \p^2$ and $X_5\cong\C \p^{2\vee}$, i.e.~the union of all segments $\subset \partial\Lambda(\mathbf{M})$ that join a configuration $K\in X_2$ and a configuration $K'\in X_5$ that contains $K$.

Let $Z$ be the coherent join of $\C \p^2$ and $\C \p^{2\vee}$. The homology of $Z$ can be calculated in the same way as for the ordinary join. There are open subsets $U,V\subset Z$ such that $U\sim \C\p^2,V\sim \C \p^{2\vee}$ and $U\cap V\sim F(3;1,2)$. The resulting Mayer-Vietoris sequence reads
\begin{equation}\label{mvcubicslastcolumn}
\cdots \to H_*(F(3;1,2),\Q)\to H_*(\C\p^2,\Q)\oplus H_*(\C\p^{2\vee},\Q)\to H_*(Z,\Q)\\ \to H_{*-1}(F(3;1,2),\Q)\to\cdots.
\end{equation}

Recall that we denote the projections $F(3;1,2)\to\C\p^2$ and $F(3;1,2)\to\C\p^{2\vee}$ by $p$, respectively $p^\vee$, see columns 2 and 3. We claim that the map $$p_*+p^\vee_*:H_*(F(3;1,2),\Q)\to H_*(\C\p^2,\Q)\oplus H_*(\C\p^{2\vee},\Q)$$ is an isomorphism in degrees 2 and 4. To prove this it would suffice to prove the same for cohomology. We will use the description of $H^*(F(3;1,2),\Q)$ given in lemma~\ref{alt_flag_variety}. In the proof of lemma~\ref{diff_2nd_column_plane_cubics} we saw already that $p^*\big(c_1(\OO_{\C\p^2}(-1))\big)=\epsilon_0$. Let us calculate $(p^\vee)^*\big(c_1(\OO_{\C\p^{2\vee}}(-1))\big)$. We have $\C\p^{2\vee}=\SL_3/P'$ where $P'\subset\SL_3(\C)$ is the stabiliser of the line $l_0$ through $(1:0:0)$ and $(0:1:0)$. The bundle $\OO_{\C\p^{2\vee}}(-1)$ has total space $$\{(f,l)\in \Pi_{1,2}\times\C\p^{2\vee}\mid f|_l=0\}.$$ The corresponding character of $P'$ takes a matrix to its bottom right entry. Restricting to $T$ and taking the differential we get $\epsilon_2$, which by lemma~\ref{alt_flag_variety} is $(p^\vee)^*\big(c_1(\OO_{\C\p^{2\vee}}(-1))\big)$. A straightforward check shows that $\epsilon_0$ and $\epsilon_2$ are linearly independent in $\Q[\epsilon_0,\epsilon_1,\epsilon_2]/I$, and so are $\epsilon_0^2$ and $\epsilon_2^2$, cf.\ the proofs of lemmas~\ref{diff_2nd_column_plane_cubics} and~\ref{diff_3d_column_plane_cubics}. This proves our claim about $p_*+p^\vee_*$, and we conclude that column 6 of the spectral sequence is as claimed in theorem~\ref{spseqcvector}, which is now proved.$\clubsuit$

\smallskip

We have now proved the first part of theorem~\ref{mhcubp2}, namely the assertion about the additive structure of $H^*(\Pi_{3,2}\setminus\N_{3,2},\Q)$.

\subsection{The space of curves and the cup product}

We will now calculate the cohomology groups of the space of nodal cubics in $\C \p^2$ and prove the remaining statements of theorem~\ref{mhcubp2}. 

\begin{Prop}
\label{nodstart}
Suppose $d\geq 2, n\geq 1$. The first few rational cohomology groups of the space $\nod_{d,n}$ (see Introduction) are as follows:
$$H^0(\nod_{d,n})=\Q(0),\ H^2(\nod_{d,n})=\Q(-1),\ H^1(\nod_{d,n})=H^3(\nod_{d,n})=0.$$
\end{Prop}

\textbf{Proof.} Set $N=\dim_\C \Pi_{d,n}$. Observe that $\N_{d,n}$ is an irreducible subvariety of $\Pi_{d,n}$ of codimension 2, and similarly, the projectivisation $\p\N_{d,n}$ is an irreducible subvariety of $\p(\Pi_{d,n})$ of codimension 2. So the top relative homology groups $H_*( \p(\Pi_{d,n}),\p\N_{d,n},\Q)$ are $\Q(N)$ in degree $2N$, $\Q(N-1)$ in degree $2N-2$ and 0 in degrees $2N-1$ and $2N-3$. The proposition now follows from the Poincar\'e duality $$H^*(\p(\Pi_{d,n})\setminus \p\N_{d,n},\Q)\cong H_{2N-*}(\p\N_{d,n}, \p(\Pi_{d,n}),\Q).$$$\clubsuit$

\medskip

Recall that if a topological group $G$ acts continuously on a space $X$, then the {\it homotopy quotient} or {\it Borel quotient} is the space $$X_{hG}=(X\times EG)/G.$$ Here $EG$ is a contractible space of which $G$ acts continuously and freely on the right in such a way that the map $EG\to EG/G$ is the projection of a locally trivial bundle, and the action of $G$ on $X\times EG$ is given by $$g\cdot (x,e)=(gx,eg^{-1}),g\in G,x\in X,e\in EG.$$

The {\it classifying space} $BG$ is $EG/G$. Note that there is a fibre bundle
\begin{equation}\label{bundlehoquot}
X\to X_{hG}\to BG,
\end{equation}
see section~\ref{sec_equi}, where we review the facts about equivariant cohomology that we will use in this paper. We note in particular that if the stabiliser $G_x$ of every $x\in X$ is finite, then, under some technical but not very onerous assumptions there is an isomorphism
\begin{equation}\label{isofinstab}
H^*(X_{hG},\Q)\cong H^*(X/G,\Q),
\end{equation}
see proposition~\ref{rational_equi_quotient}.

\smallskip

Let us now take $X=\Pi_{3,2}\setminus \N_{3,2}$ and $G=\GL_1(\C)=\C^*$. Using the cohomology groups of $X$ we computed in the section~\ref{eq32} we get the second page of the spectral sequence of $X_{hG}\to BG=\C \p^\infty$:

\begin{equation}
\label{BG32}
\begin{array}{c||c|c|c|c|c|c|c|c}
\hline
11&\Q(-7)&&\Q(-8)&&\Q(-9)&&\Q(-10)&\cdots\\
\hline
10&\Q(-6)&&\Q(-7)&&\Q(-8)&&\Q(-9)&\cdots\\
\hline
9&&&&&&&&\\
\hline
8&\Q(-5)&&\Q(-6)&&\Q(-7)&&\Q(-8)&\cdots\\
\hline
7&&&&&&&&\\
\hline
6&&&&&&&&\\
\hline
5&\Q(-3)&&\Q(-4)&&\Q(-5)&&\Q(-6)&\cdots\\
\hline
4&&&&&&&&\\
\hline
3&\Q(-2)&&\Q(-3)&&\Q(-4)&&\Q(-5)&\cdots\\
\hline
2&&&&&&&&\\
\hline
1&&&&&&&&\\
\hline
0&\Q(0)&&\Q(-1)&&\Q(-2)&&\Q(-3)&\cdots\\
\hline
\hline
&0&1&2&3&4&5&6&\cdots\\
\end{array}
\end{equation}

By proposition~\ref{nodstart}, the total space has zero cohomology in degree 3, so the differential $d_4^{0,3}$ is non zero and therefore all differentials $d_4^{2n,3}$ are non zero. The cohomology of the quotient space $X_{hG}$ is finite dimensional, so we have two alternatives: 
\begin{enumerate}
\item The differentials $d^2_{0,11}$ and $d^4_{0,8}$ are non-zero. Then row 11 cancels with a tail of row 10, and row 8 with a tail of row 5.
\item We have $d^2_{0,11}=0$, but $d^4_{0,11}\neq 0$ and $d^6_{0,10}\neq 0$. Then row 11 cancels with a tail of row 8, and row 10 with a tail of row 5.
\end{enumerate}
We now rule out alternative 2. It implies that second term of the spectral sequence of the quotient map $p:\Pi_{3,2}\setminus \N_{3,2}\to\nod_{3,2}$ is as follows:
$$
\begin{array}{c||c|c|c|c|c|c|c|c|c|c|c|}
\hline
1& \Q(-1)&&\Q(-2)&&&\Q(-4)&&\Q(-5)&\Q(-6)&\Q(-6)&\Q(-7)\\
\hline
0& \Q(0)&&\Q(-1)&&&\Q(-3)&&\Q(-4)&\Q(-5)&\Q(-5)&\Q(-6)\\
\hline
\hline
& 0&1&2&3&4&5&6&7&8&9&10\\
\end{array}$$

Let $a_i\in E_2^{i,0}\cong H^i(\nod_{3,2},\Q), i=0,2,5,7,\ldots,10$ and $b\in E_2^{0,1}$ be non-zero elements. We have $d_2(b)=a_2$ because we know (see theorem~\ref{spseqcvector}) that $H^1(\Pi_{3,2}\setminus \N_{3,2},\Q)=0$. From the multiplicative structure of spectral sequence~(\ref{BG32}) we see that after rescaling $a_7$ if necessary we have $a_2a_5=a_7$, and so $d_2(ba_7)=a_2a_7=a_2^2a_5=0$. We obtain that $H^8(\Pi_{3,2}\setminus \N_{3,2},\Q)$ has dimension 2, whereas by theorem~\ref{spseqcvector} we have $H^8(\Pi_{3,2}\setminus \N_{3,2},\Q)\cong\Q(-5)$.

We deduce that alternative 1 holds, so the $E_2$ term of the spectral sequence of $p:\Pi_{3,2}\setminus \N_{3,2}\to\nod_{3,2}$ is
\begin{equation}
\label{spsec_cubics_eq_to_curves}
\begin{array}{c||c|c|c|c|c|c|c|c|c|c|c|}
\hline
1& \Q(-1)&&\Q(-2)&&&\Q(-4)&&\Q(-5)& & &\Q(-7)\\
\hline
0& \Q(0)&&\Q(-1)&&&\Q(-3)&&\Q(-4)& & &\Q(-6)\\
\hline
\hline
& 0&1&2&3&4&5&6&7&8&9&10\\
\end{array}
\end{equation}
Note that this implies the assertion of theorem~\ref{mhcubp2} about the additive structure of $H^*(\nod_{3,2},\Q)$. As above, we let $a_i\in E_2^{i,0}\cong H^i(\nod_{3,2},\Q), i=0,2,5,7,10$ and $b\in E_2^{0,1}$ be non-zero elements and deduce that $a_2a_5\neq 0$, which implies that for $\nod_{3,2}$ the cup product $H^2\otimes H^5\to H^7$ is non zero. We also have $d_2(a_2b)=d_2(a_5)=0$ and $a_2ba_5\neq 0$, so the product $H^3\otimes H^5\to H^8$ for $\Pi_{3,2}\setminus \N_{3,2}$ is also non zero. The product $H^3\otimes H^8\to H^{11}$ for $\Pi_{3,2}\setminus \N_{3,2}$ is zero because the square of an element of $H^3$ is zero.

It follows from spectral sequence (\ref{spsec_cubics_eq_to_curves}) that the rational cohomology map $p^*$ induced by $p$ is an isomorphism in degrees $0,5$ and $10$; in all other degrees $p^*=0$ for dimension reasons. We also note that $H_2(\nod_{3,2},\Q)$ is spanned by the fundamental class of any line $l\subset\p\Pi_{3,2}$ that does not intersect $\p\N_{3,2}$. By taking $l$ such that the map $q:\nod_{3,2}\to\nod_{3,2}//\PGL_3(\C)$ is not constant on $l$ we see that $q$ induces a surjective map in rational homology.

It remains to calculate the cohomology maps induced by the orbit maps $\GL_3(\C)\to\Pi_{3,2}\setminus\N_{3,2}$ and $\SL_3(\C)\to\nod_{3,2}$. Recall that all orbit maps constructed using a continuous action of a topological group on a path connected space are homotopic. We use the polynomial $f_0=\mathrm{x}_0\mathrm{x}_1\mathrm{x}_2\in\Pi_{3,2}\setminus\N_{3,2}$ and the corresponding curve $C_0\in\C\p^2$, which is a union of three lines that do not pass through a single point. The stabiliser of $f_0$, respectively of $C_0$ in $\GL_3(\C)$ contains a maximal torus $T$ of $\SL_3(\C)$, respectively a maximal torus $T'$ of $\GL_3(\C)$. We have a commutative diagram
$$
\begin{tikzcd}
\GL_3(\C)\arrow[r]\arrow{d}[swap]{=} & \GL_3(\C)/T \arrow[d]\arrow[r] & \Pi_{3,2}\setminus\N_{3,2}\arrow{d}{p}\\
\GL_3(\C) \arrow[r]& \GL_3(\C)/T' \arrow[r] & \nod_{3,2}
\end{tikzcd}
$$
where the horizontal rows are the orbit maps of $f_0$ and $C_0$.
Note that $\GL_3(\C)/T'$ is homotopy equivalent to the flag variety $F(3;2,1)$, which has cohomology in degrees $0,2,4,6$, so $H^{>7}(\GL_3(\C)/T,\Q)=0$. Using these observations and the fact that $p^*:H^5(\nod_{3,2},\Q)\to H^5(\Pi_{3,2}\setminus\N_{3,2},\Q)$ is an isomorphism we conclude from the diagram that all classes in $H^{>0}(\Pi_{3,2}\setminus\N_{3,2},\Q)$ and $H^{>0}(\nod_{3,2},\Q)$ pull back to $0$ except possibly a generator $a\in H^3(\Pi_{3,2}\setminus\N_{3,2},\Q)$.

To calculate the pullback of $a$ we use the orbit map of $f_0$ for the action of $SL_2(\C)\subset\GL_3(\C)$ on $\Pi_{3,2}\setminus\N_{3,2}$, which factorises as $\SL_2(\C)\to \SL_2(\C)/T''\to \Pi_{3,2}\setminus\N_{3,2}$ where $T''=T\cap\SL_2(\C)$. The space $\SL_2(\C)/T''$ is homotopy equivalent to $\C\p^1$, so $a$ pulls back to $0$ in $H^3(\SL_2(\C),\Q)$ and hence also in $H^3(\GL_3(\C),\Q)\cong H^3(\SL_2(\C),\Q)$. This completes the proof of theorem~\ref{mhcubp2}.$\clubsuit$


\section{Nodal quartics in \texorpdfstring{$\C \p^2$}{}}\label{quartics}
In this section we prove theorem~\ref{mhquartp2}.

\subsection{The space of equations}
\label{eq42}
We begin by applying the construction of section \ref{conres} to calculate the rational cohomology of the space of nodal quartics in $\C \p^2$. We set $V$, respectively $\Sigma$, to be $\Pi_{4,2}$, respectively, $\N_{4,2}$. Define $\mathbf{M}$ and $\Sing(f),f\in\Sigma$ 
as described in the beginning of section~\ref{plane_cubics} (see also the beginning of section~\ref{nod}). In particular, we have
$$\mathbf{M}=F(3;2,1)=\{(x,l)\mid\mbox{$x\in\C \p^2$ and $l\subset\C \p^2$ is a line through $x$}\}.$$

\begin{Prop}\label{sing_plane_quart}
The subspaces $X_1,\ldots,X_{\ref{qlast}}\subset2^\mathbf{M}$ whose generic elements are described below satisfy conditions \ref{first}-\ref{xxx} from list \ref{conditions}, the numbers in square brackets being the $d_i$'s from condition~\ref{three}.
\begin{enumerate}
\item\label{q1} A one-element subset of $\mathbf{M}$ [10].
\newline
\includegraphics[scale=0.3]{one.pdf}
\item\label{q2} A subset of $\mathbf{M}$ formed by all $(x_0,l)\in\mathbf{M}$ where $x_0$ is a fixed point [9].
\newline
\includegraphics[scale=0.3]{fullone.pdf}
\item\label{q3} A subset of $\mathbf{M}$ formed by all $(x,l_0)\in\mathbf{M}$ where $l_0$ is a fixed line [6].
\newline
\includegraphics[scale=0.3]{line.pdf}
\item\label{q4} A subset of $\mathbf{M}$ formed by all $(x,l)$ such that $x=x_0$ or $l=l_0$ where $x_0$ and $l_0\ni x_0$ are a fixed point and a fixed line respectively [5].
\newline
\includegraphics[scale=0.3]{oneline.pdf}
\item\label{q5} A two-element subset $\{(x_1,l_1),(x_2,l_2)\}\subset\mathbf{M}$ such that $x_1,x_2$ and $l_1\cap l_2$ are three distinct points [5].
\newline
\includegraphics[scale=0.3]{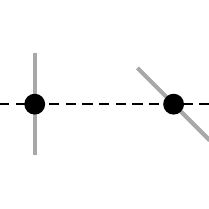}
\item\label{q6} A subset of $\mathbf{M}$ formed by all $(x,l)$ such that $x=x_1$ or $x=x_2$ or $l=l_0$ where $x_1$ and $x_2$ are fixed points and $l_0\ni x_1,x_2$  is a fixed line [4].
\newline
\includegraphics[scale=0.3]{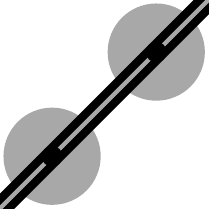}
\item\label{q7} A subset of $\mathbf{M}$ formed by all $(x,l)$ such that $x\in l_0$ for some fixed line $l_0$ [3].
\newline
\includegraphics[scale=0.3]{lineline.pdf}
\item\label{q8} A subset of $\mathbf{M}$ formed by all $(x,l)$ such that $l=l_1$ or $l=l_2$ or $x=l_1\cap l_2$ where $l_1$ and $l_2$ are fixed lines [1].
\newline
\includegraphics[scale=0.3]{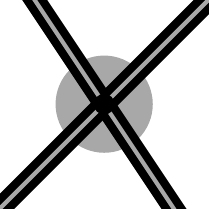}
\item\label{q9} A subset of $\mathbf{M}$ formed by all $(x,l)$ such that $x$ lies on a fixed smooth conic $Q$ and $l$ is tangent to $Q$ at $x$ [1].
\newline
\includegraphics[scale=0.3]{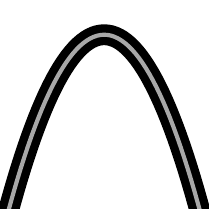}
\item\label{q10} A three-point subset $\{(x_1,l_1),(x_2,l_2),(x_3,l_3)\}\subset\mathbf{M}$ such that $x_1,x_2,x_3$ are not on a line and the intersection $l_1\cap l_2\cap l_3$ is a point different from either of the points $x_1,x_2,x_3$ [1].
\newline
\includegraphics[scale=0.3]{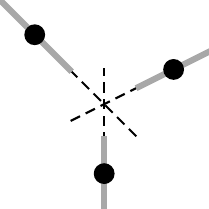}
\item\label{qlast} The whole $\mathbf{M}$ [0].
\newline
\includegraphics[scale=0.3]{all.pdf}
\end{enumerate}
\end{Prop}

{\bf Remark.} Same as for plane cubics, the space $X_3$ does not intersect the closure of $\{\Sing(f)\mid f\in\N_{4,2}\}$ in $2^{\mathbf{M}}$.

\begin{figure}[h]
\centering
\includegraphics[scale=1]{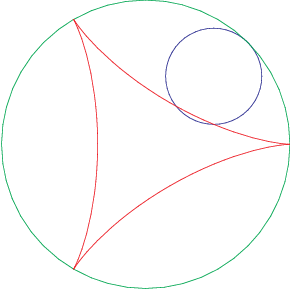}
\caption{A hoop inside a barrel: the tricuspidal quartic.}
\label{tricusp}
\end{figure}

{\bf Proof.} Conditions \ref{first} and \ref{three} can be checked using the following observations.
\begin{itemize}
\item The kernel of the Hessian of $f\in\Sigma$ at a point where $df=0$ is a vector subspace. So if $\Sing(f)\supset\{(x,l_1),(x,l_2)\}$ and $l_1\neq l_2$, then for every $l\ni x$ we have $(x,l)\in\Sing(f)$. 
\item Suppose $\Sing(f)\supset\{(x_1,l_1),(x_2,l_2)\}, x_1\neq x_2$, and $x_2\in l_1$. Then for every $x\in l_1$ we have $(x,l_1)\in\Sing(f)$; in particular, the quartic defined by $f$ contains $l_1$ as a component of multiplicity $\geq 2$.
\item There exists exactly one irredicuble plane quartic with three cusps, modulo projective transformations. (By the Pl\"ucker formulas, such a quartic is rational; one then applies e.g.\ \cite[2.2]{namba}; as a side remark, one can obtain the trucuspidal quartic as the curve traced by a point on a hoop of radius 1 rolling inside a barrel of radius 3, see figure \ref{tricusp}).
\end{itemize}

Conditions \ref{sec}-\ref{four}, \ref{last} and \ref{xxx} are easy to check directly.

In order to check condition \ref{five} we note that we have isomorphisms $X_{\ref{q2}}\cong\C \p^2,X_{\ref{q3}}\cong X_{\ref{q7}}\cong\C \p^{2\vee},X_{\ref{q4}}\cong\mathbf{M},X_{\ref{q8}}\cong B(\C \p^{2\vee},2)$. Using these we describe in the following table the geometrisation of every element of $\left(\bigcup \bar X_i\right)\setminus \left(\bigcup X_i\right)$, with the following exception: we give $\Geom(K)$ only for those $K\in\bar X_{\ref{q10}}\setminus X_{\ref{q10}}$ that consist of 3 elements; two-element configurations from $\bar X_{\ref{q10}}\setminus X_{\ref{q10}}$ can be handled in the same way as elements of $\bar X_{\ref{q5}}\setminus X_{\ref{q5}}$.
\pagebreak[3]

\begin{longtable}{|c|c|c|}
\caption{Geometrisation of limit configurations: plane quartics.}
\label{geom_table_quartics}\\
\hline
\begin{tabular}{c}$i$ such that\\ $K\in \bar X_i\setminus X_i$\end{tabular} & $K$ & $\Geom(K)$\\ 
\hline
\ref{q5} & \begin{tabular}{c}\rule{0pt}{12pt} $\{(x,l_1),(x,l_2)\}$\\ such that $l_1\neq l_2$\\
\includegraphics[scale=0.3]{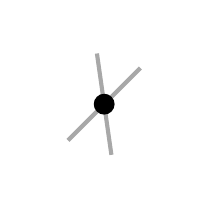}\end{tabular} & \begin{tabular}{c}the element of $X_{\ref{q2}}$\\that corresponds to $x\in\C \p^2$\\
\includegraphics[scale=0.3]{fullone.pdf}\end{tabular}\\
\hhline{~--}
&\begin{tabular}{c}$\{(x_1,l),(x_2,l)\}$\\
such that $x_1\neq x_2$\\
\includegraphics[scale=0.3]{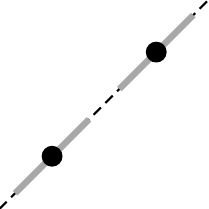}\end{tabular}& \begin{tabular}{c}\rule{0pt}{12pt}the element of $X_{\ref{q3}}$\\that corresponds to $l\in\C \p^{2\vee}$\\
\includegraphics[scale=0.3]{line.pdf}\end{tabular}\\
\hhline{~--}
\rule{0pt}{12pt}&
\begin{tabular}{c}\rule{0pt}{12pt}
$\{(x_1,l_1),(x_2,l_2)\}$\\ 
such that 
$x_1\neq x_2,$\\$l_1\neq l_2,x_2\in l_1$\\
\includegraphics[scale=0.3]{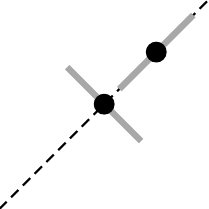}
\end{tabular}& 
\begin{tabular}{c}the element of $X_{\ref{q4}}$\\that corresponds to $(x_2,l_1)\in\mathbf{M}$\\
\includegraphics[scale=0.3]{oneline.pdf}\end{tabular}\\
\hline
\rule{0pt}{12pt}\ref{q9}&\begin{tabular}{c}\rule{0pt}{12pt}$\bigcup_{x\in l_1}\{(x,l_1)\}$\\$\cup\bigcup_{x\in l_2}\{(x,l_2)\}$\\for some
fixed lines $l_1\neq l_2$\\
\includegraphics[scale=0.3]{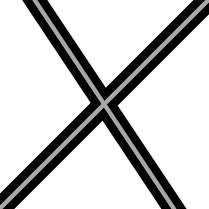}\end{tabular}&\begin{tabular}{c} the element of $X_{\ref{q8}}$ that corresponds\\ to $\{l_1,l_2\}\in B(\C \p^{2\vee},2)$\\
\includegraphics[scale=0.3]{onelineline.pdf}\end{tabular}\\
\hline
\ref{q10}&\begin{tabular}{c}\rule{0pt}{12pt}$\{(x_1,l),(x_2,l),(x_3,l)\}$\\ where $x_1,x_2$ and $x_3$ are\\distinct points on a line $l$\\
\includegraphics[scale=0.3]{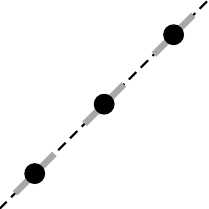}\end{tabular}&\begin{tabular}{c}the element of $X_{\ref{q3}}$\\that corresponds to $l$\\
\includegraphics[scale=0.3]{line.pdf}\end{tabular}\\
\hhline{~--}
&\begin{tabular}{c}\rule{0pt}{12pt}$\{(x_1,l_1),(x_2,l_1),(x_3,l_2)\}$\\ where $x_1\neq x_2$ are points\\ on a line $l_1$,
and $x_3\not\in l_1$ is\\ a point on a line $l_2\neq l_1$;\\ $(l_1\cap l_2)\cap\{x_1,x_2\}$ may or\\may not be empty\\
\includegraphics[scale=0.3]{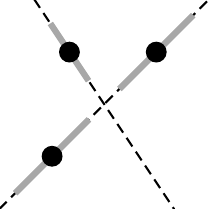}\end{tabular}&\begin{tabular}{c}the element of $X_{\ref{q8}}$\\that corresponds to\\ $\{l_1,l_2\}\in B(\C \p^{2\vee},2)$\\
\includegraphics[scale=0.3]{onelineline.pdf}\end{tabular}\\
\hhline{~--}
&\begin{tabular}{c}\rule{0pt}{12pt}$\{(x_1,l_1),(x_2,l_1),(x_3,l_2)\}$\\ where $x_1\neq x_2$ are points\\
on a line $l_1$, $l_2\neq l_1$ is
\\ another line and $l_1\cap l_2=\{x_3\}$;\\ $x_3$ may or may not \\ coincide with $x_1$ or $x_2$\\
\includegraphics[scale=0.3]{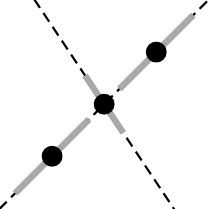}\end{tabular}&\begin{tabular}{c}the element of $X_{\ref{q4}}$\\that corresponds to\\ $(x_3,l_1)$\\
\includegraphics[scale=0.3]{oneline.pdf}\end{tabular}\\
\hhline{~--}
%
%
&
\begin{tabular}{c}
\rule{0pt}{12pt}$\{(x,l_1),(x,l_2),(x,l_3)\}$\\ where $l_1,l_2,l_3$ are\\ distinct lines through $x$\\
\includegraphics[scale=0.3]{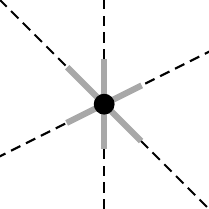}
\end{tabular}
&\begin{tabular}{c}the element of $X_{\ref{q2}}$\\that corresponds to $x\in\C \p^2$\\
\includegraphics[scale=0.3]{fullone.pdf}\end{tabular}\\
\hhline{~--}
&\begin{tabular}{c}\rule{0pt}{12pt}$\{(x_1,l_1),(x_1,l_2),(x_2,l_3)\}$\\ where $l_1,l_2,l_3$ are\\ distinct lines through $x_1$,\\
and $x_2\neq x_1$ is a point on $l_3$\\
\includegraphics[scale=0.3]{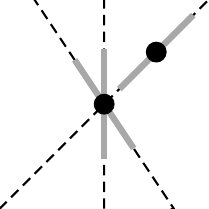}\end{tabular}&\begin{tabular}{c}the element of $X_{\ref{q4}}$\\that corresponds to $(x_1,l_3)\in\mathbf{M}$\\
\includegraphics[scale=0.3]{oneline.pdf}\end{tabular}\\
\hhline{~--}
&\begin{tabular}{c}\rule{0pt}{12pt}$\{(x_1,l_1),(x_2,l_2),(x_3,l_3)\}$\\ where $x_1,x_2,x_3$ are distinct\\ points on a line $l$,\\
$l_1,l_2,l_3$ are distinct lines\\ through a point $x$, $x\notin l$,\\and $x_i\in l_i,i=1,2,3$\\
\includegraphics[scale=0.3]{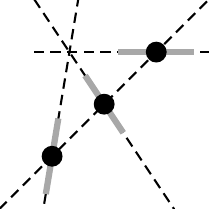}\end{tabular}&\begin{tabular}{c}the element of $X_{\ref{q7}}$\\ that corresponds to $l$\\
\includegraphics[scale=0.3]{lineline.pdf}\end{tabular}\\
\hhline{~--}
&\begin{tabular}{c}\rule{0pt}{12pt}$\{(x_1,l_1),(x_2,l_2),(x_3,l_3)\}$\\ where $l_1,l_2,l_3$ are\\ distinct lines through $x_1$,\\and $x_i$ is a point on $l_i\setminus l_1$\\for $i=2,3$\\
\includegraphics[scale=0.3]{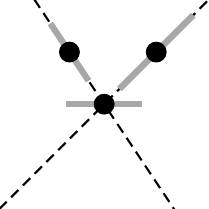}\end{tabular}&
\begin{tabular}{c}the element of $X_{\ref{q8}}$\\that corresponds to\\ $\{l_2,l_3\}\in B(\C \p^{2\vee},2)$\\
\includegraphics[scale=0.3]{onelineline.pdf}\end{tabular}\\
\hhline{~--}
\hline
\end{longtable}
%

Checking condition \ref{five} is now straightforward.$\clubsuit$

Now we apply the method described in section \ref{conres} to construct a conical resolution $\sigma$ of $\Sigma$ and a filtration $$\varnothing=F_0\subset F_1\subset\cdots\subset F_{\ref{qlast}}.$$ The first assertion of theorem \ref{mhquartp2} follows immediately from

\begin{theorem}\label{spseqqvector}
Let $(E^r_{p,q}, d^r_{p,q})$ be the resulting spectral sequence that converges to $\bar H_*(\sigma,\Q)\cong\bar H_*(\Sigma,\Q)$. The groups $(E^1_{p,q})$ are given by table \ref{spseqquart}.
\begin{table}
\caption{Spectral sequence for $\sigma$ in the case of plane nodal quartics}
\label{spseqquart}
$$\begin{array}{c|c|c|c|c|c|c|c|c|c|c|c|}
&1&2&3&4&5&6&7&8&9&10&11\\
\hline
25&\Q(13)&&&&&&&&&&\\
\hline
24&&&&&&&&&&&\\
\hline
23&\Q(12)\oplus\Q(12)&\Q(12)&&&&&&&&&\\
\hline
22&&&&&&&&&&&\\
\hline
21&\Q(11)\oplus\Q(11)&\Q(11)&&&&&&&&&\\
\hline
20&&&&&&&&&&&\\
\hline
19&\Q(10)&\Q(10)&&&&&&&&&\\
\hline
18&&&&&&&&&&&\\
\hline
17&&&&&&&&&&&\\
\hline
16&&&\Q(9)&&\Q(10)&&&&&&\\
\hline
15&&&&&&&&&&&\\
\hline
14&&&\Q(8)&&\Q(9)&\Q(9)&&&&&\\
\hline
13&&&&&&&&&&&\\
\hline
12&&&\Q(7)&&\Q(8)&\Q(8)&&&&&\\
\hline
11&&&&&&&&&&&\\
\hline
10&&&&&&\Q(7)&&&&&\\
\hline
9&&&&&&&&&&&\\
\hline
8&&&&&&&&&&&\\
\hline
7&&&&&&&&&&&\\
\hline
6&&&&&&&&\Q(6)&&&\\
\hline
5&&&&&&&&&&&\\
\hline
4&&&&&&&&\Q(5)&&&\\
\hline
3&&&&&&&&&&&\\
\hline
2&&&&&&&&\Q(4)&&&\Q(5)\\
\hline
1&&&&&&&&&&&\\
\hline
0&&&&&&&&&&&\Q(4)\\
\hline
\end{array}$$
\end{table}	
The differentials $d^1_{2,23}$, $d^1_{2,21}$, $d^1_{2,19}$, $d^3_{11,0}$, $d^3_{11,2}$, $d^1_{6,12}$, $d^1_{6,14}$ and $d_{6,10}^3$ are non-zero.
\end{theorem}

{\bf Proof.} We will now calculate all groups $E^1_{p,q}$ and all differentials except $d_{6,10}^3$; we will prove that $d_{6,10}^3\neq 0$ in the next subsection.

\textbf{Column 1.} We have $\Phi_{\ref{q1}}=X_{\ref{q1}}=\mathbf{M}$, and $F_{\ref{q1}}$ is a complex bundle over $\Phi_{\ref{q1}}$ of rank 10. Hence
the groups $E^1_{\ref{q1},i}\cong\bar H_{\ref{q1}+i}(F_{\ref{q1}},\Q)$ are as shown in table \ref{spseqquart}.

\medskip

\textbf{Column 2.} For $K\in X_{\ref{q2}}$ we have $$\partial\Lambda(K)=\tilde\partial\Lambda(K)=\bigcup_{\begin{array}{c}\scriptstyle{K'\in X_{\ref{q1}},}\\ \scriptstyle{K'\subset K}\end{array}}\Lambda(K')\cong\C \p^1.$$ So the space $\Phi_{\ref{q2}}\setminus\Phi_{\ref{q1}}$ is a fibre bundle
over $X_{\ref{q2}}$ with fibre homeomorphic to the open cone $Cone^\circ(\C \p^1)$. The space $F_{\ref{q2}}\setminus F_{\ref{q1}}$ is a rank 9 complex vector bundle over $\Phi_{\ref{q2}}\setminus\Phi_{\ref{q1}}$. This implies the statement of the theorem concerning the groups $E^1_{\ref{q2},*}$. 

\medskip

\textbf{Column 3.} The groups $E^1_{\ref{q3},*}$ can be calculated in the same way as $E^1_{\ref{q2},*}$.

\medskip

\textbf{Column 4.} For $K\in X_{\ref{q4}}$ the space $\partial\Lambda(K)$ is equal to $\tilde\Lambda(K_1)\cup\tilde\Lambda(K_2)$, where $K_1$ and $K_2$ are the maximal subconfigurations of
$K$ that belong to $X_{\ref{q2}}$ and $X_{\ref{q3}}$ respectively. Notice that $K_1\cap K_2$ is a configuration of type \ref{q1} and $\tilde\Lambda(K_1)\cap\tilde\Lambda(K_2)=\tilde \Lambda(K_1\cap K_2)=\Lambda(K_1\cap K_2)$, see lemma~\ref{cup_cap}. Since $\tilde\Lambda(K_1),\tilde\Lambda(K_2)$ and $\Lambda(K_1\cap K_2)$ are all contractible, all groups $E^1_{\ref{q4},i}$ are zero.

\medskip

\textbf{Column 5.} For $K\in X_{\ref{q5}}$ we have $\tilde\partial\Lambda(K)=\partial\Lambda(K)=\Lambda(K_1)\sqcup\Lambda(K_2)$ with $K_1,K_2\in X_{\ref{q1}}$.
So $\bar H_*(\Phi_{\ref{q5}}\setminus\Phi_{\ref{q4}},\Q)\cong\bar H_{*-1}(X_{\ref{q5}},\pm\Q)$ (recall that $X_{\ref{q5}}\subset B(\mathbf{M},2)$). The groups $\bar H_*(\Phi_{\ref{q5}}\setminus\Phi_{\ref{q4}},\Q)$ can be obtained from the fact that $\bar H_i(B(\C \p^2,2),\pm\Q)=\Q(i/2)$, if $i=2,4,6$ and is zero otherwise~\cite{vas3}.

\medskip

\textbf{Column 6.} The space $\partial\Lambda(K),K\in X_{\ref{q6}}$ is
$\partial\Lambda(K)=\tilde\Lambda(K_1)\cup\tilde\Lambda(K_2) \cup A,$
where $K_1$ and $K_1$ are maximal subconfigurations of $K$ of type \ref{q4}, and $A$ is the union of $\Lambda(K')\setminus\partial\Lambda(K')$ for all $K'\subset K$ of type \ref{q5}. The space $\tilde\Lambda(K_1)\cup\tilde\Lambda(K_2)$ has the homology of a point. Moreover, $(\tilde\Lambda(K_1)\cup\tilde\Lambda(K_2))\cap A=\varnothing$, so $$\bar H_p(\Lambda(K)\setminus\partial\Lambda(K),\Z)=\tilde H_{p-1}(\partial\Lambda(K),\Q)=\bar H_{p-1}(A,\Q).$$ The result is $\Q(2)$ if $p=6$ and 0 otherwise. The groups $\bar H_*(\Phi_{\ref{q6}}\setminus\Phi_{\ref{q5}},\Q)$ can now be obtained in the same way as $\bar H_*(\Phi_{\ref{q5}}\setminus\Phi_{\ref{q4}},\Q)$.

\medskip

\textbf{Column 7.} For any $K\in X_{\ref{q7}}$ we have a filtration on $\partial\Lambda(K)$ induced by the intersections with
$\Phi_i$, which we denote by $\partial_i \Lambda(K)$. Similarly to the 
above we get $\partial_{\ref{q4}} \Lambda(K)\sim pt$. Note also that $\partial_{\ref{q6}} \Lambda(K)\setminus \partial_{\ref{q4}} \Lambda(K)$ is fibred over the space
$$X'=\{K'\in X_6\mid K'\subset K\}.$$ Let $F=\Lambda (K')\setminus \partial_4 \Lambda (K')$ be the fibre over $K'\in X_6$. We have 
$\partial_4 \Lambda (K')\sim pt$, so $\bar H_*(F,\Q)=0$, which implies that $\tilde H_*(\partial\Lambda(K),\Q)$ is zero, and hence so is the \ref{q7}-th column of table \ref{spseqquart}.

\medskip

%
%
\textbf{Column 8.} The \ref{q8}-th column of table \ref{spseqquart} is calculated in exactly the same way as the \ref{q6}-th.

\medskip

\textbf{Column 9.} If $K\in X_{\ref{q9}}$, then $\partial\Lambda (K)$ is homeomorphic to the second self-join $(\C \p^1)^{*2}\cong (S^2)^{*2}$ of $\C \p^1$. It will follow that the \ref{q9}-th column of table \ref{spseqquart} is zero once we have proved the following lemma.

\begin{lemma}\label{selfjoincp1}
For every $k\geq 2$ the space $(\C\p^1)^{*k}$ has the rational homology of a point.
\end{lemma}

{\bf Proof of lemma \ref{selfjoincp1}.} The direct limit $(\C\p^1)^{*\infty}=\varinjlim (\C\p^1)^{*k}$ is contractible. There is a filtration
\begin{equation}\label{filtrinfselfjoin}
\C\p^1=(\C\p^1)^{*1}\subset (\C\p^1)^{*2}\subset (\C\p^1)^{*3}\subset\cdots =(\C\p^1)^{*\infty}.
\end{equation}
Every space $(\C\p^1)^{*k}\setminus (\C\p^1)^{*(k-1)}$ is fibred over $B(\C\p^1,k)$, and we have $$H_*((\C\p^1)^{*k}\setminus (\C\p^1)^{*(k-1)},\Q)\cong H_{*+k-1}(B(\C\p^1,k),\pm \Q).$$ If $k>1$, then all groups $H^i(B(\C\p^1,k),\pm \Q)$ are zero apart from $H^2(B(\C\p^1,2),\pm \Q)\cong\Q$, see \cite{vas3}. Applying the Poincar\'e duality we see that the only non-zero group $H_i(B(\C\p^1,k),\pm \Q)$ for $k>1$ is $H_2(B(\C\p^1,2),\pm \Q)\cong\Q$. So the first term of the rational homology spectral sequence obtained from (\ref{filtrinfselfjoin}) is
\begin{equation*}
\begin{array}{c|c|c|}
\hline
1&\Q&\Q\\
\hline
0 &&\\
\hline
-1 &\Q &\\
\hline
 & 1& 2
\end{array}
\end{equation*}
The only potentially non-zero differential in this sequence is non zero, and the lemma follows.$\clubsuit$

\medskip

\textbf{Column 10.} We now consider the \ref{q10}-th column of table \ref{spseqquart}. As in the case of column \ref{q5}, the groups $\bar H_*(F_{\ref{q10}}\setminus F_{\ref{q9}},\Q)$ are isomorphic up to a shift to the groups $\bar H_*(X_{\ref{q10}},\pm\Q)$ (we use the fact that $X_{\ref{q10}}\subset B(\mathbf{M},3)$). The space $X_{\ref{q10}}$ is fibred over $\C \p^2$ by taking a $K\in X_{\ref{q10}}$ to the intersection point of the three lines involved in the definition of $K$. Let $F$ be the fibre. The groups $\bar H_*(F,\pm\Q)$ are isomorphic to shifted groups $\bar H_*(B(\C \p^1,3),\pm\Q)$, which are zero by \cite{vas3}, see the previous paragraph. We conclude that $E^1_{\ref{q10},*}=0$.

\medskip

\textbf{Column 11.} Finally, to obtain the last column of table \ref{spseqquart} we need to compute the rational Borel-Moore homology of the open cone over the space

\begin{equation}\label{lambdaq}
\partial\Lambda(\mathbf{M})=\tilde\partial\Lambda(\mathbf{M})=\bigcup_{i\leq\ref{q10}}\bigcup_{K\in X_i}\tilde\Lambda(K).
\end{equation}
This space can be filtered as follows:
\begin{equation}\label{filtr}
\Phi_{\ref{q2}}\subset\Phi_{\ref{q7}}\subset\Phi_{\ref{q8}}\subset\Phi_{\ref{q9}}\subset\Phi_{\ref{q10}}.
\end{equation}
In the first term of the spectral sequence
that corresponds to this filtration, the $i$-th column for $i=3,4,5$ is equal to the $i+5$-th column of table \ref{spseqquart}, with all weights and degrees shifted down by $2d_{i+5}$. The rational Borel-Moore homology of $\Phi_{\ref{q2}}$ is isomorphic to $H_*(\C \p^2,\Q)$, and the space $\Phi_{\ref{q7}}\setminus\Phi_{\ref{q2}}$ is fibred over $\C \p^2$ with fibre
$$F=\tilde\Lambda(K)\setminus\bigsqcup_{\begin{array}{c}\scriptstyle{K'\in X_{\ref{q2}},}\\ \scriptstyle{K'\subset K}\end{array}}\tilde\Lambda(K').$$

The rational homology of the disjoint union in the right hand side of this formula is isomorphic to $H_*(\C\p^1,\Q)$, so $\bar H_i(F,\Q)=\Q(1)$ for $i=3$ and 0 otherwise, which implies that $\bar H_i(\Phi_{\ref{q7}}\setminus\Phi_{\ref{q2}},\Q)=\Q(\frac{i-1}{2})$ for $i=3,5,7$ and 0 otherwise. We conclude that only the first three columns of the first term of the spectral sequence that corresponds to (\ref{filtr}) are non zero, and they are as follows.
\begin{equation}\label{lastcolq}
\begin{array}{c|c|c|c|}
\hline
9&&&\Q(5)\\
\hline
8&&&\\
\hline
7&&&\Q(4)\\
\hline
6&&&\\
\hline
5&&\Q(3)&\Q(3)\\
\hline
4&&&\\
\hline
3&\Q(2)&\Q(2)&\\
\hline
2&&&\\
\hline
1&\Q(1)&\Q(1)&\\
\hline
0&&&\\
\hline
-1&\Q(0)&&\\
\hline
&1&2&3\\
\end{array}
\end{equation}

Comparing this with what we already know about the spectral sequence of table~\ref{spseqquart}, we see that if any of the potentially non-zero differentials of (\ref{lastcolq}) were zero, then $\Pi_{4,2}\setminus\Sigma$ would have a non-zero rational cohomology group in degree $\geq 21$. Using corollary~\ref{cor_groupactions} we conclude that this is impossible because the group $\GL_3(\C)$ acts on $\Pi_{4,2}\setminus\Sigma$ with finite stabilisers, and the geometric quotient exists and is a non-complete complex algebraic variety of dimension 6.

So the non-zero Borel-Moore homology groups of the space $\partial\Lambda(\mathbf{M})$ are $\Q(0),\Q(4)$ and $\Q(5)$ in degrees 0, 10 and 12 respectively, which implies that the groups~$E^1_{\ref{qlast},*}$ are as shown in table~\ref{spseqquart}.

\bigskip

We have thus finished computing the groups~$E^1_{*,*}$. To complete the proof of theorem \ref{spseqqvector} it remains to compute the differentials $d_{p,q}^r$. This will be done in the next few lemmas.

\begin{lemma}\label{difq1}
The differentials $d^1_{2,23}, d^1_{2,21}$ and $d^1_{2,19}$ are non zero.
\end{lemma}

{\bf Proof.} The proof is almost identical to that of lemma~\ref{diff_2nd_column_plane_cubics}, the only difference being as follows. 
We replace the vector bundles $\vartheta_{3,2}$ and $\xi_{3,2}$ with $\vartheta_{4,2}$ and $\xi_{4,2}$ respectively. The fibre of $\vartheta_{4,2}/\xi_{4,2}$ over $(x_0,l_0)$ is spanned by $\mathrm{x}_0^2\mathrm{x}_2^2$. The corresponding character of $T$ is $T(s,t)\mapsto s^2$, and the differential is $2\epsilon_1$. None of the elements $2\epsilon_1,2\epsilon_0\epsilon_1,2\epsilon_0^2 \epsilon_1\in\Q[\epsilon_0,\epsilon_1,\epsilon_2]/I$ is zero.$\clubsuit$

\begin{lemma}\label{difq2}
The differentials $d^3_{11,0}$ and $d^3_{11,2}$ are non zero.
\end{lemma}

{\bf Proof.} We will compute the rational Borel-Moore homology groups of $F_{\ref{qlast}}\setminus F_{\ref{q7}}$ using a filtration that is different from the one induced by $\varnothing=F_0\subset F_1\subset\cdots\subset F_{\ref{qlast}}$, namely the filtration
\begin{equation}\label{difq2filtr}
X\subset X\cup (F_{\ref{q8}}\setminus F_{\ref{q7}})\subset X\cup (F_{\ref{q9}}\setminus F_{\ref{q7}})\subset X\cup (F_{\ref{q10}}\setminus F_{\ref{q7}})=F_{\ref{qlast}}\setminus F_{\ref{q7}},
\end{equation}
where $X=\Phi_{\ref{qlast}}\setminus\Phi_{\ref{q7}}$, and we identify $X$ with a subspace of $\Pi_{4,2}\times\tilde\Lambda(\mathbf{M})$ via the zero section embedding.

We have $\bar{H}_{*+1}(X)=H_{*+1}(Cone(\Phi_{10}),\Phi_7)=\tilde H_{*}(\Phi_7)$, and the latter groups are given by the first two columns of spectral sequence~(\ref{lastcolq}). Hence $\bar H_i(X,\Q)=\Q(3)$ for $i=8$ and $0$ for $i\neq 8$. The
third and fourth columns of the rational Borel-Moore homology spectral sequence that corresponds to (\ref{difq2filtr}) are zero since, as we saw above,
$$\bar H_*(\Phi_{\ref{q9}}\setminus\Phi_{\ref{q8}},\Q)=\bar H_*(\Phi_{\ref{q10}}\setminus\Phi_{\ref{q9}},\Q)=0.$$ To determine the second column
we need to know the groups $\bar H_*((F_{\ref{q8}}\setminus F_{\ref{q7}})\setminus X,\Q)$

The space
$(F_{\ref{q8}}\setminus F_{\ref{q7}})\setminus X$ coincides with $F_{\ref{q8}}\setminus F_{\ref{q7}}$ minus the zero section of the bundle
$F_{\ref{q8}}\setminus F_{\ref{q7}}\stackrel{p}\to\Phi_{\ref{q8}}\setminus\Phi_{\ref{q7}}$. Recalling the way we computed above the Borel-Moore homology of $F_{\ref{q8}}\setminus F_{\ref{q7}}$, we see that $$\bar H_*((F_{\ref{q8}}\setminus F_{\ref{q7}})\setminus X,\Q)\cong \bar H_{*-6}(Y,\pm\Q)\otimes\Q(2)$$ where $Y$ is the space of all couples $$(f,K)\in \left(\Pi_{4,2}\setminus\{0\}\right)\times B(\C \p^{2\vee},2)$$ such that $f$ is divisible by the square of an equation of each of the lines $l_1,l_2$ such that $K=\{l_1,l_2\}$.

Let $\tilde Y$ be the total space of the pullback of the vector bundle $Y\to B(\C \p^{2\vee},2)$ to $F(\C \p^{2\vee},2)$. The group $S_2$ acts on $\tilde Y$ in a natural way, and
$\bar H_*(Y,\pm\Q)$ is isomorphic to the $S_2$-anti-invariant part of $\bar H_*(\tilde Y,\Q)$. Let us compute $\bar H_*(\tilde Y,\Q)$ together with the action of $S_2$ on it. By the
Poincar\'e duality, it suffices to consider $H^*(\tilde Y,\Q)$.

\begin{lemma}
\label{ring}
The ring $H^*(F(\C \p^{2\vee},2),\Q)$ is generated by the pullbacks $a_1$ and $a_2$ of the standard generator $c_1(\OO_{\C\p^{2\vee}}(1))$ of $H^2(\C \p^{2\vee},\Q)$ under the
projections $F(\C \p^{2\vee},2)\rightrightarrows\C \p^{2\vee}$, subject to the relations $a_1^3=a_2^3=a_1^2+a_2^2+a_1a_2=0$. The generator of $S_2$ interchanges $a_1$ and $a_2$.
\end{lemma}

{\bf Proof of lemma~\ref{ring}.} There is a general recipe for calculating the rational cohomology algebra of ordered configuration spaces of smooth compact complex algebraic varieties, see \cite{totaro}. In our case it gives the following answer: $H^*(F(\C \p^{2\vee},2),\Q)$ is the cohomology of $\C \p^{2\vee}\times \C \p^{2\vee}$ modulo the ideal generated by the Poincar\'e dual of the diagonal. The standard formula that expresses the class of the diagonal in terms of the Poincar\'e duality shows that this class is $a_1^2+a_2^2+a_1a_2$.$\clubsuit$

The line bundle over $\C \p^{2\vee}$ with total space
$$\{(f,l)\in\Pi_{2,2}\times\C \p^{2\vee}\mid\mbox{$f$ is the square of an equation of $l$}\}$$ is isomorphic to $\mathcal{O}_{\C \p^{2\vee}}(-2)$, so $\tilde Y$ can be viewed
as the space of non-zero vectors of the restriction of $\xi_1\otimes\xi_2$ to $F(\C \p^{2\vee},2)$ where $\xi_1$ and $\xi_2$ are the pullbacks of $\OO_{\C \p^{2\vee}}(-2)$ under the projections $\C \p^{2\vee}\times\C \p^{2\vee}\rightrightarrows\C \p^{2\vee}$. The first Chern class of $\xi_1\otimes\xi_2$ is $-2a_1-2a_2$. Using proposition~\ref{euler-thom} we deduce that $H^*(\tilde Y,\Q)$ is $S_2$-equivariantly isomorphic to the cohomology of the algebra $(H^*(F(\C\p^{\vee},2),\Q)\otimes\Lambda[b],d)$ where $b$ has degree $1$, the action of $S_2$ on $\Lambda[b]$ is trivial, $d(b)=2a_1+2a_2$, and $d=0$ on $H^*(F(\C\p^{\vee},2),\Q)$. 

The quotient algebra $H^*(F(\C\p^{2\vee},2),\Q)/(a_1+a_2)$ is two dimensional over $\Q$ and is spanned by the classes of $1$ and $a_1-a_2$. An easy check shows now that the $S_2$-anti-invariant part of $H^*(\tilde Y,\Q)$ is $\Q(-1)$ in degree 2, $\Q(-4)$ in degree 7 and zero in all other degrees. This implies $\bar H_8(Y,\pm\Q)=\Q(4),
\bar H_3(Y,\pm\Q)=\Q(1)$, and $\bar H_i(Y,\pm\Q)=0$ for $i\neq 3,8$. So $(F_{\ref{q8}}\setminus F_{\ref{q7}})\setminus X$ has non-zero rational Borel-Moore homology groups
only in degrees 9 and 14, and these are respectively $\Q(3)$ and $\Q(6)$.

Now we can write the $E^1$ page of the spectral sequence obtained from (\ref{difq2filtr}) that converges to $\bar H_*(F_{\ref{qlast}}\setminus F_{\ref{q7}},\Q)$.
$$
\begin{array}{c|c|c|}
\hline
12&&\Q(6)\\
\hline
11&&\\
\hline
10&&\\
\hline
9&&\\
\hline
8&&\\
\hline
7&\Q(3)&\Q(3)\\
\hline
&1&2\\
\end{array}
$$
Comparing this with the last four columns of table \ref{spseqquart}, which show that $\bar H_8(F_{\ref{qlast}}\setminus F_{\ref{q7}},\Q)=0$, we conclude that $\bar H_{14}(F_{\ref{qlast}}\setminus F_{\ref{q7}},\Q)=\Q(6)$, all other groups
$\bar H_*(F_{\ref{qlast}}\setminus F_{\ref{q7}},\Q)$ are zero, and so the differentials $d^3_{11,0}$ and $d^3_{11,2}$ of the spectral sequence of
theorem \ref{spseqqvector} are non zero. Lemma \ref{difq2} is proved.$\clubsuit$

{\bf Remark.} An alternative way to prove lemma~\ref{difq2} would be to use lemma~\ref{alt_flag_variety} and the transitive action of $\GL_3(\C)$ on $F(\C\p^{2\vee},2)$, cf.\ the proofs of lemmas~\ref{diff_2nd_column_plane_cubics} and~\ref{difq1}.

\medskip

Lemmas \ref{difq1} and \ref{difq2} will suffice (together with table \ref{spseqquart}) to determine the cohomology of the quotient space $(V\setminus\Sigma)/\GL_3(\C)$. However, to complete the proof of theorem \ref{spseqqvector} we need to compute the remaining differentials $d_{p,q}^r$. In the case of $d^1_{6,12}$ and $d^1_{6,14}$, the argument is analogous to the one used above to prove lemmas \ref{difq1} and \ref{difq2}. 
There doesn't seem to be an easy direct way to show that the differential $d_{6,10}^3$ is non zero. We will calculate $d_{6,10}^3$ alongside the cohomology of the moduli space in section~\ref{moduli42}.

\begin{lemma}\label{difq3}
The differentials $d^1_{6,12}$ and $d^1_{6,14}$ are non zero.
\end{lemma}

{\bf Proof.} We will use lemma~\ref{aux_lemma_diff} to calculate the connecting homomorphism of the pair $(F_{\ref{q6}}\setminus F_{\ref{q4}},F_{\ref{q5}}\setminus F_{\ref{q4}})$. The space $F_{\ref{q6}}\setminus F_{\ref{q4}}$ is the union of $F_{\ref{q5}}\setminus F_{\ref{q4}}$, which is a rank 5 vector bundle over $\Phi_{\ref{q5}}\setminus \Phi_{\ref{q4}}$, and the space $A$ of all couples $(f,x)$ such that $x\in\tilde\Lambda(K)\setminus\Phi_{\ref{q4}}$ for some $K\in X_{\ref{q6}}$ and $f\in L(K)$, which is a rank 4 vector bundle over $\Phi_{\ref{q6}}\setminus \Phi_{\ref{q4}}$.

Observe that $\Phi_{\ref{q6}}\setminus \Phi_{\ref{q4}}$ is a fibre bundle over $X_{\ref{q6}}$ with fibre $\tilde\Lambda(K)\setminus (\tilde\Lambda(K_1)\cup\tilde\Lambda(K_2))$, where $K\in X_{\ref{q6}}$ and $K_1\neq K_2$ are the configurations from $X_{\ref{q4}}$ that are contained in $K$. Since $\tilde\Lambda(K_1),\tilde\Lambda(K_2)$ and $\tilde\Lambda(K_1)\cap \tilde\Lambda(K_2)$ are all contractible, we have $\bar H_*(\Phi_{\ref{q6}}\setminus \Phi_{\ref{q4}},\Q)=0$, so all connecting homomorphisms of the pair $(\Phi_{\ref{q6}}\setminus \Phi_{\ref{q4}},\Phi_{\ref{q5}}\setminus \Phi_{\ref{q4}})$ are isomorphisms.

To apply lemma~\ref{aux_lemma_diff} we also need to calculate the Euler class of the quotient of the vector bundle $F_{\ref{q5}}\setminus F_{\ref{q4}}\to \Phi_{\ref{q5}}\setminus \Phi_{\ref{q4}}$ by the restriction of $A\to \Phi_{\ref{q6}}\setminus \Phi_{\ref{q4}}$ to $\Phi_{\ref{q5}}\setminus \Phi_{\ref{q4}}$, as well as the action of the Euler class on $\bar H_*(\Phi_{\ref{q5}}\setminus \Phi_{\ref{q4}},\Q)$. Recall that there is a homotopy equivalence $\Phi_{\ref{q5}}\setminus \Phi_{\ref{q4}}\to X_{\ref{q5}}$. Observe that the bundle we are after is the pullback of $\xi/\eta$ from $X_{\ref{q5}}$ where 
$$\tot(\xi)=\{(K,f)\in X_{\ref{q5}}\times\Pi_{4,2}\mid f\in L(K)\}, \mbox{ and }$$
$$\tot(\eta)=\{(K,f)\in X_{\ref{q5}}\times\Pi_{4,2}\mid f\in L\ (\mbox{the unique configuration from $X_{\ref{q6}}$ containing $K$})\}.$$

Let $\tilde X$ be the space of \emph{ordered} couples $((x',l'),(x'',l''))$ such that $\{(x',l'),(x'',l'')\}\in X_{\ref{q5}}$. Note that $S_2$ acts on $\tilde X$ and $H^*(X_5,\Q)$ is the invariant part $H^*(\tilde X,\Q)^+$ of $H^*(\tilde X,\Q)$. Note also that
$$\bar H_*(\Phi_{\ref{q5}}\setminus \Phi_{\ref{q4}},\Q)\cong \bar H_{*-1}(X_{\ref{q5}},\pm\Q)\cong H^{12-(*-1)}(X_{\ref{q5}},\pm\Q)\cong H^{13-*}(\tilde X,\Q)^-$$
as $H^*(X_{\ref{q5}},\Q)\cong H^*(\tilde X,\Q)^+$-modules, where upper $-$ denotes the $S_2$-anti-invariant part. Using lemma~\ref{aux_lemma_diff} we conclude that lemma~\ref{difq3} will follow once we have proved
\begin{lemma}\label{difq3_aux}
Let $\tilde\xi$ and $\tilde\eta$ be the pullbacks of $\xi$, respectively $\eta$ to $\tilde X$. Then $e(\tilde\xi/\tilde\eta)\in H^2(\tilde X,\Q)^+$ and the maps $H^2(\tilde X,\Q)^-\to H^4(\tilde X,\Q)^-$ and $H^4(\tilde X,\Q)^-\to H^6(\tilde X,\Q)^-$ given by multiplication by $e(\tilde\xi/\tilde\eta)$ are non zero.
\end{lemma}

{\bf Proof of lemma~\ref{difq3_aux}.} The group $\SL_3(\C)$ acts transitively on $\tilde X$. Let $(x_0,l_0)$ be as in the proof of lemma~\ref{diff_2nd_column_plane_cubics}. Set $x_1=(0:0:1)$ and $l_1$ equal the projectivisation of the plane through $(0,1,0)$ and $(0,0,1)$. The stabiliser of the couple $((x_0,l_0),(x_1,l_1))$ is then the standard maximal torus $T\subset\SL_3(\C)$, see lemma~\ref{alt_flag_variety}. The description of $H^*(\SL_3(\C)/B,\Q)$ given in that lemma is also valid for $H^*(\SL_3(\C)/T,\Q)$. The fibre of $\tilde\xi/\tilde\eta$ over $((x_0,l_0),(x_1,l_1))$ is spanned by $\mathrm{x}_0^2\mathrm{x}_2^2$. As we saw in the proof of lemma~\ref{difq1}, the differential of the corresponding character of $T$ is $2\epsilon_1=c_1(\tilde\xi/\tilde\eta)=e(\tilde\xi/\tilde\eta)$.

Let us now calculate the action of $S_2$ on $H^*(\tilde X,\Q)$. Let $p', p'':\tilde X\to\C\p^2$ be the projections given by $((x',l'),(x'',l''))\mapsto x'$, respectively by $((x',l'),(x'',l''))\mapsto x''$. The $SL_3(\C)$-stabiliser of $x_0$ is the subgroup $P$ from the proof of lemma~\ref{diff_2nd_column_plane_cubics}. So the map $p'$ is the projection $\SL_3(\C)/T\to \SL_3(\C)/P$. As in the proof of lemma~\ref{diff_2nd_column_plane_cubics} we find that $(p')^*\big(c_1(\OO_{\C\p^2}(-1))\big)=\epsilon_0$, and a similar argument shows that $(p'')^*\big(c_1(\OO_{\C\p^2}(-1))\big)=\epsilon_2$. So the generator of $S_2$ interchanges $\epsilon_0$ and $\epsilon_2$ (and preserves $\epsilon_1=-\epsilon_0-\epsilon_2$). The elements $$(\epsilon_0-\epsilon_2), 2\epsilon_1(\epsilon_0-\epsilon_2), 2\epsilon_1^2(\epsilon_0-\epsilon_2)$$ belong to the $S_2$ anti-symmetric part of $\Q[\epsilon_0,\epsilon_1,\epsilon_2]/I$, and none of them is zero. This proves lemma~\ref{difq3_aux}, and also lemma~\ref{difq3}.$\clubsuit$

\begin{lemma}\label{difq4}
The differential $d^3_{6,10}$ is non zero.
\end{lemma}

\textbf{Proof.} The claim follows from computation of cohomology of the moduli space $(\Pi_{4,2}\setminus \N_{4,2})/\GL_3(\C)$ and is presented in section~\ref{moduli42}.$\clubsuit$

\smallskip

This completes the calculation of the spectral sequence given in table~\ref{spseqquart} and concludes the proof of Theorem~\ref{spseqqvector}.$\clubsuit$

\subsection{The moduli space}
\label{moduli42}

It follows from proposition~\ref{examples_quotients} that there exists a geometric quotient $\mathcal{M}^\nod_{4,2}=(\Pi_{4,2}\setminus \N_{4,2})/\GL_3(\C)$. To obtain the cohomology groups of $\mathcal{M}^\nod_{4,2}$ we compute the Leray spectral sequence of the projection $r:{\Pi_{4,2}\setminus \N_{4,2}\rightarrow \mathcal{M}^\nod_{4,2}}$.

\begin{lemma}
\label{mod42lemma}
The second page of the Leray spectral sequence $(E^{a,b}_c,d_c)$ of the projection $\Pi_{4,2}\setminus \N_{4,2}\stackrel{r}{\longrightarrow} \mathcal{M}^\nod_{4,2}$ is as follows:
\end{lemma}

\begin{equation}
\label{mod42}
\begin{array}{c||c|c|c|c|c|c|c|}
\hline
9&\Q (-6)&&\Q (-7)&&\Q (-8)&&\Q (-9)  \\
\hline
8&\Q (-5)&&\Q (-6)&&\Q (-7)&&\Q (-8)  \\
\hline
7&&&&&&&  \\
\hline
6&\Q (-4)&&\Q (-5)&&\Q (-6)&&\Q (-7)  \\
\hline
5&\Q (-3)&&\Q (-4)&&\Q (-5)&&\Q (-6)  \\
\hline
4&\Q (-3)&&\Q (-4)&&\Q (-5)&&\Q (-6)  \\
\hline
3&\Q (-2)&&\Q (-3)&&\Q (-4)&&\Q (-5)  \\
\hline
2&&&&&&&  \\
\hline
1&\Q (-1)&&\Q (-2)&&\Q (-3)&&\Q (-4)  \\
\hline
0&\Q (0)&&\Q (-1)&&\Q (-2)&&\Q (-3)  \\
\hline
\hline
&0&1&2&3&4&5&6 \\ 
\end{array}
\end{equation}

\textbf{Proof.} The cohomology groups $H^*(\Pi_{4,2}\setminus \N_{4,2},\Q)$ of the total space were partially computed in section~\ref{eq42} and are given by the following table:

\begin{table}[h]
\caption{Rational cohomology of $\Pi_{4,2}\setminus \N_{4,2}$.}
\label{gray}
\begin{equation*}
\begin{array}{c||c|c|c|c|c|c|c|c|c|c|c|c|c|c|c|c|c|}
\hline
\Q (0)&&&\Q (-2)&&\Q (-3)&&\Q (-4)&\Q (-5)&& \Q (-6)&&\Q (-7)& \Q(-8)? & \Q(-8)?&\Q (-9) \\
\hline
\hline
0&1&2&3&4&5&6&7&8&9&10&11&12&13&14&15\\ 
\end{array}
\end{equation*}
\end{table}

{\itshape The question marks in this table refer to the groups $\Q(-8)$ which are present if and only if the differential $d^3_{6,10}$ from the spectral sequence in table~\ref{spseqquart} is zero. We will shortly see that $H^*(\Pi_{4,2}\setminus \N_{4,2},\Q)$ is actually zero in degrees 13 and 14, which will complete the calculation of the spectral sequence in table~\ref{spseqquart}.}

\smallskip

We obtain from corollary~\ref{cor_groupactions} that $$E_2^{ab}=H^a(\mathcal{M}^\nod_{4,2},H^b(\GL_3(\C),\Q)).$$ The cohomology ring of $\GL_3(\C)$ is $\Lambda(\xi_1,\xi_3,\xi_5)$, where $\xi_{2i-1}$ spans $\Q(-i)$ and lives in degree $2i-1$. Together with table~\ref{gray} these facts ensure that the zeroth and the last columns of spectral sequence~(\ref{mod42}) are as shown. Since the first non-trivial cohomology group of the total space is located in degree 3 and is of weight 4, the first column of~(\ref{mod42}) is zero and the second column is as shown, moreover, $d_2^{0,1}$ is non zero. 

It remains to determine the groups $H^i(\mathcal{M}^\nod_{4,2},\Q),i=3,4,5$. Our first claim is that all these groups are 0 or pure of weight 4. Suppose the contrary. Then one of these groups contains an element of weight $>4$ or an element of weight $<4$. Let us consider the first possibility. Let $w\geq 5$ be the maximum weight, and let $m$ be the maximum degree in which elements of weight $w$ occur. Then in the group $E_2^{m,9}$ there are elements of weight $w+12$. By inspection, such an element cannot cancel with anything. The second possibility can be handled in a similar way, which proves the claim.

Our second claim is that $H^i(\mathcal{M}^\nod_{4,2},\Q)=0$ for $i=3,5$. If $E_2^{3,0}\neq 0$, then $E_2^{3,1}$ is non zero and pure of weight~6. By table~\ref{gray}, this group cannot survive in $E_\infty$, but there are no differentials that can kill it. Similarly, if $E_2^{5,0}\neq 0$, then it survives until $E_\infty$, but there are no classes of weight $4$ in $H^5(\Pi_{4,2}\setminus \N_{4,2},\Q)$ by table~\ref{gray}. Finally, if $\dim H^4(\mathcal{M}^\nod_{4,2},\Q)\neq 1$, then one of the groups $H^j(\Pi_{4,2}\setminus \N_{4,2},\Q),j=3,4$ will have the wrong dimension.$\clubsuit$

\smallskip

We are now ready to prove lemma~\ref{difq4}.

{\bf Proof of lemma~\ref{difq4}.} Suppose the differential $d^3_{6,10}$ of spectral sequence of table~\ref{spseqquart} is non zero. Then the groups $H^j(\Pi_{4,2}\setminus \N_{4,2},\Q)\cong\Q(-8)$ for $j=13,14$, and the differential $d_2^{4,9}$ of spectral sequence~(\ref{mod42}) is zero. Using multiplication in spectral sequence~(\ref{mod42}) and the fact that $d^{1,0}_2$ is an isomorphism we deduce that the rational cohomology cup product $H^2\otimes H^4\to H^6$ for $\mathcal{M}^\nod_{4,2}$ is zero. This in turn implies that the differential $d_2^{4,6}=0$. Note also that $d_2^{0,9}$ is an isomorphism, again using multiplication in the spectral sequence. The last two observations imply that $H^{10}(\Pi_{4,2}\setminus \N_{4,2},\Q)$ is two dimensional, which contradicts table~\ref{gray}.$\clubsuit$

\smallskip

\begin{lemma}\label{quartics_ring_moduli}
We have $H^*(\mathcal{M}^\nod_{4,2},\Q)\cong\Q[a]/a^4$ where $a\in H^2(\mathcal{M}^\nod_{4,2},\Q)$ is any non-zero element.
\end{lemma}

{\bf Proof.} We use the notation of lemma~\ref{mod42lemma}. Let $a_i\in H^i(\mathcal{M}^\nod_{4,2},\Q), i=2,4,6$ be a non-zero element. The argument we have already used in the proof of lemma~\ref{difq4} shows that $a_2 a_4\neq 0$. If $a_2^2=0$, then the differential $d_2^{2,1}$ of spectral sequence~(\ref{mod42}) is $0$. Comparing this with table~\ref{gray} we see that the differential $d_4^{3,0}$ must then be non-zero, which implies $E_\infty^{0,3}=0$.

Recall that the zeroth column of the $E_\infty$ page is in fact isomorphic to the image of $O^*$ where $O:\GL_3(\C)\to\Pi_{4,2}\setminus\N_{4,2}$ is an orbit map. So $E_\infty^{0,3}=0$ is impossible, because by~\cite{kolya}, the map $O$ induces a surjective map of $H^3(-,\Q)$.$\clubsuit$

\subsection{The space of equations and the space of curves}
We will now compute cohomology of the space $\nod_{4,2}=\p\Pi_{4,2}\setminus \p\N_{4,2}$. Let us summarise the results of section~\ref{moduli42}. 

\begin{itemize}
\item The rational cohomology spectral sequence for the projection $\Pi_{4,2}\setminus \N_{4,2}\stackrel{r}{\longrightarrow} \mathcal{M}^\nod_{4,2}$ is given by~(\ref{mod42}); 
\item the cohomology of the total space $\Pi_{4,2}\setminus \N_{4,2}$ is given by table~\ref{gray} with the items $\Q(-8)?$ replaced by $0$'s;
\item the cohomology of the base is given by $H^*(\mathcal{M}^\nod_{4,2},\Q)=\Q[a]/a^4$;
\item the cohomology of the fibre is isomorphic to $H^*(\GL_3(\C),\Q)\cong \Lambda (\xi_1,\xi_3,\xi_5)$;
\item identifying $\Lambda (\xi_1,\xi_3,\xi_5)$ with the zeroth column of~(\ref{mod42}), we have $d_2(\xi_1)= a$, $d_{\geq 2}(\xi_3)=0$, $d_{\geq 2}(\xi_5)=0$.
\end{itemize}

\begin{lemma}\label{quartics_ring_equations}
We have an isomorphism $H^*(\Pi_{4,2}\setminus \N_{4,2},\Q)\cong \Lambda (\zeta_3,\zeta_5',\zeta_7')$, where $\zeta_i,i=3,5,7$ is a generator of $H^i(\Pi_{4,2}\setminus \N_{4,2},\Q)$. Moreover, if $O:\GL_3(\C)\to \Pi_{4,2}\setminus\N_{4,2}$ is an orbit map, we have $O^*(\zeta_3)=\xi_3$, $O^*(\zeta_5')=\xi_5$, $O^*(\zeta_7')=0$.
\end{lemma}

\textbf{Proof.} Given the information in the above list, a straightforward check shows that the page $E_3=E_\infty$ of spectral sequence~(\ref{mod42}) is

$$
{\begin{array}{c||c|c|c|c|c|c|c|}
\hline
9&&&&&&&\Q (-9)  \\
\hline
8&\Q (-5)&&&&&&  \\
\hline
7&&&&&&&  \\
\hline
6&&&&&&&\Q (-7)  \\
\hline
5&\Q (-3)&&&&&&  \\
\hline
4&&&&&&&\Q (-6)  \\
\hline
3&\Q (-2)&&&&&&  \\
\hline
2&&&&&&&  \\
\hline
1&&&&&&&\Q (-4)  \\
\hline
0&\Q (0)&&&&&&  \\
\hline
\hline
&0&1&2&3&4&5&6 \\ 
\end{array}}
$$

Moreover, as an algebra $E_\infty$ is freely generated by elements of $E_\infty^{0,3}, E_\infty^{0,5}$ and $E_\infty^{6,1}$. We lift these generators to elements $\zeta_3,\zeta_5', \zeta_7'\in {H^*(\Pi_{4,2}\setminus \N_{4,2},\Q)}$ with $\deg \zeta_3=3,\deg \zeta_5'=5,\deg \zeta_7'=7$. After rescaling $\zeta_3$ and $\zeta_5'$ if necessary we have then $O^*(\zeta_3)=\xi_3$, $O^*(\zeta_5')=\xi_5$ and $O^*(\zeta_7')=0$. (The first two of these equalities also follow from~\cite{kolya}.) Moreover, the map $\Lambda (\zeta_3,\zeta_5',\zeta_7')\rightarrow H^*(\Pi_{4,2}\setminus \N_{4,2},\Q)$ is an isomorphism.$\clubsuit$

\begin{lemma}\label{quartics_ring_curves}
We have an isomorphism $H^*(\nod_{4,2},\Q)\cong\Lambda[\zeta_5,\zeta_7]\otimes\Q[a']/{a'}^2$. The elements $\zeta_5$ and $\zeta_7$ pull back to $\zeta_5'$, respectively $\zeta_7'$ under $p:\Pi_{4,2}\setminus \N_{4,2}\rightarrow \nod_{4,2}$, and $a'=q^*(a)$ where $q:\nod_{4,2}\to \mathcal{M}^\nod_{4,2}$ is the map into the quotient and $a$ is a multiplicative generator of $H^*(\mathcal{M}^\nod_{4,2},\Q)$, see lemma~\ref{quartics_ring_moduli}.
\end{lemma}

{\bf Proof.} We consider the spectral sequence of the fibration $X_{hG}=(X\times EG)/G\rightarrow BG$ (see section~\ref{sec_equi}) for $X=\Pi_{4,2}\setminus \N_{4,2}$ and $G=\GL_1(\C)=\C^*$. Note that since the geometric quotient $\nod_{4,2}=X/G$ exists, we have $H^*(X_{hG},\Q)\cong H^*(X/G,\Q)=H^*(\nod_{4,2},\Q)$ as mixed Hodge structures by propositions~\ref{rational_equi_quotient} and~\ref{prop_equi}. Recall that $BG\cong\C \p^\infty$. The beginning of the spectral sequence is as follows.

\begin{equation*}
\label{BG42}
\begin{array}{c||c|c|c|c|c|c|c|c}
\hline
15&\Q (-9)&&\Q (-10)&&\Q (-11)&&\Q (-12)&\cdots\\
\hline
14&&&&&&&&\cdots\\
\hline
13&&&&&&&&\cdots\\
\hline
12&\Q (-7)&&\Q (-8)&&\Q (-9)&&\Q (-10)&\cdots\\
\hline
11&&&&&&&&\\
\hline
10&\Q (-6)&&\Q (-7)&&\Q (-8)&&\Q (-9)&\cdots\\
\hline
9&&&&&&&&\\
\hline
8&\Q (-5)&&\Q (-6)&&\Q (-7)&&\Q (-8)&\cdots\\
\hline
7&\Q (-4)&&\Q (-5)&&\Q (-6)&&\Q (-7)&\cdots\\
\hline
6&&&&&&&&\\
\hline
5&\Q (-3)&&\Q (-4)&&\Q (-5)&&\Q (-6)&\cdots\\
\hline
4&&&&&&&&\\
\hline
3&\Q (-2)&&\Q (-3)&&\Q (-4)&&\Q (-5)&\cdots\\
\hline
2&&&&&&&&\\
\hline
1&&&&&&&&\\
\hline
0&\Q (0)&&\Q (-1)&&\Q (-2)&&\Q (-3)&\cdots\\
\hline
\hline
&0&1&2&3&4&5&6&\cdots\\
\end{array}
\end{equation*}

Identifying the zeroth column with $H^*(\Pi_{4,2}\setminus\N_{4,2},\Q)$ described in lemma~\ref{quartics_ring_equations}, we have $d_2=0$ on the generators $\zeta_3,\zeta_5',\zeta_7'$, so the differential $d_2$ is zero, and hence so is $d_3$. By proposition~\ref{nodstart} we have $H^3(\nod_{4,2})=0$, therefore $d_4^{0,3}$ is an isomorphism. This implies that $d_4$ is non zero whenever this is allowed by the weights, and as a consequence, the page $E_5=E_\infty$ is as follows.

\begin{equation*}
\begin{array}{c||c|c|c|}
\hline
12&\Q (-7)&&\Q (-8)\\
\hline
11&&&\\
\hline
10&&&\\
\hline
9&&&\\
\hline
8&&&\\
\hline
7&\Q (-4)&&\Q (-5)\\
\hline
6&&&\\
\hline
5&\Q (-3)&&\Q (-4)\\
\hline
4&&&\\
\hline
3&&&\\
\hline
2&&&\\
\hline
1&&&\\
\hline
0&\Q (0)&&\Q (-1)\\
\hline
\hline
&0&1&2\\
\end{array}
\end{equation*}

The zeroth column of this page is isomorphic to the image of $H^*(X_{hG},\Q)\cong H^*(\nod_{4,2},\Q)$ in $H^*(X,\Q)=H^*(\Pi_{4,2}\setminus\N_{4,2},\Q)$. We set $\zeta_5, \zeta_7\in H^*(\nod_{4,2},\Q)$ to be lifts of the generators $\zeta_5'\in E_\infty^{0,5}$ and $\zeta_7'\in E_\infty^{0,7}$ respectively. Set also $a'\in H^2(\nod_{4,2},\Q)$ to be the unique lift of a generator of $E_\infty^{2,0}$. We have $a'^2=0$ from the spectral sequence, and so we obtain a ring isomorphism $$\Lambda[\zeta_5,\zeta_7]\otimes\Q[a']/{a'}^2\rightarrow H^*(\nod_{4,2},\Q)$$
such that $p^*$ takes the image of $\zeta_i,i=5,7$ to $\zeta_i'$.

To finish the proof of lemma~\ref{quartics_ring_curves} we need to determine the map $q^*$. Using what we already know we calculate the $E_\infty$ page of the Leray spectral sequence of $q:\nod_{4,2}\rightarrow \mathcal{M}^\nod_{4,2}$ and obtain the following result.

$$
\begin{array}{c||c|c|c|c|c|c|c|}
\hline
8&&&&&\Q (-7)&&\Q (-8)  \\
\hline
7&&&&&&&  \\
\hline
6&&&&&&&  \\
\hline
5&\Q (-3)&&\Q (-4)&&&&  \\
\hline
4&&&&&&&  \\
\hline
3&&&&&\Q (-4)&&\Q (-5)  \\
\hline
2&&&&&&&  \\
\hline
1&&&&&&&  \\
\hline
0&\Q (0)&&\Q (-1)&&&&  \\
\hline
\hline
& 0 & 1 & 2 & 3 & 4 & 5 & 6\\
\end{array}
$$

We conclude that $q$ induces an isomorphism of $H^2(-,\Q)$, which completes the proof of lemma~\ref{quartics_ring_curves}.$\clubsuit$

Theorem~\ref{mhquartp2} follows from lemmas~\ref{quartics_ring_moduli}, \ref{quartics_ring_equations} and~\ref{quartics_ring_curves}.$\clubsuit$


\section{Nodal cubics in \texorpdfstring{$\C \p^3$}{}}\label{cubics}

In this section we prove theorem~\ref{mhcubp3}. 

\subsection{The moduli space}

We consider the space $\Hyp_{3,3}=\p(\Pi_{3,3})$, the projectivisation of the space of homogeneous polynomials of degree~3 on $\mathbb{C}^4$. 
For a cubic polynomial $Q\in\Pi_{3,3}$, we denote the corresponding point of $\Hyp_{3,3}$ by $|Q|$. The group $\SL_4(\C)$ acts naturally on the argument space $\mathbb{C}^4$, on $\Hyp_{3,3}$, and on the line bundle $\OO(1)$ over $\Hyp_{3,3}$. In the rest of the section we consider (non) stable and semistable elements of $\Hyp_{3,3}$ with respect to $\OO(1)$ and its natural linearisation.

Recall (see Introduction) that we denote the subspace of $\Hyp_{3,3}$ formed by nodal cubics by $\nod_{3,3}=\p\Pi_{3,3}\setminus \p\N_{3,3}$. Our goal here is to describe the corresponding moduli space, namely $\mathcal{M}^\nod_{3,3}=\text{Nod}_{3,3}/\SL_4(\C)$.

\begin{lemma}\label{lemma_moduli_cubics_p3}
Let $Q\in\Pi_{3,3}$. The cubic $|Q|$ is stable if and only if $Q\in\Pi_{3,3}\setminus\N_{3,3}$ and semistable but not stable if and only if the closure of the $\SL_4(\C)$-orbit of $|Q|$ intersects the closure of the orbit of the cubic $|Q_s|$ defined in the proof (see~(\ref{qs})).

The categorical quotient $\Hyp_{3,3}^{ss}//\SL_4(\C)$ is isomorphic to the weighted projective space $\C\p(1,2,3,4,5)$. The geometric quotient $$\mathcal{M}^\nod_{3,3}=\nod_{3,3}/\SL_4(\C)=\Hyp_{3,3}^{s}/\SL_4(\C)$$ is isomorphic to $\C\p(1,2,3,4,5)\setminus\{(0:0:0:0:1)\}$, and we have
$$H^*(\mathcal{M}^\nod_{3,3},\Q)\cong \Q[a]/a^4, \deg a=2,\mbox{\normalfont{ and }} P_{mH}(\mathcal{M}^\nod_{3,3})=1+t^2w^1+t^4w^2+t^6w^3.$$
\end{lemma}

{\bf Proof.} In the proof $X_1,X_2, X_3$ are coordinates in $\C^3$. We use $$O(|(\mbox{some variables})|^m\cdot\mbox{a monomial $M$})$$ to denote the product of $M$ and a not necessarily homogeneous polynomial of degree $\geqslant m$ in the variables. We also write $O(|X|^m \cdot M)$ instead of $O(|(X_1,X_2,X_3)|^m \cdot M)$. For example, $O(|(X_1,X_2)|^2\cdot X_3)$ stands for the product of $X_3$ and a degree $\geqslant 2$ polynomial in $X_1,X_2$, and $O(|X|^m)$ is a polynomial in $X_1,X_2,X_3$ of degree $\geqslant m$.

Using~\cite[Chapter 4.2]{mum}, we see that for $Q\in\Pi_{3,3}$, the cubic $|Q|$ 
is {\itshape not stable} ($|Q|\notin \Hyp_{3,3}^s$) if and only if it has a singular point $x_0$ of one of the following types:
\begin{equation}
\label{sing1}
\begin{cases}
Q= X_1X_2+O(|X|^3) & \text{ (binode, generic \includegraphics[scale=0.1]{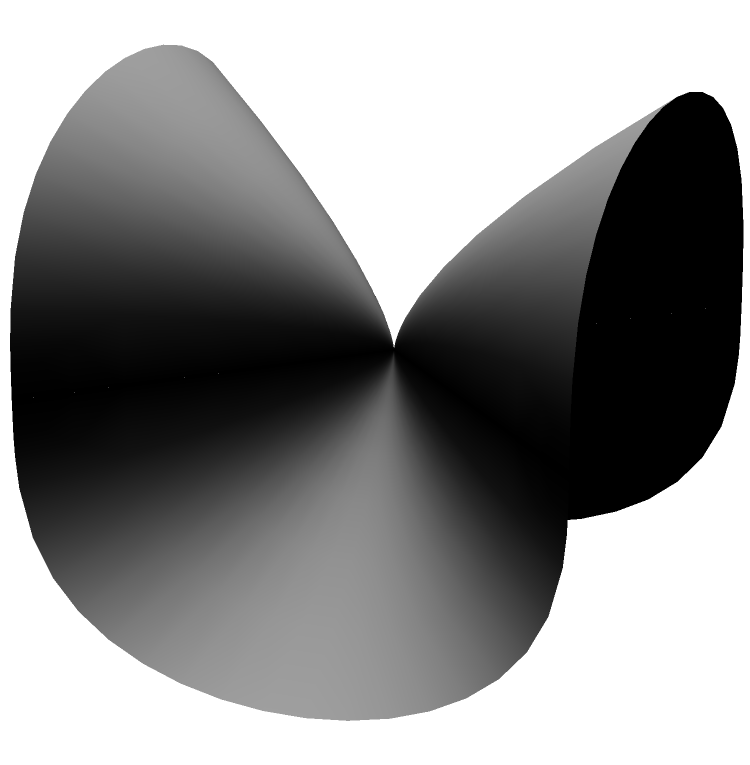})},\\
Q= X_1^2+O(|X|^3) & \text{ (unode, generic \includegraphics[scale=0.1]{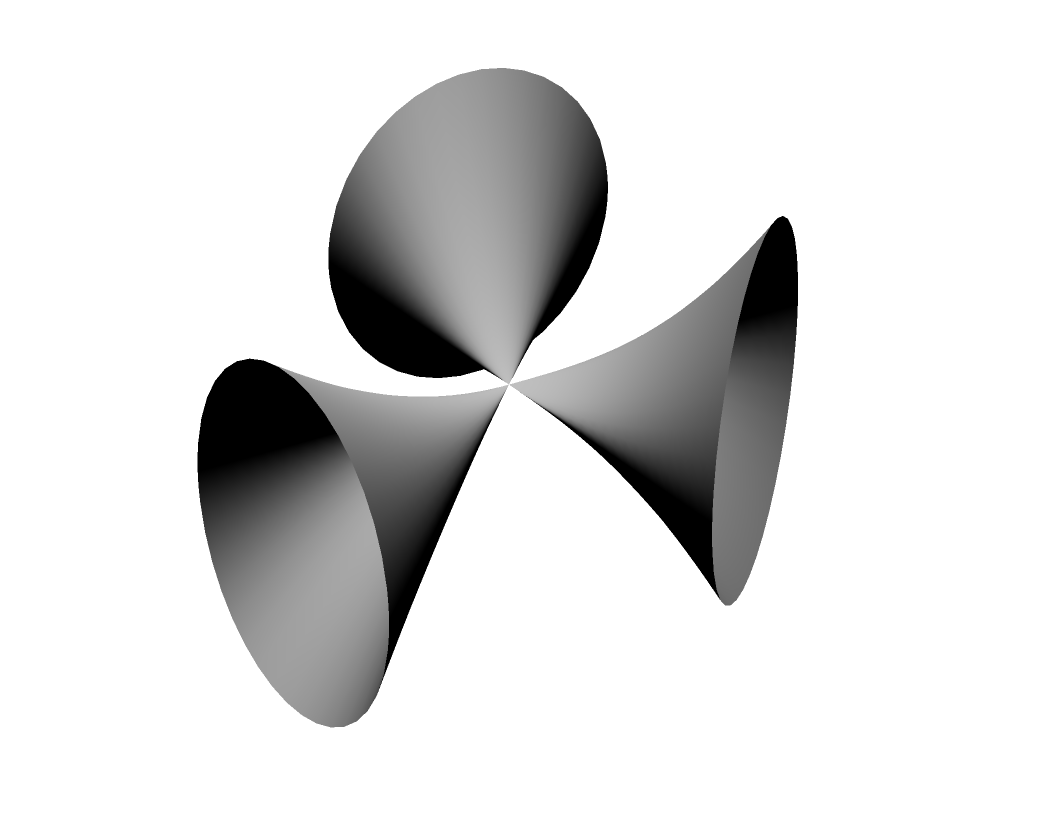})},\\
Q= O(|X|^3) & \text{ (triple point, generic \includegraphics[scale=0.1]{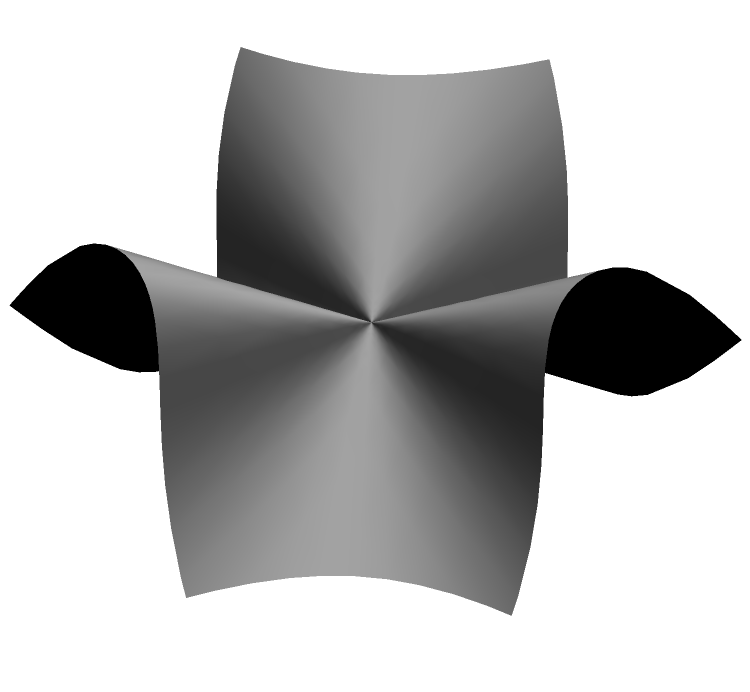})};
\end{cases}   
\end{equation}
or equivalently, if and only if it has a critical point $x_0$ such that the quadratic part of the expansion of $Q$ at $x_0$ is degenerate. So $|Q|\in \Hyp_{3,3}$ is stable if and only if every singular point of $Q$ is non degenerate, i.e.\ Morse, or in other words, if and only if $Q$ is nodal.

Similarly, using~\cite[Chapter 4.2]{mum} we conclude that a cubic $|Q|$ in $\Hyp_{3,3}$ is {\itshape not semistable} (${|Q|\notin \Hyp_{3,3}^{ss}}$), if and only if it has a singular point $x_0$ of one of the following types:
\begin{equation}
\label{sing2}
\begin{cases}
Q= X_1X_2+O(|X|^2\cdot X_1)+O(|X|^2\cdot X_2) & \text{ (the $X_3^3$ term is absent; binode of higher type, generic \includegraphics[scale=0.1]{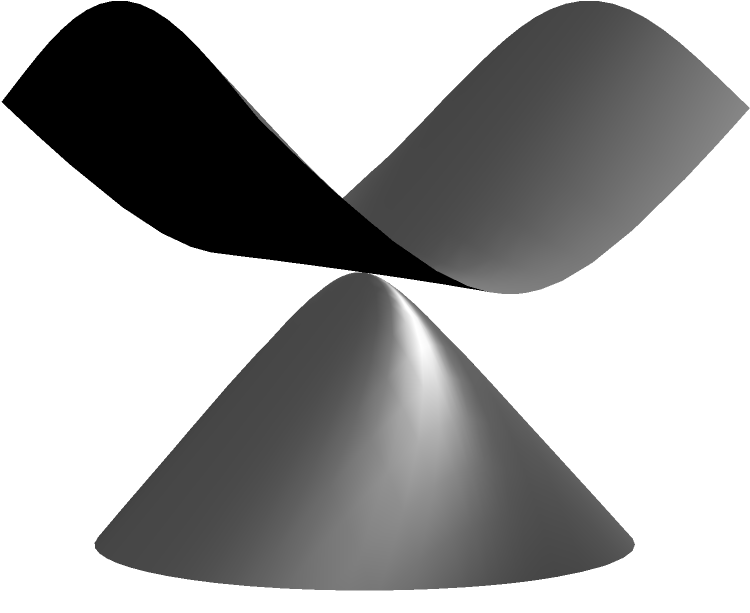})}, \\
Q= X_1^2+O(|X|^3) & \text{ (unode, generic  \includegraphics[scale=0.1]{tritrumpet.pdf})},\\
Q= O(|X|^3) & \text{ (triple point, generic  \includegraphics[scale=0.1]{triple.pdf})}.
\end{cases}
\end{equation}
Here a binode of higher type is a binode $x_0$ which is not resolved by a single blow-up at $x_0$. Comparing~\eqref{sing1} and~\eqref{sing2}, one deduces that all semistable but not stable cubics $Q$ in $\mathbb{C}\p^3$ have a simple binode, namely a point $x_0$ such that in some coordinates $(X_1,X_2,X_3)$ centred at $x_0$ the equation of the cubic takes the form
$$Q=X_1X_2+X_3^3+O(|X|^2\cdot X_1)+O(|X|^2 \cdot X_2),$$ or equivalently, in coordinates $(\mathrm{x}_0,\mathrm{x}_1,\mathrm{x}_2,\mathrm{x}_3)$ in $\mathbb{C}^4$ such that $x_0=(0:0:1:0)$ we have
\begin{equation}
\label{q}
Q=\mathrm{x}_0\mathrm{x}_1\mathrm{x}_2+\mathrm{x}_3^3+O(|(\mathrm{x}_0,\mathrm{x}_1,\mathrm{x}_3)|^2\cdot \mathrm{x}_0)+O(|(\mathrm{x}_0,\mathrm{x}_1,\mathrm{x}_3)|^2\cdot \mathrm{x}_1).
\end{equation}
Let now $Q_s$ be given by
\begin{equation}
\label{qs}
Q_s(\mathrm{x}_0,\mathrm{x}_1,\mathrm{x}_2,\mathrm{x}_3)=\mathrm{x}_0\mathrm{x}_1\mathrm{x}_2+\mathrm{x}_3^3.
\end{equation}
If any cubic $Q$ is semistable but not properly stable, and therefore can be put in the form~\eqref{q}, then $|Q_s|\in \overline{\SL_4(\C)(|Q|)}$ since 
\begin{gather}
\lim\limits_{t\rightarrow \infty}\textbf{diag}(t,t,t^{-2},1)\cdot(\mathrm{x}_0\mathrm{x}_1\mathrm{x}_2+\mathrm{x}_3^3+O(|(\mathrm{x}_0,\mathrm{x}_1,\mathrm{x}_3)|^2\mathrm{x}_0)+O(|(\mathrm{x}_0,\mathrm{x}_1,\mathrm{x}_3)|^2\mathrm{x}_1))=\notag\\
=\mathrm{x}_0\mathrm{x}_1\mathrm{x}_2+\mathrm{x}_3^3=Q_s(\mathrm{x}_0, \mathrm{x}_1, \mathrm{x}_2, \mathrm{x}_3).\notag
\end{gather}
(Here we consider the {\it left} action of $\SL_4(\C)$ on polynomials; in particular, for $g\in SL_4(\C)$ and $\mathrm{x}=(\mathrm{x}_0, \mathrm{x}_1, \mathrm{x}_2, \mathrm{x}_3)^T$ we have $(g\cdot Q)(\mathrm{x})=Q(g^{-1}\mathrm{x})$.) So in the categorical quotient $\Hyp_{3,3}//\SL_4(\C)$ all non stable cubics are represented by a single point, $|Q_s|$. The corresponding cubic surface is shown below.

\includegraphics[scale=0.25]{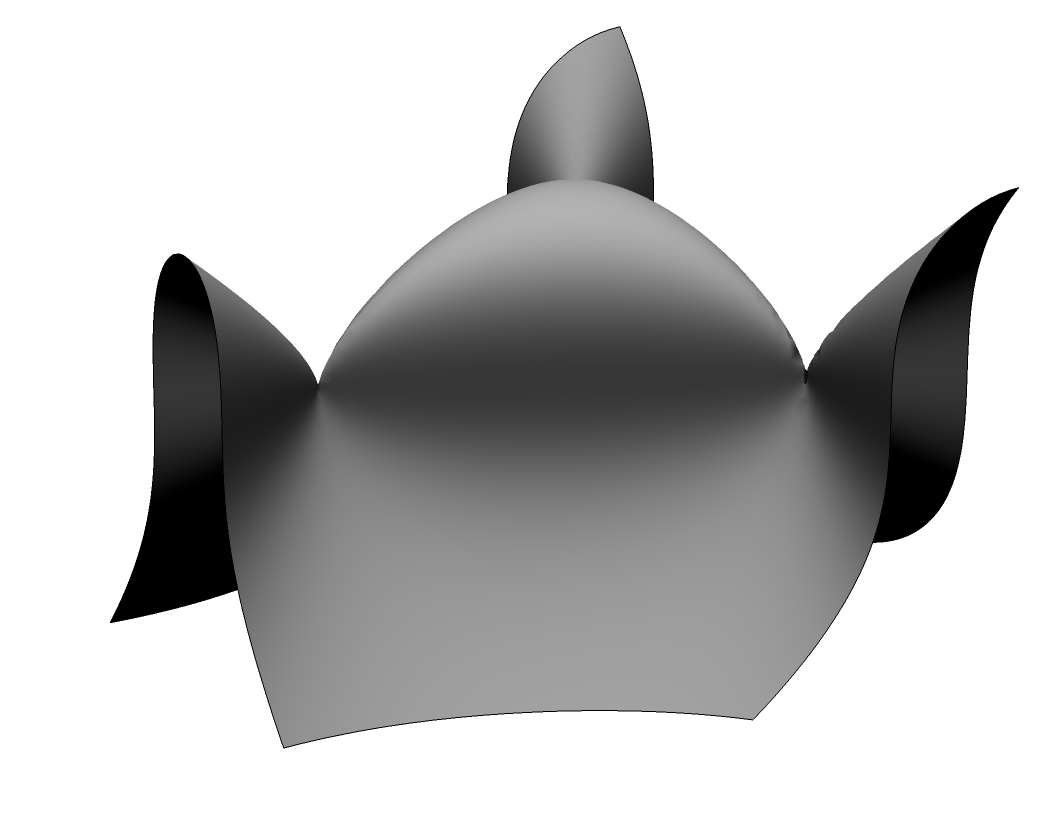}

We conclude that
\begin{equation}
\label{1st}
\mathcal{M}^\nod_{3,3}=\nod_{3,3}/\SL_4(\C)\cong \Hyp_{3,3}^s/\SL_4(\C)\cong (\Hyp_{3,3}//\SL_4(\C))\setminus \{|Q_s|\}
\end{equation}  

Now we need to understand $C_{3,3}=\Hyp_{3,3}^{ss}//\SL_4(\C)$ and locate the point $|Q_s|$ there. For these purposes we use the description of the quotient $C_{3,3}$ provided in~\cite[Chapter 10.4]{dolgachev}, which expresses $C_{3,3}$ as a weighted projective space: 
$$C_{3,3}\cong\C\p(1,2,3,4,5)\cong \text{Projm} (k[I_8,I_{16},I_{24},I_{32},I_{40}]),$$ 
where $I_j, j=8,16,24,32,40$ is a certain $\SL_4(\C)$-invariant polynomial on $\Pi_{3,3}$ of degree $j$. An explicit construction of ${I_n}$ is provided in~\cite[Chapter XV]{sal} (the invariants $I_8,I_{16},I_{24},I_{32},I_{40}$ are denoted there $A,B,C,D,E$ respectively). Computing their values on a cubic $Q$ is particularly easy if one is able to find coordinates $(\mathrm{x},\mathrm{y},\mathrm{z},\mathrm{v})$ in $\mathbb{C}^4$ such that $Q$ becomes a sum of five cubes:
\begin{gather}
\label{5}
Q(\mathrm{x},\mathrm{y},\mathrm{z},\mathrm{v})=a\mathrm{x}^3+b\mathrm{y}^3+c\mathrm{z}^3+d\mathrm{v}^3+e\mathrm{w}^3\\
\text{where }\mathrm{w}=-(\mathrm{x}+\mathrm{y}+\mathrm{z}+\mathrm{v}).\notag
\end{gather}
Then, letting $\sigma_i$ be the $i$-th elementary symmetric polynomial in the variables $(a,b,c,d,e)$ defined by~\eqref{5}, we get the following values of the invariants at $Q$:
\begin{equation}
\label{sigma}
I_{8}=\sigma_4^2-4\sigma_3\sigma_5,\ \ I_{16}=\sigma_1\sigma_5^3,\ \ I_{24}=\sigma_4\sigma_5^4,\ \ I_{32}=\sigma_2\sigma_5^6,\ \ I_{40}=\sigma_5^8,
\end{equation}
see e.g.~\cite[formula (9.59)]{dolgachev1}. 
Instead of computing the invariants for $Q_s$ directly, we compute them for a family of semistable cubics $Q(t)\stackrel{t\rightarrow 0}{\longrightarrow} Q_s$ which are in the form~(\ref{5}).  
We let $Q(t)$ be the cubic obtained by substituting
\begin{equation}\label{subs}
\mathrm{x}=\mathrm{x}_0-\mathrm{x}_1-\mathrm{x}_2, \mathrm{y}=\mathrm{x}_1-\mathrm{x}_2-\mathrm{x}_0, \mathrm{z}=\mathrm{x}_2-\mathrm{x}_0-\mathrm{x}_1, \mathrm{v}=\mathrm{x}_0+\mathrm{x}_1+\mathrm{x}_2-t\mathrm{x}_3, \mathrm{w}=t\mathrm{x}_3
\end{equation}
in 
\begin{gather}
\label{t}
\mathrm{x}^3+\mathrm{y}^3+\mathrm{z}^3+\mathrm{v}^3+\frac{1}{T}\mathrm{w}^3
\text{ where }T=\frac{t^3}{24}.
\end{gather}
Note that the sum of the right hand sides of (\ref{subs}) is $0$, so expression~\eqref{t} is of the form~\eqref{5}. 

A direct computation shows that $24Q_s=\lim\limits_{t\rightarrow 0} Q(t)$. 
Since
$$(a,b,c,d,e)(Q(t))=\left(1,1,1,1,\frac{1}{T}\right),$$ 
we see that the elementary symmetric polynomials $\sigma_i(a,b,c,d,e)$ are given by
$$
\sigma_1=\frac{1+4T}{T},\ \sigma_2=\frac{4+6T}{T},\ \sigma_3=\frac{6+4T}{T},\ \sigma_4=\frac{4+T}{T},\ \sigma_5=\frac{1}{T}.
$$
Substituting these in~\eqref{sigma} we get
\begin{gather}
(I_{8},I_{16},I_{24},I_{32},I_{40})(Q(t))=\left(\sigma_4^2-4\sigma_3\sigma_5,\sigma_1\sigma_5^3,\sigma_4\sigma_5^4,\sigma_2\sigma_5^6,\sigma_5^8\right)=\notag\\
=\left(\frac{T^2-8T-8}{T^2},\frac{1+4T}{T^4},\frac{4+T}{T^5},\frac{4+6T}{T^7},\frac{1}{T^8}\right).\notag
\end{gather}
We conclude that every cubic $|Q(t)|,t\in\C^*$ is semistable, and the map $\Hyp_{3,3}^{ss}\to\Hyp^{ss}_{3,3}//\SL_4(\C)$ takes $|Q(t)|$ to the point
\begin{equation}\label{quot_point_family}
\left(\frac{T^2-8T-8}{T^2}:\frac{1+4T}{T^4}:\frac{4+T}{T^5}:\frac{4+6T}{T^7}:\frac{1}{T^8}\right)\in\C\p (8,16,24,32,40)\cong\C\p (1,2,3,4,5).
\end{equation}
(In fact, all surfaces $|Q(t)|$ are smooth and hence stable, but we will not need this.) As $t\to 0$, the limit of (\ref{quot_point_family}) is $(0:0:0:0:1)$, which is the point of the categorical quotient that corresponds to $|Q_s|$.$\clubsuit$

\subsection{The space of equations}
\label{eq33}
Recall that we have $H^*(\GL_4(\C),\mathbb{Q})\cong \Lambda(\xi_1,\xi_3,\xi_5,\xi_7)$, where $\xi_{2i-1}$ spans $\Q(-i)$ and has degree $2i-1$. Recall also that we know by lemma~\ref{lemma_moduli_cubics_p3} that there exists a geometric quotient $\nod_{3,3}/\SL_4(\C)$, or equivalently, a geometric quotient $(\Pi_{3,3}\setminus \N_{3,3})/\GL_4(\C)$.

\begin{lemma}\label{lemma_equations_cubics_p3}
We have $$H^*(\Pi_{3,3}\setminus\N_{3,3},\Q)\cong\Lambda[\zeta'_3,\zeta'_5,\zeta'_7,b']$$ where $\zeta_i',i=3,5,7$ has degree $i$ and is Tate of weight $i+1$ and $b'$ is Tate of weight $8$. If $O:\GL_4(\C)$ is an orbit map, we have $O^*(\zeta'_i)=\xi_i$ for $i=3,5,7$ and $O^*(b')=0$. We have
\begin{equation*}
P_{\mathrm{mH}}(\Pi_{3,3}\setminus\N_{3,3})=(1+t^3 w^2)(1+t^5 w^3)(1+t^7 w^4)^2.
\end{equation*}

\end{lemma}

{\bf Proof.} We obtain using corollary~\ref{cor_groupactions} that the second page of the rational cohomology Leray spectral sequence of the map $$(\Pi_{3,3}\setminus \N_{3,3})\rightarrow (\Pi_{3,3}\setminus \N_{3,3})/\GL_4(\C)=\mathcal{M}^\nod_{3,3}\cong \C\p(1,2,3,4,5)\setminus \{(0:0:0:0:1)\}$$ is given by
$$E_2^{ab}=H^a(\mathcal{M}^\nod_{3,3},H^b(\GL_4(\C),\Q)).$$
Using lemma~\ref{lemma_moduli_cubics_p3} we get the following result.
$$
\begin{array}{c||c|c|c|c|c|c|c|}
&0&1&2&3&4&5&6\\
\hline
\hline
16&\Q(-10)&&\Q(-11)&&\Q(-12)&&\Q(-13)\\
\hline
15&\Q(-9)&&\Q(-10)&&\Q(-11)&&\Q(-12)\\
\hline
14&&&&&&&\\
\hline
13&\Q(-8)&&\Q(-9)&&\Q(-10)&&\Q(-11)\\
\hline
12&\Q(-7)&&\Q(-8)&&\Q(-9)&&\Q(-10)\\
\hline
11&\Q(-7)&&\Q(-8)&&\Q(-9)&&\Q(-10)\\
\hline
10&\Q(-6)&&\Q(-7)&&\Q(-8)&&\Q(-9)\\
\hline
9&\Q(-6)&&\Q(-7)&&\Q(-8)&&\Q(-9)\\
\hline
8&\Q(-5)\oplus\Q(-5)&&\Q(-6)\oplus\Q(-6)&&\Q(-7)\oplus\Q(-7)&&\Q(-8)\oplus\Q(-8)\\
\hline
7&\Q(-4)&&\Q(-5)&&\Q(-6)&&\Q(-7)\\
\hline
6&\Q(-4)&&\Q(-5)&&\Q(-6)&&\Q(-7)\\
\hline
5&\Q(-3)&&\Q(-4)&&\Q(-5)&&\Q(-6)\\
\hline
4&\Q(-3)&&\Q(-4)&&\Q(-5)&&\Q(-6)\\
\hline
3&\Q(-2)&&\Q(-3)&&\Q(-4)&&\Q(-5)\\
\hline
2&&&&&&&\\
\hline
1&\Q(-1)&&\Q(-2)&&\Q(-3)&&\Q(-4)\\
\hline
0&\Q(0)&&\Q(-1)&&\Q(-2)&&\Q(-3)\\
\hline
\end{array}
$$ 

Let us identify column $0$ with $H^*(\GL_4(\C);\mathbb{Q})\cong \Lambda(\xi_1,\xi_3,\xi_5,\xi_7)$. Note that the differential $d_2^{0,1}$ is an isomorphism because otherwise $H^1(\Pi_{3,3}\setminus \N_{3,3},\Q)\neq 0$, which is impossible since $\N_{3,3}$ has codimension 2 in $\Pi_{3,3}\cong \mathbb{C}^{20}$. Using the cup product structure of $H^*(\mathcal{M}^\nod_{3,3},\Q)$, which we know by lemma~\ref{lemma_moduli_cubics_p3}, we deduce that all differentials $d_2^{a,b}$ which may be non zero for weight reasons are in fact non zero. So the $E_3$ page, and also the $E_\infty$ page is as follows.

\begin{equation}\label{e_infty_equations_cubic_surfaces}
\begin{array}{c||c|c|c|c|c|c|c|}
&0&1&2&3&4&5&6\\
\hline
\hline
16&&&&&&&\Q(-13)\\
\hline
15&\Q(-9)&&&&&&\\
\hline
14&&&&&&&\\
\hline
13&&&&&&&\Q(-11)\\
\hline
12&\Q(-7)&&&&&&\\
\hline
11&&&&&&&\Q(-10)\\
\hline
10&\Q(-6)&&&&&&\\
\hline
9&&&&&&&\Q(-9)\\
\hline
8&\Q(-5)&&&&&&\Q(-8)\\
\hline
7&\Q(-4)&&&&&&\\
\hline
6&&&&&&&\Q(-7)\\
\hline
5&\Q(-3)&&&&&&\\
\hline
4&&&&&&&\Q(-6)\\
\hline
3&\Q(-2)&&&&&&\\
\hline
2&&&&&&&\\
\hline
1&&&&&&&\Q(-4)\\
\hline
0&\Q(0)&&&&&&\\
\hline
\end{array}
\end{equation}

Let us pick elements $\zeta'_i\in H^i(\Pi_{3,3}\setminus\N_{3,3},\Q), i=3,5,7$ that project to generators of $E_\infty^{0,i}$. Recall that the zeroth column of $E_\infty$ is the image of $O^*$, so by rescaling $\zeta'_i$'s if necessary we have $O^*(\zeta'_i)=\xi_i\in H^i(\GL_3(\C),\Q)$ for $i=3,5,7$ (cf.~\cite{kolya}). Let us also choose a non-zero element $b'\in H^7(\Pi_{3,3}\setminus\N_{3,3},\Q)$ that comes from $E_\infty^{6,1}$. We have then $O^*(b')=0$, and the natural map $\Lambda(\zeta'_3,\zeta'_5,\zeta'_7,b')\rightarrow H^*(\Pi_{3,3}\setminus \N_{3,3},\Q)$ is an isomorphism.$\clubsuit$

When computing $H^*(\nod_{3,3},\Q)$ in section~\ref{sec_space_cub_surfaces} we will show that $b'=p^*(b)$, $\zeta'_5=p^*(\zeta_5)$ and $\zeta'_7=p^*(\zeta_7)$ for $p:\Pi_{3,3}\setminus \N_{3,3}\rightarrow \nod_{3,3}$ equal the natural projection, and certain classes $b,\zeta_5,\zeta_7\in H^*(\nod_{3,3},\Q)$. 

\subsection{The space of surfaces}\label{sec_space_cub_surfaces}

Here we calculate the cohomology of $\nod_{3,3}$, the space of nodal cubics in $\C \p^3$. We proceed as in section~\ref{eq33}. Recall that we have $H^*(\PGL_4(\C),\Q)\cong \Lambda(\xi_3,\xi_5,\xi_7)$ where $\xi_i,i=3,5,7$ spans $\Q(-i)$ and has degree $2i-1$. Note that we use the same notation for the generators of $H^*(\PGL_4(\C),\Q)$ as we used for the generators of $H^*(\GL_4(\C),\Q)$ in section~\ref{eq33}: these generators can be chosen to correspond to each other under the pullback map induced by $\GL_4(\C)\to\PGL_4(\C)$.

\begin{lemma}\label{lemma_surfaces_cubics_p3}
We have $$H^*(\nod_{3,3},\Q)\cong\Lambda[\zeta_5,\zeta_7,b]\otimes\Q[a']/a'{}^2$$ where $\zeta_i,i=5,7$ has degree $2i-1$ and is Tate of weight $i$, the element $a'$ has degree $2$ and is Tate of weight 2, and $b$ has degree $7$ and is Tate of weight $8$. We have
\begin{equation*}
P_{\mathrm{mH}}(\Pi_{3,3}\setminus\N_{3,3})=(1+t^2 w)(1+t^5 w^3)(1+t^7 w^4)^2.
\end{equation*}
If $O:\PGL_4(\C)\to\nod_{3,3}$ is an orbit map, we have $O^*(\zeta_i)=\xi_i, i=5,7$ and $O^*(b)=O^*(a')=0$. If $p:\Pi_{3,3}\setminus\N_{3,3}\to\nod_{3,3}$ and $q:\Pi_{3,3}\setminus\N_{3,3}\to \mathcal{M}^\nod_{3,3}$ are the natural projections, we have $p^*(\zeta_i)=\zeta'_i$ for $i=5,7$, $p^*(b)=b'$ and $q^*(a)=a'$ where $\zeta'_5$ and $\zeta'_7$ are the generators of $H^*(\Pi_{3,3}\setminus\N_{3,3},\Q)$ defined in lemma~\ref{lemma_equations_cubics_p3}, and $a\in H^2(\mathcal{M}^\nod_{3,3},\Q)$ is the generator from lemma~\ref{lemma_moduli_cubics_p3}.
\end{lemma}

{\bf Proof.} As in section~\ref{eq33} we find that the $E_2$ page of the rational cohomology Leray spectral sequence of the quotient map $$\nod_{3,3}\rightarrow \nod_{3,3}/\PGL_4(\C)\cong \C\p(1,2,3,4,5)\setminus \{(0:0:0:0:1)\}$$
is
$$
\begin{array}{c||c|c|c|c|c|c|c|}
&0&1&2&3&4&5&6\\
\hline
\hline
15&\Q(-9)&&\Q(-10)&&\Q(-11)&&\Q(-12)\\
\hline
14&&&&&&&\\
\hline
13&&&&&&&\\
\hline
12&\Q(-7)&&\Q(-8)&&\Q(-9)&&\Q(-10)\\
\hline
11&&&&&&&\\
\hline
10&\Q(-6)&&\Q(-7)&&\Q(-8)&&\Q(-9)\\
\hline
9&&&&&&&\\
\hline
8&\Q(-5)&&\Q(-6)&&\Q(-7)&&\Q(-8)\\
\hline
7&\Q(-4)&&\Q(-5)&&\Q(-6)&&\Q(-7)\\
\hline
6&&&&&&&\\
\hline
5&\Q(-3)&&\Q(-4)&&\Q(-5)&&\Q(-6)\\
\hline
4&&&&&&&\\
\hline
3&\Q(-2)&&\Q(-3)&&\Q(-4)&&\Q(-5)\\
\hline
2&&&&&&&\\
\hline
1&&&&&&&\\
\hline
0&\Q(0)&&\Q(-1)&&\Q(-2)&&\Q(-3)\\
\hline
\end{array}
$$ 

Identifying column 0 with $H^*(\PGL_4(\C),\Q)\cong \Lambda(\xi_3,\xi_5,\xi_7)$ we deduce that $d_2=d_3=0$. By proposition~\ref{nodstart} we have $H^3(\nod_{3,3},\Q)=0$, so $d_4^{0,3}$ is an isomorphism.
As in the proof of lemma~\ref{lemma_equations_cubics_p3} we deduce from $d_4^{0,3}\neq 0$ that $d_4^{a,b}\neq 0$ whenever the source and target have the same weight, and as a result, that the $E_5=E_\infty$ page is as follows.

$$
\begin{array}{c||c|c|c|c|c|c|c|}
&0&1&2&3&4&5&6\\
\hline
\hline
15&&&&&\Q(-11)&&\Q(-12)\\
\hline
14&&&&&&&\\
\hline
13&&&&&&&\\
\hline
12&\Q(-7)&&\Q(-8)&&&&\\
\hline
11&&&&&&&\\
\hline
10&&&&&\Q(-8)&&\Q(-9)\\
\hline
9&&&&&&&\\
\hline
8&&&&&\Q(-7)&&\Q(-8)\\
\hline
7&\Q(-4)&&\Q(-5)&&&&\\
\hline
6&&&&&&&\\
\hline
5&\Q(-3)&&\Q(-4)&&&&\\
\hline
4&&&&&&&\\
\hline
3&&&&&\Q(-4)&&\Q(-5)\\
\hline
2&&&&&&&\\
\hline
1&&&&&&&\\
\hline
0&\Q(0)&&\Q(-1)&&&&\\
\hline
\end{array}
$$ 

By comparing this with~(\ref{e_infty_equations_cubic_surfaces}) we see that one can pick generators $\zeta_i,i=5,7$ of $H^i(\nod_{3,3},\Q)$ such that $p^*(\zeta_i)=\zeta'_i\in H^i(\Pi_{3,3}\setminus \N_{3,3},\Q)$ and hence $O^*(\zeta_i)=\xi_i\in H^i(\PGL_4(\C),\Q)$. Also set $a'=q^*(a)\in H^2(\nod_{3,3},\Q)\cong E_\infty^{2,0}$. Let $b\in H^7(\nod_{3,3},\Q)$ be a non-zero element that comes from $E_\infty^{4,3}$. We conclude from the spectral sequence that resulting map $$\Lambda(\zeta_5,\zeta_7,b)\otimes\Q[a']/a'{}^2\rightarrow H^*(\nod_{3,3})$$ is an isomorphism, and that $O^*(b)=O^*(a')=0$ (recall that all elements of $H^*(\nod_{3,3},\Q)$ that map to $E_\infty^{>0,*}$ go to $0$ under $O^*$).

It remains to calculate $p^*(b)$. From the spectral sequence of $p:\Pi_{3,3}\setminus \N_{3,3}\to\nod_{3,3}$ we see that $\ker p^*\subset H^*(\nod_{3,3},\Q)$ is the ideal generated by $a'$. This ideal does not contain $b$, so $p^*(b)\neq 0$. On the other hand, $b\in\ker O^*$, so $p^*(b)$ is in the kernel of the map $O'{}^*$ induced by an orbit map $O':\GL_4(\C)\to\Pi_{3,3}\setminus\N_{3,3}$. It follows from~(\ref{e_infty_equations_cubic_surfaces}) that the kernel of $O'{}^*$ in degree 7 is spanned by $b'$, so by rescaling $b$ if necessary we get $p^*(b)=b'$.$\clubsuit$


Theorem~\ref{mhcubp3} follows from lemmas~\ref{lemma_moduli_cubics_p3}, \ref{lemma_equations_cubics_p3} and~\ref{lemma_surfaces_cubics_p3}.$\clubsuit$

\appendix

\section{Quotients, slices, and equivariant cohomology}
\label{b}

In the main text of the paper we often encountered the following problem: suppose $X$ is a space on which a connected complex algebraic group $G$ acts with finite stabilisers. One would like to deduce from this that every higher direct image $R^i p_* \underline{\Q}_X$ under the projection $p:X\to X/G$ is constant of rank $\dim_\Q H^i(G,\Q)$. Here $X/G$ is equipped with the quotient topology induced from $X$.
In our examples $X$ is a complex algebraic variety and the action of $G$ on $X$ is algebraic. Here we collect several facts about algebraic and topological quotients and explain how together they imply the required result about $R^i p_* \underline{\Q}_X$ in our examples (theorem~\ref{groupactions} and corollary~\ref{cor_groupactions}).

The material we discuss here is standard or folklore (or both), but not all of it is easy to find in the literature. When we were unable to do so, we provided proofs. We refer the reader to~\cite[Section 2]{kirwan_et_al} for a brief review of geometric invariant theory, and to~\cite{vinpop}, \cite{newstead} and~\cite{dolgachev} for a more detailed account. The original sources for most of the appendix are~\cite{mum} and~\cite{luna}.

In section~\ref{sec_geo_quot} we prove propoition~\ref{examples_quotients}, in which we summarise the tools which allow one to show the existence of a geometric quotient in all examples we consider in this paper, and also a technical result about principal and associated bundles (proposition~\ref{assoc_bund}). In section~\ref{sec_slices} we review slices and prove that under some conditions, the direct images of a locally constant sheaf under the map into the quotient are again locally constant (theorem~\ref{groupactions}). We also compare the complex analytic topology on a geometric quotient with the quotient topology (theorem~\ref{slices_geom_quot}). In section~\ref{sec_equi} we review equivariant cohomology. In proposition~\ref{rational_equi_quotient} we give a sufficient condition for the equivariant cohomology to be isomorphic to the cohomology of the quotient space, and in proposition~\ref{prop_equi} we show that in the algebraic case equivariant cohomology has a functorial mixed Hodge structure.

{\bf Conventions.} All algebraic varieties and groups will be over $\C$, and all algebraic groups will moreover be affine. All group actions are assumed algebraic unless otherwise stated. In the sequel $G$ and $X$ will be an algebraic group and an algebraic $G$-variety respectively. We use the definitions of the categorical and geometric quotients given in~\cite[Definition 2.1.4]{kirwan_et_al}; in particular, we require the map from $X$ to the geometric quotient by $G$ to be affine. If a $G$-invariant regular map $X\to Y$ has all properties of the geometric quotient given in~\cite[Definition 2.1.4]{kirwan_et_al} except possibly being affine, we call it the {\it weak geometric quotient}.

We consider left group actions unless specified otherwise, but we put the group on the right to denote a quotient. So for example, $A/H$ is the quotient of a space $A$ by a continuous action of a topological group $H$ regardless whether the action is on the left or on the right. Similarly, if there exists the geometric, respectively categorical quotient of $X$ by $G$, we will denote it $X/G$, and respectively $X//G$. 

If a group $H$ acts on $A$ on the right and on $B$ on the left, we introduce a right action of $H$ on $A\times B$ by setting $(a,b)h=(ah, h^{-1}b)$ for $a\in A,b\in B,h\in H$, and we use $A\times_H B$ to denote the resulting quotient (topological, geometric or categorical). 
For an $a\in A$ we denote the orbit of $a$ by $H(a)$ and the stabiliser of $a$ by $H_a$.

We denote the constant sheaf on a space $X$ with stalk an abelian group $M$ by $\underline{M}_X$, but we will often write $H^*(X,M)$ instead of $H^*(X,\underline{M}_X)$.

\subsection{Geometric quotients}\label{sec_geo_quot}

The following result (see~\cite[Theorem 1.10]{mum}) is one of the main tools for constructing geometric quotients.

\begin{theorem}\label{main_thm_git}
Suppose $G$ is reductive, and $L$ is a $G$-equivariant line bundle over $X$. Let $X^s\subset X$ be the stable locus with respect to $L$. Then there exists a geometric quotient $X^s\to X^s/G$.
\end{theorem}

$\clubsuit$

The following proposition allows one to prove the existence of a geometric quotient in most examples we consider in this paper.
\begin{Prop}\label{examples_quotients}
Suppose $V$ is a finite-dimensional complex vector space.
\begin{itemize}
\item Let $U$ be an open subset of the space of $k$-frames in $V$ which is invariant under the action of $GL_k(\C)$, and let $\bar U$ be the image of $X$ in $\Gr(k,V)$. Then the natural map $U\to \bar U$ is a geometric quotient. 
\item Suppose we are given a representation $\GL_n(\C)\to\GL(V)$. Let $U\subset\p(V)^s\subset \p(V)$ be an open $SL_n(\C)$-invariant subset of the set of points which are stable with respect to the bundle $\OO_{\p(V)}(1)$ on $\p(V)$, and let $U'$ be the preimage of $U$ in $V\setminus\{0\}$. Then there exist geometric quotients $U'/\C^*, U'/\GL_n(\C), U'/\SL_n(\C)$ and $U/\PGL_n(\C)$.
\end{itemize}
\end{Prop}

{\bf Proof.} First, recall that if $p:X\to Y$ is a geometric quotient, then $p$ is open (\cite[Chapter 0, \S 2, Remark 4]{mum}) and if $W\subset Y$ is open, then the map $p^{-1}(W)\to W$ induced by $p$ is again a geometric quotient (\cite[Proposition 3.10]{newstead}). So it suffices to prove part 1 for $U$ equal the whole space of frames and part 2 for $U=\p(V)^s$. Then part~1 becomes straightforward.

Let us prove part 2. For $U/\SL_n(\C)$ we use theorem~\ref{main_thm_git}, and the case of $U/\PGL_n(\C)$ follows. The quotient $U'/\SL_n(\C)$ can be handled in a similar way using the fact that every element of $U'$ is stable with respect to the natural action of $\SL_n(\C)$ on the trivial line bundle on $V\setminus\{0\}$. Finally, the case of $U'/\GL_n(\C)$ follows from the existence of the geometric quotients $U'/\C^*$ and $U/\PGL_n(\C)=(U'/\C^*)/(\GL_n(\C)/\C^*)$.$\clubsuit$
%

The remaining examples we consider in this paper are principal and associated bundles. In section~\ref{sec_geo_quot} {\it principal $G$-bundles} will be understood to be locally trivial with respect to the \'etale topology, see e.g.~\cite[Definition 2.1.4]{kirwan_et_al}. A typical example is the map from the first part of proposition~\ref{examples_quotients} (in which case $G=\GL_k(\C)$). We note that if $X\to Y$ is a principal $G$-bundle, then $Y$ is automatically the geometric quotient of $X$ by $G$.

A variety $Z$ with an action of an algebraic group {\it has property} $(*)$ if every finite subset of $Z$ is contained in an affine open subvariety.

\begin{Prop}\label{assoc_bund}
Suppose $P\to Y$ is a principal $G$-bundle and let $Z$ be an arbitrary $G$-variety that has property~$(*)$. Then there exists a geometric quotient $P\times Z\to P\times_G Z$, which is a principal $G$-bundle. The resulting regular map $P\times_G Z\to P/G=Y$ is an \'etale locally trivial fibre bundle with fibre $Z$.
\end{Prop}

We will call $P\times_G Z\to Y$ the {\it bundle with fibre $Z$ associated to the principal bundle $P\to Y$}.

{\bf Remark.} One can get rid of the assumption that $Z$ has property $(*)$ (and of the same assumption on $X$ and $Y$ in proposition~\ref{prop_equi}) by considering algebraic spaces rather than algebraic varieties. Nevertheless, we note that e.g.\ all quasi-projective varieties have this property, which suffices for our examples, so we omit the details.

{\bf Proof.} See~\cite{serre_efa}, Proposition 4 and example c.\ in section 3.2.$\clubsuit$

\begin{Prop}\label{principal_bundle_restr}
Suppose $P\to Y$ is a principal $G$-bundle and let $H\subset G$ be an algebraic subgroup. Then there exists a geometric quotient $P\to P/H$, which is a principal $H$-bundle. The resulting regular map $P/H\to P/G=Y$ is an \'etale locally trivial fibre bundle with fibre $G/H$.
\end{Prop}

{\bf Proof.} We apply~\cite{serre_efa}, Proposition 8 and again example c., section 3.2; note that ``locally trivial'' in that proposition means Zariski locally trivial.$\clubsuit$
\subsection{Slices}\label{sec_slices}

We will now review slices. Let $Y$ be topological space on which a topological group $H$ acts continuously. A subspace $S\subset Y$ is called a {\it topological slice} at $y\in S$ for the action of $H$ (see~\cite[Definition 2.1.1]{palais}) if and only if $H(S)$ is open in $Y$, and there is an $H$-equivariant continuous map $H(S)\to H/H_y$ such that the preimage of $H_y$ under this map is $S$. Here $H(S)$ is the union of the orbits of all points of $S$. 

\begin{lemma}\label{slice_lemma}
\begin{enumerate}
\item Suppose $y\in Y$ and the projection $q:H\to H/H_y$ has a continuous section over a neighbourhood of every point of $H/H_y$. (This will be the case e.g.\ if $H$ is a Lie group.) Let $S$ be a slice through $y$. Then the map $(h,s)\mapsto hs$ from $H\times S$ to $H(S)$ induces an $H$-equivariant homeomorphism $g:H \times_{H_y} S\to H(S)$. 
\item If $S$ is an $H_y$-invariant subspace through $y\in Y$ such that $H(S)$ is open in $Y$ and the map $g:H \times_{H_y} S\to H(S)$ defined as in the previous part is a homeomorphism, then $S$ is a slice at $y$. Moreover, under the same assumptions $S$ is locally closed and the image of $y$ in $Y/H$ has an open neighbourhood homeomorphic to $S/H_y$.
\end{enumerate}
\end{lemma}

{\bf Proof.} We will only prove part 1, part 2 being straightforward. The map $g$ is clearly well defined, surjective and continuous. Let $f:H(S)\to H/H_y$ be an $H$-equivariant continuous map such that $f^{-1}(H_y)=S$. Let us show that $g$ is injective. Suppose $h_1,h_2\in H$ and $s_1,s_2\in S$ are such that $h_1 s_1=h_2 s_2$. Multiplying both sides by $h_2^{-1}$ and applying $f$, we get $h_2^{-1} h_1\in H_y$, so $(h_1,s_1)=(h_2,s_2)\cdot h_2^{-1} h_1$, which means that $(h_1,s_1)$ and $(h_2,s_2)$ are in the same $H_y$-orbit.

Let us now show that $g$ is open. Let $W\subset H/H_y$ be an open set over which there is a continuous section $\varphi:W\to H$ of $q$. Set $W'\subset H\times_{H_y} S$ to be the image of $q^{-1}(W)\times S\subset H\times S$. This is an open subset of $H\times_{H_y} S$. We claim that $g(W')\subset f^{-1}(W)$. Indeed, suppose $(h,s)\in q^{-1}(W)\times S$. We have $f(hs)=h f(s)=q(h)$ where the element in the middle is the result of applying $h\in H$ to $f(s)=H_y\in H/H_y$. So $f(g([h,s]))=q(h)\in W$ where we use the square brackets to denote the image of an element of $H\times S$ in $H\times_{H_y} S$. This proves the claim.

We now define a continuous map $g_1:f^{-1}(W)\to W'$ by setting $g_1(w)=[\varphi(f(w)), \varphi(f(w))^{-1} w]$. To see that $g_1$ is well defined we observe that if $w\in f^{-1}(W)$, we have $q\big(\varphi(f(w))\big)=f(w)\in W$ and $f(\varphi(f(w))^{-1} w)= \varphi(f(w))^{-1} f(w)=H/H_y$, so $\varphi(f(w))\in q^{-1}(W)$ and $\varphi(f(w))^{-1} w\in S$. Moreover, we have $g(g_1(w))=w$ for all $w\in f^{-1}(W)$, and $g_1$ is continuous. Combining this with the above we see that $g_1$ and the map $W'\to f^{-1}(W)$ induced by $g$ are mutually inverse homeomorphisms. This shows that $g|_{W'}$ is open, which implies that $g$ is open and hence a homeomorphism.$\clubsuit$

There are straightforward analogues of slices for smooth manifolds and complex analytic varieties. In these cases slices are required to be locally closed submanifolds and subvarieties respectively. The algebraic version is more subtle. Namely, suppose $X,G$ are as in section~\ref{sec_geo_quot} and there is a categorical quotient $X//G$. An affine subvariety $S\subset X$ is an {\it \'etale slice} at $x\in S$ (see~\cite{luna}) if
\begin{enumerate}
\item $S$ is $G_x$-invariant;
\item the categorical quotient $(G\times_{G_x} S)//G$ exists, and the natural morphism $\varphi: (G\times_{G_x} S)//G\to X//G$ is \'etale;
\item the natural morphism from $G\times_{G_x} S$ to the pullback $((G\times_{G_x} S)//G)\times_{X//G} X$ is an isomorphism.
\end{enumerate}
Here the geometric quotient $G\times_{G_x} S$ exists by~proposition~\ref{assoc_bund}. Note that the categorical quotient $(G\times_{G_x} S)//G$ exists if and only if $S//G_x$ does, and if both quotients exist, they are isomorphic.

D.~Luna proves in~\cite{luna} that if $X$ is affine and $G$ is reductive, then there is an \'etale slice through every $x\in X$ such that $G(x)$ is closed. We note that a neighbourhood of such an $x$ in $S$ is a complex analytic (and hence a topological) slice though $x$, see~\cite[Section 6.6]{vinpop}. 
Also note that if $G$ is reductive, $X, S$ are affine and $G(x)$ is closed, then the categorical quotients $(G\times_{G_x} S)//G$ and $S//G_x$ automatically exist as $G_x$ in this case is reductive by the Matsushima criterion (see~\cite{richardson} and references therein).

\begin{theorem}\label{slices_geom_quot}
If a weak geometric quotient $p:X\to X/G$ exists, then the following is true.
\begin{enumerate}
\item Suppose $G$ is reductive and $p:X\to X/G$ is a genuine geometric quotient. Then for every $x\in X$ there is an affine neighbourhood $U$ of $x$ and an \'etale slice through $x$ in $U$. In particular, there is a complex analytic (and hence a topological) slice through every $x\in X$.
\item 
The complex analytic topology on $X/G$ coincides with the quotient of the complex analytic topology on $X$ by the equivalence relation induced by the action of $G$.
\end{enumerate}
\end{theorem}

{\bf Proof.} Let us prove the first part. By~\cite[Proposition 3.10]{newstead} and the definition of a good quotient, every point of $X/G$ has an affine neighbourhood $\bar U$ such that $U=p^{-1}(\bar U)$ is affine and $U\to \bar U$ is the geometric quotient of the action of $G$ on $U$. Since all $G$-orbits are closed in $X$, the first part follows from Luna's theorem and the remarks above about analytic slices. 
To prove part 2 we note that by~\cite[Chapter 0, \S2, Remark 4]{mum}, $p$ is universally open, which implies that $p$ is open in the complex analytic topology, see~\cite{moret_bailly} and references therein.$\clubsuit$

\begin{theorem}
\label{groupactions}
Suppose that $H$ is a connected Lie group that acts continuously with finite stabilisers on a topological space $Y$, and let $p:Y\to Y/H$ be the projection. Suppose $\eu L$ is a local system on $Y$, and that for every $y\in Y$ the following holds:
\begin{enumerate}
\item The stabiliser $H_y$ of $y$ is finite.
\item There is a topological slice $S$ though $y$ such that $S$ is $H_y$-equivariantly homotopy equivalent to $\{y\}$.
\item The pullback of $\eu L$ under every orbit map $O:H\to Y$ is constant.
\end{enumerate}
Then each sheaf $R^i p_*{\eu L}$ on $Y/H$ is isomorphic to $p_*\eu L\otimes H^i(H,\Q)$. Moreover, if $\eu L=\underline{\Q}_Y$, then $p_*\eu L\cong\underline{\Q}_{Y/H}$.
\end{theorem}

{\bf Remark.} Suppose $Y$ is a smooth manifold and the action $H\times Y\to Y$ is smooth. Then one can construct a slice $S$ at $y\in Y$ that $H_y$-equivariantly deformation retracts onto $y$ as follows. We equip $Y$ with an $H_y$-invariant Riemannian metric and take $S$ to be the exponential of the $\varepsilon$-neighbourhood of $0$ in $(T_y H(y))^\perp\subset T_y Y$ for $\varepsilon$ sufficiently small. If moreover we already have a slice $S'$ at $y$ that is smooth near $y$, then by using the same idea but exponentiating in $S'$ rather than in $Y$ we obtain a smaller smooth slice $S''\subset S'$ at $y$ that $H_y$-equivariantly deformation retracts onto~$y$.

{\bf Proof.} The proof is very similar to that of Theorem 3.8 in \cite{quintics}. If $\eu{L}=\underline{\Q}_Y$, one can use a set-theoretical (i.e.\ not necessarily continuous) section $Y/H\to Y$ and the orbit maps $O:H\to Y$ to construct, for each $y\in Y$, an isomorphism between $H^*(H,\Q)$ and the stalk of $R^i p_*{\eu L}$ over $y$. Then one uses slices to show that these isomorphisms trivialise $R^i p_*{\eu L}$.

To prove the theorem in the general case we use the cup product map $R^0 p_*\eu{L}\otimes R^i p_*\underline{\Q}_Y\to R^i p_*\eu{L}$. Using slices again we see that on the stalks over the orbit of $y\in Y$ this map is the cup product
\begin{equation}\label{iso_groupactions}
H^0(H(y),\eu{L})\otimes H^i(H(y),\Q))\to H^i(H(y),\eu{L}).
\end{equation}
By our assumptions, $\eu{L}$ trivialises over $H$, so the action of $\pi_1(H(y))$ on the stalk of $\eu{L}$ factorises through $\pi_1(H(y))\to H_y$. In particular, this action is completely reducible.

The covering $q:H\to H(y)$ induces an isomorphism in rational cohomology, and $q_*(\underline{\Q}_H)$ is the local system obtained by composing $\pi_1(H(y))\to H_y$ with the action of $H_y$ on its group algebra $\Q[H_y]$. (Recall that the latter contains all irreducible $H_y$-modules.) So if $\eu{L}'$ is an irreducible local system on $H(y)$ that is obtained from a representation of $H_y$, we have $H^*(H(y),\eu{L}')\cong H^*(H,\Q)$ if $\eu{L}'$ is constant and $H^*(H(y),\eu{L}')=0$ otherwise. We conclude that~(\ref{iso_groupactions}) is an isomorphism.$\clubsuit$

\begin{corollary}\label{cor_groupactions}
Suppose $X$ is a complex algebraic variety and $G$ is a complex algebraic group that acts algebraically on $X$. If $X$ is smooth, $G$ is connected, every stabiliser $G_x,x\in X$ is finite and there exists a geometric quotient $p:X\to X/G$, then the second page of the Leray spectral sequence of $p$ (that converges to $H^*(X,\Q)$) is $E_2^{p,q}\cong H^p(X/G,\Q)\otimes H^q(G,\Q)$. \end{corollary}

{\bf Proof.} By~\cite[Theorem 6.5]{vinpop}, the analytic slice at $x\in S$ one obtains from the first part of theorem~\ref{slices_geom_quot} may be assumed smooth. Using the remark above we see that we may assume furthermore that this slice $G_x$-equivariantly deformation retracts onto $x$. We now apply theorem~\ref{groupactions} to $\eu{L}=\underline{\Q}_X$.$\clubsuit$

{\bf Remark.} It seems likely that corollary~\ref{cor_groupactions} can be generalised to the case when $X$ is not necessarily smooth.

\subsection{Equivariant cohomology}\label{sec_equi}

Recall that if $H$ is a Lie group, then the {\it universal $H$-bundle} $EH$ is a contractible space with a free continuous right action of $H$ such that $EH\to EH/H$ is a locally trivial fibre bundle. The quotient $BH=EH/H$ 
is called the {\it classifying space} of $H$. There are several ways to choose an $EH$, but they all result in the same $BH$ up to a preferred homotopy equivalence. More precisely, suppose $E'H$ and $E''H$ are two choices. Then the projections $E'H\times_H E''H\to E'H/H$ and $E'H\times_H E''H\to E''H/H$ are homotopy equivalences, so they induce a homotopy equivalence $e^{E''}_{E'}: E'H/H\to E''H/H$. Moreover, $e^{E'}_{E'}$ is homotopic to the identity, and if $E'''H$ is yet another choice, then $e^{E'''}_{E'}$ and $e^{E'''}_{E''}\circ e^{E''}_{E'}$ are homotopic.

If $Y$ is a topological space with a continuous action of $H$, then the {\it equivariant cohomology} $H_H^*(Y, R)$ of $Y$ with coefficients in a commutative ring $R$ is the cohomology $H^*(EH\times_H Y, R)$ of the {\it homotopy quotient} (also called the {\it equivariant quotient} or {\it Borel quotient}) $EH\times_H Y$, which we denote $Y_{hH}$. Same as $BH$, the equivariant quotient is uniquely defined up to a preferred homotopy equivalence. The map $EH \times_H Y\to BH$ is the projection of a locally trivial bundle with fibre $Y$; moreover, if $H$ is path connected, then $BH$ is simply connected, and we obtain a spectral sequence
\begin{equation}\label{sp_seq_equi}
(E^{p,q}_r,d_r)\Rightarrow H^{p+q}(EH \times_H Y,R),E^{p,q}_2=H^p(BH,H^q(Y,R)).
\end{equation}

Let us now compare the equivariant cohomology with the cohomology of the quotient space $Y/H$. Let $f:EH\times_H Y\to Y/H$ be the natural map. 

\begin{Prop}\label{rational_equi_quotient}
Suppose for every $y\in Y$ the stabiliser $H_y$ is finite and every neighbourhood $U$ of $y$ in $Y$ contains a slice $S$ through $y$ that $H_y$-equivariantly deformation retracts onto $y$. (Note that this will be the case when $Y$ is a smooth manifold and $H$ is a Lie group that acts smoothly on $Y$.) Then $f$ induces an isomorphism in rational cohomology.
\end{Prop}

{\bf Proof.} By lemma~\ref{slice_lemma}, $S/H_y$ is an open neighbourhood of $H(y)$ in $Y/H$, and we have
$$f^{-1}(S/H_y)=EH\times_H H(S)\cong EH\times_H(H\times_{H_y} S)\cong (EH \times_H H)\times_{H_y} S\cong EH\times_{H_y} S\sim B H_y,$$ which has the rational cohomology of a point. This implies that $R^i f_*\underline{\Q}_{EH\times_H Y}$ is $\underline{\Q}_{Y/H}$ if $i=0$ and $0$ otherwise.$\clubsuit$

Let now $X,G$ be as in section~\ref{sec_geo_quot}. We wish to put a mixed Hodge structure on (\ref{sp_seq_equi}) for $Y=X$ and $H=G$. There are several equivalent ways to do this. In~\cite{th_hodge_3} P.~Deligne uses the bar construction to introduce a mixed Hodge structure on $H^*(BG,\Q)$, and this approach generalises to equivariant cohomology. Here we describe another construction which uses principal bundles and finite-dimensional approximations. It will allow us to quickly prove everything we need (proposition~\ref{prop_equi}). From now on we suppose that $X$ has property $(*)$ (e.g. that\ $X$ is quasi-projective). 
We call a sequence $P_0\to P_1\to\cdots $ of \'etale locally trivial principal $G$-bundles and closed $G$-embeddings {\it approximating} if for all $j$, the group $\pi_j(P_i)$ is trivial for all sufficiently large $i$.

\begin{lemma}
Approximating sequences exist.
\end{lemma}

{\bf Proof.} For $G=\GL_k(\C)$ we can take $P_i,i\geq k$ to be the Stiefel manifold of all $k$-frames in $\C^i$ (see proposition~\ref{examples_quotients}) and $P_{<k}=G$. The same recipe works for $G\subset \GL_k(\C)$ by proposition~\ref{principal_bundle_restr}.$\clubsuit$

Let us take any approximating sequence $P_0\to P_1\to\cdots $. In the rest of the appendix we consider all varieties in the complex analytic topology. The space $P=\varinjlim P_i$ is a model for $EG$, so $P/G$ is a model for $BG$; moreover, we have $P/G\cong\varinjlim (P_i/G)$. 
%
We get an isomorphism
\begin{equation}\label{equi_hodge}
H^*_G(X,\Q)=H^*(EG \times_G X,\Q)\cong \varprojlim H^*(P_i\times_G X,\Q).
\end{equation}
The geometric quotients $P_i\times_G X$ exist by proposition~\ref{assoc_bund}. The group $H^*(P_i\times_G X,\Q)$ carries a natural mixed Hodge structure and eventually stabilises in each degree as $i\to\infty$, as do the spectral sequences of the projections $P_i\times_G X\to P_i/G$, which stabilise to (\ref{sp_seq_equi}) for $H=G,Y=X$ if $G$ is path connected. For every $i$ the maps $X\to P_i\times_G X$ and $P_i\times_G X\to P_i/G$ induces maps of mixed Hodge structures in rational cohomology, and so does the map $f_i:P_i\times_G X\to X//G$ if there exists at least a categorical quotient $X//G$.

Furthermore, if $P_0\to P_1\to\cdots $ and $Q_0\to Q_1\to\cdots $ are approximating sequences, by proposition~\ref{assoc_bund} so is $P_0\times Q_0\to P_1\times Q_1\to\cdots $, and $\varinjlim(P_i\times Q_i)$ is naturally homeomorphic to $\varinjlim P_i\times\varinjlim Q_i$ by~\cite[18.5]{whitehead}. This allows one to identify isomorphisms (\ref{equi_hodge}) obtained by taking different approximating sequences in the same way as we showed above that the classifying space is unique up to a preferred homotopy equivalence. Moreover, the identification is compatible with the mixed Hodge structures on the right hand sides. 

Finally, suppose $g:H\to G$ is a homomorphism of algebraic groups, and $P_0\to P_1\to\cdots $ is an approximating sequence for $G$. If $Q_0\to Q_1\to\cdots $ is an approximating sequence for $H$, then so is $P_0\times Q_0\to P_1\times Q_1\to\cdots $ with the diagonal action of $H$, again by proposition~\ref{assoc_bund}, and the projections $P_i\times Q_i\to P_i$ are $g$-equivariant. Moreover, if $Y$ is an $H$-variety that has property $(*)$, and $h:Y\to X$ is a $g$-equivariant regular map, then the following diagram commutes for every $i$:
$$
\begin{tikzcd}
(P_i\times Q_i)\times_H Y\arrow[r]\arrow[d] & P_i\times_G X\arrow[d]\\
P_i\times_H Q_i\arrow[r] & P_i/G.
\end{tikzcd}
$$

Summarising, we get the following proposition, where we write MHS for ``mixed Hodge structure''.

\begin{Prop}\label{prop_equi} Suppose $G$ is an algebraic group and $X$ is a $G$-variety with property $(*)$. Then the following is true.
\begin{enumerate}
\item There is an MHS on $H^*_G(X,\Q)=H^*(EG \times_G X,\Q)$ given by the isomorphism~(\ref{equi_hodge}), which is independent of the choice of an approximating sequence and the choice of a model for $EG$.
\item There is an MHS on the spectral sequence of $EG \times_G X\to BG$, and the spectral sequence converges as an MHS to $H^*_G(X,\Q)$ (equipped with the MHS from the previous part).
\item The maps $X\to EG\times_G X$ and $EG \times_G X\to BG$ induce maps of MHSs in rational cohomology, and so does the map $EG\times_G X\to X//G$ if the target is well defined (i.e.\ if the categorical quotient exists).
\item More generally, the MHSs from the first two parts are functorial in the following sense. Suppose $g:H\to G$ is homomorphism of algebraic groups, $Y$ is an $H$-variety with property $(*)$, and $h:Y\to X$ is a $g$-equivariant regular map. Then the resulting map $H^*_G(X,\Q)\to H^*_H(Y,\Q)$, and also the map of the spectral sequences induced by
$$
\begin{tikzcd}
EH\times_H Y\arrow[r]\arrow[d] & EG\times_G X\arrow[d]\\
BH\arrow[r] & BG
\end{tikzcd}
$$
are maps of MHSs. 
\end{enumerate}
\end{Prop}

$\clubsuit$

{\bf Remark.} As an example, by applying part 4 to $Y=X=pt$, we see that the map $H^*(BH,\Q)\to H^*(BG,\Q)$ induced by $g$ is a map of mixed Hodge structures. In particular, we can take $g$ to be the inclusion $T\to G$ of a maximal torus. 
As an application we obtain that $H^*(BG,\Q)$ is zero or Tate of weight equal degree in even degrees and zero in odd degrees, see~\cite[Theorem 9.1.1]{th_hodge_3}. Note however that $H^2(BG,\Q)=0$ if $G$ has no non-trivial characters, so $BG$ can not in general be ``approximated'' by smooth compact varieties.

\section{Spectral sequences and mixed Hodge structures (by Alexey Gorinov and Nikolay Konovalov)}\label{app_mhs_spectral}

Here we review and/or prove some facts about Borel-Moore homology, mixed Hodge structures and spectral sequences. Unless stated otherwise, all spaces in this appendix are assumed Hausdorff, locally compact and paracompact, the ground ring $R$ is a principal ideal domain, and all sheaves are sheaves of $R$-modules. As in appendix~\ref{b}, we will use $\underline{M}_X$ to denote the constant sheaf on a space $X$ with stalk an $R$-module $M$, and we will often abbreviate $H^*(X,\underline{M}_X)$ to $H^*(X,M)$, and similarly for cohomology with compact support. We use $C_*(-,R)$ and $C^*(-,R)$ to denote the singular (co)chain functors.

In section~\ref{sec_bm_hom} we recall the definition of the Borel-Moore homology for sheaves. For constant coefficients the result is the same as the homology of the complex of locally finite chains (proposition~\ref{homology_iso_sheaf_singular}). In section~\ref{borel_moore} we review mixed Hodge structures on the cohomology and Borel-Moore homology of semi-simplicial (in particular cubical) complex algebraic varieties (propositions~\ref{mhs_coh} and~\ref{mhs_compact_support_bm}).

In section~\ref{thom_iso_ccavb} we prove the Thom isomorphism theorem for cubical complex algebraic vector (and more generally, affine) bundles (theorem~\ref{thom_ccavb}). The proof is surprisingly indirect as we first need to establish the Hodge theoretic versions of the K\"unneth isomorphism and the relative cup product map for general simplicial, respectively semi-simplicial complex algebraic varieties (theorems~\ref{simplicial_kunneth_mhs} and~\ref{cup_product_mhs}). Finally, in section~\ref{sec_mhs_sseq} we equip two spectral sequences with mixed Hodge structures. In proposition~\ref {mhs_sseq_filtred} we consider the filtration of a proper semi-simplicial variety by closed subvarieties, and in proposition~\ref{mhs_sseq_leray} the Leray spectral sequence of the augmentation map $|X_\bullet|\to X_\varnothing$ for a proper cubical complex algebraic variety $X_\bullet$.

\subsection{Borel-Moore homology}\label{sec_bm_hom}

For a topological space $X$ we let $J_X$ be a functorial injective resolution of the constant sheaf $\underline{R}_X$. Here ``functorial'' means that for each continuous map $f:X\to Y$ we have a map $J_Y\to f_*(J_X)$, and if $g:Y\to Z$ is another continuous map, then the diagram
$$
\begin{tikzcd}
g_* f_* J_X & g_* J_Y\ar[l]\\
& J_Z\ar[u]\ar[ul],
\end{tikzcd}
$$
commutes, where the horizontal arrow is $g_*(J_Y\to f_* J_X)$. Let us also choose an injective resolution $R\to I^*R$ of $R$ as a module over itself. By our assumption on $R$, we may and will assume that the length of $I^*R$ is $\leq 1$. We define the {\it dualising complex} $D_X$ of a space $X$ as the complex of sheaves that takes an open $U\subset X$ to
$$
\Hom(\Gamma_c(U,J_X^*),I^*R),
$$
see~\cite[V.3]{bredon}; recall that if $C_1^*, C_2^*$ are two cochain complexes, then $\Hom^*(C_1^*,C_2^*)$ is another cochain complex with $$\Hom^n(C_1^*,C_2^*)=\prod_{-p+q=n} \Hom^*(C_1^p,C_2^p);$$ so $D_X$ is a complex of sheaves concentrated in degree $\leq 1$.

For a sheaf $S$ on $X$, we define the {\it Borel-Moore homology of $X$ with coefficients in $S$} to be the groups $$\bar H_*(X,S)=H^{-*}(\Gamma(X,D_X\otimes S)),$$ see~ibid.
Similarly, the {\it homology of $X$ with coefficients in $S$} is set to be $$H_*(X,S)=H^{-*}(\Gamma_c(X,D_X\otimes S)).$$

These groups are functorial, see~ibid., V.4: 
\begin{Prop}
If $f:X\to Y$ is proper and $S$ is a sheaf on $Y$, we have an induced map $f_*:\bar H_*(X,f^{-1}(S))\to \bar H_*(Y,S)$ of Borel-Moore homology groups. If $g:Y\to Z$ is another proper map, then $(g\circ f)_*=g_*\circ f_*$ for every sheaf on $Z$. The analogues of these statements for homology hold without assuming the maps proper.
\end{Prop}

$\clubsuit$

Let $i:A\to X$ be the inclusion of a closed subset and let $S$ be a sheaf on $X$. Then the induced maps $\Gamma(A,D_A\otimes i^{-1}(S))\to \Gamma(X,D_X\otimes S))$ and $\Gamma_c(A,D_A\otimes i^{-1}(S))\to \Gamma_c(X,D_X\otimes S))$ are injective~\cite[V.5]{bredon}. One then defines relative chains by taking the quotient complexes, and one obtains {\it relative Borel-Moore homology groups} $\bar H_*(X,A,S)$ and {\it relative homology groups} $H_*(X,A,S)$ by taking homology. We have
$$\bar H_*(X,A,S)\cong H_*(X,A,S)$$ if $X$ is compact. Moreover, if $S$ is a locally constant sheaf on $X$ of finite rank, we have
\begin{equation}\label{bm_open_subset}
\bar H_*(X,A,S)\cong \bar H_*(X\setminus A,S|_{X\setminus A}),
\end{equation}
see~ibid., Corollary~V.5.10, so in this case we have a long exact sequence
$$\cdots \to \bar H_*(A,S|A)\to \bar H_*(X,S)\to \bar H_*(X\setminus A,S|_{X\setminus A})\to \bar H_{*-1}(A,S|_A)\to\cdots.$$
Short exact sequences of sheaves result in long exact sequences of homology groups by ibid., V.3. The same applies to complexes of sheaves. 
We have an exact sequence
\begin{equation}\label{uct_bm_hc}
0\to\mathop{\mathrm{Ext}}(H^{n+1}_c(X,R),R)\to \bar H_n(X,R)\to \Hom(H^n_c(X,R),R)\to 0
\end{equation}
for every space $X$ and $n\in\Z$, see ibid., sequence (9) on p.~292. Moreover, if $X$ is a topological manifold of dimension $m$ without boundary and $S$ is a sheaf on $X$, we have the Poincar\'e isomorphisms 
$$H^*(X,O_X\otimes S)\cong\bar H_{m-*}(X,S)\mbox{ and }  H^*_c(X,O_X\otimes S)\cong H_{m-*}(X,S),$$ where $O_X$ is the orientation sheaf of $X$, see~ibid., Theorem V.9.2.

\medskip

Another approach to homology is via singular chains. For a space $X$ and an $R$-module $M$ we set $\bar C_*(X,M)$ to be the complex of locally finite singular chains of $X$ with coefficients in $M$. Since we assume $X$ locally compact, we have
$$\bar C_*(X,M)\cong\varprojlim C_*(X,X\setminus K,M)$$ where the limit is taken with respect to compact $K\subset X$. These complexes are functorial with respect to proper maps.

If $A\subset X$ is a closed subspace, then $\bar C_*(A,M)\subset \bar C_*(X,M)$ is a subcomplex. We set $\bar C_*(X,A,M)=C_*(X,M)/C_*(A,M)$. The homology groups of the latter complex are the {\it singular Borel-Moore homology groups of $X$ modulo $A$ with coefficients in $M$,} denoted $\bar H_*(X,A,M)$.

\begin{Prop}\label{homology_iso_sheaf_singular}
Let $A\subset X$ be a closed subspace, and assume $X$ and $A$ are finite dimensional and locally contractible. Then we have natural isomorphisms
$$\bar H_*(X,A,\underline{M}_X)\cong \bar H_*(X,A,M)\mbox{ and }H_*(X,A,\underline{M}_X)\cong H_*(X,A,M)$$
for any $R$-module $M$.
\end{Prop}

{\bf Proof.} The statement follows from~\cite[Corollary V.13.6 and the remarks at the end of V.1.19]{bredon}. For completeness we will now explain in more detail how these remarks extend from the case $M=R$ to arbitrary $M$. For simplicity we assume that $A=\varnothing$.

Note that all components of the source and target complexes of the map~(2) and hence of the map~(3) in ibid., V.1.19 are flat. 
Let us replace the homological grading of the source and target of the latter map by the cohomological one by setting $\triangle^*(X,R)=\triangle_{-*}(X,R), \mathcal{S}^*(X,R)=\mathcal{S}_{-*}(X,R)$. Tensoring the resulting map $\triangle^*(X,R)\to\mathcal{S}^*(X,R)$ by identity of the constant sheaf $\underline{M}_X$ we get a quasi-isomorphism, as can be seen from the diagram
$$
\begin{tikzcd}
\triangle^*(X,R)\otimes\underline{M}_X \ar[d]& \triangle^*(X,R)\otimes\underline{F}_X^*\ar[d]\ar[l]\\
\mathcal{S}^*(X,R)\otimes\underline{M}_X & \mathcal{S}^*(X,R)\otimes\underline{F}_X^*\ar[l]
\end{tikzcd}
$$
where $F^*\to M$ is a flat resolution of $M$ (it lives in non-positive degrees). In the diagram, the horizontal arrows and the arrow on the right are quasi-isomorphisms because bounded above complexes of flat sheaves are K-flat by~\cite[tag 06YD]{stacks}, so the arrow on the left is also a quasi-isomorphism.

Each component of $\mathcal{S}^*(X,R)$ is c-soft by construction (\cite[V.1.18]{bredon}), hence soft by~ibid., II.16.5 and our assumptions on $X$. Since all these components are also flat, the product $\mathcal{S}^*(X,R)\otimes\underline{M}_X$ is degree-wise c-soft (hence also degree-wise soft) by ibid., Corollary II.16.31. All components of $\triangle^*(X,R)$, hence of $\triangle^*(X,R)\otimes\underline{M}_X$ are homotopically fine by ibid., Exercise II.32. Applying ibid., Theorem IV.2.2 we see that $$\triangle^*(X,R)\otimes\underline{M}_X\to\mathcal{S}^*(X,R)\otimes\underline{M}_X$$ induces an isomorphism of $\Gamma(-)$ and $\Gamma_c(-)$. The global sections of $\triangle^*(X,R)\otimes\underline{M}_X$ can be identified with $\bar C_{-*}(X,M)$ by ibid., solution of Exercise I.12. Under this isomorphism compactly supported sections are identified with $C_{-*}(X,M)$.$\clubsuit$

{\bf Remark.} One can show that the complexes $\Gamma(X,D_X\otimes S)$ and $\Gamma_c(X,D_X\otimes S)$ form a chain theory for sheaf theoretic homology. Namely, if $f:X\to Y$ is a proper map and $S$ is a sheaf on $Y$, then we have a chain map $f_\#: \Gamma(X,D_X\otimes f^{-1}(S))\to \Gamma(Y,D_Y\otimes S)$ which induces $f_*:\bar H_*(X,f^{-1}(S))\to \bar H_*(Y,S)$. These maps are functorial: if $g:Y\to Z$ is another proper map, then $(g\circ f)_\#=g_\#\circ f_\#$ for every sheaf on $Z$. For homology we do not need to assume the maps proper.

Moreover, one can show that under the assumptions of proposition~\ref{homology_iso_sheaf_singular} the functors $X\mapsto\Gamma(X,D_X\otimes \underline{M}_X)$ and $X\mapsto\Gamma_c(X,D_X\otimes \underline{M}_X)$ are quasi-isomorphic to $X\mapsto \bar C_*(X,M)$, respectively $X\mapsto C_*(X,M)$, and similarly for pairs of spaces. We will not need these results, and the details are somewhat technical, so they will be presented elsewhere.

\subsection{Borel-Moore homology of cubical varieties}\label{borel_moore}

We refer the reader e.g.\ to~\cite[5.1.1]{stepet} for general information about semi-simplicial objects. Recall that starting from an $n$-cubical space $X_\bullet=(X_I,\varphi_{IJ})$ one can construct an augmented semi-simplicial one, which we will denote by $S(X_\bullet)_\bullet$, by setting for $k\leq n-1$
\begin{equation}\label{semi_simp_from_cub}
S(X_\bullet)_k=\bigsqcup_{|I|=k+1} X_I
\end{equation}
and $S(X_\bullet)_{\geq n}=\varnothing$, and using $\varphi_{IJ}$ to define the structure maps and the augmentation: Let $\beta$ be a strictly increasing map $\{0,\ldots,l\}\to \{0,\ldots,k\}$. Writing a subset $I\subset\{1,\ldots, n\},|I|=k+1$ as $I=\{i_0<\cdots< i_k\}$, set $J=\{i_{\beta(0)},\ldots,i_{\beta(l)}\}$, and define the map $S(X_\bullet)_k\to S(X_\bullet)_l$ induced by $\beta$ to be $\varphi_{JI}$ on $X_I$. The augmentation towards $X_\varnothing$ is defined via the maps $\varphi_{\varnothing I}$. See~\cite[Remark 5.4.2 and Figure 5.3]{stepet}. This construction is clearly functorial. So from now on we will mainly focus on the semi-simplicial case.

We will call a(n augmented) semi-simplicial topological space {\it proper} if all its structure maps, including the augmentation in the augmented case, are proper, and similarly for (augmented) semi-simplicial complex algebraic varieties. Note that this is not the same as requiring the map from $X_\bullet$ to a point, viewed as a constant semi-simplicial space, to be proper.

A collection $A_\bullet=(A_i)_{i\in\Z_{\geq 0}}$ of subspaces of a semi-simplicial space $X_\bullet$ is a {\it semi-simplicial subspace} if the structure maps of $X_\bullet$ preserve $A_\bullet$, and similarly for semi-simplicial complex algebraic varieties. The notions of compactness, open subspaces, dense subspaces etc.\ are transferred to semi-simplicical topological spaces and varieties component-wise.

We will now define (co)homology for semi-simplicial topological spaces. For simplicity we consider only the case of constant coefficients. For a semi-simplicial space $X_\bullet$, an $R$-module $M$ and $n\in\Z_{\geq 0}$ we define the {\it alternating face differential}
$$\delta:C^*(X_{n-1},M)\to C^*(X_{n},M), \delta=\sum_{i=0}^n (-1)^i (d_i^{X_\bullet})^\#$$ where $d_i^{X_\bullet}:X_{n}\to X_{n-1}$ is induced by the unique face $\{0,\ldots, n-1\}\to\{0,\ldots, n\}$ that misses $i$. We have $\delta\circ\delta=0$, and $\delta$ commutes with the cochain differentials on the source and the target. By taking the total cochain complex of the resulting double complex and then taking the cohomology we obtain the {\it cohomology groups of $X_\bullet$ with coefficients in $M$}, denoted $H^*(X_\bullet,M)$. 

Using the same construction we define the {\it homology groups} $H_*(X_\bullet, M)$ and also the {\it relative homology} and {\it cohomology groups} $H^*(X_\bullet, A_\bullet, M)$ and $H_*(X_\bullet, A_\bullet, M)$ where $A_\bullet\subset X_\bullet$ is a semi-simplicial subspace. The long exact sequences for couples and triples of spaces have straightforward semi-simplicial analogues.

Suppose now $X_\bullet$ is a proper semi-simplicial space and $A_\bullet\subset X_\bullet$ is a closed semi-simplical subspace. Recall that we use $\bar C_*(-,M)$ to denote the complex of locally finite singular chains with coefficients in $M$. Similarly, let $$C^*_c(X,M)=\varinjlim_{K} C^*(X,X\setminus K,M)$$ be the complex of cochains with compact support and coefficients in $M$. Here the limit is taken with respect to compact $K\subset X$. Both $\bar C_*(-,M)$ and $C^*_c(-,M)$ are functorial with respect to proper maps, and we use these functors to the {\it (relative) Borel-Moore homology} and {\it compactly supported cohomology groups of $X_\bullet$ modulo $A_\bullet$ with coefficients in $M$}, which we denote $\bar H_*(X_\bullet,A_\bullet,M)$, respectively $H^*_c(X_\bullet, A_\bullet, M)$. We again have long exact sequences of couples and triples of semi-simplicial spaces.

{\bf Remark.} Let $X_\bullet$ be a semi-simplicial space and let $S^\bullet$ be a sheaf on $X_\bullet$, see~\cite[5.1.2]{stepet}. One can define the cohomology groups $H^*(X_\bullet, S^\bullet)$ using Godement resolutions. One can show that if all components of $X$ are locally contractible and $S^\bullet$ is constant with stalk $M$, then the resulting groups coincide with $H^*(X_\bullet, M)$ constructed above. The same remarks apply to compactly supported cohomology if $X_\bullet$ is proper.

In a similar way, if a semi-simplicial space $X_\bullet$ is proper and augmented towards a space $X$ (see e.g.~\cite[5.1.1]{stepet}), and each $S^n$ is the pullback of the same sheaf $S$ on $X$, then one can use the Godement resolutions of $D_{X_n}\otimes S^n$ to define the homology groups $H_*(X_\bullet, S^*)$ and, if $X_\bullet$ is proper, the Borel-Moore homology groups $\bar H_*(X_\bullet, S^*)$. These coincide with the groups $H_*(X_\bullet, M)$ and $\bar H_*(X_\bullet, M)$ defined above if $X$ and all components of $X_\bullet$ are finite-dimensional and locally contractible, and $S=\underline{M}_X$. We will not use these results in the sequel, so we omit the details.

\smallskip

The construction of the geometric realisation extends from cubical to arbitrary semi-simplicial spaces, see e.g.~\cite[5.1.1]{stepet}. Geometric realisation commutes with (co)homology: if $(X_\bullet, A_\bullet)$ is a closed pair of semi-simplicial spaces and $M$ is an $R$-module, we have 
\begin{equation}\label{geom_realis_hom_coh}
H_*(X_\bullet,A_\bullet, M)\cong H_*(|X_\bullet|,|A_\bullet|, M); H^*(X_\bullet, A_\bullet, M)\cong H^*(|X_\bullet|,|A_\bullet|, M).
\end{equation}
(This can be proved in the same way as~\cite[Theorem 3.1]{carlson}.) If in addition $X_\bullet$ is proper, we have
\begin{equation}\label{geom_realis_bm_hc}
\bar H_*(X_\bullet,A_\bullet, M)\cong \bar H_*(|X_\bullet|,|A_\bullet|, M); H^*_c(X_\bullet, A_\bullet, M)\cong H^*_c(|X_\bullet|,|A_\bullet|, M).
\end{equation}
If $X_\bullet$ is compact, we have $\bar H_*(X_\bullet, A_\bullet,M)=H_*(X_\bullet, A_\bullet,M)$ and $H^*_c(X_\bullet, A_\bullet,M)=H^*(X_\bullet, A_\bullet,M)$.

\smallskip

A {\it compactification} of a complex (semi-simplicial) algebraic variety $X$ is a compact complex (semi-simplicial) algebraic variety $X'$ that contains $X$ as an open dense subvariety.

\begin{lemma}\label{open_dense_cubical}
\begin{enumerate}
\item Let $f:A\to B$ be a continuous map such that $A, B$ are Hausdorff and $B$ is locally compact, and let $V\subset A, W\subset B$ be subspaces such that $V$ is dense in $A$, $W$ is open in $B$, and $f(V)\subset W$. If the map $g:V\to W$ induced by $f$ is proper, then $f(A\setminus V)\subset B\setminus W$.
\item Let $U_\bullet$ be an open dense semi-simplicial subspace of a locally compact Hausdorff semi-simplicial space $X_\bullet$. If $U_\bullet$ is proper, then the collection $(X_I\setminus U_I)$ is a (closed) semi-simplicial subspace, which we denote $X_\bullet\setminus U_\bullet$. If in addition $X_\bullet$ and $U_\bullet$ are augmented towards a space $X$ and an open subset $U\subset X$ respectively, then $X_\bullet\setminus U_\bullet$ is augmented towards $X\setminus U$.
\end{enumerate}
\end{lemma}

$\clubsuit$

Let $X_\bullet$ be a proper semi-simplicial space, and let $A_\bullet\subset X_\bullet$ be a closed subspace. Suppose $X'_\bullet$ is a compact semi-simplicial space that contains $X_\bullet$ as an open subspace, and suppose $\partial X'_\bullet=X'_\bullet\setminus X_\bullet$ is a semi-simplicial subspace of $X'_\bullet$. (By lemma~\ref{open_dense_cubical}, the latter is automatically the case if $X_\bullet$ is dense in $X'_\bullet$.) Let $A'_\bullet$ be the closure of $A_\bullet$ in $X'_\bullet$. Using~(\ref{geom_realis_bm_hc}), (\ref{uct_bm_hc}) and~(\ref{geom_realis_hom_coh}) we then have
\begin{equation}\label{homk}
\bar H_*(X_\bullet,A_\bullet,\mathbf{k})\cong \Hom_\mathbf{k} (H^{*}_c(X_\bullet,A_\bullet,\mathbf{k}),\mathbf{k})\cong \Hom_\mathbf{k}(H^*(X'_\bullet, A'_\bullet\cup \partial X'_\bullet,\mathbf{k}),\mathbf{k})
\end{equation}
%
where $\mathbf{k}$ is a field. We will now use formula (\ref{homk}) to put a mixed Hodge structure on Borel-Moore homology groups of a proper semi-simplicial algebraic variety.

\begin{lemma}\label{compactification}
\begin{enumerate}
\item Suppose $f:X\to Y$ is a proper morphism of complex algebraic varieties.
Let $X'$, respectively $Y'$ be a compactification of $X$, respectively $Y$. The closure $Z$ of the graph of $f$ in $X'\times Y'$ is a compactification of $X$, and $f$ extends to a morphism $g:Z\to Y'$ that takes the boundary $X'\setminus X$ to $Y'\setminus Y$.
\item For every proper (augmented) semi-simplicial complex algebraic variety $X_\bullet$ 
there exists a(n augmented) compactification $X'_\bullet$.
\end{enumerate}
\end{lemma}

{\bf Proof.} The first part is straightforward. Let us prove the second part. We choose a compactification $X_0'$ of $X_0$ arbitrarily. Then using part 1 we construct for every $n>0$ a compactification $X_n'$ of $X_n$ to which all face maps $X_n\to X_{n-1}$ extend. A similar argument applies in the augmented case.$\clubsuit$

Suppose now $X_\bullet$ is a semi-simplicial complex algebraic variety, and $A_\bullet\subset X_\bullet$ is a semi-simplicial subvariety. Using~\cite[IV 2]{gnpp} (cf.\ also~\cite[5.5.1]{stepet}) we equip the relative cohomology groups $H^*(X_\bullet, A_\bullet, \Q)$ with a mixed Hodge structure with the following properties. 
Note that in~\cite{gnpp} the authors use {\it strict simplicial objects} for what we call semi-simplicial objects, see ibid., p. 89.

\begin{Prop}\label{mhs_coh}
\begin{enumerate}
\item The mixed Hodge structure is functorial in $(X_\bullet, A_\bullet)$.
\item All maps in the cohomology long exact sequence of the couple $(X_\bullet, A_\bullet)$ are maps of mixed Hodge structures.
\end{enumerate}
\end{Prop}

{\bf Proof.} We use~\cite[IV 2.13 and IV 2.16]{gnpp}.$\clubsuit$

{\bf Remark.} We note that as a consequence, if $X_\bullet, A_\bullet$ are as above and $B_\bullet\subset A_\bullet$ is a semi-simplicial subvariety, then all maps in the long exact sequence of the triple $(X_\bullet, A_\bullet,B_\bullet)$ are also maps of mixed Hodge structures as they can be represented as compositions of maps from parts 1 and 2 of the proposition.

\smallskip

Let now $X_\bullet$ be a proper semi-simplicial variety, and let $A_\bullet\subset X_\bullet$ be a closed subvariety. Let  $X'_\bullet$ be a compactification of $X_\bullet$. Set $A'_\bullet$ to be the closure of $A_\bullet$ in $X'_\bullet$. By lemma~\ref{open_dense_cubical}, $\partial X'_\bullet=X'_\bullet\setminus X_\bullet$ is a semi-simplicial subvariety. We have then
$H^*_c(X_\bullet,A_\bullet,\Q)\cong H^*(X'_\bullet, A'_\bullet\cup\partial X'_\bullet,\Q)$, and we use this formula to put a mixed Hodge structure on the left hand side and, after dualising, on $\bar H_*(X_\bullet,A_\bullet,\Q)$.

Suppose now $f:Y_\bullet\to X_\bullet$ is a proper morphism of proper semi-simplicial complex algebraic varieties, $B_\bullet\subset Y_\bullet$ is a closed semi-simplicial subvariety, and $f(B_\bullet)\subset A_\bullet$. Then (cf.\ the proof of lemma~\ref{compactification}) there is a compactification $Y'_\bullet$ of $Y_\bullet$ such that $f$ extends to a morphism $Y'_\bullet\to X'_\bullet$, which we denote again by $f$. Using lemma~\ref{open_dense_cubical} again we see that $\partial Y'_\bullet=Y'_\bullet\setminus Y_\bullet$ is a semi-simplicial subvariety of $Y'_\bullet$. Since $Y_\bullet$ is dense in $Y'_\bullet$ (as $Y'_\bullet$ is a compactification of $Y_\bullet)$, and $f:Y_\bullet\to X_\bullet$ is proper, we have $f(\partial Y'_\bullet)\subset\partial X'_\bullet$. Set $B'_\bullet$ to be the closure of $B_\bullet$ in $Y'_\bullet$. We have $f(B'_\bullet)\subset A'_\bullet$, so $f$ induces a map of pairs $(Y'_\bullet, B'_\bullet\cup\partial Y'_\bullet)\to (X'_\bullet, A'_\bullet\cup\partial X'_\bullet)$. This implies that the maps $f^*:H^*_c(X_\bullet, A_\bullet,\Q)\to H^*_c(Y_\bullet, B_\bullet,\Q)$ and $f_*:\bar H_*(Y_\bullet,B_\bullet,\Q)\to \bar H_*(X_\bullet,A_\bullet,\Q)$ induced by $f$ are maps of mixed Hodge structures.

Finally, suppose $X_\bullet,X'_\bullet,\partial X'_\bullet, A_\bullet$ and $A'_\bullet$ are as above and $C_\bullet\subset A_\bullet$ is a closed subvariety. We set $C'_\bullet=$ the closure of $C_\bullet$ in $X'_\bullet$. The long exact sequence for compactly supported cohomology of the triple $(X_\bullet,A_\bullet,C_\bullet)$ is isomorphic to the cohomology long exact sequence of $(X'_\bullet,A'\cup \partial X'_\bullet,C'\cup\partial X'_\bullet)$, so all maps in both sequences are maps of mixed Hodge structures, and similarly for Borel-Moore homology. 
Summarising the last three paragraphs we get the following result.

\begin{Prop}\label{mhs_compact_support_bm}
If $X_\bullet$ is a proper semi-simplicial (in particular cubical) complex algebraic variety, and $A_\bullet\subset X_\bullet$ is a closed semi-simplicial subvariety, then $H^*_c(X_\bullet, A_\bullet, \Q)$ and $\bar H_*(X_\bullet, A_\bullet, \Q)$ have natural mixed Hodge structures, which are functorial with respect to proper morphisms.

If $C_\bullet\subset A_\bullet$ is a closed semi-simplicial subvariety, then the resulting long exact sequences of the cohomology groups with compact support, and of the Borel-Moore homology groups of the triple $(X_\bullet,A_\bullet,C_\bullet)$ are sequences of mixed Hodge structures. 
\end{Prop}
$\clubsuit$

\smallskip

{\bf Remark.} Another approach to the mixed Hodge structure on Borel-Moore homology of complex algebraic varieties is via M.\ Saito's theory of mixed Hodge modules, see e.g.~\cite[Corollary 14.9]{stepet}. There does not seem to be an easy way however to extend this theory to semi-simplicial or cubical varieties as the rationalisation functor (ibid., Axiom A in 14.1.1) is only defined on the level of derived categories. So we proceed via mixed Hodge complexes of sheaves and vector spaces, see~\cite{th_hodge_3}; more precisely, we use the version described in~\cite{gnpp} and~\cite[Chapter 5]{stepet}.

Moreover, it is not clear to us how to construct an explicit mixed Hodge complex using 
locally finite chains, or cochains with compact support (or a variant), which is why we use compactifications.

\subsection{K\"unneth theorem, cup product, and Thom isomorphism}\label{thom_iso_ccavb}

We will now prove the Thom isomorphism theorem for semi-simplicial complex algebraic vector bundles over proper semi-simplicial complex algebraic varieties (theorem~\ref{thom_ccavb}). This covers the case of CCAVB's, see p.~\pageref{ccavb}. There does not seem to be a reference or a shortcut for this, so we first need to define the relative cup product for (semi-)simplicial spaces and then prove that in the complex algebraic case the result is compatible with the mixed Hodge structures (theorem~\ref{cup_product_mhs}). In sections~\ref{products_etc} and~\ref{cones} we review and/or prove some preliminary facts about (semi-)simplicial spaces and varieties, in particular the K\"unneth theorem (theorem~\ref{simplicial_kunneth_mhs}). All spaces in this section are assumed compactly generated and Hausdorff. We equip all product spaces with the compactly generated topology induced by the product topology. In this section and the next we will sometimes abbreviate ``mixed Hodge structure'' to MHS.

\subsubsection{Products, geometric realisation, and the K\"unneth theorem}\label{products_etc} Given a semi-simplicial space $X_\bullet$ one can construct a simplicial one $EX_\bullet$ by adding degeneracies so that we get a left adjoint functor $E:X_\bullet \mapsto EX_\bullet$ to the forgetful functor $F$ from simplicial spaces to semi-simplicial spaces, see~\cite[Section 1.1]{ebert_randal_williams}. If $X_\bullet$ is a semi-simplicial complex algebraic variety, then $EX_\bullet$ is a simplicial complex algebraic variety. The components
\begin{equation}\label{comps_of_unit}
X_\bullet \to F(E(X_\bullet))
\end{equation}
of the unit of the adjunction induce weak equivalences of the geometric realisations for every semi-simplicial space~$X_\bullet$ (ibid., Lemma~2.6). 

{\bf Definition.} We define the {\it thin geometric realisation} $\|X_\bullet\|$ of a simplicial space $X_\bullet$ by the formula
$$\|X_\bullet\| = \coprod_{n\geq 0} X_n\times \Delta^n/(x,\theta_*y)\sim (\theta^*x,y), $$
where $\Delta^n$ is the standard simplex in $\R^{n+1}$ and $\theta\colon [m]\to[n] $ is a map from the category of finite ordered sets $\Delta$, see~\cite[Appendix A]{segal} or~\cite[Section 1.2]{ebert_randal_williams}.

{\bf Definition.} We set the {\it fat geometric realisation} $|X_\bullet|$ of a simplicial space $X_\bullet$ to be $|F(X_\bullet)|$. 

{\bf Remark.} We note that the fat geometric realisation does not take degeneracies into account, and so it may not coincide with the thin realisation even up to weak equivalence. However, for \emph{good} simplicial spaces in the sense of G.~Segal, i.e.\ for those where each degeneracy map is closed and has the homotopy extension property, the thin and fat realisations are weakly homotopy equivalent, see~\cite[Appendix A]{segal}. Also, simplicial complex algebraic varieties are good as in this case all degeneracy maps are inclusions of retracts, hence embeddings of closed subvarieties. In this subsection we will use both realisations, although in the rest of the paper we use the fat geometric realisation unless otherwise stated.

If $X_\bullet,Y_\bullet$ are simplicial spaces, then the {\it product} $(X\times Y)_\bullet$ is the simplicial space defined by $(X\times Y)_n=X_n\times Y_n$ with the product structure maps. By~\cite[Proposition A.1]{segal} (see also~\cite[Theorem 7.2]{ebert_randal_williams}), we have a weak equivalence
\begin{equation*}
|X_\bullet\times Y_\bullet|\xrightarrow{\simeq} |X_\bullet|\times |Y_\bullet|.
\end{equation*}
We also have a homeomorphism
\begin{equation*}
\|X_\bullet\times Y_\bullet\|\xrightarrow{\cong} \|X_\bullet\|\times \|Y_\bullet\|.
\end{equation*}
Both are induced by the projections $X_\bullet\times Y_\bullet\to X_\bullet$ and $X_\bullet\times Y_\bullet\to Y_\bullet$. The {\it product} $(X\times Y)_\bullet$ of two semi-simplicial spaces $X_\bullet,Y_\bullet$ is defined as in the simplicial case, i.e.\ level-wise. Note however that this construction does not commute with geometric realisation, either on the nose or up to homotopy.

\smallskip

The {\it (co)homology} of a simplicial space is defined as the (co)homology of the corresponding semi-simplicial space, see section~\ref{borel_moore}. Let $X_\bullet, Y_\bullet$ be simplicial spaces and $R$ a commutative ring. Using fat geometric realisations one can construct the K\"unneth map 
\begin{equation}\label{kunneth_map}
H^*(X_\bullet,R)\otimes H^*(Y_\bullet,R)\to H^*(X_\bullet\times Y_\bullet,R),
\end{equation}
which is an isomorphism under the usual assumptions, e.g.\ if $R$ is a field and at least one of the graded modules $H^*(X_\bullet,R), H^*(Y_\bullet,R)$ has finite rank in each degree.

A cochain level description of the K\"unneth map is as follows. Using singular or other cochains of the components of $X_\bullet$ and $Y_\bullet$ one obtains co-simplicial cochain complexes, which we denote $C^{*,*}_{X_\bullet}$, respectively $C^{*,*}_{Y_\bullet}$, with the first star being the cochain degree and the second the co-simplicial degree. Let $C^{*,*}(X_\bullet)$ and $C^{*,*} (Y_\bullet)$ be the resulting double complexes, see section~\ref{borel_moore}; the resulting total complexes $sC^*(X_\bullet)$ and $sC^* (Y_\bullet)$ compute the cohomology of $X_\bullet$, respectively $Y_\bullet$ with coefficients in $R$.

Starting from the co-simplicial cochain complexes $C^{*,*}_{X_\bullet}$ and $C^{*,*}_{Y_\bullet}$ one can form the {\it exterior tensor product} $C^{*,*}_{X_\bullet}\boxtimes C^{*,*}_{Y_\bullet}$, which is a bi-co-simplicial cochain complex, cf.~\cite[Definition 1.4]{ebert_randal_williams}, where the authors use $\otimes$ instead of $\boxtimes$. Associated to any bi-co-simplicial cochain complex $C$ are its {\it total (triple) complex} $TC$ and the {\it diagonal double complex} $DC$; the corresponding simple complexes $sTC$ and $sDC$ are quasi-isomorphic by the Cartier-Eilenberg-Zilber theorem, see~\cite[Section 2]{dold_puppe} and also~\cite{hanna} for a proof of the bi-simplicial analogue. For $C=C^{*,*}_{X_\bullet}\boxtimes C^{*,*}_{Y_\bullet}$ the complex $sTC$ is isomorphic to $sC^*(X_\bullet)\otimes sC^* (Y_\bullet)$, and $sDC$ is the cochain complex that computes $H^*(X_\bullet\times Y_\bullet,R)$. The K\"unneth map is
$$H^*(X_\bullet,R)\otimes H^*(Y_\bullet,R)\cong H^*(sC^*(X_\bullet))\otimes H^*(sC^* (Y_\bullet))\to H^*(
sC^*(X_\bullet)\otimes sC^* (Y_\bullet))\cong H^*(X_\bullet\times Y_\bullet,R).$$

We will abbreviate ``simplicial complex algebraic variety'' to SCAV. A subvariety $U_\bullet\subset X_\bullet$ of a smooth SCAV $X_\bullet$ is the {\it complement of a normal crossing divisor} if for every $n$, $X_n\setminus U_n$ intersected with each component is a divisor with normal crossings, the whole component, or empty. Note that $(X_n\setminus U_n)$ need not be a simplicial subvariety of $X_\bullet$. 
\begin{theorem}\label{simplicial_kunneth_mhs}
\begin{enumerate}
\item For every SCAV $X_\bullet$ there are a smooth compact SCAV $Y_\bullet$, the complement $U_\bullet\subset Y_\bullet$ of a normal crossing divisor and a SCAV map $U_\bullet\to X_\bullet$ that induces an isomorphism of $H^*(-,\Q)$.
\item If $X_\bullet$ and $Y_\bullet$ are simplicial complex algebraic varieties, the K\"unneth isomorphism
$$H^*(X_\bullet,\Q)\otimes H^*(Y_\bullet,\Q)\to H^*(X_\bullet\times Y_\bullet,\Q)$$ is an isomorphism of mixed Hodge structures.
\end{enumerate}
\end{theorem}

{\bf Proof.} To prove the first part we use~\cite[6.4.4]{th_hodge_3} to construct bi-simplicial complex algebraic varieties $X_{\bullet,\bullet}\subset\bar X_{\bullet,\bullet}$ such that each $\bar X_{m,n}$ is smooth and compact; each $X_{\bullet,n}$ is the complement of a normal crossing divisor in $\bar X_{\bullet,n}$; $X_{\bullet,\bullet}$ is 1-augmented towards $X_\bullet$; and for every $n$, the augmentation $X_{\bullet,n}\to X_n$ is of cohomological descent. (For the latter part see ibid., 5.3.5 (IV); 
see also~\cite[Theorem 7.9]{conrad_descent_notes} for a detailed proof that a proper hypercovering is of cohomological descent.) Then we take $U_\bullet$ to be the diagonal simplicial variety, i.e.\ we set $U_n=X_{n,n}$ for all $n$, and apply~\cite[6.4.2.2]{th_hodge_3}.

Let us now prove the second part. If $X_\bullet$ and $Y_\bullet$ are the complements of normal crossing divisors in smooth compact SCAV's, the statement we are after is ibid., 8.1.25. Applying part 1 we obtain the general case.$\clubsuit$

\subsubsection{Cones}\label{cones}

For a semi-simplicial space $X_\bullet$ and $n\in\Z_{\geq 1},i\in\Z_{\geq 0}, i\leq n$ let $d_i^{X_\bullet}:X_n\to X_{n-1}$ be the face map induced by the strictly increasing map $\{0,\ldots,n-1\}\to\{0,\ldots, n\}$ that misses $i$. 
The {\it cone} of $X_\bullet$ is the semi-simplicial space $Cone(X_\bullet)$ defined as follows
\begin{equation}\label{semi_simplicial_cone_deg_n}
Cone(X_\bullet)_n = \{(x,q)\mid x\in X_{n-q}, \; 0 \leq q\leq n+1\} =X_n\sqcup X_{n-1}\sqcup \ldots \sqcup X_{-1}, \;\; X_{-1}=pt,
\end{equation}
with face maps $d_i^{Cone(X_\bullet)}\colon Cone(X_\bullet)_n \to Cone(X_\bullet)_{n-1}$, $0\leq i\leq n$ given by
\begin{equation}\label{semi_simplicial_cone_face_maps}
d^{Cone(X_\bullet)}_i(x,q)=\left\{
\begin{array}{ll}
(d^{X_\bullet}_i(x),q), & \mbox{$0 \leq i\leq n-q$, $q\leq n$, $x\in X_{n-q},$}\\
(x,q-1), & \mbox{$n-q<i\leq n$, $q\leq n$, $x\in X_{n-q}$,} \\
(pt,n), & \mbox{$q=n+1$.}
\end{array}
\right.
\end{equation}
Note that there is a natural embedding $X_\bullet \hookrightarrow Cone(X_\bullet)$ given by $x\mapsto (x,0)$, $x\in X_n$, $n\geq 0$.

Given a map $f_\bullet: X_\bullet\to Y_\bullet$ of two semi-simplicial topological spaces, we define the {\it mapping cone} of $f_\bullet$ as the following pushout
$$Cone(f_\bullet) = Y_\bullet \sqcup_{X_\bullet} Cone(X_\bullet), $$
compare with~\cite[6.3]{th_hodge_3} and~\cite[IV.1.C]{gnpp}. 
(In~\cite[IV.1.C]{gnpp} the authors also describe a more economical construction of the cone and prove that it is homologous to the one given above, which in turn has the advantage that it extends to the simplicial case.) 

We can extend the definition of the cone (respactively the mapping cone) to simplicial spaces (respectively to maps of simplicial spaces). Namely, for a simplicial space $X_\bullet$ we let $s_i^{X_\bullet}:X_n\to X_{n+1}$ be the degeneracy map induced by the increasing map $\{0,\ldots, n+1\}\to \{0,\ldots,n\}$ that repeats $i$. Given a simplicial space $X_\bullet$, we define degeneracy maps $s^{Cone(X_\bullet)}_i\colon Cone(X_\bullet)_n \to Cone(X_\bullet)_{n+1}$, $0\leq i\leq n$ as follows
\begin{equation}\label{simplicial_cone_degeneracy_maps}
s^{Cone(X_\bullet)}_i(x,q)=\left\{
\begin{array}{ll}
(s^{X_\bullet}_i(x),q), & \mbox{$0 \leq i\leq n-q$, $q\leq n$, $x\in X_{n-q}$,}\\
(x,q+1), & \mbox{$n-q<i\leq n$, $q\leq n$, $x\in X_{n-q}$,} \\
(pt,n+2), & \mbox{$q=n+1$.}
\end{array}
\right.
\end{equation}
The face maps are defined as in the semi-simplicial case, see~\eqref{semi_simplicial_cone_deg_n} and~\eqref{semi_simplicial_cone_face_maps}. Finally, if $f_\bullet:X_\bullet\to Y_\bullet$ is a map of simplicial spaces, we again define the simplicial mapping cone $Cone(f_\bullet)$ as the pushout 
$$Cone(f_\bullet) = Y_\bullet \sqcup_{X_\bullet} Cone(X_\bullet).$$

\begin{lemma}\label{lemma_cone_is_contractible}
The (semi-)simplicial cone $Cone(X_\bullet)$ is contractible for any (semi-)simplicial space $X_\bullet$.
\end{lemma}

{\bf Proof.} We set $Cone(X_\bullet)_{-1}=X_{-1}=pt$, see~\eqref{semi_simplicial_cone_deg_n}. Then there is a unique (semi-)simplicial map $$\epsilon_\bullet \colon Cone(X_\bullet)\to Cone(X_\bullet)_{-1}$$ to the constant (semi-)simplicial space $Cone(X_\bullet)_{-1}$, and so the cone $Cone(X_\bullet)$ can be considered as an augmented (semi-)simplicial space, see~\cite[Definition~1.3]{ebert_randal_williams}. By~\cite[Lemma~1.12]{ebert_randal_williams}, it suffices to construct extra degeneracy maps $s_{n+1}\colon Cone(X_\bullet)_{n} \to Cone(X_\bullet)_{n+1}$, $n\geq -1$ such that
$$d^{Cone(X_\bullet)}_{n+1}s_{n+1} = \id_{Cone(X_\bullet)_n}, \;\; n\geq -1;$$
$$d^{Cone(X_\bullet)}_{i}s_{n+1}=s_{n}d^{Cone(X_\bullet)}_i, \;\; 0\leq i<n+1;$$
$$\epsilon_0s_0 =\id_{Cone(X_\bullet)_{-1}}; $$
and if $X_\bullet$ is a simplicial space, then we also require 
$$s^{Cone(X_\bullet)}_{i}s_{n+1}=s_{n+2}s^{Cone(X_\bullet)}_i, \;\; 0\leq i\leq n+1. $$
In the notation of~\eqref{semi_simplicial_cone_deg_n}, we define the extra degeneracy maps $s_{n+1}\colon Cone(X_\bullet)_{n} \to Cone(X_\bullet)_{n+1}$, $n\geq -1$ by $s_{n+1}(x,q)=(x,q+1)$, $x\in X_{n-q}$, $0 \leq q\leq n+1$. A routine check using formulas~\eqref{semi_simplicial_cone_face_maps} and~\eqref{simplicial_cone_degeneracy_maps} shows that the maps $s_{n+1}$ has all required properties.$\clubsuit$

The (semi-)simplicial $Cone(f_\bullet)$ is {\it pointed}, i.e.\ there is a map $pt_\bullet\to Cone(f_\bullet)$ where $pt_\bullet$ denotes the constant (semi-)simplicial point. For pointed (semi-)simplicial spaces one can define the {\it wedge sum} in a straightforward way, and in the simplicial case also the {\it smash product}. If $X_\bullet$ is a pointed (semi-)simplicial space, we define the {\it reduced cohomology} $\tilde H^*(X_\bullet)$ of $X_\bullet$ as $\ker( H^*(X_\bullet)\to H^*(pt_\bullet))$. In the simplicial case, the K\"unneth map~(\ref{kunneth_map}) takes $H^*(\tilde X_\bullet,R)\otimes \tilde H^*(Y_\bullet,R)$ to $\tilde H^*(X_\bullet\times Y_\bullet,R)$. 
If $(X_\bullet,A_\bullet)$ is a pair of (semi-)simplicial spaces, we define the {\it relative cohomology} $H^*(X_\bullet, A_\bullet)$ to be $\tilde H^*(Cone(A_\bullet\to X_\bullet))$. In the semi-simplicial case we get the same result as above in section~\ref{borel_moore}.

In the sequel we wish to study the cup product on the relative cohomology groups $H^*(X_\bullet, A_\bullet)$. The mapping cone $Cone(f_\bullet)$ is not quite well suited for this purpose. Therefore, we will introduce a different model for the mapping cone.

Given a (semi-)simplicial space $X_\bullet$, we define a (semi-)simplicial space $C(X_\bullet)$ as follows
\begin{equation}\label{semi_simplicial_naive_cone_deg_n}
C(X_\bullet)_n = \{(x,q): x\in X_{n}, \; 0 \leq q\leq n\} \sqcup X_{-1} =X_n\sqcup X_{n}\sqcup \ldots \sqcup X_{n} \sqcup X_{-1}, \;\; X_{-1}=pt,
\end{equation}
with face maps $d_i^{C(X_\bullet)}\colon C(X_\bullet)_n \to C(X_\bullet)_{n-1}$, $0\leq i\leq n$ given by
\begin{equation}\label{semi_simplicial_naive_cone_face_maps}
d^{C(X_\bullet)}_i(x,q)=\left\{
\begin{array}{ll}
(d^{X_\bullet}_i(x),q), & \mbox{$0 \leq i\leq n-q$, $q\leq n$, $x\in X_{n}$}\\
(d^{X_\bullet}_i(x),q-1), & \mbox{$n-q<i\leq n$, $q\leq n$, $x\in X_{n}$,}
\end{array}
\right.
\end{equation}
and in the simplicial case, with degeneracy maps $s^{C(X_\bullet)}_i\colon C(X_\bullet)_n \to C(X_\bullet)_{n+1}$, $0\leq i\leq n$ given by
\begin{equation}\label{simplicial_naive_cone_degeneracy_maps}
s^{C(X_\bullet)}_i(x,q)=\left\{
\begin{array}{ll}
(s^{X_\bullet}_i(x),q), & \mbox{$0 \leq i\leq n-q$, $q\leq n$, $x\in X_{n}$,}\\
(s^{X_\bullet}_i(x),q+1), & \mbox{$n-q<i\leq n$, $q\leq n$, $x\in X_{n}$.}
\end{array}
\right.
\end{equation}
Finally, all face and degeneracy maps preserve the last component $X_{-1} \subset C(X_\bullet)_n$. 

Given a map $f_\bullet: X_\bullet\to Y_\bullet$ of two (semi-)simplicial topological spaces, we define a (semi-)simplicial space $C(f_\bullet)$ as the following pushout
$$C(f_\bullet) = Y_\bullet \sqcup_{X_\bullet} C(X_\bullet), $$
where $X_\bullet \hookrightarrow C(X_\bullet)$ is a natural embedding given by $x\mapsto (x,0)$, $x\in X_n$, $n\geq 0$. 

{\bf Remark.} We note that the (semi-)simplicial space $C(X_\bullet)$ is isomorphic to the quotient (semi-)simplicial space $$X_\bullet \times \Delta^1/X_\bullet \times \{1\},$$ where
$\Delta^1$ is the simplicial interval, i.e.\ $\Delta^1$ is the simplicial set such that $\Delta^1$ has exactly two non-degenerate $0$-simplices, one non-degenerate $1$-simplex, and the set of $n$-simplices $|\Delta^1_{n}|$ has cardinality $n+1$. In particular, $C(X_\bullet)$ is contractible. 

{\bf Remark.} Since the thin geometric realisation commutes with products, the cone operation $C(-)$ commutes with thin geometric realisation, i.e.\ the thin geometric realisation of the cone of a map is homeomorphic to the cone of the induced map of the thin geometric realisations.

\begin{lemma}\label{lemma_geometric_realisation_hep}
Let $f_\bullet\colon A_\bullet \to X_\bullet$ be a map of simplicial spaces such that both $A_\bullet, X_\bullet$ are good and the level-wise maps $f_n\colon A_n \to X_n$, $n\geq 0$ are closed and have the homotopy extension property. Then the induced map $\|f_\bullet\|$ of the thin geometric realisations also has the homotopy extension property.
\end{lemma}

{\bf Proof.} We endow the category $\mathcal{T}$ of topological spaces with the Str\o m model structure. Recall that in this case the weak equivalences are the homotopy equivalences and the cofibrations are the closed maps which have the homotopy extension property, see~\cite{strom}. Next, we endow the category $s\mathcal{T}$ of simplicial spaces with the Reedy model structure~\cite[Chapter~15]{hirschhorn}. By applying Lillig's union theorem~\cite{lillig}, we can observe that the map $f_\bullet\colon A_\bullet \to X_\bullet$ is a Reedy cofibration provided $f_\bullet$ satisfies the assumptions of the lemma, see e.g. the introduction of~\cite{angelini-knoll-salch}. Finally, by e.g.~\cite[Proposition~3.45]{horel}, the thin geometric realisation $\|f_\bullet\|$ is a cofibration in the Str\o m model structure, i.e. it is closed and has the homotopy extension property.$\clubsuit$

\begin{lemma}\label{lemma_cone_and_contraction}
Let $f_\bullet\colon A_\bullet \to X_\bullet$ be a map of simplicial spaces such that both $A_\bullet, X_\bullet$ are good and the level-wise maps $f_n\colon A_n \to X_n$, $n\geq 0$ are closed and have the homotopy extension property. Then the natural map $C(f_\bullet) \to X_\bullet/A_\bullet$ induces a homotopy equivalence after taking the thin geometric realisation.
\end{lemma}

{\bf Proof.} Consider the commutative square
$$
\begin{tikzcd}
\|C(f_\bullet)\| \arrow{r}{\cong}\arrow{d} & C(\|f_\bullet\|) \arrow{d}{\simeq}\\
\|X_\bullet/A_\bullet\|\arrow{r}{\cong} & \|X_\bullet\|/\|A_\bullet\|.
\end{tikzcd}
$$
Here the horizontal maps are homeomorphisms and the left vertical map is a homotopy equivalence by lemma~\ref{lemma_geometric_realisation_hep} and~\cite[Proposition~6.6]{switzer}. The lemma follows.$\clubsuit$

\begin{lemma}\label{lemma_two_cones_give_same_model}
Let $X_\bullet$ be a (semi-)simplicial space. Then there is a homotopy equivalence 
$$\varphi_\bullet \colon C(X_\bullet) \to Cone(X_\bullet) $$
which commutes with the natural embeddings $X_\bullet \hookrightarrow C(X_\bullet)$, $X_\bullet \hookrightarrow Cone(X_\bullet)$. In particular, if $f_\bullet: X_\bullet\to Y_\bullet$ is a map of (semi-)simplicial spaces and $X_\bullet$ is good, then the induced map $\|C(f_\bullet)\| \to \|Cone(f_\bullet)\|$ is a homotopy equivalence.
\end{lemma}

{\bf Proof.} Since both (semi-)simplicial spaces $C(X_\bullet)$, $Cone(X_\bullet)$ are contractible (see lemma~\ref{lemma_cone_is_contractible} and the first remark before lemma~\ref{lemma_geometric_realisation_hep}), any map $\varphi_\bullet \colon C(X_\bullet) \to Cone(X_\bullet)$ will be a homotopy equivalence.

In the notation of~\eqref{semi_simplicial_cone_deg_n} and~\eqref{semi_simplicial_naive_cone_deg_n}, we define maps $\varphi_{n}\colon C(X_\bullet)_{n} \to Cone(X_\bullet)_{n}$, $n\geq 0$ by the following rule 
$$\varphi_{n}(x,q)=\left\{
\begin{array}{ll}
(x,0), & \mbox{$q=0$, $x\in X_n$,}\\
(d^{X_\bullet}_{n-q+1}\circ \ldots \circ d^{X_\bullet}_{n}(x),q), & \mbox{$1 \leq q\leq n$, $x\in X_{n}$.}
\end{array}
\right.
$$
Finally, we set $\varphi_n(X_{-1})=X_{-1}\subset Cone(X_\bullet)_n$. A routine check with the simplicial identities and formulas~\eqref{semi_simplicial_cone_face_maps},~\eqref{semi_simplicial_naive_cone_face_maps},~\eqref{simplicial_cone_degeneracy_maps}, and~\eqref{simplicial_naive_cone_degeneracy_maps} shows that the maps $\varphi_{n}$, $n\geq 0$ indeed define a map of (semi-)simplicial spaces $\varphi_{\bullet}$ that commutes with the embeddings of $X_\bullet$.

For the second part, we note that if $X_\bullet$ is a good simplicial space, then both cones $C(X_\bullet)$ and $Cone(X_\bullet)$ are also good and the natural embeddings $X_\bullet \hookrightarrow C(X_\bullet)$, $X_\bullet \hookrightarrow Cone(X_\bullet)$ satisfy the assumptions of lemma~\ref{lemma_geometric_realisation_hep}. This implies that the homotopy equivalence $\|C(X_\bullet)\| \to \|Cone(X_\bullet)\|$ also induces a homotopy equivalence $\|C(f_\bullet)\| \to \|Cone (f_\bullet)\|$ of the mapping cones.$\clubsuit$



\begin{lemma}\label{lemma_product_of_cones}
If $A_\bullet\subset X_\bullet$ and $B_\bullet\subset Y_\bullet$ are pairs of simplicial spaces, then we have a map of simplicial spaces
\begin{equation}\label{product_of_cones}
\psi_{\bullet}\colon C (A_\bullet\to X_\bullet)\times C (B_\bullet\to Y_\bullet)\to C((A_\bullet\times Y_\bullet)\cup (X_\bullet\times B_\bullet)\to X_\bullet\times Y_\bullet)
\end{equation}
such that the following diagram commutes
$$
\begin{tikzcd}[column sep=huge]
C (A_\bullet\to X_\bullet)\times C (B_\bullet\to Y_\bullet) \arrow{r}{\psi_\bullet}\arrow{d} & C((A_\bullet\times Y_\bullet)\cup (X_\bullet\times B_\bullet)\to X_\bullet\times Y_\bullet) \arrow{d}\\
X_\bullet/A_\bullet \times Y_\bullet/B_\bullet \arrow{r} & X_\bullet \times Y_\bullet /((A_\bullet\times Y_\bullet)\cup (X_\bullet\times B_\bullet)).
\end{tikzcd}
$$

Moreover, if $A_\bullet\subset X_\bullet$ is a pair of SCAV's, then the cone $C(A_\bullet\to X_\bullet)$ is also a SCAV. If $B_\bullet\subset X_\bullet$ is another pair of SCAV's, then~(\ref{product_of_cones}) can be chosen to be a morphism of SCAV's.
\end{lemma}

{\bf Proof.} By formula~\eqref{semi_simplicial_naive_cone_deg_n}, the space of $n$-simplices of the left hand side of~(\ref{product_of_cones}) consists of quadruples $(x,q_1,y,q_2)$, where $q_1,q_2 \in \mathbb{Z}_{\geq 0}$, $0\leq q_1,q_2\leq n+1$, 
$$x \in \left\{
\begin{array}{ll}
X_n, & \mbox{$q_1=0$,}\\
A_n, & \mbox{$1 \leq q_1\leq n$,} \\
A_{-1}=pt, & \mbox{$q_1=n+1$,}
\end{array}
\right.
$$
and similarly,
$$y \in \left\{
\begin{array}{ll}
Y_n, & \mbox{$q_2=0$,}\\
B_n, & \mbox{$1 \leq q_2\leq n$,} \\
B_{-1}=pt, & \mbox{$q_2=n+1$.}
\end{array}
\right.
$$
Analogously, the space of $n$-simplices of the right hand side of~\eqref{product_of_cones} consists of triples $(x,y,q)$, where $q\in \mathbb{Z}_{\geq 0}$, $0\leq q\leq n+1$, and
$$(x,y) \in \left\{
\begin{array}{ll}
X_n\times Y_n, & \mbox{$q=0$,}\\
A_n \times Y_n \cup X_n\times B_n, & \mbox{$1 \leq q\leq n$,} \\
pt, & \mbox{$q=n+1$.}
\end{array}
\right.
$$

We now define maps $\psi_n$, $n\geq 0$ as follows:
$$\psi_n(x,q_1,y,q_2) = \left\{
\begin{array}{ll}
(x,y,\max(q_1,q_2)), & \mbox{$q_1,q_2\leq n$,}\\
(pt,n+1), & \mbox{$\max(q_1,q_2)=n+1$.}
\end{array}
\right.
$$
A routine check with formulas~\eqref{semi_simplicial_naive_cone_face_maps} and~\eqref{simplicial_naive_cone_degeneracy_maps} shows that the maps $\psi_{n}$, $n\geq 0$ indeed define a map of simplicial spaces $\psi_{\bullet}$ and the desired diagram is commutative.$\clubsuit$ 

A pair $\{A_1,A_2\}$ of subspaces of a topological space $X$ is {\it excisive} if the inclusion of the complexes of singular chains $C_*(A_1,\Z)+C_*(A_2,\Z)\subset C(A_1\cup A_2,\Z)$ is a quasi-isomorphism. This definition extends level-wise to (semi-)simpliclal spaces. Recall that an excisive pair of subspaces induces a Mayer-Vietoris long exact sequence.

\begin{corollary}\label{corollary_cones_and_kunneth}
Let $A_\bullet \subset X_\bullet$, $B_\bullet\subset Y_\bullet$ be pairs of simplicial spaces, and assume that $\{A_\bullet\times Y_\bullet,X_\bullet\times B_\bullet\}$ is an excisive pair of simplicial subspaces of $X_\bullet\times Y_\bullet$. Then, for an arbitrary commutative ring $R$, the map $\psi^*_\bullet$ of $\tilde H^*(-,R)$ induced by~\eqref{product_of_cones} is injective, and its image contains the image of $$\tilde H^*(C(A_\bullet\to X_\bullet),R)\otimes \tilde H^*(C(B_\bullet\to Y_\bullet),R)$$ under the K\"unneth map~(\ref{kunneth_map}).
\end{corollary}

{\bf Proof.} If all spaces $A_\bullet, X_\bullet, B_\bullet, Y_\bullet$ are good and the inclusions $A_\bullet \subset X_\bullet$, $B_\bullet\subset Y_\bullet$ are level-wise closed embeddings with the homotopy extension property, we use lemmas~\ref{lemma_cone_and_contraction} and~\ref{lemma_product_of_cones}. To prove the general case we apply level-wise the functor $\|Sing(-)\|$ to all spaces involved where $Sing(-)$ denotes the singular simplicial set functor. Note that $\|Sing(-)\|$ preserves products; moreover, if $A\subset X, B\subset Y$ are pairs of spaces and $\{A\times Y,X\times B\}$ is an excisive pair of subspaces of $X\times Y$, then there are cohomology isomorphisms
$$H^*(A\times Y\cup X\times B) \cong H^*(\|Sing(A\times Y\cup X\times B)\|)\cong H^*(\|Sing(A)\|\times \|Sing(Y)\|\cup \|Sing(X)\|\times \|Sing (B)\|)$$ with any coefficients. This allows one to reduce the general case to the one we have already considered above.$\clubsuit$

\subsubsection{Relative cup product and the Thom isomorphism}
\begin{theorem}\label{cup_product_mhs}
Let $X_\bullet$ be a (semi-)simplicial topological space, and let $A_\bullet, B_\bullet\subset X_\bullet$ be subspaces. Assume that $\{A_\bullet\times Y_\bullet,X_\bullet\times B_\bullet\}$ is an excisive pair of (semi-)simplicial subspaces of $X_\bullet\times Y_\bullet$. Then for any commutative ring~$R$, there is a cup product map $$H^*(X_\bullet,A_\bullet,R)\otimes H^*(X_\bullet,B_\bullet,R)\to H^*(X_\bullet,A_\bullet\cup B_\bullet,R),$$ which is moreover a map of mixed Hodge structures in the complex algebraic case for $R=\Q$.
\end{theorem}

{\bf Proof.} To begin with, note that the semi-simplicial case follows from the simplicial one: The functor $E$ (see the beginning of section~\ref{products_etc}) preserves unions, and the maps~(\ref{comps_of_unit}) induce cohomology isomorphisms. 
(In particular, for semi-simplicial spaces $X_\bullet\supset A_\bullet, B_\bullet$, if $\{A_\bullet\times Y_\bullet,X_\bullet\times B_\bullet\}$ is excisive in $X_\bullet\times Y_\bullet$, then so is $\{EA_\bullet\times EY_\bullet,EX_\bullet\times EB_\bullet\}$ in $EX_\bullet\times EY_\bullet$.)

In the simplicial case the cup product map is constructed as follows:
\begin{multline*}
H^*(X_\bullet,A_\bullet)\otimes H^*(X_\bullet,B_\bullet)\cong \tilde H^*(Cone (A_\bullet\to X_\bullet))\otimes \tilde H^*(Cone(B_\bullet\to X_\bullet))\\
\stackrel{K}{\to} \tilde H^*(Cone (A_\bullet\to X_\bullet)\times Cone(B_\bullet\to X_\bullet))\to \tilde H^*(Cone((A_\bullet\times X_\bullet)\cup (X_\bullet\times B_\bullet)\to X_\bullet\times X_\bullet))\\
\stackrel{D}{\to} \tilde H^*(Cone(A_\bullet\cup B_\bullet\to X_\bullet))\cong H^*(X_\bullet,A_\bullet\cup B_\bullet).
\end{multline*}
Here all cohomology is with coefficients in $R$; the first and last isomorphisms follow by definition; the map marked $K$ is the K\"unneth map~(\ref{kunneth_map}); the arrow in the middle of the second line is obtained using lemmas~\ref{lemma_two_cones_give_same_model},~\ref{lemma_product_of_cones} and corollary~\ref{corollary_cones_and_kunneth}, and the arrow marked $D$ is obtained from the diagonal embedding of $X_\bullet$, which induces the commutative diagram
$$
\begin{tikzcd}
A_\bullet\cup B_\bullet \ar[r]\ar[d] & X_\bullet\ar[d]\\
(A_\bullet\times X_\bullet)\cup (X_\bullet\times B_\bullet)\ar[r] & X_\bullet\times X_\bullet.
\end{tikzcd}
$$

In the complex algebraic case and taking $R=\Q$, the K\"unneth map is compatible with mixed Hodge structures by theorem~\ref{simplicial_kunneth_mhs}, and by construction so are all other ingredients in the definition of the cup product and the cohomology map induced by~(\ref{comps_of_unit}).$\clubsuit$

\smallskip

We are now ready to prove the Thom isomorphism theorem. Instead of vector bundles we consider the slightly more general case of affine ones. A {rank $r$ \it affine algebraic bundle} over an algebraic variety $X$ is a Zariski locally trivial bundle with fibre $A_\C^r$ and structure group equal the group of affine automorphisms. 
(Semi)-simplicial and cubical complex algebraic affine bundles are defined in a straightforward way.

\begin{theorem}\label{thom_ccavb}
Let $E_\bullet$ be a semi-simplicial complex algebraic affine bundle of rank $r$ over a proper semi-simplicial complex algebraic variety $X_\bullet$ augmented towards a complex algebraic variety $X$. (This includes the case when $E_\bullet$ is a CCAVB over a proper cubical complex algebraic variety.) Then we have $\bar H_*(E_\bullet,\Q)\cong \bar H_*(X_\bullet,\Q)\otimes \Q(r)[-2r]$ where $\Q(r)[-2r]$ denotes the Tate Hodge structure $\Q(r)$ placed in degree~$2r$.
\end{theorem}

{\bf Proof.} Note that $|E_\bullet|$ is a topological vector bundle over $|X_\bullet|$. So we have the topological Thom isomorphism
\begin{equation}\label{thom_iso_comp_sup_cubical}
H^*_c(X_\bullet,\Q)\to H^{*+2r}_c(E_\bullet,\Q),
\end{equation}
and to prove theorem~\ref{thom_ccavb} it would suffice to show that~(\ref{thom_iso_comp_sup_cubical}) is a map of mixed Hodge structures. 

Using lemma~\ref{compactification} we construct an augmented algebraic compactification $X'_\bullet$ of 
$X_\bullet$. Let $E\to X$ be the component of $E_\bullet$ over the augmentation variety $X$. We compactify the bundle projection $E_\bullet\to X_\bullet$ by taking $\bar E$ to be a complex algebraic variety over $X$ that contains $E$ as an open dense subvariety, and then setting $\bar E_i=\bar E\times_X X_i, i\geq 0$. 
We then construct an algebraic compactification $Y'_\bullet$ of $\bar E_\bullet$ such that the proper morphism $\bar E_\bullet\to X_\bullet$ extends to a (proper) morpshim $f:Y'_\bullet\to X'_\bullet$ of augmented semi-simplicial varieties. Let $X'$ and $Y'$ be the target varieties of the augmentation morphisms $X_\bullet'\to X', Y_\bullet'\to Y'$. 
Note that $Y'_\bullet$ contains $\bar E_\bullet$ and hence $E_\bullet$ as an open dense subvariety. So by lemma~\ref{open_dense_cubical}, $B'_\bullet=Y'_\bullet\setminus E_\bullet$ is a closed cubical subvariety of $Y'_\bullet$. It is augmented, and we let $B'=Y'\setminus E$ denote the target of the augmentation $B'_\bullet \to B'$. The topological Thom class of $E_\bullet$ lives in $H^{2r}_c(E_\bullet,\Q)\cong H^{2r}(Y'_\bullet, B'_\bullet,\Q)$.

Applying lemma~\ref{open_dense_cubical} again, we see that $\partial X'_\bullet=X'_\bullet\setminus X_\bullet$ is a closed cubical subvariety of $X'_\bullet$. Let $\partial Y'_\bullet=f^{-1}(\partial X'_\bullet)$. Observe that $\partial Y'_\bullet\subset B'_\bullet$. 
The Thom isomorphism~(\ref{thom_iso_comp_sup_cubical}) is then the composition
$$
\begin{tikzcd}
H^*_c(X_\bullet)\ar[r,"\cong"] & H^*(X'_\bullet,\partial X'_\bullet)\ar[r,"f^*"] &  H^*(Y'_\bullet,\partial Y'_\bullet)\ar[r,"-\smile u"] & H^{*+2r}(Y'_\bullet,B'_\bullet) \ar[r,"\cong"] & H^{*+2r}_c(E_\bullet)
\end{tikzcd}
$$
where all coefficients are rational and $u\in H^{2r}(Y'_\bullet,B'_\bullet)\cong H^{2r}(|Y'_\bullet|,|B'_\bullet|)$ is the topological Thom class of $E'_\bullet$. 

Let us show that $u$ is Tate of weight $2r$. Let $E_{\bullet,const}, X'_{\bullet,const},Y'_{\bullet,const},B'_{\bullet,const}$ be the constant semi-simplicial varieties that correspond to $E, X', Y',B'$ respectively. We have morphisms $X'_{\bullet}\to X'_{\bullet,const}, Y'_{\bullet}\to Y'_{\bullet,const}$ of semi-simplicial complex algebraic varieties, and the class $u$ is the image of the Thom class $u'$ of $E_{\bullet,const}$ under
$$H^{2r}_c(E_{\bullet, const})\cong H^{2r}(Y'_{\bullet,const},B'_{\bullet,const})\to H^{2r}(Y'_{\bullet},B'_{\bullet})\cong H^{2r}_c(E_\bullet).$$ The class $u'$ is Tate of weight $2r$, as can be seen e.g.\ by pulling back $E_{\bullet, const}$ to a ``semi-simplicial point'' $pt_\bullet\to  X_{\bullet,const}$. So $u$ is also Tate of weight $2r$. We finish the proof of theorem~\ref{thom_ccavb} by applying proposition~\ref{mhs_coh} and theorem~\ref{cup_product_mhs}.$\clubsuit$

{\bf Remark.} One can show that if $X$ is a scheme that is separated and of finite type over a Noetherian scheme $S$, then any locally free sheaf $E$ on $X$ of rank $r<\infty$ extends to a locally free sheaf of rank $r$ on some compactification of $X$ over $S$, cf.~\cite[Lemma 1]{diaz}. 
So if in theorem~\ref{thom_ccavb}, $E_\bullet$ is a semi-simplicial vector bundle, and not just an affine one, then one has an explicit candidate for the semi-simplicial variety $Y'_\bullet$ from the proof of the theorem. Namely, one can extend $E_\bullet$ to a semi-simplicial algebraic vector bundle $E'_\bullet$ on some compactification $X'_\bullet$ of $X_\bullet$, and then take $Y'_\bullet$ to be $\p(E'_\bullet\oplus \OO_{X_\bullet'})$.

\subsection{Spectral sequences}\label{sec_mhs_sseq}

Recall that we write MHS for ``mixed Hodge structure''. A homological spectral sequence $(E_{pq}^r,d^r)$ of rational vector spaces that converges to a graded rational vector space $H_*$ {\it has a natural mixed Hodge structure starting from $E^r$} if $H_*$ and all $E_{p,q}^r,\infty\geq r\geq 1$ are MHSs, all differentials $d^{\geq r}$ are maps of MHSs, the resulting filtration on $H_*$ is a filtration by MHSs, and $E^\infty$ is isomorphic to the associated graded as an MHS. There is a straightforward analogue of this in the cohomological case.

\begin{Prop}\label{mhs_sseq_filtred}
Let $X_\bullet$ be a proper semi-simplicial (for example cubical) complex algebraic variety. Suppose $X_\bullet^i\subset X_\bullet,i=-1,\ldots, N$ are closed semi-simplicial subvarieties such that  $X^{-1}_\bullet=\varnothing$, $X^N_\bullet=X_\bullet$ and $X^{i-1}_\bullet\subset X^{i}_\bullet$ for $i=0,\ldots, N$. The corresponding Borel-Moore homology spectral sequence $(E_{p,q}^r,d^r)\Rightarrow \bar H_*(X_\bullet,\Q)$ has a natural MHS starting from $E^1$. 
\end{Prop}

{\bf Proof.} 
The proposition follows from proposition~\ref{mhs_compact_support_bm}: one of the ways to construct the spectral sequence $(E_{p,q}^r,d^r)$ is to use the exact couple obtained by weaving together the individual long exact sequences of the pairs $(X^{i}_\bullet,X^{i-1}_\bullet),i=0,\ldots, N$.$\clubsuit$
\smallskip

Let now $X_\bullet$ be a proper semi-simplicial (for example cubical) complex algebraic variety, and let $A_\bullet\subset X_\bullet$ be a closed semi-simplicial subvariety. Suppose $X_\bullet$ is augmented towards a complex algebraic variety $X$. Set $U=|X_\bullet|\setminus |A_\bullet|$, and let $j:U\to |X_\bullet|$ be the inclusion. We have
$$H^*_c(|X_\bullet|,|A_\bullet|,\Q)\cong H^*_c(|X_\bullet|,j_!(\underline{\Q}_U)), \bar H_*(|X_\bullet|,|A_\bullet|,\Q)\cong H^{-*}(|X_\bullet|,Rj_* (D_U)).$$ (Recall that we use $D_X$ to denote the dualising complex of $X$, see section~\ref{sec_bm_hom}. The second formula follows from~(\ref{bm_open_subset}); note that $D_U$ is quasi-isomorphic to a bounded complex of sheaves.) The Leray spectral sequences for the augmentation map $\varepsilon:|X_\bullet|\to X$ and the complexes $j_!(\underline{\Q}_U)$ and $Rj_* (D_U)$ converge to $H^*_c(|X_\bullet|,|A_\bullet|,\Q)$, respectively $\bar H_*(|X_\bullet|,|A_\bullet|,\Q)$.

\begin{Prop}\label{mhs_sseq_leray}
Assume $X$ is quasi-projective. Then each of the Leray spectral sequences has an MHS starting from the second page.
\end{Prop}

{\bf Proof.} The Borel-Moore homology spectral sequence is the linear dual of the one for compactly supported cohomology from the second page on, 
so it suffices to consider the latter, in which case the argument used in the proof of~\cite[Corollary 4.4]{ara} works with minor modifications. Let us review the main steps.

We may assume that $X$ is projective and $X_\bullet$ is compact by lemma~\ref{compactification}, cf.\ also the proof of proposition~\ref{mhs_compact_support_bm}.
We use the Jouanolou trick to construct an affine bundle $m:Y\to X$, and we let $Y_\bullet$ be the pullback of $X_\bullet$ to $Y$. Set $k:|Y_\bullet|\to |X_\bullet|$ to be the map of the geometric realisations. The augmentation map $\varepsilon$ is proper, and for every sheaf $S$ on $X$ the map $H^*(X,S)\to H^*(Y,m^{-1}(S))$ is an isomorphism. So the Leray spectral sequence for $\varepsilon:|X_\bullet|\to X$ and the sheaf $j_!(\underline{\Q}_U))$ is isomorphic to the Leray spectral sequence of the augmentation $\epsilon:|Y_\bullet|\to Y$ and the sheaf $T=k^{-1}(j_!(\underline{\Q}_U))=j'_!(\underline{\Q}_V))$ where $V=k^{-1}(U)$ and $j':V\to |Y_\bullet|$ is the inclusion (cf.~ibid, Lemma~4.2). 

All direct images $R^q\epsilon_* T$ are algebraically constructible. Using ibid., Lemma~3.7 we construct a filtration $Y^{-1}=\varnothing\subset Y^0\subset\cdots\subset Y^N=Y$ of $Y$ by closed subvarieties which is {\it cellular} with respect to $T_\oplus=\bigoplus_q R^q\epsilon_* T$, meaning that $\dim Y_i=i$ for all $i=0,\ldots, N$, and the relative cohomology groups $H^p(Y^{i},Y^{i-1},T_\oplus)$ are zero except maybe for $p=i$. Set $Y^i_\bullet=\epsilon^{-1}(Y_i)$. The long exact sequence of the triple $(|Y_\bullet|, |Y^i_\bullet|, |Y^{i-1}_\bullet|)$ with coefficients in $T$ is isomorphic to the long exact sequence of the triple $(|Y_\bullet|, |Y^i_\bullet\cup A_\bullet|, |Y^{i-1}_\bullet\cup A_\bullet|)$ with $\Q$ coefficients. The latter are sequences of MHSs by the remark after proposition~\ref{mhs_coh}, and formula~(\ref{geom_realis_hom_coh}). 

By combining these long exact sequences we get an exact couple and hence a spectral sequence that converges to $H^*(|Y_\bullet|,T)$, see~ibid., proof of Lemma~3.8. The Leray spectral sequence for $\epsilon$ and $T$ maps into this spectral sequence (ibid., Lemma 3.13), and since the filtration on $Y$ is cellular, the map of the spectral sequences is an isomorphism from page 2 onwards (ibid., proof of Theorem~3.1).$\clubsuit$

{\bf Remark.} We use the notation of the proof of the proposition. Let $j^X:W^X\to X, j^Y:W^Y\to Y$ be open embeddings such that $\varepsilon(W^X)\subset U,\epsilon (W^Y)\subset V$. Then, since $\varepsilon$ and $\epsilon$ are proper, we have
$$(R^q \epsilon_*) j'_!(\underline{\Q}_V)= (R^q \epsilon_!) j'_!(\underline{\Q}_V)=j^Y_! R^q\epsilon_!(\underline{\Q}_V),\mbox{ and } (R^q \varepsilon_*) j_!(\underline{\Q}_U)= (R^q \varepsilon_!) j_!(\underline{\Q}_U)=j^X_! R^q\varepsilon_!(\underline{\Q}_U).$$

{\bf Remark.} The assumption that $X$ should be quasi-projective is technical and it is only needed in order to apply the Jouanolou trick. We expect that this assumption can be removed.

Aleksandr Berdnikov: School of Mathematics, Institute for Advanced Study, 1 Einstein Dr., Princeton, NJ 08540, USA, email: \url{aberdnik@ias.edu}, \url{beerdoss@mail.ru}

Alexey Gorinov: Faculty of Mathematics, Higher School of Economics, 6 Usacheva ulitsa, Moscow, Russia 119048, email: \url{agorinov@hse.ru}, \url{gorinov@mccme.ru}

Nikolay Konovalov: Department of Mathematics, University of Notre Dame, 255 Hurley Hall, Notre Dame, IN 46556, USA, email: \url{nkonoval@nd.edu}, \url{ nikolay.konovalov.p@gmail.com}
\end{document}